\documentclass[a4paper,12pt]{article}
\usepackage[all]{xy}
\usepackage[T1]{fontenc}
\usepackage[dvips]{graphicx}
\usepackage[final]{pdfpages}
\usepackage{amsmath,amssymb,amsfonts,amsthm,latexsym,mathrsfs,textcomp,verbatim, enumerate, tikz, tikz-cd, extpfeil, extpfeil}

\begin{document}

\title{Denseness conditions, morphisms \\and equivalences of toposes}
\author{Olivia Caramello}

\bgroup           % fake a titlepage 
\let\footnoterule\relax  % no rule above thanks footnotes 
\date{August 3, 2020}

\maketitle

\begin{abstract}
We systematically investigate morphisms and equivalences of toposes from multiple points of view. We establish a dual adjunction between morphisms and comorphisms of sites, introduce the notion of weak morphism of toposes and characterize the functors which induce such morphisms. In particular, we examine continuous comorphism of sites and show that this class of comorphisms notably includes all fibrations as well as morphisms of fibrations. We also establish a characterization theorem for essential geometric morphisms and locally connected morphisms in terms of continuous functors, and a relative version of the comprehensive factorization of a functor.  

Then we prove a general theorem providing necessary and sufficient explicit conditions for a morphism of sites to induce an equivalence of toposes. This stems from a detailed analysis of arrows in Grothendieck toposes and denseness conditions, which yields results of independent interest. We also derive site characterizations of the property of a geometric morphism to be an inclusion (resp. a surjection, hyperconnected, localic), as well as site-level descriptions of the surjection-inclusion and hyperconnected-localic factorizations of a geometric morphism.
\end{abstract} 
\vspace{5cm}
\egroup

%MACROS-----------------------------------------------------------------------------------------------------------------------

%	European dates ``19 April 1990'' not ``April 19, 1990''
\def\Monthnameof#1{\ifcase#1\or
   January\or February\or March\or April\or May\or June\or
   July\or August\or September\or October\or November\or December\fi}
\def\today{\number\day~\Monthnameof\month~\number\year}

%===========================================================================
%	END OF PROOF BOX
%
%
%  The complexity of the macro necessary to get a little box on the
%  right-hand-side at the end of a proof is amazing.  It really does
%  have to be this long!  Otherwise you're liable to get it at the
%  beginning of the next line, or even on the next page.
%
\def\pushright#1{{%        set up
   \parfillskip=0pt            % so \par doesnt push \square to left
   \widowpenalty=10000         % so we dont break the page before \square
   \displaywidowpenalty=10000  % ditto
   \finalhyphendemerits=0      % TeXbook exercise 14.32
  %
  %                 horizontal
   \leavevmode                 % \nobreak means lines not pages
   \unskip                     % remove previous space or glue
   \nobreak                    % don't break lines
   \hfil                       % ragged right if we spill over
   \penalty50                  % discouragement to do so
   \hskip.2em                  % ensure some space
   \null                       % anchor following \hfill
   \hfill                      % push \square to right
   {#1}                        % the end-of-proof mark (or whatever)
  %
  %                   vertical
   \par}}                      % build paragraph

% prefer proofs with statements, also space after
\def\qed{\pushright{$\square$}\penalty-700 \smallskip}

\newtheorem{theorem}{Theorem}[section]

\newtheorem{proposition}[theorem]{Proposition}

\newtheorem{scholium}[theorem]{Scholium}

\newtheorem{lemma}[theorem]{Lemma}

\newtheorem{corollary}[theorem]{Corollary}

\newtheorem*{theorem*}{Theorem}

\newtheorem{conjecture}[theorem]{Conjecture}

\newenvironment{proofs}%
 {\begin{trivlist}\item[]{\bf Proof }}%
 {\qed\end{trivlist}}

  \newtheorem{rmk}[theorem]{Remark}
\newenvironment{remark}{\begin{rmk}\em}{\end{rmk}}

  \newtheorem{rmks}[theorem]{Remarks}
\newenvironment{remarks}{\begin{rmks}\em}{\end{rmks}}

  \newtheorem{defn}[theorem]{Definition}
\newenvironment{definition}{\begin{defn}\em}{\end{defn}}

  \newtheorem{eg}[theorem]{Example}
\newenvironment{example}{\begin{eg}\em}{\end{eg}}

  \newtheorem{egs}[theorem]{Examples}
\newenvironment{examples}{\begin{egs}\em}{\end{egs}}

%%%%%%%%%%%%%%%%%%%%%%%%%%%%%%%%%%%%%%%%%%%%%%%%%%%%%%%%%%%%%%%%%%%%%%

\makeatletter
\newcommand*{\doublerightarrow}[2]{\mathrel{
		\settowidth{\@tempdima}{$\scriptstyle#1$}
		\settowidth{\@tempdimb}{$\scriptstyle#2$}
		\ifdim\@tempdimb>\@tempdima \@tempdima=\@tempdimb\fi
		\mathop{\vcenter{
				\offinterlineskip\ialign{\hbox to\dimexpr\@tempdima+1em{##}\cr
					\rightarrowfill\cr\noalign{\kern.5ex}
					\rightarrowfill\cr}}}\limits^{\!#1}_{\!#2}}}
\newcommand*{\triplerightarrow}[1]{\mathrel{
		\settowidth{\@tempdima}{$\scriptstyle#1$}
		\mathop{\vcenter{
				\offinterlineskip\ialign{\hbox to\dimexpr\@tempdima+1em{##}\cr
					\rightarrowfill\cr\noalign{\kern.5ex}
					\rightarrowfill\cr\noalign{\kern.5ex}
					\rightarrowfill\cr}}}\limits^{\!#1}}}
\makeatother

\newcommand{\comma}[2]% comma object
 {\mbox{$(#1\!\downarrow\!#2)$}}

\newcommand{\cod}% codomain
 {{\rm cod}}

\newcommand{\comp}% composition
 {\circ}

\newcommand{\Cont}% category of continuous G-sets
 {{\bf Cont}}
 
\newcommand{\ac}
{`}

\newcommand{\dom}% domain
 {{\rm dom}}

\newcommand{\empstg}% empty string
 {[\,]}

\newcommand {\acc}
{`}

\newcommand{\epi}% epimorphism
 {\twoheadrightarrow}

\newcommand{\hy}% hyphen (in math mode)
 {\mbox{-}}

\newcommand{\im}% image
 {{\rm im}}

\newcommand{\imp}% implication
 {\!\Rightarrow\!}

\newcommand{\Ind}[1]% ind-completion of #1
 {{\rm Ind}\hy #1}

\newcommand{\mono}% monomorphism 
 {\rightarrowtail}

\newcommand{\ob}% class of objects
 {{\rm ob}}
 
 \newcommand{\Hom}% class of objects
 {{\rm Hom}}

\newcommand{\op}% opposite category
 {^{\rm op}}

\newcommand{\Set}% category of sets
 {{\bf Set }}

\newcommand{\Sh}% category of sheaves
 {{\bf Sh}}

\newcommand{\sh}% category of sheaves
 {{\bf sh}}

\newcommand{\Sub}% subobject lattice
 {{\rm Sub}}

\tableofcontents

\section{Introduction}\label{intro}

In this work we establish a number of results around the theme of morphisms and equivalences of toposes and their characterizations in terms of sites. Its contents can be summarized as follows.

In section \ref{sec:arrows}, after establishing a basic result allowing to describe the subobjects of the sheafification of a certain object $A$ in terms of subobjects of $A$ which are closed with respect to the associated closure operation, we show that arrows in a topos $\Sh({\cal C}, J)$ between objects of the form $l_{J}^{{\cal C}}(c)$ for $c\in {\cal C}$ (where $l_{J}^{{\cal C}}$ is the functor ${\cal C}\to \Sh({\cal C}, J)$ given by the composite of the Yoneda embedding $y_{\cal C}:{\cal C}\to [{\cal C}^\textup{op}, \Set]$ with the associated sheaf functor $a_{J}:[{\cal C}^\textup{op}, \Set] \to \Sh({\cal C}, J)$) can all be locally represented in terms of arrows coming from the site. Then we investigate, more generally, arrows in $\Sh({\cal C}, J)$ between objects of the form $a_{J}(P)$, and obtain an explicit characterization for them in terms of $J$-functional relations, also describing the operation of composition of such relations which corresponds to composition of the corresponding arrows. Further, we discuss the specialization of this characterization in the case of arrows between objects of the form $l_{J}^{{\cal C}}(c)$ for $c\in {\cal C}$ and its relationship with the alternative description previously obtained. Finally, we introduce a notion of cofinal functor relative to a Grothendieck topology, study its properties and derive a number of applications which will be useful in the subsequent sections, including an explicit characterization of colimits in a Grothendieck topos in terms of generalized elements. 

In section \ref{sec:morphismsandcomorphisms} we investigate the notions of morphism and comorphism of sites. We work in the setting of \emph{small-generated} sites (in the sense of \cite{Shulman} -- these are called \emph{essentially small} sites in \cite{El}) rather than in the usual, restricted context of small sites since all the fundamental results for morphisms of sites, flat functors and comorphisms of sites actually hold at this higher level of generality and their formulations in the extended setting provide a much greater flexibility both theoretically and in connection with applications. In the context of our general study of comorphisms of sites, we focus in particular on fibrations, establishing a number of results about them which will be useful in the following sections; we show in particular that the smallest Grothendieck topology making a fibration towards a category endowed with a Grothendieck topology a comorphism of sites admits an explicit description involving cartesian arrows. Then we investigate the relationship between morphisms and comorphisms of sites, introducing constructions allowing to naturally turn a morphism of sites into a comorphism of sites inducing the same geometric morphism, and conversely. We show that these constructions actually yield a dual adjunction between a category of morphisms under a given site and a category of comorphisms over that site. Lastly, we investigate another fundamental construction, which we call the \emph{fibration of generalized elements} of a functor; this notably allows one to represent a geometric morphism induced by an arbitrary comorphism of sites as a morphism induced by the fibration thus associated with it.      

In section \ref{sec:essentialcomorphisms} we introduce the notion of \emph{weak morphism} of toposes; by this we simply mean a pair of adjoint functors, without the requirement that the inverse image (that is, the left adjoint) preserve finite limits. Then we recall, in the context of small-generated sites $({\cal C}, J)$ and $({\cal D}, K)$, the concept of $(J, K)$-continuous (or simply continuous, when $J$ and $K$ can be unambigously inferred from the context) functor ${\cal C}\to {\cal D}$. This notion generalizes that of morphism of sites $({\cal C}, J)\to ({\cal D}, K)$; in fact, we show that $(J, K)$-continuous functors are precisely the functors ${\cal C}\to {\cal D}$ which induce a weak morphism of toposes $\Sh({\cal D}, K)\to \Sh({\cal C}, J)$. This results from a general equivalence theorem between the category of weak morphisms from a Grothendieck topos ${\cal E}$ to a topos $\Sh({\cal C}, J)$ of sheaves on a small-generated site $({\cal C}, J)$  and the a category of continuous functors ${\cal C}\to {\cal E}$. We obtain a number of equivalent characterizations of continuous functors as well as of continuous comorphisms of sites $({\cal C}, J)\to ({\cal D}, K)$, including an explicit characterization involving relative cofinality conditions. We also derive from the above theorem an equivalence result between essential geometric morphism $\Sh({\cal C}, J)\to {\cal E}$ and $(J, J^{\textup{can}}_{\cal E})$-continuous comorphism of sites $({\cal C}, J)\to ({\cal E}, J^{\textup{can}}_{\cal E})$, and characterize the essential geometric morphisms $\Sh({\cal C}, J)\to \Sh({\cal D}, K)$ that are induced by a $(J, K)$-continuous comorphism of sites $({\cal C}, J)\to ({\cal D}, K)$ in terms of properties of the essential images of such morphisms. In particular, we show that, if the objects coming from the two sites can be characterized by means of suitable generalized compactness conditions, then the categories of such objects will be preserved by the essential images of such morphisms. This is relevant, for instance, in connection with the method for constructing dualities by means of \ac functorializing' topos-theoretic \ac bridges' introduced in \cite{OC10}. In order to describe an arbitrary essential geometric morphism as a morphism induced by a comorphism of sites, we investigate the operation  extending a comorphism of sites to a comorphism between the associated presheaf toposes, endowed with the presheaf \ac liftings' of the given Grothendieck topologies. We also discuss essential morphisms between localic toposes. 

Then we prove a theorem asserting that every fibration gives rise to a continuous comorphism of sites. This result is complemented by the identification of a general context where the continuity of a comorphism of sites can be \ac transported' along a fibration. This analysis leads to a number of corollaries, asserting in particular the continuity of morphisms of fibrations and that of all \ac localizations' of a continuous comorphism of sites. Further, we investigate locally connected morphisms and establish characterizations for them, including a criterion providing necessary and sufficient conditions on a continuous comorphism of sites for the associated geometric morphism to be locally connected. Lastly, we introduce the notion of \emph{terminally connected} geometric morphism. This notion generalizes that of connected and locally connected morphism to the setting of arbitrary essential geometric morphisms. We show that terminally connected morphisms are orthogonal to local homeomorphisms in the $2$-category of toposes, and prove that every essential geometric morphism can be (uniquely up to isomorphism) factored as a terminally connected morphism followed by a local homeomorphism. Then we introduce a relative version of the comprehensive factorization of a functor and show that it corresponds, in the case of a continuous comorphism of sites, to the above-mentioned factorization of the geometric morphism induced by it. In passing, we establish an adjunction between the category of sheaves on a small-generated site and that of functors towards the category underlying the site. 

In section \ref{sec:denseness} we recall the notion of a dense morphism of sites and show that, if the target site is subcanonical, it corresponds precisely to the property of the associated geometric morphism to be an equivalence. In order to generalize this result to the setting of arbitrary sites, we introduce the notion of a weakly dense morphism of sites, giving an explicit characterization of it, and discuss an example of a weakly dense morphism which is not dense. The resulting general theorem providing necessary and sufficient conditions for a morphism of sites to induce an equivalence of toposes represents a vast extension of Grothendieck's Comparison Lemma \cite{grothendieck}. Next we deduce from this result a criterion for a $J$-continuous flat functor ${\cal C}\to {\cal E}$ to induce, via Diaconescu's equivalence, an equivalence of toposes ${\cal E}\simeq \Sh({\cal C}, J)$, and a criterion for the toposes of sheaves on two small-generated sites to be equivalent. Lastly, we introduce some notions of local faithfulness, local fullness and local surjectivity and show that they are naturally related to our denseness conditions. 

In section \ref{sec:propgeometricmorphisms}, by applying results obtained in the previous section, we explicitly characterize the morphisms of sites whose corresponding geometric morphism is a surjection (resp. an inclusion, hyperconnected, localic); this applies in particular to continuous flat functors, giving necessary and sufficient conditions for the geometric morphisms corresponding to them to satisfy such properties. This analysis notably leads to a characterization of the property of a geometric morphism to be an inclusion entirely in terms of its inverse image functor. We then give site-level descriptions of the surjection-inclusion and hyperconnected-localic factorizations, and derive alternative criteria for a morphism of sites (resp. a continuous flat functor) to induce an equivalence of toposes. In this section, we also identify a most general framework for defining induced Grothendieck topologies (which subsumes the classical notion of Grothendieck topology induced on a dense subcategory), and investigate a notion of image of a Grothendieck topology under a functor. 

In section \ref{sec:comorphismsofsites}, we investigate the properties of geometric morphisms to be a surjection (resp. an inclusion, hyperconnected, localic, an equivalence) in the context of morphisms induced by comorphisms of sites, obtaining site-theoretic characterizations for them. Some of these results generalize, in a natural but non-trivial way, a number of known criteria. Lastly, we investigate functors that are both (weak) morphisms and comorphisms of sites, and show in particular that every fully faithful functor which is both a morphism and a comorphism of sites induces a coadjoint retraction of toposes, defined by two morphisms, one local and the other one essential.

\section{Arrows in a Grothendieck topos}\label{sec:arrows}

In this section we make a detailed site-theoretic analysis of arrows in a Grothendieck topos, with the purpose of preparing the ground for the results in the following sections of the paper.

Recall that a \emph{small-generated site} (in the sense of \cite{Shulman}) (also called an \emph{essentially small site} in \cite{El}) is a site $({\cal C}, J)$ such that $\cal C$ is locally small and has a small $J$-dense subcategory. Notice that, for any Grothendieck topos $\cal E$, the site $({\cal E}, J_{\cal E}^{\textup{can}})$, where $J_{\cal E}^{\textup{can}}$ is the canonical topology on $\cal E$, is small-generated.

Even though a Grothendieck topos is defined as a category equivalent to the category of sheaves on a \emph{small} site, it is technically convenient, in the study of morphisms of toposes, to consider, more generally, small-generated sites, in order to be able in particular to consider a topos as a site (namely, its canonical site).

\subsection{The closure operation on subobjects}

For any small-generated site $({\cal C}, J)$, we shall denote by $l_{J}^{{\cal C}}$ (or simply by $l$, when there is no risk of ambiguity) the functor ${\cal C}\to \Sh({\cal C}, J)$ given by the composite of the Yoneda embedding $y_{\cal C}:{\cal C}\to [{\cal C}^\textup{op}, \Set]$ with the associated sheaf functor $a_{J}:[{\cal C}^\textup{op}, \Set] \to \Sh({\cal C}, J)$. We shall say that a family of arrows in $\cal C$ with common codomain is $J$-covering if the sieve generated by it is. The unit of the adjunction between the inclusion functor $i$ of $\Sh({\cal C}, J)$ into $[{\cal C}^\textup{op}, \Set]$ and the associated sheaf functor $a_{J}$ will be denoted by $\eta_{J}$ (or simply by $\eta$).

Notice that we have a closure operation $c_{J}$ on subobjects in the presheaf category $[{\cal C}^\textup{op}, \Set]$, which admits the following explicit description: for any subobject $A\mono E$ in $[{\cal C}^\textup{op}, \Set]$, we have
\[
c_{J}(A)(c)=\{x\in E(c) \mid \{f:d\to c \mid E(f)(x)\in A(d)\}\in J(c)\}
\]   
for any $c\in {\cal C}$. Note also that, for any small-generated site $({\cal C}, J)$, the canonical inclusion functor $\Sh({\cal C}, J)\hookrightarrow [{\cal C}^{\textup{op}}, \Set]$ admits a left adjoint $a_{J}$, which we call, by analogy with the classical case, the associated sheaf functor, and which is (colimit and) finite-limit-preserving. Indeed, if $\cal D$ is a $J$-dense small subcategory of $\cal C$ then, by the Comparison Lemma, we have an equivalence $c:\Sh({\cal D}, J|_{\cal D})\to \Sh({\cal C}, J)$ which is the restriction, along the inclusions $\Sh({\cal D}, J|_{\cal D})\hookrightarrow [{\cal D}^{\textup{op}}, \Set]$ and $\Sh({\cal C}, J)\hookrightarrow [{\cal C}^{\textup{op}}, \Set]$, of the functor $[{\cal D}^{\textup{op}}, \Set]\to [{\cal C}^{\textup{op}}, \Set]$ given by the right Kan extension along the embedding $i$ of ${\cal D}^{\textup{op}}$ into ${\cal C}^{\textup{op}}$, whence the composite functor $c\circ a_{J|_{\cal D}}\circ (-\circ i):[{\cal C}^{\textup{op}}, \Set] \to \Sh({\cal C}, J)$ yields a left adjoint to the inclusion $\Sh({\cal C}, J)\hookrightarrow [{\cal C}^{\textup{op}}, \Set]$.

The closure operation $c_{J}$ can be notably used for obtaining site-theoretic descriptions of many properties and constructions on the topos $\Sh({\cal C}, J)$ which can be expressed in terms of subobjects in it; for instance, we have the following result:

\begin{lemma}\label{lemmalift}
	Let $({\cal C}, J)$ be a small-generated site and $\alpha:F\to G$ an arrow in the presheaf category $[{\cal C}^\textup{op}, \Set]$. Then
	\begin{enumerate}[(i)]
		\item $a_{J}(\alpha)$ is a monomorphism in $\Sh({\cal C}, J)$ if and only if for every $c\in {\cal C}$ and any elements $x, x'\in F(c)$ such that $\alpha(c)(x)=\alpha(c)(x')$, the sieve $\{f:d\to c \mid F(f)(x)=F(f)(x')\}$ is $J$-covering.
		
		\item $a_{J}(\alpha)$ is an epimorphism in $\Sh({\cal C}, J)$ if and only if for every $c\in {\cal C}$ and $x\in G(c)$, the sieve $\{f:d\to c \mid G(f)(x)\in \textup{Im}(\alpha(d)) \}$ is $J$-covering. 
		
		\item $a_{J}(\alpha)$ is an isomorphism in $\Sh({\cal C}, J)$ if and only if both conditions (ii) and (iii) are satisfied.
	\end{enumerate}
\end{lemma}

\begin{proofs}
	(i) An arrow $A\to B$ in a category with pullbacks is a monomorphism if and only if the diagonal monomorphism $A\to A\times_{B} A$ is an isomorphism. We can thus conclude that $a_{J}(\alpha)$ is a monomorphism in $\Sh({\cal C}, J)$ if and only if the diagonal arrow $F\to F\times_{G} F$ is $c_{J}$-dense, that is if and only if for every $c\in {\cal C}$ and any elements $x, x'\in F(c)$ such that $\alpha(c)(x)=\alpha(c)(x')$, the sieve $\{f:d\to c \mid F(f)(x)=F(f)(x')\}$ is $J$-covering.  
	
	(ii) Clearly, $\alpha$ is sent by $a_{J}$ to an epimorphism if and only if its image $\textup{Im}(\alpha)$ is $c_{J}$-dense, that is if and only if, for every $c\in {\cal C}$ and $x\in G(c)$, the sieve $\{f:d\to c \mid G(f)(x)\in \textup{Im}(\alpha(d)) \}$ is $J$-covering. 
	
	(iii) This follows from the fact that every topos is a balanced category.
\end{proofs}

\begin{remarks}
	\begin{enumerate}[(a)]
		\item Let $A_{c, x}:=\{f:d\to c \mid G(f)(x)\in \textup{Im}(\alpha(d)) \}$ (for each $c\in {\cal C}$ and $x\in G(c)$). Then for any arrow $\xi:c\to c'$ and any elements $x\in G(c)$ and $x'\in G(c')$ such that $G(\xi)(x')=x$, $A_{c, x}=\xi^{\ast}(A_{c', x'})$. So, if the codomain of $\alpha$ is a representable $y_{\cal C}(c)$, condition (ii) becomes equivalent to the requirement that the sieve $\{f:d\to c \mid f\in \textup{Im}(\alpha(d))\}$ be $J$-covering.
		
		\item The arrows of $[{\cal C}^{\textup{op}}, \Set]$ whose image under the associated sheaf functor $a_{J}$ is an isomorphism will be called \emph{$J$-bicovering}, in comformity with the terminology \ac\ac \emph{bicouvrante}'' introduced in section II.5 of \cite{grothendieck} for qualifying this type of arrows.  
	\end{enumerate}
	
\end{remarks}

\subsection{Subsheaves and closed subobjects}

The following result is probably well-known, but we did not find a proof of it in the literature, so we give it here.

\begin{proposition}\label{propclosedsubobjects}
	Let $({\cal C}, J)$ be a small-generated site. Then for any object $F$ of $[{\cal C}^{\textup{op}}, \Set]$, denoting by $\textup{ClSub}_{[{\cal C}^{\textup{op}}, \Set]}(F)$ the sub-lattice of $\Sub_{[{\cal C}^{\textup{op}}, \Set]}(F)$ consisting of the $c_J$-closed subobjects, we have a lattice isomorphism
	\[
	\Sub_{\Sh({\cal C}, J)}(a_{J}(F))\simeq \textup{ClSub}_{[{\cal C}^{\textup{op}}, \Set]}(F)
	\]
	which sends a subobject in $\textup{ClSub}_{[{\cal C}^{\textup{op}}, \Set]}(F)$ to its image under $a_{J}$ and a subobject of $a_{J}(F)$ to the pullback of it under the unit arrow $F\to a_{J}(F)$.  
\end{proposition}

\begin{proofs}
	Given a subobject $m:A\to a_{J}(F)$ in $\Sh({\cal C}, J)$, let us associate with it the subobject $n:A'\mono F$ of $F$ defined by the following pullback diagram in $[{\cal C}^{\textup{op}}, \Set]$:
	$$
	\xymatrix{
		A' \ar[d]^{n} \ar[r] & \ar[d]^{m} A  \\
		F \ar[r]^{\eta_{F}\:\:}  & a_{J}(F)
	}
	$$
	Recall that the $c_{J}$-closure $c_{J}(k):A''\mono F$ of a subobject $k:A'\mono F$ in $[{\cal C}^{\textup{op}}, \Set]$ is characterized by the following pullback diagram:
	$$
	\xymatrix{
		A'' \ar[d]^{c_J(k)} \ar[r] & \ar[d]^{a_{J}(k)} a_J(A')  \\
		F \ar[r]^{\eta_{F}\:\:}  & a_{J}(F)
	}
	$$ 
	It thus follows that $n$ is $c_J$-closed, since by applying the pullback-preserving functor $a_J$ to the former pullback square we obtain that $m\cong a_J(n)$.
	Conversely, we associate with any $c_J$-closed subobject $n$ the subobject $a_J(n)$ of $\Sh({\cal C}, J)$. It is immediate to see that these two operations are inverse to each other, and that they are order-preserving; so they are lattice isomorphisms.
\end{proofs}

\begin{remarks}\label{remimageclosedsubobject}
	\begin{enumerate}[(a)]
		\item Given an arrow $\alpha:P\to Q$ in $[{\cal C}^{\textup{op}}, \Set]$, the $c_J$-closed subobject $Q_{\alpha}$ of $Q$ corresponding via Proposition \ref{propclosedsubobjects} to the image of the arrow $a_{J}(\alpha)$ is given by $$Q_{\alpha}(c)=\{x\in Q(c) \mid \{f:d\to c \mid Q(f)(x)\in \textup{Im}(\alpha(d))\}\in J(c) \}$$
		for any $c\in {\cal C}$.
		
		\item Proposition \ref{propclosedsubobjects} clearly generalizes to the setting of elementary toposes, replacing $[{\cal C}^{\textup{op}}, \Set]$ by an arbitrary elementary topos $\cal E$, $J$ by a local operator on it and $c_{J}$ by the associate closure operation on subobjects in $\cal E$.
	\end{enumerate}
	
\end{remarks}

\subsection{Site-theoretic description of arrows between objects coming from the site}\label{sec:arrowscomingfromthesite}

Given a site $({\cal C}, J)$, for two arrows $h,k:c\to d$ in $\cal C$ we shall write $h\equiv_{J} k$ for \emph{$J$-local equality}, that is, to mean that there exists a $J$-covering sieve $S$ on $c$ such that $h\circ f=k\circ f$ for every $f\in S$. Notice that $l(h)=l(k)$ if and only if $h\equiv_{J} k$.

\begin{proposition}\label{propexplicit}
	Let $({\cal C}, J)$ be a small-generated site. 
	
	\begin{enumerate}[(i)]
		\item Then for any arrow $\xi:l(c)\to l(d)$ in $\Sh({\cal C}, J)$ there exists a family of arrows $\{f_{u}:c_{u}\to c, g_{u}:c_{u}\to d \mid u\in U\}$ such that $\{f_{u}:c_{u}\to c \mid u\in U\}$ generates a $J$-covering sieve, for any object $e$ and arrows $h:e\to c_{u}$ and $k:e\to c_{u'}$ such that $f_{u}\circ h=f_{u'}\circ k$ we have $g_{u}\circ h\equiv_{J} g_{u'}\circ k$, and $\xi \circ l(f_{u})=l(g_{u})$ for every $u\in U$.  
		
		\item Conversely, any family of arrows ${\cal F}:\{f_{u}:c_{u}\to c, g_{u}:c_{u}\to d \mid u\in U\}$ such that $\{f_{u}:c_{u}\to c \mid u\in U\}$ generates a $J$-covering sieve and for any object $e$ and arrows $h:e\to c_{u}$ and $k:e\to c_{u'}$ such that $f_{u}\circ h=f_{u'}\circ k$ we have $g_{u}\circ h\equiv_{J} g_{u'}\circ k$, determines a unique arrow $\xi_{\cal F}:l(c)\to l(d)$ in $\Sh({\cal C}, J)$ such that $\xi_{\cal F}\circ l(f_{u})=l(g_{u})$ for every $u\in U$. 
		
		\item Two families ${\cal F}=\{f_{u}:c_{u}\to c, g_{u}:c_{u}\to d \mid u\in U\}$ and ${\cal F}'=\{f'_{v}:e_{v}\to c, g'_{v}:e_{v}\to d \mid v\in V\}$ as in (ii) determine the same arrow $l(c)\to l(d)$ (i.e. $\xi_{\cal F}=\xi_{{\cal F}'}$) if and only if there exist a $J$-covering family $\{a_{k}:b_{k}\to c \mid k\in K\}$ and factorizations of it through both of them by arrows $x_{k}:b_{k}\to c_{u(k)}$ and $y_{k}:b_{k}\to e_{v(k)}$ (i.e. $f_{u(k)}\circ x_{k}=a_{k}=f'_{v(k)}\circ y_{k}$ for every $k\in K$) such that $g_{u(k)}\circ x_{k}\equiv_{J} g'_{v(k)}\circ y_{k}$ for every $k\in K$.
		
		\item Given two families ${\cal F}=\{f_{u}:c_{u}\to c, g_{u}:c_{u}\to d \mid u\in U\}$ and ${\cal G}=\{h_{v}:d_{v}\to d, k_{v}:d_{v}\to e \mid v\in V\}$, the composite arrow $\xi_{\cal G}\circ \xi_{\cal F}:l(c)\to l(e)$ is induced as in (ii) by the family $\{f_{u}\circ x:\textup{dom}(x)\to c, k_{v}\circ y:\textup{dom}(y)\to e \mid (u,v,x,y )\in Z\}$, where $Z=\{(u, v, x, y) \mid  u\in U, v\in V, \textup{dom}(x)=\textup{dom}(y), \textup{cod}(x)=c_{u}, \textup{cod}(y)=d_{v}, h_{v}\circ y=g_{u}\circ x\}$. 
	\end{enumerate} 
\end{proposition}

\begin{proofs}
	(i) Consider the following pullback square in $[{\cal C}^\textup{op}, \Set]$:
	$$
	\xymatrix{
		R \ar[d]^{r} \ar[rr]^t & & y_{\cal C}(d) \ar[d]^{\eta_{y_{\cal C}(d)}}  \\
		y_{\cal C}(c) \ar[rr]^{\xi\circ \eta_{y_{\cal C}(c)}} & & l(d)
	}
	$$ 
	and the sieve $S=\{f:e\to c  \mid \exists g:e\to d \textup{ with } (f, g)\in R\}$ on $c$ given by the image in $[{\cal C}^\textup{op}, \Set]$ of the arrow $r:R\to y_{\cal C}(c)$. This sieve is $J$-covering since it is sent by $a_{J}$ to an isomorphism. Let us show that for any $(f,g)\in R$, $\xi\circ l(f)=l(g)$. Consider the adjunction between $i$ and $a_{J}$: the arrow $\eta_{y_{\cal C}(d)} \circ y_{\cal C}(g)$ is the transpose of $l(g)$ and the arrow $\xi\circ \eta_{y_{\cal C}(c)} \circ y_{\cal C}(f)$ is the transpose of $\xi \circ l(f)$, so it is equivalent to verify that $\eta_{y_{\cal C}(d)} \circ y_{\cal C}(g)=\xi\circ \eta_{y_{\cal C}(c)} \circ y_{\cal C}(f)$. But this follows by the commutativity of the above square, since $(f, g)\in R$. To check that the family of arrows $\{(f, g) \mid (f,g)\in R\}$ satisfies the required condition, we have to check that for any $(f, g), (f',g')\in R$ and arrows $h, h'$ such that $f\circ h=f'\circ h'$, $g\circ h\equiv_{J} g'\circ h'$, equivalently $l(g\circ h)=l(g'\circ h')$. Since $(f\circ h, g\circ h)\in R$ as $R$ is functorial, we have $\xi\circ l(f\circ h)=l(g\circ h)$; similarly, we have $\xi \circ l(f'\circ h')=l(g'\circ h')$. So $l(g\circ h)=l(g'\circ h')$, as required. 
	
	(ii) Any sieve $S$ on an object $c$ of $\cal C$, regarded as a presheaf on ${\cal C}$, is the colimit of the canonical diagram defined on the full subcategory $\int S$ of ${\cal C}\slash c$ on the objects which are elements of $S$. If $S$ is $J$-covering then the canonical momonomorphism $S\mono y_{\cal C}(c)$ is sent by the associated sheaf functor $a_{J}$ to an isomorphism, so, since $a_{J}$ is colimit-preserving, $l(c)$ is the colimit of the functor $a_{S}:\int S \to \Sh({\cal C}, J)$ sending any object $f:d\to c$ of $\int S$ to $l(d)$ and any arrow $h:(f:d\to c)\to (f':d'\to c)$ in $\int S$ to $l(h):l(d')\to l(d)$. Therefore, giving an arrow $l(c)\to l(d)$ in $\Sh({\cal C}, J)$ amounts to giving a cocone on $a_{S}$ with vertex $l(d)$. Any family of arrows $\{f_{u}:c_{u}\to c, g_{u}:c_{u}\to d \mid u\in U\}$ such that $\{f_{u}:c_{u}\to c \mid u\in U\}$ generates a $J$-covering sieve and for any object $e$ and arrows $h:e\to c_{u}$ and $k:e\to c_{u'}$ such that $f_{u}\circ h=f_{u'}\circ k$ we have $g_{u}\circ h\equiv_{J} g_{u'}\circ k$, defines such a cocone. Indeed, by taking $S$ equal to the sieve generated by the family $\{f_{u}:c_{u}\to c\}$ we obtain, thanks to the property that for any object $e$ and arrows $h:e\to c_{u}$ and $k:e\to c_{u'}$ such that $f_{u}\circ h=f_{u'}\circ k$ we have $g_{u}\circ h\equiv_{J} g_{u'}\circ k$, a well-defined family of arrows $l(\textup{dom}(f))\to l(d)$ (for $f\in S$) satisfying the commutativity conditions for a cocone. It is clear that the arrow $l(c)\to l(d)$ thus defined satisfies the required property.     
	
	(iii) Let us first prove the `if' direction. Since the sieve generated by the family  $\{a_{k}:b_{k}\to c \mid k\in K\}$ is $J$-covering, the arrows $l(a_{k})$ to $l(c)$ are jointly epimorphic and hence $\xi_{\cal F}=\xi_{{\cal F}'}$ if and only if $\xi_{\cal F}\circ l(a_{k})=\xi_{{\cal F}'}\circ l(a_{k})$ for every $k\in K$. But $\xi_{\cal F}\circ l(a_{k})=\xi_{\cal F} \circ l(f_{u(k)})\circ l(x_{k})= l(g_{u(k)})\circ l(x_{k})=l(g_{u(k)}\circ x_{k})=l(g'_{v(k)}\circ y_{k})=l(g'_{v(k)})\circ l(y_{k})=\xi_{{\cal F}'}\circ l(f'_{v(k)})\circ l(y_{k})=\xi_{{\cal F}'}\circ l(a_{k})$, where the equality $l(g_{u(k)}\circ x_{k})=l(g'_{v(k)}\circ y_{k})$ follows from the fact that $g_{u(k)}\circ x_{k}\equiv_{J} g'_{v(k)}\circ y_{k}$. It thus remains to prove the `only if' direction. Suppose that $\xi_{\cal F}=\xi_{{\cal F}'}$. Let us define the family $\{a_{k}:b_{k}\to c \mid k\in K\}$ to be the intersection of the sieve generated by the family $\{f_{u}:c_{u}\to c \mid u\in U\}$ with that generated by the family $\{f'_{v}:e_{v}\to c \mid v\in V\}$; since these two sieves are $J$-covering by our hypotheses, their intersection is $J$-covering as well. Let us choose factorizations $x_{k}:b_{k}\to c_{u(k)}$ and $y_{k}:b_{k}\to e_{v(k)}$ of it through the arrows of these families, so that $f_{u(k)}\circ x_{k}=a_{k}=f'_{v(k)}\circ y_{k}$ for every $k\in K$. We want to prove that, for every $k\in K$, $g_{u(k)}\circ x_{k}\equiv_{J} g'_{v(k)}\circ y_{k}$, equivalently that $l(g_{u(k)}\circ x_{k})=l(g'_{v(k)}\circ y_{k})$. But $l(g_{u(k)}\circ x_{k})=l(g_{u(k)})\circ l(x_{k})=\xi_{\cal F}\circ l(f_{u(k)})\circ l(x_{k})=\xi_{{\cal F}'}\circ l(f'_{v(k)}\circ y_{k})=l(g'_{v(k)}\circ y_{k})$, as required.
	
	(iv) First of all, we have to show that the family of arrows $\{f_{u}\circ x \mid (u,v,x,y )\in Z\}$ is $J$-covering. To this end, consider the pullback
	$$
	\xymatrix{
		P \ar[d] \ar[r] & H \ar[d]  \\
		G \ar[r]  & y_{\cal C}(d)
	} 
	$$
	where $H$ is the sieve on $d$ generated by the arrows $\{h_{v}:d_{v}\to d \mid v\in I\}$, $G$ is the sieve on $d$ generated by the arrows $\{g_{u}:c_{u}\to d \mid u\in U\}$ and the arrows are the canonical inclusion monomorphisms.
	
	Since $H$ is $J$-covering, the canonical monomorphism $H\to y_{\cal C}(d)$ is $c_{J}$-dense and hence the monomorphism $P\mono G$ in the above square is $c_{J}$-dense too; this means that for any $u\in U$, the sieve $\{\xi:\textup{dom}(\xi)\to c_{u} \mid g_{u}\circ \xi\in P\}$ is $J$-covering. But $P=H\cap G$, so this means that for any $u\in U$, the sieve $\{x:\textup{dom}(x)\to c_{u} \mid \exists v\,\exists y\,(u,v,x,y)\in Z\}$ is $J$-covering. Since this sieve is contained in the pullback along $f_{u}$ of the sieve $\{f_{u}\circ x \mid (u,v,x,y )\in Z\}$, it follows from the transitivity axiom for Grothendieck topologies that this latter sieve is $J$-covering, as desired. 
	
	To conclude the proof of (iv), we have to verify that for every $(u,v,x,y)\in Z$, $\xi_{\cal G}\circ \xi_{\cal F}\circ l(f_{u})\circ l(x)=l(k_{v})\circ l(y)$. Recalling the definitions of $\xi_{\cal F}$ and $\xi_{\cal G}$, we obtain that $\xi_{\cal G}\circ \xi_{\cal F}\circ l(f_{u})\circ l(x)=\xi_{\cal G}\circ l(g_{u})\circ l(x)=\xi_{\cal G}\circ l(h_{v})\circ l(y)=l(k_{v})\circ l(y)=l(k_{v}\circ y)$, as required.   
\end{proofs}

\begin{remarks}\label{rempropexplicit}
	\begin{enumerate}[(a)]
		\item Proposition \ref{propexplicit} could alternatively have been derived from the description of the topos $\Sh({\cal C}, J)$ as a completion of the site $({\cal C}, J)$ provided by Theorem 8.22 \cite{Shulman}.
		
		\item If we assume the axiom of choice, given a family of arrows $\{f_{u}:c_{u}\to c, g_{u}:c_{u}\to d \mid u\in U\}$ such that $\{f_{u}:c_{u}\to c \mid u\in U\}$ generates a $J$-covering sieve and for any object $e$ and arrows $h:e\to c_{u}$ and $k:e\to c_{u'}$ such that $f_{u}\circ h=f_{u'}\circ k$ we have $g_{u}\circ h=g_{u'}\circ k$, as in points (i) and (ii) of Proposition \ref{propexplicit}, we can clearly suppose without loss of generality $U$ to be the $J$-covering sieve on $c$ generated by the family $\{f_{u}:c_{u}\to c \mid u\in U\}$. So such a family can be equivalently identified with a family of arrows  $g_{f}:\textup{dom}(f)\to d$ indexed by the arrows $f$ of a $J$-covering sieve $S$ on $c$ with the property that for any arrow $h$ composable with $f$, $g_{f\circ h}\equiv_{J} g_{f}\circ h$. 
		
		\item Under the identification of Remark \ref{rempropexplicit}(b), the family $\{a_{k}:b_{k}\to c \mid k\in K\}$ in point (iii) of Proposition \ref{propexplicit} can be assumed to be the intersection of the sieves $U$ and $V$.
		
	\end{enumerate} 
\end{remarks}

\begin{proposition}
	Let $({\cal C}, J)$ be a small-generated site and $h:c\to d, k:c\to e$ arrows in $\cal C$. Then there exists an arrow $\xi:l(e)\to l(d)$ such that $\xi \circ l(h)=l(k)$ if and only if there exists a family of arrows $g_{f}:\textup{dom}(f)\to d$ indexed by the arrows $f$ of a $J$-covering sieve $S$ on $e$ such that $g_{f\circ z}=g_{f}\circ z$ whenever $z$ is composable with $f$, and such that for any $\chi\in h^{\ast}(S)$, $k\circ \chi\equiv_{J} g_{h\circ \chi}$.
\end{proposition}

\begin{proofs}
	By Proposition \ref{propexplicit}, giving an arrow $\xi:l(e)\to l(d)$ amounts to specifying a family of arrows $g_{f}:\textup{dom}(f)\to d$ indexed by the arrows $f$ of a $J$-covering sieve $S$ on $e$ such that $g_{f\circ z}=g_{f}\circ z$ whenever $z$ is composable with $f$, which satisfies $\xi\circ l(f)=l(g_{f})$ for each $f\in S$. Now, $\xi \circ l(h)=l(k)$ if and only if $\xi \circ l(h)\circ l(\chi)=l(k)\circ l(\chi)$ for every $\chi\in h^{\ast}(S)$ (since, $h^{\ast}(S)$ being $J$-covering as $S$ is, the arrows $l(\chi)$ are jointly epimorphic). But $\xi \circ l(h)\circ l(\chi)=\xi \circ l(h\circ \chi)=l(g_{h\circ \chi})$, while $l(k)\circ l(\chi)=l(k\circ \chi)$, so the condition $\xi \circ l(h)=l(k)$ is equivalent to $k\circ \chi\equiv_{J} g_{h\circ \chi}$, as required.   	
\end{proofs}

\begin{proposition}\label{proplifttosite}
	Let $({\cal C}, J)$ be a small-generated site and $\xi:A\to l(c)$ an arrow of $\Sh({\cal C}, J)$. Then the family of arrows $\chi:l(d)\to A$, where $d\in {\cal C}$, such that $\xi \circ \chi=l(f)$ for some arrow $f:d\to c$ in $\cal C$ is epimorphic on $A$.	
\end{proposition}

\begin{proofs}
	The family $R$ of arrows to $A$ whose domain is of the form $l(e)$ for some $e\in {\cal C}$ is clearly epimorphic on $A$. Let us denote by $T$ the family of arrows $\chi:l(d)\to A$, where $d\in {\cal C}$, such that $\xi \circ \chi=l(f)$ for some arrow $f:d\to c$ in $\cal C$. To show that $T$ is epimorphic, it clearly suffices, by transitivity, to show that for any  $\alpha\in R$, $\alpha^{\ast}(T)$ is epimorphic. But this follows from Proposition \ref{propexplicit}(i) applied to the arrow $\xi \circ \alpha$, since that implies that there exists a $J$-covering family $\{f_{i} \mid i\in I\}$ of arrows on $\textup{dom}(\alpha)$ and for each $i\in I$ an arrow $g_{i}:\textup{dom}(f_{i})\to c$ such that $\xi\circ \alpha \circ l(f_{i})=l(g_{i})$; indeed, the family of arrows $\{l(f_{i}) \mid i\in I\}$ is epimorphic and contained in $\alpha^{\ast}(T)$. 
\end{proofs}

\begin{remark}\label{remrelativetopology}
	If $A=a_{J}(P)$ for a presheaf $P$ on $\cal C$ then the collection of objects $(d, x)$ of of the category $\int P$ of elements of $P$ such that, denoting by $\lambda_{(d, x)}:l(d)\to a_{J}(P)=A$ the canonical colimit arrow, $\xi\circ \lambda_{(d, x)}$ is of the form $l(f)$ for some $f:d\to c$ defines a full $J_{P}$-dense subcategory of $\int P$ (where $J_{P}$ is the Grothendieck topology on $\int P$ whose covering sieves are precisely those sent by the canonical projection functor ${\int P}\to {\cal C}$ to $J$-covering sieves).  
\end{remark}

\subsection{Locally functional relations}

Let us now consider the problem of explicitly describing arrows in a sheaf topos $\Sh({\cal C}, J)$. Given two objects $P$ and $Q$ of $[{\cal C}^{\textup{op}}, \Set]$, we shall describe the arrows $a_{J}(P)\to a_{J}(Q)$ in $\Sh({\cal C}, J)$ by explicitly characterizing their graphs. Notice that a relation $r:R\mono A\times B$ in a topos is the graph of an arrow $A\to B$ if and only if the composite of the canonical projection $A\times B \to A$ with $r$ is an isomorphism. By Proposition \ref{propclosedsubobjects}, the relations $r:R\mono a_{J}(P)\times a_{J}(Q)\cong a_{J}(P\times Q)$ in $\Sh({\cal C}, J)$ can be identified with the $c_{J}$-closed subobjects of $P\times Q$, via the correspondence sending $r$ to its pullback $r':R'\mono P\times Q$ along the arrow $\eta_{P\times Q}:P\times Q \to a_{J}(P\times Q)$; so $R\cong a_{J}(R')$. Notice that a subobject $R'\mono P\times Q$ is $c_{J}$-closed if and only if for any $c\in {\cal C}$ and any $(x, y)\in P(c)\times Q(c)$, if the sieve $\{f:d\to c \mid\ (P(f)(x), Q(f)(x))\in R(d)\}$ is $J$-covering then $(x,y)\in R(c)$. The condition for $R$ to be the graph of an arrow $a_{J}(P)\to a_{J}(Q)$ in $\Sh({\cal C}, J)$ can thus be reformulated as the requirement that the arrow $\pi_{P}\circ r'$ should be sent by $a_{J}$ to an isomorphism. In light of Lemma \ref{lemmalift}, this motivates the following definition.

\begin{definition}\label{defJfunctionalrelation}
	In a presheaf topos $[{\cal C}^{\textup{op}}, \Set]$, a relation $R\mono P\times Q$ (that is, an assignment $c\to R(c)$ to each object $c$ of $\cal C$ of a subset $R(c)$ of $P(c)\times Q(c)$ which is \emph{functorial} in the sense that for any arrow $f:c\to c'$ in $\cal C$, $P(f)\times Q(f)$ sends $R(c')$ to $R(c)$), is said to be \emph{$J$-functional from $P$ to $Q$} if it satisfies the following properties:
	\begin{enumerate}[(i)]
		\item for any $c\in {\cal C}$ and any $(x, y)\in P(c)\times Q(c)$,\\ if $\{f:d\to c \mid\ (P(f)(x), Q(f)(y))\in R(d)\}\in J(c)$ then $(x,y)\in R(c)$;
		
		\item for any $c\in {\cal C}$ and any $(x, y), (x', y')\in R(c)$, if $x=x'$ then $\{f:d\to c \mid Q(f)(y)=Q(f)(y')\}\in J(c)$;
		
		\item for any $c\in {\cal C}$ and any $x\in P(c)$, $\{f:d\to c \mid \exists y\in Q(d) \,(P(f)(x), y)\in R(d) \}\in J(c)$.
	\end{enumerate}
\end{definition}

\begin{remarks}\label{remfunctionalrelations}
	\begin{enumerate}[(a)]
		\item By the above discussion, the $J$-functional relations $R\mono P\times Q$ from $P$ to $Q$ are precisely the $c_{J}$-closed relations $R$ on $[{\cal C}^{\textup{op}}, \Set]$ whose image under $a_{J}$ is the graph of an arrow $a_{J}(P)\to a_{J}(Q)$ in $\Sh({\cal C}, J)$, that is the subobjects corresponding to the functional relations from $a_{J}(P)$ to $a_{J}(Q)$ under the bijection of Proposition \ref{propclosedsubobjects}.
		
		\item Any relation $R\mono P\times Q$ on $[{\cal C}^{\textup{op}}, \Set]$ satisfying conditions (ii) and (iii) of Definition \ref{defJfunctionalrelation} admits a \emph{$J$-closure} $\overline{R}^{J}$, consisting of all the pairs $(x, y)\in P(c)\times Q(c)$ such that $\{f:d\to c \mid\ (P(f)(x), Q(f)(y))\in R(d)\}\in J(c)$, which satisfies all the conditions of the definition. 
		
		\item The $J$-functional relation $R_{\xi}$ corresponding to an arrow $\xi:a_{J}(P)\to a_{J}(Q)$ can be concretely described as follows: for any $(x,y)\in P(c)\times Q(c)$, $(x, y)\in R_{\xi}$ if and only if $(\xi\circ \eta_{P})(c)(x)=(\eta_{Q}(c))(y)$, where $\eta_{P}:P\to a_{J}(P)$ and $\eta_{Q}:Q\to a_{J}(Q)$ are the unit arrows.			
	\end{enumerate}
\end{remarks} 

The following result shows that $J$-functional relations from $P$ to $Q$ can be identified with certain ways of assigning (possibly empty) $\equiv_{J}$-equivalence classes of elements of $Q$ to elements of $P$.

\begin{proposition}\label{propfuncionalassignment}
	We have a natural bijection between the $J$-functional relations $R$ from $P$ to $Q$ and the functions $f$ assigning to each element $x\in P(c)$ a $\equiv_{J}$-equivalence class of elements of $Q(c)$ (that is, a subset $A$ of $Q(c)$ such that for any $y\in A$ and $y'\in Q(c)$, $y\equiv_{J} y'$ if and only if $y'\in A$) which is \emph{functorial} in the sense that for any arrow $g:d\to c$ in $\cal C$, if $y\in f(x)$ then $Q(g)(y)\in f(P(g)(x))$, \emph{$J$-closed} in the sense that if $Q(g)(y)\in f(P(g)(x))$ for all the arrows $g$ of a $J$-covering sieve then $y\in f(x)$ and \emph{$J$-locally pointed} in the sense that for any $x\in P(c)$, $\{g:d\to c \mid f(P(g)(x))\neq \emptyset \}\in J(c)$. This correspondence sends 
	\begin{itemize}
		\item a function $f$ satisfying the above properties to the relation $R_{f}\mono P\times Q$ consisting of the pairs $(x, y)$ such that $y\in f(x)$;
		
		\item a $J$-functional relation $R$ from $P$ to $Q$ to the function $f_{R}$ given by:
		\[
		f_{R}(x)=\{y\in Q(c) \mid (x, y)\in R(c)\}
		\]
		for any $x\in P(c)$.
	\end{itemize}
\end{proposition}

\begin{proofs}
	First, let us show that the correspondence is well-defined. Given $f$, let us verify that $R_{f}$ is a $J$-functional relation. The functoriality of $R_{f}$ follows from that of $f$. The relation $R_{f}$ satisfies condition (i) by the $J$-closedness of $f$, condition (ii) since $f$ takes values in $\equiv_{J}$-equivalence classes and condition (iii) since $f$ is $J$-locally pointed. Conversely, let us show that, given a $J$-functional relation $R$ from $P$ to $Q$, the function $f_{R}$ takes values in $\equiv_{J}$-equivalence classes, is functorial, $J$-closed and $J$-locally pointed. Given $x\in P(c)$ and $y, y'\in Q(c)$ such that $y\in f_{R}(x)$, let us prove that $y\equiv_{J} y$ if and only if $y'\in f_{R}(x)$. The `if' direction follows at once from the fact that $R$ satisfies condition (ii), so it remains to prove the `only if' one. If $y\equiv_{J} y$ then there is a $J$-covering sieve $T$ on $c$ such that for any $t\in T$, $Q(t)(y)=Q(t)(y')$. Now, since $(x, y)\in R(c)$, by the functoriality of $R$ we have $(P(t)(x), Q(t)(y))\in R(\textup{dom}(t))$ for each $t\in T$; condition (i) thus implies that $(x, y')\in R(c)$, as required. The fact that $f_{R}$ is functorial (resp. $J$-closed,  $J$-locally pointed) follows from the fact that $R$ is functorial (resp. satisfies condition (i), condition (iii)). Now that we have proved that the assignments $R \to f_{R}$ and $f \to R_{f}$ are well-defined, it remains to show that they are inverse to each other. The equality $f=f_{R_{f}}$ follows from the fact that for any $x$, $f_{R_{f}}=\{y \mid (x, y)\in R_{f}\}=\{y \mid y\in f(x)\}=f(x)$, while the equality $R=R_{f_{R}}$ follows from the fact that $R_{f_{R}}=\{(x, y) \mid y\in f_{R}(x)\}=R$. 
\end{proofs}

Now that we have seen that we can naturally represent arrows $a_{J}(P)\to a_{J}(Q)$ in $\Sh({\cal C}, J)$ in terms of $J$-functional relations from $P$ to $Q$, it is natural to consider how composition of arrows in $\Sh({\cal C}, J)$ can be described in terms of an operation of composition of such relations in $[{\cal C}^{\textup{op}}, \Set]$. The $J$-functional relation corresponding to the composite of two arrows in $\Sh({\cal C}, J)$ induced by $J$-functional relations $R$ and $S$ is the $c_{J}$-closure of the relation given by the (relational) composite of $R$ and $S$ in $[{\cal C}^{\textup{op}}, \Set]$. Indeed, the operation of composition of relations is clearly preserved by geometric functors; but $a_{J}$ is such a functor and, by Proposition \ref{propclosedsubobjects}, there is just one $c_{J}$-closed subobject whose image under $a_{J}$ is a given subobject. Notice that the composition of relations in a presheaf topos is computed pointwise.  

Summarizing, we have the following result:

\begin{theorem}\label{thmfunctionalrelations}
	Let $({\cal C}, J)$ be a small-generated site. Then, for any presheaves $P,Q\in [{\cal C}^{\textup{op}},\Set]$, the arrows $a_{J}(P)\to a_{J}(Q)$ in $\Sh({\cal C}, J)$ are in natural bijection with the $J$-functional relations from $P$ to $Q$ in $[{\cal C}^{\textup{op}},\Set]$.
	
	Moreover, under this bijection, the composition of arrows in $\Sh({\cal C}, J)$ corresponds to the $c_{J}$-closure of the composition in $[{\cal C}^{\textup{op}},\Set]$ of the associated $J$-functional relations. This operation, which we shall denote by the symbol $\ast$, admits the following explicit description: given a $J$-functional relation $R$ from $P$ to $Q$ and a $J$-functional relation $S$ from $Q$ to $Z$, $S\ast R$ is given by the formula
	\begin{equation*}\begin{split}
	(S\ast R)(c)&=\{(x, z)\in P(c)\times Z(c) \mid \{f:d\to c \mid \exists y\in Q(d)\, \\
	& \quad\quad(P(f)(x), y)\in R(d)\textup{ and } (y, Z(f)(z))\in S(d)\}\in J(c) \}
	\end{split}\end{equation*}
	for any $c\in {\cal C}$.
\end{theorem}

Theorem \ref{thmfunctionalrelations} can be notably applied to the description of the full subcategory $a_{J}({\cal C})$ of $\Sh({\cal C}, J)$ on the objects of the form $l(c)$ for $c\in {\cal C}$.

\begin{corollary}\label{corfullsubcategorytopos}
	Let $({\cal C}, J)$ be a small-generated site. Then the full subcategory $a_{J}({\cal C})$ of $\Sh({\cal C}, J)$ on the objects of the form $l(c)$ for $c\in {\cal C}$ is equivalent to the category ${\cal C}_{J}$ defined as follows: 
	\begin{itemize}
		\item the objects of ${\cal C}_{J}$ are the objects of $\cal C$; 
		
		\item the arrows $c\to d$ are the collections $R$ of pairs of arrows $(f:e\to c,g:e\to d)$ where $e$ varies among the objects of $\cal C$ (we write $(f, g)\in R(e)$ to mean that $(f, g)\in R$ and $\textup{dom}(f)=\textup{dom}(g)=e$) satisfying the following properties:
		
		\begin{enumerate}[(i)]
			\item for any arrow $k:e'\to e$, if $(f,g)\in R$ then $(f\circ k, g\circ k)\in R$; 
			
			\item for any $e\in {\cal C}$ and any arrows $x:e\to c, y:e\to d$, if $\{f:e'\to e \mid\ (x\circ f, y\circ f)\in R\}\in J(e)$ then $(x, y)\in R$; 
			
			\item for any $e\in {\cal C}$ and any $(x, y), (x', y')\in R(e)$, if $x=x'$ then $\{f:e'\to e \mid y\circ f=y'\circ f\}\in J(e)$;
			
			\item for any $e\in {\cal C}$ and any $x:e\to c$, $\{f:e'\to e \mid \exists y:e'\to d \text{ such that } (x\circ f, y)\in R(e') \}\in J(e)$;
		\end{enumerate}
		
		\item the composite $S\ast R$ of two arrows $R:c\to d$ and $S:d\to a$ is given by the following formula: 
		\begin{equation*}\begin{split}
		(S\ast R)(e)&=\{(x:e\to c, z:e\to a) \mid \{f:e'\to e \mid \exists y:e'\to d\, \\
		& \quad\quad(x\circ f, y)\in R(e')\textup{ and } (y, z\circ f)\in S(e')\}\in J(e) \}
		\end{split}\end{equation*}
		for any $e\in {\cal C}$. 
	\end{itemize}
	
	One half of the equivalence ${\cal C}_{J}\simeq a_{J}({\cal C})$ is given by the functor ${\cal C}_{J}\to a_{J}({\cal C})$ acting on objects by sending any $c\in {\cal C}_{J}$ to $l(c)$ and any arrow $R:c\to c'$ in ${\cal C}_{J}$ to the arrow $l(c)\to l(c')$ in $a_{J}({\cal C})$ whose graph is $a_{J}(R)$.   
	
\end{corollary}\qed

\begin{remark}\label{remrelationstwoapproaches}
	Let us clarify the relationship between the description of arrows in $a_{J}({\cal C})$ given by Proposition \ref{propexplicit} with that provided by Corollary 	\ref{corfullsubcategorytopos} (modulo the equivalence $a_{J}({\cal C})\simeq {\cal C}_{J}$). We stress that, unlike the former description (as equivalence classes of ways of representing them), the latter is \emph{canonical}; in particular, the composition operation, not being defined by using representatives, does not require any choices. Indeed, $J$-functional relations \emph{are} equivalence classes themselves.
	
	Given an arrow $R:c\to d$ in ${\cal C}_{J}$, by property (iv) we have (by taking $e=c$ and $x$ equal to the identity arrow on $c$) that $S:=\{f:e\to c \mid \exists y:e\to d \text{ such that } (f, y)\in R \}\in J(c)$. By using the axiom of choice, we can therefore choose, for each $f\in S$, an arrow $y_{f}:e\to d$ such that $(f, y_{f})\in R$, and hence the family $\{f:e\to c, y_{f}:e\to d \mid f\in S\}$ presents the arrow $l(c)\to l(d)$ corresponding to $R$ in the sense of Proposition \ref{propexplicit}. In the converse direction, we claim that, given a family ${\cal F}=\{f:\textup{dom}(f)\to c, g_{f}:\textup{dom}(f) \to d \mid f\in S\}$ inducing an arrow $\xi:l(c)\to l(d)$ in the sense of Proposition \ref{propexplicit} (we present $\cal F$ in the form of Remark \ref{rempropexplicit}(b) for simplicity), the $J$-functional relation $R$ corresponding to the arrow $\xi$ admits the following explicit description:
	$$
	(h,k)\in R \quad\text{if and only if}\quad \text{ for all } \chi\in h^{\ast}(S), k\circ \chi\equiv_{J} g_{h\circ \chi}
	$$
	By definition of $R$ as the pullback of the graph of the arrow $\xi:l(c)\to l(d)$ along the arrow $\eta_{y_{\cal C}(c)}\times \eta_{y_{\cal C}(d)}$, $R$ is characterized by the following universal property: for any $e\in {\cal C}$, $(h,k)\in R(e)$ if and only if $\xi\circ l(h)=l(k)$. Now, since the sieve $S$ is $J$-covering, for any $h$, the sieve $h^{\ast}(S)$ is $J$-covering as well and hence the family of arrows $\{l(\chi) \mid \chi\in h^{\ast}(S)\}$ is epimorphic on $l(e)$. So $\xi\circ l(h)=l(k)$ if and only if $\xi\circ l(h)\circ l(\chi)=l(k)\circ l(\chi)$ for every $\chi\in h^{\ast}(S)$. But $\xi\circ l(h)\circ l(\chi)=\xi\circ l(h\circ \chi)=l(g_{h\circ \chi})$, while $l(k)\circ l(\chi)=l(k\circ \chi)$, so $\xi\circ l(h)\circ l(\chi)=l(k)\circ l(\chi)$ if and only if $l(g_{h\circ \chi})=l(k\circ \chi)$, i.e. if and only if $k\circ \chi\equiv_{J} g_{h\circ \chi}$, as required.   
\end{remark}

The following corollary of Theorem \ref{thmfunctionalrelations} will be instrumental in section \ref{sec:hyperconnfact} for obtaining a site-level description of the hyperconnected-localic factorization of a geometric morphism.

\begin{corollary}\label{corhyperconnectedsieves}
	Let $({\cal C}, J)$ be a small-generated site. Then the full subcategory of $\Sh({\cal C}, J)$ on the objects of the form $a_{J}(S)$ for $c\in {\cal C}$ and $S$ a ($J$-closed) sieve on $c$ is equivalent to the category ${\cal C}_{J}^{s}$ defined as follows: 
	\begin{itemize}
		\item the objects of ${\cal C}_{J}^{s}$ are the pairs $(c, S)$ consisting of an object $c$ of $\cal C$ and a ($J$-closed) sieve $S$ on $c$; 
		
		\item the arrows $(c, S)\to (d, T)$ are the collections $R$ of pairs of arrows $(x:e\to c,y:e\to d)$ where $e$ varies among the objects of $\cal C$ (we write $(x, y)\in R(e)$ to mean that $(x, y)\in R$ and $\textup{dom}(x)=\textup{dom}(y)=e$) satisfying the following properties:
		
		\begin{enumerate}[(i)]
			\item for any $(x,y)\in R$, $x\in S$ and $y\in T$;
			
			\item for any arrow $f:e'\to e$, if $(x,y)\in R$ then $(x\circ f, y\circ f)\in R$; 
			
			\item for any $e\in {\cal C}$ and any $(x:e\to c, y:e\to d)\in S\times T$,  if $\{f:e'\to e \mid\ (x\circ f, y\circ f)\in R\}\in J(e)$ then $(x,y)\in R$; 
			
			\item for any $e\in {\cal C}$ and any $(x, y), (x', y')\in R(e)$, if $x=x'$ then $\{f:e'\to e \mid y\circ f=y'\circ f\}\in J(e)$;
			
			\item for any $e\in {\cal C}$ and any $x:e\to c$ in $S$, $\{f:e'\to e \mid \exists y:e'\to d \text{ such that } (x\circ f, y)\in R(e') \}\in J(e)$;
		\end{enumerate}
		
		\item the composite $R'\ast R:(c, S)\to (a, Z)$ of two arrows $R:(c, S)\to (d, T)$ and $R':(d, T)\to (a, Z)$ is given by the following formula: 
		\begin{equation*}\begin{split}
		(R'\ast R)(e)&=\{(x:e\to c, z:e\to a)\in S\times Z \mid \{f:e'\to e \mid \exists y:e'\to d\, \\
		& \quad\quad(x\circ f, y)\in R(e')\textup{ and } (y, z\circ f)\in R'(e')\}\in J(e) \}
		\end{split}\end{equation*}
		for any $e\in {\cal C}$. 
	\end{itemize} 	
\end{corollary}\qed

\begin{remark}\label{remsubcanonicalclosedsubobjects}
	If $\cal C$ is a geometric category (i.e. a well-powered category which has finite limits, images of arrows which are stable under pullback and arbitrary unions of subobjects that are stable under pullback) and $J$ is the geometric topology on it (whose covering sieves are those which contain small families of arrows the union of whose images is the given object) then $\cal C$ is closed under subobjects in $\Sh({\cal C}, J)$, so the category ${\cal C}_{J}^{s}$ is equivalent to $\cal C$. 
\end{remark}

The following proposition characterizes the $J$-functional relations $R$ from $P$ to $Q$ which induce a monomorphism (resp. an epimorphism) $a_{J}(P)\to a_{J}(Q)$:

\begin{proposition}\label{propfunctionalrelationmonoepi}
	Let $R$ be a $J$-functional relation from $P$ to $Q$ in $[{\cal C}^{\textup{op}}, \Set]$ and $\alpha_{R}:a_{J}(P)\to a_{J}(Q)$ the arrow in $\Sh({\cal C}, J)$ induced by $R$. Then 
	\begin{enumerate}[(i)]
		\item  $\alpha_{R}$ is a monomorphism if and only if for any $c\in {\cal C}$ and any elements $(x,y), (x',y')\in R(c)$ such that $y=y'$, the sieve $\{f:d\to c \mid P(f)(x)=P(f)(x')\}$ is $J$-covering.
		
		\item  $\alpha_{R}$ is an epimorphism if and only if for any $c\in {\cal C}$ and $y\in Q(c)$, the sieve $\{f:d\to c \mid (\exists x\in P(d))((x, Q(f)(y))\in R(d))\}$ is $J$-covering.
	\end{enumerate}
\end{proposition}

\begin{proofs}
	Given an arrow $\alpha:A\to B$ in a topos and its graph $R_{\alpha}\mono A\times B$, $\alpha$ is a monomorphism (resp. an epimorphism) if and only if $\pi_{B}:R_{\alpha}\to B$ is a monomorphism (resp. an epimorphism).
	But by Lemma \ref{lemmalift}, $\alpha_{R}$ is a monomorphism if and only if 
	for every $c\in {\cal C}$ and any elements $(x,y), (x',y')\in R(c)$ such that $y=y'$, the sieve $\{f:d\to c \mid P(f)(x)=P(f)(x')\}$ is $J$-covering, while $\alpha_{R}$ is an epimorphism if and only if for any $c\in {\cal C}$ and $y\in Q(c)$, the sieve $\{f:d\to c \mid (\exists x\in P(d))((x, Q(f)(y)\in R(d))) \}$ is $J$-covering.  
\end{proofs}

It is interesting to describe the arrows $l(c)\to a_{J}(P)$ in a topos $\Sh({\cal C}, J)$ of sheaves on a small-generated site $({\cal C}, J)$ both from the point of view of their \ac local' representation (already investigated in section \ref{sec:arrowscomingfromthesite} in the context of arrows $l(c)\to l(d)$ between objects coming from the site) and from the point of view of $J$-functional relations.

\begin{proposition}\label{propsheafification}
	Let $({\cal C}, J)$ be a small-generated site. Then 
	\begin{enumerate}[(i)]
		\item An arrow $\xi:l(c)\to a_{J}(P)$ in $\Sh({\cal C}, J)$ (equivalently, an element of $a_{J}(P)(c)$) can be identified with an equivalence class of families $\{x_{f}\in P(\textup{dom}(f)) \mid f\in S\}$ of elements of $P$ indexed by the arrows $f$ of a $J$-covering sieve $S$ on $c$ which are \emph{locally matching} in the sense that for any arrow $g$ composable with an arrow $f\in S$, $x_{f\circ g}\equiv_{J} P(g)(x_{f})$, modulo the equivalence which identifies two such families when they are locally equal on a common refinement. 
		
		\item Any such family yields a local representation of $\xi$ in the sense that $\xi\circ l(f)=r_{x_{f}}$ for each $f\in S$, where $r_{x_{f}}$ is the image under $a_{J}$ of the arrow $y_{\cal C}(\textup{dom}(f))\to P$ corresponding to the element $x_{f}\in P(\textup{dom}(f))$ via the Yoneda lemma. 
		
		\item Alternatively, the elements of $a_{J}(P)(c)$ can be identified with functions sending to each arrow $f:d\to c$ in $\cal C$ a (possibly empty) $\equiv_{J}$-equivalence class $C_{f}\subseteq P(\textup{dom}(f))$ of elements of $P(\textup{dom}(f))$ in such a way that 
		\begin{itemize}
			\item for any $g$ composable with $f$, $P(g)(C_{f})\subseteq C_{f\circ g}$;
			
			\item for any $y\in P(\textup{dom}(f))$, if $\{g:\textup{dom}(g)\to \textup{dom}(f) \mid P(g)(y)\in C_{f\circ g}\}\in J(\textup{dom}(f))$ then $y\in C_{f}$;
			
			\item $\{g:\textup{dom}(g)\to \textup{dom}(f) \mid C_{f\circ g}\neq \emptyset \}\in J(\textup{dom}(f))$.
		\end{itemize}    
	\end{enumerate}	
\end{proposition}

\begin{proofs}
	The proof of (i) and (ii) is a straightforward generalisation of that of Proposition \ref{propexplicit}, while part (iii) follows at once from Theorem \ref{thmfunctionalrelations} and Proposition \ref{propfuncionalassignment}.  	
\end{proofs}

\begin{remark}
	Proposition \ref{propsheafification}(i) gives an explicit description of the associated sheaf functor $a_{J}(P)$ of a presheaf $P$, different from the usual construction of it by means of the double plus construction. This alternative construction of the associated sheaf functor seems to have been first discovered (albeit not published) by Eduardo Dubuc in the eighties. 
\end{remark} 

\subsection{Relative cofinality}\label{sec:relativecofinality}

In this section we introduce a notion of cofinal functor relative to a site, and show that it naturally corresponds to the property of inducing an isomorphism between colimits in the relevant toposes. 

Let $({\cal C}, J)$ be a small-generated site and $F:{\cal A}\to {\cal C}$ and $F:{\cal A}'\to {\cal C}$ two functors to $\cal C$ related by a functor $\xi:{\cal A}\to {\cal A}'$ and a natural transformation $\alpha:F\to F'\circ \xi$. We have an arrow 
\[
\tilde{\alpha}: \textup{colim}_{[{\cal C}^{\textup{op}}, \Set]}(y_{\cal C}\circ F) \to \textup{colim}_{[{\cal C}^{\textup{op}}, \Set]}(y_{\cal C}\circ F')
\]
in $[{\cal C}^{\textup{op}}, \Set]$, defined through the universal property of the colimit by setting
\[
\tilde{\alpha} \circ \lambda_{a}= \chi_{\xi(a)}\circ y_{\cal C}(\alpha(a))
\]
(for any $a\in {\cal A}$), where $\lambda_{a}:y_{\cal C}(F(a))\to \textup{colim}_{[{\cal C}^{\textup{op}}, \Set]}(y_{\cal C}\circ F)$ (for $a\in {\cal A}$) and $\chi_{a'}:y_{\cal C}(F'(a'))\to \textup{colim}_{[{\cal C}^{\textup{op}}, \Set]}(y_{\cal C}\circ F')$ (for $d'\in {\cal A}'$) are the canonical colimit arrows. 

To show that this definition is indeed well-posed, we recall that colimits in functor categories are computed pointwise; so, for each $c\in {\cal C}$, we have 
\[
\textup{colim}_{[{\cal C}^{\textup{op}}, \Set]}(y_{\cal C}\circ F)(c)=(\coprod_{a\in {\cal A}}\textup{Hom}_{\cal C}(c, F(a)))\slash R_{c},
\] 
where $R_{c}$ is the equivalence relation given by: $(x:c\to F(a), x':c\to F(b))\in R_{c}$ if and only if $x$ and $x'$ belong to the same connected component of the category $(c\downarrow F)$. Similarly, we have 
\[
\textup{colim}_{[{\cal C}^{\textup{op}}, \Set]}(y_{\cal C}\circ F')(c)=(\coprod_{a'\in {\cal A}'}\textup{Hom}_{\cal C}(c, F(a')))\slash R'_{c},
\]
where $R'_{c}$ is the equivalence relation given by: $(y:c\to F(a'), y':c\to F(b'))\in R'_{c}$ if and only if $y$ and $y'$ belong to the same connected component of the category $(c\downarrow F')$.

We thus see that the arrows $\{\chi_{\xi(a)}\circ y_{\cal C}(\alpha(a)) \mid a\in {\cal A}\}$ form indeed a cocone under the diagram $y_{\cal C}\circ F$, since they respect these equivalence relations. 

We want to understand under which conditions the arrow $\tilde{\alpha}$ is sent by the associated sheaf functor $a_{J}:[{\cal C}^{\textup{op}}, \Set] \to \Sh({\cal C}, J)$ to an isomorphism. For this, we can use Lemma \ref{lemmalift}(iii); this yields the following result: 

\begin{proposition}\label{procofinality}
	Let $({\cal C}, J)$ be a small-generated site and $F:{\cal A}\to {\cal C}$ and $F':{\cal A}'\to {\cal C}$ two functors to $\cal C$ related by a functor $\xi:{\cal A}\to {\cal A}'$ and a natural transformation $\alpha:F\to F'\circ \xi$. Then the canonical arrow 
	\[
	\tilde{\alpha}: \textup{colim}_{[{\cal C}^{\textup{op}}, \Set]}(y_{\cal C}\circ F) \to \textup{colim}_{[{\cal C}^{\textup{op}}, \Set]}(y_{\cal C}\circ F')
	\]
	defined above is sent by $a_{J}$ to an isomorphism 
	\[
	a_{J}(\tilde{\alpha}): \textup{colim}_{\Sh({\cal C}, J)}(l\circ F) \to \textup{colim}_{\Sh({\cal C}, J)}(l\circ F) 
	\]
	if and only if $(\xi, \alpha)$ satisfies the following \ac cofinality' conditions:
	\begin{enumerate}[(i)]
		\item For any object $c$ of $\cal C$ and any arrow $y:c\to F'(a')$ in $\cal C$ there are a $J$-covering family $\{f_i: c_i \to c \mid i\in I\}$ and for each $i\in I$ an object $a_{i}$ of $\cal A$ and an arrow $y_i:c_{i}\to F(a_{i})$ such that $(y\circ f_{i}, \alpha(a_{i})\circ y_{i})\in R'_{c_{i}}$.
		
		\item For any object $c$ of $\cal C$ and any arrows $x:c\to F(a)$ and $x':c\to F(b)$ in $\cal C$ such that $(\alpha(a)\circ x, \alpha(b)\circ x')\in R'_{c}$ there is a $J$-covering family $\{f_i: c_i \to c \mid i\in I\}$ such that $(x\circ f_i, x'\circ f_i)\in R_{c_{i}}$ for each $i\in I$. 
	\end{enumerate}	
\end{proposition}\qed

To better understand the conditions of Proposition \ref{procofinality}, it is natural to consider the fibrations $\pi_{\cal C}^{F}:(1_{\cal C}\downarrow F)\to {\cal C}$ and $\pi^{F'}_{\cal C}:(1_{\cal C}\downarrow F)\to {\cal C}$ given by the canonical projection functors (cf. section \ref{sec:fibrgenelements} for an analysis of this type of fibrations). Then the natural transformation $\alpha$ yields a functor $\alpha_{\textup{fib}}: (1_{\cal C}\downarrow F) \to (1_{\cal C}\downarrow F')$, sending to each object $(c, a, x:c\to F(a))$ of $(1_{\cal C}\downarrow F)$ the object $(c, \xi(a), \alpha(a)\circ x:c\to F'(\xi(a)))$ of $(1_{\cal C}\downarrow F')$, which is actually a morphisms of fibrations.  

\begin{remarks}\label{rempropcofinality}
	\begin{enumerate}[(a)]
		\item Condition (i) of Proposition \ref{procofinality} can be replaced by the single condition that for any object $u:=(c, a', y:c\to F'(a'))$ of the category $(1_{\cal C}\downarrow F')$, there is a $J$-covering family such that the category $((u\circ f_{i})\downarrow \alpha_{\textup{fib}})$ is non-empty for each $i\in I$.
		
		\item In condition (i) of Proposition \ref{procofinality}, if the claim is satisfied for an arrow $y:d\to F(a')$ then it is satisfied by any arrow of the form $y\circ z$ for an arrow $z$ composable with $y$. It follows in particular that if $\xi$ is essentially surjective and $\alpha$ is the identity then the condition is satisfied.
		
		\item Proposition \ref{procofinality} can be notably applied in the context of colimit-preserving functors $$Z:\Sh({\cal C}, J)\to \Sh({\cal D}, K)$$ between Grothendieck toposes such that $Z\circ l=l'\circ z$ for some functor $z:{\cal C}\to {\cal D}$ (such as, for instance, inverse image functors of geometric morphisms induced by a morphism of sites, or essential images of essential geometric morphisms induced by a continuous comorphism of sites, as investigated in section \ref{sec:essentialmorphismsandcomorphisms}). Indeed, the value of such a functor $Z$ on a $J$-sheaf $Q$ on $\cal C$ is given by $\textup{colim}_{\Sh({\cal D}, K)}(Z\circ l \circ \pi_{Q})=\textup{colim}_{\Sh({\cal D}, K)}(l'\circ z \circ \pi_{Q})$, where $\pi_{Q}:{\int Q} \to {\cal C}$ is the canonical projection functor from the category $\int Q$ of elements of $Q$. See section \ref{sec:locallyconnectedmorphisms} below for an application of this remark. 
		
		\item We shall see in section \ref{sec:essentialmorphismsandcomorphisms} (cf. Proposition \ref{propcofinalcontinuous}) that, under the hypotheses of Proposition \ref{procofinality}, if $p$ is a continuous comorphism of sites $({\cal C}, J)\to ({\cal D}, K)$ such that the essential image of the associated geometric morphism $C_{p}$ is conservative then $(\xi, \alpha)$ satisfies the conditions of Proposition \ref{procofinality} if and only if $(\xi, p\alpha)$ does.
	\end{enumerate}
\end{remarks}

Let us now apply Proposition \ref{procofinality} in two notable cases:
\begin{enumerate}[(1)]
	\item $F=\xi:{\cal A}\to {\cal C}$, $F'=1_{\cal C}$, $\alpha$ is the identity. In this case the equivalence relations $R'_{c}$ are the total relations, and the proposition yields the following criterion: the canonical arrow
	\[
	\textup{colim}_{[{\cal C}^{\textup{op}}, \Set]}(y_{\cal C}\circ \xi) \to \textup{colim}_{[{\cal C}^{\textup{op}}, \Set]}(y_{\cal C})=1_{[{\cal C}^{\textup{op}}, \Set]}
	\]
	is sent by $a_{J}$ to an isomorphism
	\[
	\textup{colim}_{\Sh({\cal C}, J)}(l\circ \xi) \to 1_{\Sh({\cal C}, J)}
	\]
	if and only if the following conditions are satisfied:
	\begin{enumerate}[(i)]
		\item For any object $c$ of $\cal C$ there are a $J$-covering family $\{f_i: c_i \to c \mid i\in I\}$ and for each $i\in I$ an object $a_{i}$ of $\cal A$ and an arrow $y_i:c_{i}\to F(a_{i})$.
		
		\item For any object $c$ of $\cal C$ and any arrows $x:c\to F(a)$ and $x':c\to F(b)$ in $\cal C$, there is a $J$-covering family $\{f_i: c_i \to c \mid i\in I\}$ such that $(x\circ f_i, x'\circ f_i)\in R_{c_{i}}$ for each $i\in I$. 
	\end{enumerate}	
	
	Notice that if $J$ is the trivial Grothendieck topology then the above conditions amount precisely to the requirement that for every object $c$ of $\cal C$, the category $(c\downarrow \xi)$ be non-empty and connected, which means precisely that the functor $F$ is cofinal (in the classical sense). It is therefore natural to call \emph{$J$-cofinal} a functor $F:{\cal A}\to {\cal C}$ which satisfies the two conditions above (see Definition \ref{defrelcofinality} below). 
	
	\item $F'$ is the forgetful functor $U_{c_{0}}:{\cal C}\slash c_{0}\to {\cal C}$ for an object $c_{0}$ of $\cal C$, $\xi$ is a cocone $\{\xi_{a}: F(a)\to c_{0} \mid a\in {\cal A}\}$ under the functor $F$ with vertex $c_{0}$ and $\alpha$ is the identity. In this case we have: for any $y:c \to (U_{c_{0}}([p:c'\to c_{0}]))$ and $y':c \to (U_{c_{0}}([q:c''\to c_{0}]))$, $(y, y')\in R'_{c}$ if and only if $p\circ y=q\circ y'$. So the proposition yields the following criterion: the cocone $\xi$ is sent by the canonical functor ${\cal C}\to \Sh({\cal C}, J)$ to a colimit if and only if the following conditions are satisfied:
	\begin{enumerate}[(i)]
		\item For any object $c$ of $\cal C$ and any arrow $y:c\to c_{0}$ in $\cal C$ there are a $J$-covering family $\{f_i: c_i \to c \mid i\in I\}$ and for each $i\in I$ an object $a_{i}$ of $\cal A$ and an arrow $y_i:c_{i}\to F(a_{i})$ such that $y\circ f_{i}=\xi_{a_{i}}\circ y_{i}$.
		
		\item For any object $c$ of $\cal C$ and any arrows $x:c\to F(a)$ and $x':c\to F(b)$ in $\cal C$ such that $\xi_{a}\circ x= \xi_{b}\circ x'$ there is a $J$-covering family $\{f_i: c_i \to c \mid i\in I\}$ such that $(x\circ f_i, x'\circ f_i)\in R_{c_{i}}$ for each $i\in I$. 
	\end{enumerate}	 
\end{enumerate}

\begin{definition}\label{defrelcofinality}
	Given a small site $({\cal C}, J)$, a functor $F:{\cal A}\to {\cal C}$ is said to be \emph{$J$-cofinal} if the following conditions are satisfied:
	\begin{enumerate}[(i)]
		\item For any object $c$ of $\cal C$ there are a $J$-covering family $\{f_i: c_i \to c \mid i\in I\}$ and for each $i\in I$ an object $a_{i}$ of $\cal A$ and an arrow $y_i:c_{i}\to F(a_{i})$.
		
		\item For any object $c$ of $\cal C$ and any arrows $x:c\to F(a)$ and $x':c\to F(b)$ in $\cal C$ there is a $J$-covering family $\{f_i: c_i \to c \mid i\in I\}$ such that $x\circ f_i$ and $x'\circ f_i$ belong to the same connected component of the category $(c_i \downarrow F)$ for each $i\in I$. 
	\end{enumerate}		
\end{definition}

The above discussion thus yields the following two corollaries.

\begin{corollary}\label{corcofinalcharacterization}
	Let $({\cal C}, J)$ be a small-generated site and $F:{\cal A}\to {\cal C}$	a functor. Then $F$ is $J$-cofinal if and only if the canonical arrow
	\[
	\textup{colim}_{\Sh({\cal C}, J)}(l\circ F) \to 1_{\Sh({\cal C}, J)}
	\] 
	is an isomorphism.
\end{corollary}\qed

Given a functor $D:{\cal A}\to {\cal C}$ and a cocone $\{\xi_{a}: D(a)\to c_{0} \mid a\in {\cal A}\}$ under $D$ with vertex $c_{0}$, we can lift $D$ to a functor $D_{\xi}:{\cal A}\to {\cal C}\slash c_{0}$ such that $U_{c_{0}}\circ D_{\xi}=D$. Let $J_{c_{0}}$ be the Grothendieck topology on ${\cal C}\slash c_{0}$ whose covering sieves are precisely those whose image under $U_{c_{0}}$ is $J$-covering. Then the arrow
\[
\textup{colim}_{[{\cal C}^{\textup{op}}, \Set]}(y_{\cal C}\circ D) \to \textup{colim}_{[{\cal C}^{\textup{op}}, \Set]}(y_{\cal C}\circ U_{c_{0}})\cong y_{\cal C}(c_{0})
\]
induced by $\xi$, regarded as an arrow in $[{\cal C}^{\textup{op}}, \Set] \slash y_{\cal C}(c_{0})$, corresponds to the unique arrow
\[
\textup{colim}_{[({{\cal C}\slash c_{0}})^{\textup{op}}, \Set]}(y_{({\cal C}\slash c_{0})}\circ D_{\xi}) \to 
1_{[{({\cal C}\slash c_{0})}^{\textup{op}}, \Set]} 
\] 
and hence the image of the former under $a_{J}$ corresponds to the image of the latter under $a_{J_{c_{0}}}$ via the equivalence $\Sh({\cal C}, J)\slash l(c_{0}) \simeq \Sh({\cal C}\slash c_{0}, J_{c_{0}})$. So the cocone $\xi$ is sent by the canonical functor ${\cal C}\to \Sh({\cal C}, J)$ to a colimit if and only if $D_{\xi}$ is $J_{c_{0}}$-cofinal:

\begin{corollary}\label{corcolimitintopos}
	Let $D:{\cal A}\to {\cal C}$ be a functor and $\xi$ a cocone $\{\xi_{a}: D(a)\to c_{0} \mid a\in {\cal A}\}$ under $D$ with vertex $c_{0}$. Then $\xi$ is sent by the canonical functor $l:{\cal C}\to \Sh({\cal C}, J)$ to a colimit cocone if and only if the functor $D_{\xi}$ is $J_{c_{0}}$-cofinal, equivalently if and only if the following conditions are satisfied:
	\begin{enumerate}[(i)]
		\item For any object $c$ of $\cal C$ and any arrow $y:c\to c_{0}$ in $\cal C$ there are a $J$-covering family $\{f_i: c_i \to c \mid i\in I\}$ and for each $i\in I$ an object $a_{i}$ of $\cal A$ and an arrow $y_i:c_{i}\to D(a_{i})$ such that $y\circ f_{i}=\xi_{a_{i}}\circ y_{i}$.
		
		\item For any object $c$ of $\cal C$ and any arrows $x:c\to D(a)$ and $x':c\to D(b)$ in $\cal C$ such that $\xi_{a}\circ x= \xi_{b}\circ x'$ there is a $J$-covering family $\{f_i: c_i \to c \mid i\in I\}$ such that $x\circ f_i$ and $x'\circ f_i$ belong to the same connected component of the category $(c_{i} \downarrow D)$ for each $i\in I$. 
	\end{enumerate}		
\end{corollary}\qed

Corollary \ref{corcolimitintopos} can be applied in particular to characterize colimits in a Grothendieck topos in terms of generalized elements:

\begin{corollary}
	Let $D:{\cal A}\to {\cal E}$ be a functor from a small category $\cal A$ to a Grothendendieck topos $\cal E$ and $\xi$ a cocone $\{\xi_{a}: D(a)\to e_{0} \mid a\in {\cal A}\}$ under $D$ with vertex $e_{0}$. Then $\xi$ is a colimit cocone if and only if the functor $D_{\xi}$ is $({J^{\textup{can}}_{\cal E}})_{e_{0}}$-cofinal, equivalently if and only if the following conditions are satisfied:
	\begin{enumerate}[(i)]
		\item For any object $e$ of $\cal E$ and any arrow $y:e\to e_{0}$ in $\cal E$ there are an epimorphc family $\{f_i: e_i \to e \mid i\in I\}$ in $\cal E$ and for each $i\in I$ an object $a_{i}$ of $\cal A$ and an arrow $y_i:e_{i}\to D(a_{i})$ such that $y\circ f_{i}=\xi_{a_{i}}\circ y_{i}$.
		
		\item For any object $e$ of $\cal E$ and any arrows $x:e\to D(a)$ and $x':e\to D(b)$ in $\cal E$ such that $\xi_{a}\circ x= \xi_{b}\circ x'$ there is an epimorphic family $\{f_i: e_i \to e \mid i\in I\}$ in $\cal E$ such that $x\circ f_i$ and $x'\circ f_i$ belong to the same connected component of the category $(e_{i} \downarrow D)$ for each $i\in I$. 
	\end{enumerate}		
\end{corollary} 

The following proposition expresses a natural relation between the notion of $J$-cofinality and the conditions of Proposition \ref{procofinality} (see section \ref{sec:essentialcomorphisms} for the notions of continuous functor and weak morphism of toposes):

\begin{proposition}\label{propcontinuitycofinality}
	Let $A:{\cal A}\to {\cal C}$ be a functor and $v:({\cal C}, J)\to ({\cal C}', J')$ a continuous functor (that is, a weak morphism of sites). If $A$ is $J$-cofinal then it satisfies, together with the identical natural transformation on $v\circ F$, the conditions of Proposition \ref{procofinality}. Moreover, the converse holds if $\Sh(v)^{\ast}$ reflects isomorphisms.  
\end{proposition}

\begin{proofs}
	Since $v$ is continuous, the following square commutes:	
	$$
	\xymatrix{
		{\cal C} \ar[r]^{v} \ar[d]_{l_{{\cal C}}} & {\cal C}' \ar[d]^{l_{\cal C}'} \\
		\Sh({\cal C}, J) \ar[r]^{\Sh(v)^{\ast}} & \Sh({\cal C}', J')
	}
	$$   
	By Corollary \ref{corcofinalcharacterization}, $F$ is $J$-cofinal if and only if the unique arrow  
	\[
	\textup{colim}_{\Sh({\cal C}, J)}(l_{\cal C}\circ F) \to 1_{\Sh({\cal C}, J)}
	\] 
	is an isomorphism. Since $(C_{v})_{!}$ is colimit-preserving and the above square commutes, we have 
	\[
	\Sh(v)^{\ast}(\textup{colim}_{\Sh({\cal C}, J)}(l_{\cal C}\circ F))\cong \textup{colim}_{\Sh({{\cal C}'}, J')}(l_{{\cal C}'}\circ v\circ F) 
	\] 
	and
	\[
	\Sh(v)^{\ast}(1_{\Sh({\cal C}, J)})=\textup{colim}_{\Sh({{\cal C}'}, J')}(l_{{\cal C}'}\circ v)
	\]
	So if the above arrow is an isomorphism then the canonical arrow
	\[
	\textup{colim}_{\Sh({{\cal C}'}, J')}(l_{{\cal C}'}\circ v\circ F) \to \textup{colim}_{\Sh({{\cal C}'}, J')}(l_{{\cal C}'}\circ v)
	\] 
	induced by $F$ is an isomorphisms (that is, $F$ satisfies the conditions of Proposition \ref{procofinality} together with the identical natural transformation on $v\circ F$), and the converse holds if $\Sh(v)^{\ast}$ reflects isomorphisms.  
\end{proofs}

\section{Morphisms and comorphisms of sites}\label{sec:morphismsandcomorphisms}

\subsection{Geometric morphisms and flat functors}

Recall that any functor $A:{\cal C}\to {\cal E}$ from an essentially small category $\cal C$ to a Grothendieck topos $\cal E$ induces an adjuction $(L_{A} \dashv R_{A})$, where $R_{A}:{\cal E}\to [{\cal C}^{\textup{op}}, \Set]$ is the \acc hom' functor $\textup{Hom}_{\cal E}(A(-),-)$ and $L_{A}:[{\cal C}^{\textup{op}}, \Set] \to {\cal E}$ is the unique (up to isomorphism) colimit-preserving functor which extends $A$ along the Yoneda embedding $y_{\cal C}:{\cal C}\to [{\cal C}^{\textup{op}}, \Set]$, in other words, the left Kan extension of $A$ along the Yoneda embedding $y_{\cal C}:{\cal C}\to [{\cal C}^{\textup{op}}, \Set]$, which sends each presheaf $P$ on $\cal C$ to the colimit of the diagram $A\circ \pi_{P}$, where $\pi_{P}:{\int P} \to {\cal C}$ is the canonical projection from the category ${\int P}$ of elements  of $P$ (see, for instance, section I.5 of \cite{MM}). The functor $A$ is said to be \emph{flat} if $L_{A}$ preserves finite limits. In this case, the pair $(L_{A}\dashv R_{A})$ defines a geometric morphism ${\cal E}\to [{\cal C}^{\textup{op}}, \Set]$. If $J$ is a Grothendieck topology on $\cal C$ and $A$ is flat then the geometric morphism $(L_{A}\dashv R_{A}):{\cal E}\to [{\cal C}^{\textup{op}}, \Set]$ factors through the canonical geometric inclusion $\Sh({\cal C}, J)\hookrightarrow [{\cal C}^{\textup{op}}, \Set]$ if and only if $A$ is $J$-continuous. In other words, we have a categorical equivalence, usually called \emph{Diaconescu's equivalence}
\[
\textbf{Geom}({\cal E}, \Sh({\cal C}, J))\simeq \textbf{Flat}_{J}({\cal C}, {\cal E}),
\]
where $\textbf{Geom}({\cal E}, \Sh({\cal C}, J))$ is the category of geometric morphisms from ${\cal E}$ to $\Sh({\cal C}, J)$ and $\textbf{Flat}_{J}({\cal C}, {\cal E})$ is the category of $J$-continuous flat functors ${\cal C}\to {\cal E}$.

In section VII.9 of \cite{MM}, the flat functors ${\cal C}\to {\cal E}$ are characterized precisely as the functors ${\cal C}\to {\cal E}$ that are \emph{filtering}, in the sense of the following result:

\begin{theorem}[Theorem VII.9.1 \cite{MM}]\label{filteringfunctor}
	A functor $A:{\cal C}\to {\cal E}$ is flat if and only if it is filtering in the sense that it satisfies the following conditions:	
	\begin{enumerate}[(i)]
		\item The family of arrrows $A(c)\to 1_{\cal E}$ (for $c\in {\cal C}$) is epimorphic in $\cal E$;
		
		\item For any objects $c$ and $d$ of $\cal C$, the family of arrows of the form $\langle A(u), A(v)\rangle :A(b)\to A(c)\times A(d)$, where $b$ is an object of $\cal C$ $u$ and $v$ are arrows respectively $b\to c$ and $b\to d$ in $\cal C$, is epimorphic in $\cal E$;
		
		\item For any two arrows $u, v:c\to d$ in $\cal C$, let $i:\textup{Eq}(A(u), A(v))\mono A(c)$ be the equalizer of $A(u)$ and $A(v)$ in $\cal E$. Then the family of factorizations $A(b)\to \textup{Eq}(A(u), A(v))$ through $i$ of arrows of the form $A(w):A(b)\to A(c)$, where $w:b\to c$ is an arrow in $\cal C$ such that $u\circ w=v\circ w$, is epimorphic. 		
	\end{enumerate}
\end{theorem}

\subsection{Morphisms of sites}

In \cite{grothendieck} (Expos\'e III.1), a functor $F:{\cal C}\to {\cal D}$ between the underlying categories $\cal C$ and $\cal D$ of two (small-generated) sites $({\cal C}, J)$ and $({\cal D}, K)$ is defined to be (what is now called) a \emph{morphism of sites} $({\cal C}, J)\to ({\cal D}, K)$ if it induces a geometric morphism $f:\Sh({\cal D}, K)\to \Sh({\cal C}, J)$, in the sense that there is a commutative diagram 
	$$
\xymatrix{
	{\mathcal C} \ar[r]^F \ar[d]_{l} &{\cal D} \ar[d]^{l'} \\
	\Sh({\cal C}, J) \ar[r]^{f^{\ast}} & \Sh({\cal D}, K),
}
$$
where $l$ and $l'$ are the canonical functors. 

If $F$ is a morphism of sites  $({\cal C}, J)\to ({\cal D}, K)$, the geometric morphism $\Sh({\cal D}, K)\to \Sh({\cal C}, J)$ that it induces will be denoted by $\Sh(F)$.

By Diaconescu's equivalence, the geometric morphisms 
\[
\Sh({\cal D}, K)\to \Sh({\cal C}, J)
\]
correspond precisely to the flat $J$-continuous functors ${\cal C}\to \Sh({\cal D}, K)$, and, thanks to Theorem \ref{filteringfunctor} and Proposition \ref{propexplicit}(i) (which allows us to $K$-locally represent the arrows in the image of $l'$ in terms of arrows in $\cal D$), we can explicitly rephrase the condition for a functor $F:{\cal C}\to {\cal D}$ that $l'\circ F$ be flat and $J$-continuous. This leads to the following characterization/definition of morphisms of sites:

\begin{definition}\label{defmorphismsites}
	A functor $F:{\cal C}\to {\cal D}$ is said to be a \emph{morphism of sites} $({\cal C}, J)\to ({\cal D}, K)$, where $J$ is a collection of sieves in $\cal C$ and $K$ is a Grothendieck topology on $\cal D$, if it satisfies the following conditions:	
	
	\begin{enumerate}[(i)]
		\item $F$ sends every $J$-covering family in ${\mathcal C}$ into a $K$-covering family in ${\mathcal D}$.
		
		\item Every object $d$ of ${\mathcal D}$ admits a $K$-covering family
		$$
		d_i \longrightarrow d \, , \quad i \in I \, ,
		$$
		by objects $d_i$ of ${\mathcal D}$ which have morphisms 
		$$
		d_i \longrightarrow F(c'_i)
		$$
		to the images under $F$ of objects $c'_i$ of ${\mathcal C}$.
		
		\item For any objects $c_1 , c_2$ of ${\mathcal C}$ and any pair of morphisms of ${\mathcal D}$
		$$
		g_1 : d \longrightarrow F(c_1) \, , \quad g_2 : d \longrightarrow F(c_2) \, ,
		$$
		there exists a $K$-covering family
		$$
		g'_i : d_i \longrightarrow d \, , \quad i \in I \, ,
		$$
		and a family of pairs of morphisms of ${\mathcal C}$
		$$
		f_1^i : c'_i \longrightarrow c_1 \, , \quad f_2^i : c'_i \to c_2 \, , \quad i \in I \,  ,
		$$
		and of morphisms of ${\mathcal D}$
		$$
		h_i : d_i \longrightarrow F(c'_i) \, , \quad i \in I \, ,
		$$
		making the following squares commutative:
		$$
		\xymatrix{
			d_i \ar[r]^{g'_i} \ar[d]_{h_i} & d \ar[d]^{g_1} \\
			F(c'_i) \ar[r]^{F(f_1^i)} &F(c_1)
		} 
		\qquad \qquad \qquad
		\xymatrix{
			d_i \ar[r]^{g'_i} \ar[d]_{h_i} & d \ar[d]^{g_2} \\
			F(c'_i) \ar[r]^{F(f_2^i)} &F(c_2)
		} 
		$$
		
		\item For any pair of arrows $f_1 , f_2 : c_1 \rightrightarrows c_2$ of ${\mathcal C}$ and any arrow of ${\mathcal D}$
		$$
		g : d \longrightarrow F(c_{1})
		$$
		satisfying
		$$
		F(f_1) \circ g = F(f_2) \circ g \, ,
		$$
		there exist a $K$-covering family
		$$
		g_i : d_i \longrightarrow d \, , \quad i \in I \, ,
		$$
		and a family of morphisms of ${\mathcal C}$
		$$
		k_i : c'_i \longrightarrow c_{1} \, , \quad i \in I \, ,
		$$
		satisfying
		$$
		f_1 \circ k_i = f_2 \circ k_i \, , \quad \forall \, i \in I \, ,
		$$
		and of morphisms of ${\mathcal D}$
		$$
		h_i : d_i \longrightarrow F(c'_i) \, , \quad i \in I \, ,
		$$
		making the following squares commutative:
		$$
		\xymatrix{
			d_i \ar[r]^{g_i} \ar[d]_{h_i} & d \ar[d]^{g} \\
			F(c'_i) \ar[r]^{F(k_i)} & F(c_{1})
		} 
		$$ 
	\end{enumerate} 
\end{definition}

\begin{remarks}\label{remmorphismsofsites}
	\begin{enumerate}[(a)]
		\item A functor $F:{\cal C}\to {\cal D}$ is a morphism of sites $({\cal C}, J)\to ({\cal D}, K)$ (in the sense of Definition \ref{defmorphismsites}) if and only if $l'\circ F$ is a $J$-continuous flat functor ${\cal C}\to \Sh({\cal D}, K)$, and, conversely, for any small-generated site $({\cal C}, J)$ and any Grothendieck topos $\cal E$, a $J$-continuous flat functor ${\cal C}\to {\cal E}$ is precisely a morphism of sites $({\cal C}, J)\to ({\cal E}, J_{\cal E}^{\textup{can}})$. 
		
		\item The above characterization of morphisms of sites $({\cal C}, J)\to ({\cal D}, K)$ as the functors $F:{\cal C}\to {\cal D}$ which induce a geometric morphism $f:\Sh({\cal D}, K) \to \Sh({\cal C}, J)$ actually holds for arbitrary collections $J$ of sieves in $\cal C$, if we define $\Sh({\cal C}, J):=\Sh({\cal C}, \tilde{J})$, where $\tilde{J}$ is the Grothendieck topology on $\cal C$ generated by $J$ (cf. Proposition \ref{propgeneratedmorphismofsites} below).
		
		\item Definition \ref{defmorphismsites} is equivalent to the one given in \cite{Shulman}, which specifies that a functor is a morphism of sites when it is cover-preserving and covering-flat (in the sense that for any finite diagram $D$ in $\cal C$ every cone over an object of the form $F(c)$ factors locally through the $F$-image of a cone over $D$), and also proves the characterization in (b) by using this latter definition. 
	\end{enumerate}
\end{remarks}

\begin{proposition}\label{propgeneratedmorphismofsites}
	Let $J'$ a collection of sieves on a category $\cal C$, and $J$ the Grothendieck topology on $\cal C$ generated by it. Then any morphism of sites $({\cal C}, J')\to ({\cal D}, K)$ is a morphism of sites $({\cal C}, J)\to ({\cal D}, K)$.   
\end{proposition}

\begin{proofs}
We have to show that a morphism of sites $F:({\cal C}, J')\to ({\cal D}, K)$ is cover-preserving as a morphism $({\cal C}, J)\to ({\cal D}, K)$ (the other conditions in the notion of morphism of sites being independent from the Grothendieck topology in the source site). For this, it clearly suffices to prove that the collection of sieves $S$ in $\cal C$ such that $F(S)$ generates a $K$-covering sieve is a Grothendieck topology containing $J'$. 

The maximality and transitivity axioms are clearly satisfied, so it remains to prove the pullback-stability axiom. We shall deduce this from the fact that $F:({\cal C}, J')\to ({\cal D}, K)$ satisfies conditions (ii) and (iii) in the definition of morphism of sites. We have to prove that if $S$ is a sieve on $c$ such that the sieve $\langle  F(S)\rangle$ generated by $F(S)$ is $K$-covering then for any arrow $f:c'\to c$, the sieve $\langle  F(f^{\ast}(S))\rangle  $ generated by $F(f^{\ast}(S))$ is also $K$-covering. For this, using the transitivity axiom for $K$, we are reduced to verifying that for each $\xi:d\to F(c')$ in $F(f)^{\ast}(\langle  F(S)\rangle  )$, $\xi^{\ast}(\langle  F(f^{\ast}(S))\rangle  )$ is $K$-covering. Since $\xi\in F(f)^{\ast}(\langle  F(S)\rangle  )$, there exist $s:c''\to c$ in $S$ and $\chi:d\to F(c'')$ in $\cal D$ such that $F(f)\circ \xi=F(s)\circ \chi$. By condition (ii) of Definition \ref{defmorphismsites}, there exist a $K$-covering family
$$
g'_i : d_i \longrightarrow d \, , \quad i \in I \, ,
$$
and a family of pairs of morphisms of ${\mathcal C}$
$$
f_1^i : c'_i \longrightarrow c'' \, , \quad f_2^i : c'_i \to c' \, , \quad i \in I \,  ,
$$
and of morphisms of ${\mathcal D}$
$$
h_i : d_i \longrightarrow F(c'_i) \, , \quad i \in I \, ,
$$
such that $F(f_1^i)\circ h_{i}=\chi \circ g'_i$ and $F(f_2^i)\circ h_{i}=\xi \circ g'_i$. Now, for each $i\in I$, consider the arrows $s\circ f_{1}^{i}, f\circ f_{2}^{i}: c'_{i} \rightrightarrows c$, which satisfy $F(s\circ f_{1}^{i})\circ h_{i}=F(f\circ f_{2}^{i})\circ h_{i}$. By condition (iii) of Definition \ref{defmorphismsites}, there exist a $K$-covering family
$$
g^i_j : d_j^i \longrightarrow d_{i} \, , \quad j \in J_{i} \, ,
$$
and a family of morphisms of ${\mathcal C}$
$$
k^i_j : {c'}_j^i \longrightarrow c'_{i}  \, , \quad j \in J_{i} \, ,
$$
satisfying
$$
s\circ f_{1}^{i} \circ k^i_j =f\circ f_{2}^{i} \circ k^i_j \, , \quad \forall \, j \in J_{i} \, ,
$$
and of morphisms of ${\mathcal D}$
$$
h^i_j : d^i_j \longrightarrow F({c'}_j^i) \, , \quad j \in J_{i} \, ,
$$
such that $h_{i}\circ g^i_j=F(k^i_j)\circ h^i_j$.

Therefore, the family of arrows $\{g'_i \circ g^i_j \mid i\in I, j\in J_{i}\}$ is $K$-covering and is contained in the sieve $\xi^{\ast}(\langle F(f^{\ast}(S))\rangle  )$. Indeed, $\xi\circ g'_i \circ g^i_j =F(f_2^i)\circ h_{i}\circ g^i_j=F(f_2^i \circ k^i_j)\circ h^i_j$, which belongs to $\langle F(f^{\ast}(S))\rangle  $ since the equality $s\circ f_{1}^{i} \circ k^i_j =f\circ f_{2}^{i} \circ k^i_j$ implies that $f_{2}^{i} \circ k^i_j\in f^{\ast}(S)$ and hence that $F(f_{2}^{i} \circ k^i_j)\in F(S)$. So $\xi^{\ast}(\langle F(f^{\ast}(S))\rangle  )$ is $K$-covering, as desired.	
\end{proofs}

\subsection{Comorphisms of sites}\label{sec:comorphismsofsites1}

Recall that a \emph{comorphism of sites} $({\cal D}, K)\to ({\cal C}, J)$ (where $J$ and $K$ are Grothendieck topologies respectively on $\cal C$ and $\cal D$) is a functor $F:{\cal D}\to {\cal C}$ which has the \emph{covering lifting property}, that is the property that for every $d\in {\cal D}$ and any $J$-covering sieve $S$ on $F(d)$ there is a $K$-covering sieve $R$ on $c$ such that $F(R)\subseteq S$.

As it is well-known (cf. Theorem VII.5 \cite{MM}), every comorphism of sites $F:({\cal D}, K)\to ({\cal C}, J)$ induces a $J$-continuous flat functor
\[
A_{F}:{\cal C}\to \Sh({\cal D}, K),
\]
which assigns to each object $c$ of $\cal C$ the $K$-sheaf given by $a_{K}(\textup{Hom}_{\cal C}(F(-),c))$ (and acts on the arrows in the obvious way), and hence, by Diaconescu's equivalence, a geometric morphism $C_{F}:\Sh({\cal D}, K)\to \Sh({\cal C}, J)$ (whose inverse image $C_{F}^{\ast}:\Sh({\cal C}, J) \to \Sh({\cal D}, K)$ sends a $J$-sheaf $P$ to $a_{K}(P\circ F^{\textup{op}})$). In fact, more generally, a functor $F:{\cal D} \to {\cal C}$ is a comorphism of sites $({\cal D}, K) \to ({\cal C}, J)$ if and only if the following diagram commutes:
$$
\xymatrix{
	[{\mathcal C}^{\textup{op}}, \Set] \ar[r]^{-\circ F^{\textup{op}}} \ar[d]_{a_{J}} & [{\cal D}^{\textup{op}}, \Set] \ar[d]^{a_{K}} \\
	\Sh({\cal C}, J) \ar[r]^{C_{F}^{\ast}} & \Sh({\cal D}, K).
}
$$

We shall also write $D_{F}$ for the functor $-\circ F^{\textup{op}}:[{\mathcal C}^{\textup{op}}, \Set] \to [{\cal D}^{\textup{op}}, \Set]$.

Notice that the canonical geometric inclusion $\Sh({\cal C}, J)\hookrightarrow [{\cal C}^{\textup{op}}, \Set]$ can be seen either as induced by the morphism of sites $({\cal C}, T)\to ({\cal C}, J)$ given by the identity functor $1_{\cal C}$ on $\cal C$, where $T$ is the trivial Grothendieck topology on $\cal C$, or as induced in this way by the comorphism of sites $1_{\cal C}:({\cal C}, J)\to ({\cal C}, T)$. 

Comorphisms of sites notably arise in the context of flat functors inducing equivalences of toposes. Let $F:{\cal C}\to {\cal E}$ be a flat functor on an essentially small category $\cal C$ with values in a Grothendieck topos $\cal E$, and let $J_{F}$ be the induced Grothendieck topology on $\cal C$ (in the sense of Proposition \ref{reminducedtopology}). Suppose that $F$ induces an equivalence ${\cal E} \simeq \Sh({\cal C}, J_{F})$ (equivalently, satisfies the conditions of Corollary \ref{flatequivalence}, that is, is $J_{F}$-full and the objects of the form $F(c)$ for $c\in {\cal C}$ form a separating family of objects of $\cal E$). Then $F$ is both a morphism and a comorphism of sites $({\cal C}, J_{F})\to ({\cal E}, J_{\cal E}^{\textup{can}})$ (cf. Corollary \ref{corequivalencecovliftproperty}), and the geometric morphisms
\[
\Sh(F):\Sh({\cal E}, J_{\cal E}^{\textup{can}}) \to \Sh({\cal C}, J_{F}) 
\]
and 
\[
C_{F}:\Sh({\cal C}, J_{F})  \to \Sh({\cal E}, J_{\cal E}^{\textup{can}})
\]
are quasi-inverse to each other. 

For any small-generated site $({\cal C}, J)$ the canonical functor $$l:({\cal C}, J) \to (\Sh({\cal C}, J), J^{\textup{can}}_{\Sh({\cal C}, J)})$$ is a comorphism of sites (which induces the identity geometric morphism on $\Sh({\cal C}, J)$), and for any site of definition $({\cal D}, K)$ of a Grothendieck topos $\cal E$, a comorphism of sites $({\cal C}, J)\to ({\cal E}, J^{\textup{can}}_{\cal E})$ yields a comorphism of sites $({\cal C}, J)\to ({\cal D}, K)$ if and only if it factors through the canonical functor ${\cal D}\to \Sh({\cal D}, K)$.

\subsubsection{The smallest Grothendieck topology making a functor a comorphism of sites}\label{sec:smallestgrothendiecktopologycomorphism}

Recall (see, for instance, Lemma C2.3.19(i) \cite{El}) that, given a functor $A:{\cal C}\to {\cal D}$ and a Grothendieck topology $K$ in $\cal D$, there is a smallest Grothendieck topology on $\cal C$ which makes $A$ a comorphism of sites to $({\cal D}, K)$. This topology, which we denote by $M^{A}_{K}$, is generated by the (pullback-stable) family of sieves of the form $S^{A}_{R}:=\{f:\textup{dom}(f)\to c \mid A(f)\in R\}$ for an object $c$ of $\cal C$ and a $K$-covering sieve $R$ on $A(c)$. 

\begin{lemma}\label{lemmacomorphismsofsitesgenerators}
	Let $({\cal C}, J)$ be a site and $F:{\cal C}\to {\cal D}$ a functor. Let us suppose that $\tilde{K}$ is a Grothendieck topology on $\cal D$ generated by a pullback-stable family of sieves $K$ on $\cal D$. Then, if for every sieve $R$ in $K$ on an object of the form $F(c)$ there is a $J$-covering sieve $S$ on $c$ such that $F(S)\subseteq R$, $F$ is a comorphism of sites $({\cal C}, J)\to ({\cal D}, \tilde{K})$. 
\end{lemma}

\begin{proofs}
	We can explicitly describe $\tilde{K}$ in terms of $K$ as follows: the $\tilde{K}$-covering sieves are precisely those which contain a finite multicomposition of $K$-covering sieves. But the covering-lifting property is clearly inherited by taking bigger sieves or multicompositions of sieves, whence our thesis follows.
\end{proofs}

\begin{corollary}\label{corliftcomorphismsofsites}
	Let $G:{\cal A}\to {\cal C}$ and $F:{\cal C}\to {\cal D}$ be functors such that $F\circ G$ is a comorphism of sites $({\cal A}, A)\to ({\cal D}, L)$. Then $G$ is a comorphism of sites $({\cal A}, A)\to ({\cal C}, M^{F}_{L})$.
	
	If $F$ is full and faithful, then $F:({\cal C}, M^{F}_{L})\to ({\cal D}, L)$ is cover-reflecting.
\end{corollary}

\begin{proofs}
Our thesis follows at once from Lemma \ref{lemmacomorphismsofsitesgenerators} in light of the fact that, as recalled above, the Grothendieck topology $M^{F}_{L}$ is generated by the pullback-stable collection of sieves of the form $\{f:\dom(f)\to c \mid F(f)\in R\}$ where $R$ varies among the $L$-covering sieves on $F(c)$.
	
Let us now suppose that $F$ is full and faithful. If $S$ is a sieve such that $F(S)$ is $K$-covering then the family $\{g\in {\cal C} \mid F(g)\in \langle F(S)\rangle \}$ is $M^{F}_{K}$-covering. But, since $F$ is full and faithful, $F(g)\in \langle F(S)\rangle $ if and only if $g\in \langle S\rangle $, so $S$ is $M^{F}_{K}$-covering, as required. 	
\end{proofs}

We shall now proceed to investigate the Grothendieck topologies of the form $M^{i}_{J}$, where $i:{\cal D}\to {\cal C}$ is a full and faithful functor. 

Consider the geometric morphism $C_{i}:[{\cal D}^{\textup{op}}, \Set]\to [{\cal C}^{\textup{op}}, \Set]$ induced by $i$; it is an inclusion (cf. Remark \ref{remarkinclusion}) and hence we have a bijective correspondence between the subtoposes of $[{\cal D}^{\textup{op}}, \Set]$ (equivalently, the Grothendieck topologies on $\cal D$) and the subtoposes of $[{\cal C}^{\textup{op}}, \Set]$ which are contained in $C_{i}$ (equivalently, the Grothendieck topologies $J$ on $\cal C$ which contain the rigid topology $R_{i}$ whose covering sieves on an object $c$ are the sieves containing all the arrows from objects of the form $i(d)$ to $c$). Given a Grothendieck topology $K$ on $\cal D$, the Grothendieck topology on $\cal C$ associated with it via this correspondence is precisely the topology $K^{i}$ coinduced by $K$ along $i$ (in the sense of Proposition \ref{propimagetopology}). On the other hand, for any Grothendieck topology $J$ on $\cal C$, the Grothendieck topology on $\cal D$ associated with the intersection of the subtoposes $\Sh({\cal C}, J)\hookrightarrow [{\cal C}^{\textup{op}}, \Set]$ and $C_{i}$ (that is, with the Grothendieck topology on $\cal C$ generated by $R_{i}$ and $J$) is the smallest Grothendieck topology $M^{i}_{J}$ on $\cal D$ which makes $i$ a comorphism of sites to $({\cal C}, J)$. Indeed, $M^{i}_{J}$ is characterized in the proof of Lemma C2.3.19 \cite{El} by the pullback square
\[
\xymatrix{
	\Sh({\cal D}, M^{i}_{J}) \ar[r] \ar[d]  & \Sh({\cal C}, J) \ar[d] \\
	[{\cal D}^{\textup{op}}, \Set] \ar[r]^{C_{i}} & [{\cal C}^{\textup{op}}, \Set],
} 
\]	
from which it is clear that subtopos given by the composite of $\Sh({\cal D}, M^{i}_{J})\hookrightarrow [{\cal D}^{\textup{op}}, \Set]$ with $C_{i}$ is precisely the subtopos of $[{\cal C}^{\textup{op}}, \Set]$ given by the intersection of the subtoposes $\Sh({\cal C}, J)\hookrightarrow [{\cal C}^{\textup{op}}, \Set]$ and $C_{i}$. This shows that $J\vee R_{i}=(M^{i}_{J})^{i}$.

Notice that, since the equivalence 
$$\Sh({\cal D}, M^{i}_{J})\simeq \Sh({\cal C}, J\vee R_{i})$$ makes the square
\[
\xymatrix{
	\Sh({\cal D}, M^{i}_{J}) \ar[r] \ar[d]  & \Sh({\cal C}, J\vee R_{i}) \ar[d] \\
	[{\cal D}^{\textup{op}}, \Set] \ar[r]^{C_{i}} & [{\cal C}^{\textup{op}}, \Set],
} 
\]    
commute (the morphism $\Sh({\cal D}, M^{i}_{J})\to \Sh({\cal C}, J)$ is the composite of the equivalence $\Sh({\cal D}, M^{i}_{J})\simeq \Sh({\cal C}, J\vee R_{i})$ with the canonical inclusion $\Sh({\cal C}, J\vee R_{i})\hookrightarrow \Sh({\cal C}, J)$), the functor $i$ is a comorphism of sites $({\cal D}, M^{i}_{J})\to ({\cal C}, J\vee R_{i})$ (cf. section \ref{sec:comorphismsofsites1}).

In the converse direction, for any Grothendieck topology $K$ on $\cal D$, we have a commutative square
\[
\xymatrix{
	\Sh({\cal D}, K) \ar[r] \ar[d]  & \Sh({\cal C}, K^{i}) \ar[d] \\
	[{\cal D}^{\textup{op}}, \Set] \ar[r]^{C_{i}} & [{\cal C}^{\textup{op}}, \Set],
} 
\] 
where the upper horizontal morphism is an equivalence (by the universal property of the topology $K^{i}$). This square is thus a pullback and hence from the pullback characterization of the topology $M^{i}_{K^{i}}$ it follows that $K=M^{i}_{K^{i}}$.

Summarizing, we have the following result:

\begin{theorem}\label{thmcorrespondenceinclusion}
	Let $i:{\cal D}\to {\cal C}$ be a full and faithful functor. Then the assignments $$J \mapsto M^{i}_{J}$$ and $$K\mapsto K^{i}$$ between the classes of Grothendieck topologies $J$ on $\cal C$ and $K$ on $\cal D$ satisfy the following properties:
	\begin{enumerate}[(i)]
		\item $J\vee R_{i}=(M^{i}_{J})^{i}$, where $R_{i}$ is the rigid topology on $\cal D$ associated with $i$;
		
		\item $K=M^{i}_{K^{i}}$.
	\end{enumerate}
Moreover, for any $J$ and any $K$, the functor $i$ is a comorphism of sites $({\cal D}, M^{i}_{J})\to ({\cal C}, J\vee R_{i})$ inducing an equivalence of toposes
\[
\Sh({\cal D}, M^{i}_{J})\simeq \Sh({\cal C}, J\vee R_{i}),
\]
and a comorphism of sites $({\cal D}, K)\to ({\cal C}, K^{i})$ inducing an equivalence
\[
\Sh({\cal D}, K)\simeq \Sh({\cal C}, K^{i}).
\]
\end{theorem}\qed

As remarked in the proof of Lemma C2.3.19 \cite{El}, the topology $M^{i}_{J}$ does not in general admit an easy explicit description in terms of $J$ and $i$. Nonetheless, if $i$ is $(M^{i}_{J}, J)$-continuous, it admits an alternative characterization, as a consequence of the following result: 

\begin{proposition}\label{proprecoverytopologycontinuity}
	Let $i$ be a $K$-full and $K$-faithful $(K, J)$-continuous comorphism of sites $({\cal D}, K)\to ({\cal C}, J)$. Then $K$ can be explicitly characterized in terms of $J$ as follows: for any sieve $T$ on an object $d$ of $\cal D$, $T\in K(d)$ if and only if the morphism $\textup{Lan}_{i^{\textup{op}}}(m_{T})$ is $J$-bicovering, where $m_{T}:T\mono y_{\cal D}(d)$ is the monomorphism in $[{\cal D}^{\textup{op}}, \Set]$ corresponding to the sieve $T$ and $\textup{Lan}_{i^{\textup{op}}}$ is the left Kan extension functor $[{\cal D}^{\textup{op}}, \Set]\to [{\cal C}^{\textup{op}}, \Set]$ along $i^{\textup{op}}$.  
\end{proposition} 

\begin{proofs}
	We preliminarily observe that, since $i$ is $K$-full and $K$-faithful, $C_{i}$ is an inclusion by Proposition \ref{propinclusioncomoprhismofsites}, equivalently the functor $(C_{i})_{!}$ is full and faithful. So $(C_{i})_{!}$ reflects isomorphisms. Since $i$ is $(K, J)$-continuous, by Corollary \ref{corequivcharcanonicity} we have a commutative square 
	\[
	\xymatrix{
		[{\cal D}^{\textup{op}}, \Set]  \ar[r]^{\textup{Lan}_{i^{\textup{op}}}} \ar[d]^{a_{K}}  & [{\cal C}^{\textup{op}}, \Set] \ar[d]^{a_{J}} \\
		\Sh({\cal D}, K) \ar[r]^{(C_{i})_{!}} & \Sh({\cal C}, J)
	} 
	\]
	Now, the commutativity of this diagram implies that the canonical morphism $m_{T}:T\mono y_{\cal D}(d)$ is sent by $(C_{i})_{!} \circ a_{K}$ to an isomorphism (equivalently, since $(C_{i})_{!}$ reflects isomorphisms, $T$ is $K$-covering) if and only $a_{J}$ sends $\textup{Lan}_{i^{\textup{op}}}(m_{T})$ to an isomorphism.  
\end{proofs} 

\begin{remark}
	Proposition \ref{proprecoverytopologycontinuity} asserts precisely that $K$ coincides with the Grothendieck topology induced by $i$, regarded as a functor to the site $({\cal C}, J)$, in the sense of section III-3 of \cite{grothendieck} (by Proposition III-3.2 therein -- cf. also Remark \ref{remcontinuityprelim}(c)); that is, $K$ is the largest Grothendieck topology which makes $i$ continuous.  
\end{remark}

\begin{corollary}\label{corexplicitdescriptionsmallesttopology}
	Let $({\cal C}, J)$ be a small-generated site and $i$ a full and faithful morphism of sites $({\cal D}, T_{\cal D})\to ({\cal C}, J\vee R_{i})$, where $T_{\cal D}$ is the trivial Grothendieck topology on $\cal D$. Then for any sieve $T$ on an object $d$ of $\cal D$, $T\in M^{i}_{J}(d)$ if and only if it is $i(T)$ is $(J\vee R_{i})$-covering.
\end{corollary}

\begin{proofs}
	As observed above, $i$ is a comorphism of sites $({\cal D}, M^{i}_{J})\to ({\cal C}, J\vee R_{i})$. On the other hand, by our hypothesis, $i$ is also a morphism of sites $({\cal D}, M^{i}_{J})\to ({\cal C}, J\vee R_{i})$ and hence it is $(K, J\vee R_{i})$-continuous. We can therefore apply Proposition \ref{proprecoverytopologycontinuity} to it. But, since $a_{J\vee R_{i}}\circ \textup{Lan}_{i^{\textup{op}}}$ preserves finite limits and hence monomorphisms, $\textup{Lan}_{i^{\textup{op}}}(m_{T})$ is $(J\vee R_{i})$-bicovering if and only if the flat functor $a_{J}\circ \textup{Lan}_{i^{\textup{op}}}\circ y_{\cal D}=l_{\cal C}\circ i$ sends $T$ to an epimorphic family of $\Sh({\cal C}, J\vee R_{i})$, that is, if and only if $i$ sends $T$ to a $(J\vee R_{i})$-covering family. 
\end{proofs}	

The following result shows that, under some specific hypotheses, the Grothendieck topology $J\vee R_{i}$ admits a particularly simple description.

\begin{proposition}
	Let $i$ be a full and faithful functor ${\cal D}\to {\cal C}$ and $J$ a Grothendieck topology on $\cal C$ such that for every $J$-covering sieve $T$ on an object of the form $i(d)$ for $d\in {\cal D}$, $T\cap \textup{Im}(i)\in J(i(d))$. Then the Grothendieck topology $J\vee R_{i}$ admits the following explicit description: for any sieve $T$ on an object $c$ of $\cal C$, $T\in (J\vee R_{i})(c)$ if and only if for every arrow $\xi:i(d)\to c$, $\xi^{\ast}(T)\in J(i(d))$.  
\end{proposition}

\begin{proofs}
	Let $W_{J}$ be the collection of sieves on $\cal C$ defined as follows: for any sieve $T$ on an object $c$ of $\cal C$, $T\in W_{J}(c)$ if and only if for every arrow $\xi:i(d)\to c$, $\xi^{\ast}(T)\in J(i(d))$. Then $J\subseteq W_{J}$ and $R_{i}\subseteq W_{J}$. Also, $W_{J}\subseteq J\vee R_{i}$ by the transitivity axiom for Grothendieck topologies. So, to show that $W_{J}=J\vee R_{i}$, it remains to verify that, under our assumption, $W_{J}$ is a Grothendieck topology. The fact that the maximality and pullback stability axioms are satisfied is clear, as it is the fact that any sieve containing a sieve in $W_{J}$ lies again in $W_{J}$. So it remains to show that arbitrary multicompositions of sieves in $W_{J}$ lie in $W_{J}$. So, suppose that $S\in W_{J}(c)$ and that, for each $f\in S$, $S_{f}\in W_{J}(\dom(f))$. Then the sieve $S\ast \{S_{f} \mid f\in S\}$ belongs to $W_{J}(c)$ since for any  $\xi:i(d)\to c$, $\xi^{\ast}(S\ast \{S_{f} \mid f\in S\})$ contains the multicomposition $(\xi^{\ast}(S)\cap \textup{Im}(i)) \ast \{S_{\xi\circ g} \mid g\in  (\xi^{\ast}(S)\cap \textup{Im}(i))\}$, all of whose sieves are $J$-covering (note that $S_{f}$ is $J$-covering for any arrow $f$ whose domain is of the form $i(d)$ for an object $d\in {\cal D}$). 
\end{proofs}

\subsubsection{Fibrations as comorphisms of sites}\label{subsubsecfibrationscomorphisms}

Given a fibration $p:{\cal C}\to {\cal D}$ to a category $\cal D$ endowed with a Grothendieck topology $K$, it is natural to regard it as a comorphism of sites $({\cal C}, M^{p}_{K})\to ({\cal D}, K)$. This point of view on fibrations is extensively developed in a forthcoming joint work with Riccardo Zanfa, in the context of which the following theorem, providing an explicit characterization of the Grothendieck topology $M^{p}_{K}$, has been discovered. Before stating it, we shall establish a number of technical results about fibrations. 

Recall that, given a functor $p:{\cal C}\to {\cal D}$, an arrow $\phi:c\to c'$ is said to be \emph{cartesian} (relative to $p$) if for any arrow $\psi:c'' \to c$ in $\cal C$ and $g:p(c'')\to p(c')$ such that $p(\phi)\circ g=p(\psi)$, there exists a unique arrow $\chi:c''\to c'$ such that $\psi=\phi\circ \chi$ and $p(\chi)=g$. We shall say that an arrow of $\cal C$ is said to be \emph{vertical} if its image under $p$ is an isomorphism.

Commonly, a functor $p$ is said to be a fibration if for any $d\in {\cal D}$ and any arrow $f:d\to p(c)$ in $\cal D$, there is a cartesian arrow $\phi:c'\to c$ such that $p(\phi)=f$. This definition of fibration, whilst equivalent to the original one by Grothendieck (\cite{SGA1}) and found almost everywhere in the categorical literature on fibrations, is ill-posed since it does not even include functors which are one half of an equivalence of categories. Indeed, the equality $p(\phi)=f$ implies in particular the equality \emph{of objects} $\dom(f)=p(\dom(\phi))$, which is a form of surjectivity of the functor $p$ which is too strong and unnatural from a categorical point of view. Indeed, surjectivity at the level of objects should always be replaced by essential surjectivity within a categorical framework. The good, categorically-invariant, definition, already pointed out by R. Street (see \cite{nLabstreetfibration}), is as follows: a functor $p:{\cal C}\to {\cal D}$ is a \emph{fibration} if for any arrow $f:d\to p(c)$ in $\cal D$, there is a cartesian arrow $\phi:c'\to c$ and an isomorphism $\alpha:p(c')\to d$ such that $f\circ \alpha=p(\phi)$. Notice that this definition covers in particular the case of a functor $p$ which is part of an equivalence of categories; indeed, in such a situation every arrow is cartesian and one can obtain an adjoint equivalence between $p$ (right adjoint) and its quasi-inverse $q$; this adjunction provides in particular, for any arrow $f:d\to p(c)$ in $\cal D$, an isomorphism $\alpha:p(q(d))\cong d$ given by the inverse of the unit of the adjunction at $d$, and an arrow $\phi:q(d)\to c$, given by the transpose of $f$ along the adjunction, such that $f\circ \alpha =p(\phi)$.  

Fibrations satisfy a number of elementary properties: the identity arrow on each object is cartesian, the composite of two cartesian arrow is cartesian, and any arrow $u$ of $\cal C$ can be decomposed, uniquely up to isomorphism, as the composite of a vertical arrow followed by a cartesian arrow. We shall call \ac the' cartesian arrow arising in this factorization, which is determined up to unique compatible isomorphism between the domains, the \emph{cartesian image} of $u$, and denote it by $\phi_{u}$.   

We shall call a family of arrows with common codomain a \emph{presieve}.

Given a functor $p:{\cal C}\to {\cal D}$ we have two operations $P \mapsto p(P)$ and $R \mapsto S^{p}_{R}:=\{f:\textup{dom}(f)\to c \mid p(f)\in R\}$ sending presieves on $\cal C$ to presieves on $\cal D$ and conversely. As shown by the following result, if $p$ is a fibration then these operations are naturally related to each other.

The sieve $\langle P\rangle $ generated by a presieve $P$ on an object $c$ is the collection of arrows to $c$ which factor through an arrow of $P$. We shall denote by $\overline{P}^{i}$ the collection of arrows $f$ such that $f\circ \alpha\in P$ for some isomorphism $\alpha$. If $p:{\cal C}\to {\cal D}$ is a fibration and $P$ is a presieve in $\cal C$ on an object $c$, we denote by $P^{\textup{cart}}$ the collection of cartesian images of arrows in $P$, and by $\overline{P}^{v}$ the collection of the arrows of the form $f\circ z$ where $f\in P$ and $z$ is a vertical arrow. Notice that $P\subseteq \overline{P^{\textup{cart}}}^{v}$ (since every arrow of $\cal C$ can be factored as a vertical arrow followed by its cartesian image) but the converse inclusion does not hold in general. 

\begin{lemma}\label{lemmafibrations}
	Let $p:{\cal C}\to {\cal D}$ be a fibration. Then
	\begin{enumerate}[(i)]
		\item For any presieve $P$  in $\cal C$, we have $\overline{p(\langle P\rangle )}^{i}=\langle p(P)\rangle $. 
		
		\item For any presieve $R$ on an object $p(c)$ in $\cal D$, $\langle S^{p}_{\overline{R}^{i}}\rangle =S^{p}_{\langle R\rangle }$. 
		
		\item For any presieve $P$ in $\cal C$, $S^{p}_{\overline{p(P)}^{i}}=\overline{P^{\textup{cart}}}^{v}$. In particular (cf. also (ii)), if $P$ entirely consists of cartesian arrows (whence $P^{\textup{cart}}=\overline{P}^{i}$) then $\langle S^{p}_{\overline{p(P)}^{i}}\rangle =S^{p}_{\langle p(P)\rangle }=\langle P\rangle $. 
		
		\item For any presieve $R$ on an object $p(c)$ in $\cal D$, $\overline{R}^{i}=\overline{p(S^{p}_{\overline{R}^{i}})}^{i}$; in particular, $\langle R\rangle =\langle p(S^{p}_{\overline{R}^{i}})\rangle $, so, if $R$ is a sieve, then $R=\langle p(S^{p}_{R})\rangle $. 
		
		\item For any sieve $R$ on an object $p(c)$ in $\cal D$, $S^{p}_{R}=\langle (S^{p}_{R})^{\textup{cart}}\rangle $.
		
		\item For any presieve $P$ on an object $c$ of $\cal C$ entirely consisting of cartesian arrows, and any arrow $f:c'\to c$ in $\cal C$, 
		\[
		\langle p(f^{\ast}(\langle P\rangle ))\rangle =p(f)^{\ast}(\langle p(P)\rangle ).
		\] 
		
		\item For any presieve $P$ on an object $c$ of $\cal C$ entirely consisting of cartesian arrows, the sieves of the form $f^{\ast}(\langle P\rangle )$, for an arrow $f:c'\to c$ in $\cal C$, have the property that they contain all the cartesian images of their arrows.
	\end{enumerate}	
\end{lemma}

\begin{proofs}
	(i) The inclusion $\subseteq$ is obvious (and holds in general), so we have to show the converse one. Given $r\in \langle p(P)\rangle $, $r=p(f)\circ \beta$ for some $f\in P$ and arrow $\beta$ in $\cal D$. Since $p$ is a fibration, there is an isomorphism $\alpha$ such that $\beta\circ \alpha=p(g)$ for some arrow $g$ to $\dom(f)$ in $\cal C$. So $r\circ \alpha=p(f)\circ \beta\circ \alpha=p(f)\circ p(g)=p(f\circ g)$; in other words, $r=p(f\circ g)\circ \alpha^{-1}$. Since $f\circ g\in \langle P\rangle $, $r\in \overline{p(\langle P\rangle )}^{i}$, as required. 
	
	(ii) The inclusion $\subseteq$ is obvious (and valid in general), so we only have to prove the converse one. Given $f\in S^{p}_{\langle R\rangle }$, $p(f)\in \langle R\rangle $, i.e. $p(f)=r\circ \gamma$ for some $r\in R$. Since $f$ is a fibration, there is a cartesian arrow $\phi$ and an isomorphism $\delta$ such that $r\circ \delta=p(\phi)$. Now, since $\phi$ is cartesian and $p(f)=p(\phi)\circ (\delta^{-1} \circ \gamma)$, there exists a unique arrow $\psi$ such that $p(\psi)=\delta^{-1} \circ \gamma$ and $f=\phi\circ \psi$. Now, since $p(\phi)=r\circ \delta$, $\phi\in S^{p}_{\overline{R}^{i}}$, whence $f\in \langle S^{p}_{\overline{R}^{i}}\rangle $.
	
	(iii) Let us first prove that $S^{p}_{\overline{p(P)}^{i}}\supseteq \overline{P^{\textup{cart}}}^{v}$. Let us consider \ac the' factorization $\phi_{u}\circ \gamma$ of $u$ as a vertical arrow followed by a cartesian arrow. Given $f\in \overline{P^{\textup{cart}}}^{v}$, there is an arrow $u\in P$ and a vertical arrow $s$ such that $f=\phi_{u}\circ s$, where $\phi_{u}$ is \ac the' cartesian image of $u$. Then $p(f)=p(\phi_{u})\circ p(s)$ and, since $p(s)$ is an isomorphism, $p(f)=p(u)\circ p(\gamma)^{-1}\circ p(s)\in \overline{p(P)}^{i}$, as required. Conversely, let us suppose that $f\in S^{p}_{\overline{p(P)}^{i}}$. Then $p(f)=p(u)\circ \alpha$ for some arrow $u\in P$ and some isomorphism $\alpha$. Since $\phi_{u}$ is cartesian, there is an arrow $\psi$ such that $p(\psi)=p(\gamma)\circ \alpha$ and $f=\phi_{u}\circ \psi$. So $\psi$ is a vertical arrow and hence $f\in \overline{P^{\textup{cart}}}^{v}$, as required.
	
	The second statement follows from the first in light of (ii).
	
	(iv) Let us start showing that $\overline{R}^{i}\subseteq \overline{p(S^{p}_{\overline{R}^{i}})}^{i}$, equivalently, that $R\subseteq \overline{p(S^{p}_{\overline{R}^{i}})}^{i}$. Given $r\in R$, since $p$ is a fibration, there is a cartesian arrow $\phi$ to $c$ and an isomorphism $\alpha$ such that $r=p(\phi)\circ \alpha$. Therefore $\phi\in S^{p}_{\overline{R}^{i}}$ and hence $r\in \overline{p(S^{p}_{\overline{R}^{i}})}^{i}$. Conversely, suppose that $u\in \overline{p(S^{p}_{\overline{R}^{i}})}^{i}$. Then $u=p(f)\circ \delta$ for some arrow $f\in S^{p}_{\overline{R}^{i}}$ and isomorphism $\delta$. Now, since  $f\in S^{p}_{\overline{R}^{i}}$, $p(f)=r\circ \beta$ for some arrow $r\in R$ and isomorphism $\beta$. So $u=r\circ (\alpha\circ \beta)$ and hence $u\in \overline{R}^{i}$.   
	
	The second statement clearly follows from the first, since taking the sieve generated by a presieve $S$ is the same as taking the sieve generated by the presieve $\overline{S}^{i}$. 
	
	(v) This follows from the fact that $S^{p}_{R}$ satisfies the property that for any vertical arrow $v$ and any arrow $b$ composable with it $b\circ v\in S^{p}_{R}$ if and only if $b\in S^{p}_{R}$ (since $p(v)$ is an isomorphism and $R$ is a sieve).
	
	(vi) Let us set $R:=p(f)^{\ast}(\langle p(P)\rangle )$. By (iv), we have $R=\langle p(S^{p}_{R})\rangle $, so it is equivalent to prove that $\langle p(S^{p}_{R})\rangle =\langle p(f^{\ast}(\langle P\rangle ))\rangle $. Let us first show that $\langle p(f^{\ast}(\langle P\rangle ))\rangle \subseteq \langle p(S^{p}_{R})\rangle $, or equivalently that $p(f^{\ast}(\langle P\rangle ))\subseteq \langle p(S^{p}_{R})\rangle $. This inclusion does not exploit the hypothesis that $P$ only consist of cartesian arrows. Given $s\in p(f^{\ast}(\langle P\rangle ))$, $s=p(g)$ where $g\in f^{\ast}(\langle P\rangle )$. Now, $g\in S^{p}_{R}$. Indeed, $p(g)\in R=p(f)^{\ast}(\langle p(P)\rangle )$ as $p(f)\circ p(g)=p(f\circ g)\in p(\langle P\rangle )\subseteq \langle p(P)\rangle $. So $g=p(s)\in p(S^{p}_{R})$, as required. Conversely, let us show that $p(S^{p}_{R})\subseteq \langle p(f^{\ast}(\langle P\rangle ))\rangle $. Given $s\in p(S^{p}_{R})$, we have an arrow $g\in S^{p}_{R}$ such that $s=p(g)$. Let us show that $g\in f^{\ast}(\langle P\rangle )$. Since $p(g)\in R=p(f)^{\ast}(\langle p(P)\rangle )$, $p(f)\circ p(g)=\xi\circ \gamma$ for some arrow $\xi\in P$ and some arrow $\gamma$. Now, since $\xi$ is cartesian (by our hypothesis every arrow of $P$ is cartesian), there is a unique arrow $\psi$ in $\cal C$ such that $p(\psi)=\gamma$ and $\xi\circ \psi=f\circ g$. So $g\in f^{\ast}(\langle P\rangle )$, as required.
	
	(vii) We know from (iii) that, if $P$ entirely consists of cartesian arrows then $\langle P\rangle =S^{p}_{\langle p(P)\rangle }$. We want to prove that if $g\in f^{\ast}(\langle P\rangle )$ (equivalently, $f\circ g\in \langle P\rangle $) and we have a factorization $g=g'\circ h$ where $g'$ is cartesian and $h$ is vertical, then $g'\in f^{\ast}(\langle P\rangle )$ (equivalently, $f\circ g'\in \langle P\rangle $). But this follows from (v) (take $R=\langle p(P)\rangle $, $v$ is $h$ and $b$ is $f\circ g'$).   
\end{proofs}

\begin{theorem}\label{thmfibration1}
	Let $p:{\cal C}\to {\cal D}$ be a fibration. Then a sieve $R$ is $M^{p}_{K}$-covering if and only if the collection of cartesian arrows in $R$ is sent by $p$ to a $K$-covering family.
\end{theorem}

\begin{proofs}
	Let us first show that the collection $\cal S$ of sieves $R$ such that the collection of their cartesian arrows is sent by $p$ to a $K$-covering family defines a Grothendieck topology on $\cal C$. The maximality axiom is trivially satisfied; $\cal S$ is clearly closed under taking bigger sieves and under multicomposition of arrows (since $K$ is closed under taking bigger sieves and under multicomposition and the composite of two cartesian arrows is always a cartesian arrow). So it remains to prove that it satisfies the pullback stability property. Given a sieve $R$ on an object $c$ and an arrow $f:c'\to c$ in $\cal C$, suppose that the presieve $P_{R}=\{f \mid f\in R \textup{ and } f \textup{ is cartesian}\}$ is sent by $p$ to a $K$-covering family. Then by Lemma \ref{lemmafibrations}(vii) the sieve $f^{\ast}(\langle P_{R}\rangle )$ satisfies the property that all the cartesian images of the arrows in it lie again in it. Now, by Lemma \ref{lemmafibrations}(vi), $f^{\ast}(\langle P_{R}\rangle )$ is sent by $p$ to a $K$-covering family, since the sieve generated by $p(f^{\ast}(\langle P_{R}\rangle ))$ is equal to the pullback sieve $p(f)^{\ast}(\langle p(P_{R})\rangle )$, which is $K$-covering since $\langle p(P_{R})\rangle $ is $K$-covering. It follows that the collection $(f^{\ast}(\langle P_{R}\rangle ))^{\textup{cart}}$ of cartesian images of arrows in  $f^{\ast}(\langle P_{R}\rangle )$, which is contained, as we said, in $f^{\ast}(\langle P_{R}\rangle )$ and hence \emph{a fortiori} in $f^{\ast}(R)$, is sent by $p$ to a $K$-covering family. Therefore $f^{\ast}(R)$ lies in $\cal S$, as required. 
	
	Now that we have proved that $\cal S$ is a Grothendieck topology on $\cal C$, it remains to show that it coincides with $M^{p}_{K}$. The fact that $M^{p}_{K}$ is contained in $\cal S$ follows from the fact that for any $K$-covering sieve $R$ on an object of the form $p(c)$, the sieve $S^{p}_{R}$ lies in $\cal S$ since, by Lemma \ref{lemmafibrations}(iv), $R=\langle p(S^{p}_{R})\rangle $ and, by Lemma \ref{lemmafibrations}(v), $S^{p}_{R}=\langle (S^{p}_{R})^{\textup{cart}}\rangle $, whence the collection $(S^{p}_{R})^{\textup{cart}}$, which is contained in $S^{p}_{R}$, is sent by $p$ to a $K$-covering family. To show the converse inclusion ${\cal S}\subseteq M^{p}_{K}$, it suffices to observe that, given a sieve $R$ in $\cal S$, the presieve $P_{R}$ consisting of the cartesian arrows in $R$ satisfies the property that $\langle p(P_{R})\rangle $ is $K$-covering; so the sieve $S^{p}_{\langle p(P_{R})\rangle }$ lies in $M^{p}_{K}$. But $S^{p}_{\langle p(P_{R})\rangle }=\langle P_{R}\rangle $ by Lemma \ref{lemmafibrations}(iii), so $\langle P_{R}\rangle $, and \emph{a fortiori} $R$, which contains it, is in  $M^{p}_{K}$.  
\end{proofs}

\subsection{Relating morphisms and comorphisms of sites}

In this section we shall describe a procedure for turning any morphism of sites into a comorphism of sites and conversely, by replacing the domain site with another site giving rise to the same topos. These constructions are applicable in many contexts; indeed, there are many situations where working with morphisms (resp. comorphisms) of sites is more natural than with comorphisms (resp. morphisms) of sites; think for instance of pullbacks of toposes, which admit natural descriptions in terms of comorphisms of sites, or to certain bilimits of toposes which can be effectively described in terms of colimits of morphisms of sites.  

\begin{proposition}\label{propadjointinduction}
	Let $({\cal C}, J)$ and $({\cal D}, K)$ be small-generated sites, and $F:{\cal C}\to {\cal D}\dashv G:{\cal D}\to {\cal C}$ adjoint functors. Then
	\begin{enumerate}[(i)]
		\item $F$ has the covering-lifting property if and only if $G$ is cover-preserving.
		
		\item $G$ is a morphism of sites $({\cal D}, T_{\cal D})\to ({\cal C}, T_{\cal C})$, where $T_{\cal C}$ and $T_{\cal D}$ are the trivial topologies respectively on $\cal C$ and $\cal D$, and induces an essential geometric morphism $[{\cal C}^{\textup{op}}, \Set]\to [{\cal D}^{\textup{op}}, \Set]$;
		
		\item $G$ is a morphism of sites $({\cal D}, K) \to ({\cal C}, J)$ (equivalently, by (i), $G$ is cover-preserving) if and only if $F$ is a comorphism of sites $({\cal C}, J)\to ({\cal D}, K)$.
		
		\item In the situation of (iii), the geometric morphism $C_{F}$ induced by $F$ co-incides with the geometric morphism $\Sh(G)$ induced by $G$.
	\end{enumerate}
\end{proposition}

\begin{proofs}
	(i) This is Lemma 3(ii) at p. 410 of \cite{MM}. 
	
	(ii) By definition of morphism of sites, $G$ is a morphism of sites $({\cal D}, T_{\cal D})\to ({\cal C}, T_{\cal C})$ if and only if the functor $L_{y_{\cal C}\circ G}:[{\cal D}^{\textup{op}}, \Set] \to [{\cal C}^{\textup{op}}, \Set]$ preserves finite limits. But $L_{y_{\cal C}\circ G}$ is left adjoint to $D_{G}:[{\cal C}^{\textup{op}}, \Set] \to [{\cal D}^{\textup{op}}, \Set]$. On the other hand, the adjunction $F\dashv G$ entails an adjunction $D_{F}:[{\cal D}^{\textup{op}}, \Set] \to [{\cal C}^{\textup{op}}, \Set] \dashv D_{G}:[{\cal C}^{\textup{op}}, \Set] \to [{\cal D}^{\textup{op}}, \Set]$, so, by the uniqueness (up to isomorphisms) of adjoints, we conclude that $L_{y_{\cal C}\circ G}\cong D_{F}$; in particular, it follows that $L_{y_{\cal C}\circ G}$ preserves arbitrary limits. So the geometric morphism $[{\cal C}^{\textup{op}}, \Set]\to [{\cal D}^{\textup{op}}, \Set]$ induced by $G$ is essential. The fact that the functor $y_{\cal C}\circ G$ is flat was also proved as Lemma 3(i) at p. 410 of \cite{MM}.  
	
	(iii) This follows immediately from (i) in light of (ii). Indeed, since $G$ is a morphism of sites $({\cal D}, T_{\cal D})\to ({\cal C}, T_{\cal C})$, it is a morphism of sites  $({\cal D}, K) \to ({\cal C}, J)$ if and only if it is cover-preserving.
	
	(iv) By definition of geometric morphisms induced by a morphism or by a comorphism of sites, the inverse image of $\Sh(G)$ is the functor $a_{K}\circ \textup{Lan}_{G^{\textup{op}}}\circ i_{J}$, where $\textup{Lan}_{G^{\textup{op}}}$ is the left Kan extension functor along $G^{\textup{op}}$, while the inverse image of $C_{F}$ is $a_{K}\circ D_{F}\circ i_{J}$. Now, since $F\dashv G$, $D_{F}\dashv D_{G}$, whence $\textup{Lan}_{G^{\textup{op}}}\cong D_{F}$, as required.  
\end{proofs}

\begin{remark}
	Proposition \ref{propadjointinduction}(iii)-(iv) generalize Lemma C2.5.1 and Scholium C2.5.2 \cite{El}.	   
\end{remark}

\subsubsection{From morphisms to comorphisms of sites}

The idea, for turning a morphism of sites into a comorphism of sites, is to replace the original domain site with a site related to it by a weakly dense morphism such that the composite of the given morphism of sites with it admits a left adjoint; this left adjoint will then be a comorphism of sites inducing the same geometric morphism (by Proposition \ref{propadjointinduction}).

We shall denote by $(F\downarrow G)$, for two functors $F:{\cal A}\to {\cal C}$ and $G:{\cal B}\to {\cal C}$, the comma category whose objects are the triplets $(a,b, \alpha)$ where $a\in {\cal A}$, $b\in {\cal B}$ and $\alpha$ is an arrow $F(a)\to G(b)$ in $\cal C$ (and whose arrows are defined in the obvious way).

We will consider in particular comma categories of the form $(1_{\cal C}\downarrow F)$ and $(F\downarrow 1_{\cal C})$. Given a functor $F:{\cal C}\to {\cal D}$ as in the following theorem, the objects of $(1_{\cal D}\downarrow F)$ are triplets of the form $(d, c, \alpha:d\to F(c))$ where $c\in {\cal C}$, $d\in {\cal D}$ and $\alpha$ is an arrow in $\cal D$.

\begin{theorem}\label{thmmorphismtocomorphism}
	Let $F:({\cal C}, J)\to ({\cal D}, K)$ be a morphism of small-generated sites. Let $i_{F}$ be the functor ${\cal C}\to (1_{\cal D}\downarrow F)$ sending any object $c$ of $\cal C$ to the triplet $(F(c), c, 1_{F(c)})$ (and acting on arrows in the obvious way), and $\pi_{\cal C}: (1_{\cal D}\downarrow F) \to {\cal C}$ and $\pi_{\cal D}: (1_{\cal D}\downarrow F) \to {\cal D}$ the canonical projection functors. Let $\tilde{K}$ be the Grothendieck topology on $ (1_{\cal D}\downarrow F)$ whose covering sieves are those whose image under $\pi_{\cal D}$ is $K$-covering. Then
	\begin{enumerate}[(i)]
		\item $\pi_{\cal C} \dashv i_{F}$, $\pi_{\cal D}\circ i_{F}=F$, $i_{F}$ is a morphism of sites $({\cal C}, J)\to ((1_{\cal D}\downarrow F), \tilde{K})$ and $c_{F}:=\pi_{\cal C}$ is a comorphism of sites $((1_{\cal D}\downarrow F), \tilde{K}) \to ({\cal C}, J)$; 	
		
		\item $\pi_{\cal D}: ((1_{\cal D}\downarrow F), \tilde{K}) \to ({\cal D}, K)$ is both a ($\tilde{K}$-full, $K$-dense and cover-reflecting) morphism of sites and a comorphism of sites inducing equivalences
		\[
		C_{\pi_{\cal D}}:\Sh((1_{\cal D}\downarrow F), \tilde{K}) \to \Sh({\cal D}, K)
		\]
		and
		\[
		\Sh(\pi_{\cal D}):\Sh({\cal D}, K) \to \Sh((1_{\cal D}\downarrow F), \tilde{K})
		\] 
		which are quasi-inverse to each other and make the following triangle commute:
		\begin{equation*}
		\begin{tikzcd}
		\Sh((1_{\cal D}\downarrow F), \tilde{K}) \arrow[rr, "C_{\pi_{\cal D}}", yshift=1ex] \ar[rr, phantom, "\sim"{yshift=-0.25ex}]  \arrow[dr, "C_{\pi_{\cal C}}\cong \Sh(i_{F})"{below, xshift=-5ex}]   & &  \Sh({\cal D}, K) \ar[ll, "\Sh(\pi_{\cal D})", yshift=-1ex] \ar[ld, "\Sh(F)", yshift=0.5ex, xshift=-1ex]            \\
		& \Sh({\cal C}, J)  & 
		\end{tikzcd}
		\end{equation*}	
	\end{enumerate}	
\end{theorem}

\begin{proofs}
	(i) It is immediate to see that $\pi_{\cal C} \dashv i_{F}$. Let us show that $\pi_{\cal D}: ((1_{\cal D}\downarrow F), T) \to ({\cal D}, K)$ is a morphism of sites, where $T$ is the trivial Grothendieck topology on $(1_{\cal D}\downarrow F)$. Notice that this will ensure that the collection of sieves $\tilde{K}$ defined in the statement of the theorem is indeed a Grothendieck topology, since it will be the Grothendieck topology induced by the morphism of sites $\pi_{\cal D}: ((1_{\cal D}\downarrow F), T) \to ({\cal D}, K)$ (in the sense of Proposition \ref{reminducedtopology}).  
	
	We will deduce that $\pi_{\cal D}:((1_{\cal D}\downarrow F), T) \to ({\cal D}, K)$ satisfies properties (ii), (iii) and (iv) in the definition of morphism of sites (that is, Definition \ref{defmorphismsites}) from the fact that $F$ does. 
	
	Let us start with condition (ii). Since $F:({\cal C}, J) \to ({\cal D}, K)$ satisfies it, every object $d$ of ${\mathcal D}$ admits a $K$-covering family
	$$
	d_i \longrightarrow d \, , \quad i \in I \, ,
	$$
	by objects $d_i$ of ${\mathcal D}$ which have morphisms 
	$$
	h_{i}:d_i \longrightarrow F(c'_i)
	$$
	to the images under $F$ of objects $c'_i$ of ${\mathcal C}$. Then we have, for each $i\in I$, an object $(c'_i, d_{i}, h_{i}:d_i \longrightarrow F(c'_i))$ of the category $(1_{\cal D}\downarrow F)$ and morphisms 
	\[
	d_{i}\to \pi_{\cal D}((c'_i, d_{i}, h_{i}:d_i \longrightarrow F(c'_i)))=d_{i}
	\]
	given by the identity arrows on the $d_{i}$. So $\pi_{\cal D}$ satisfies condition (ii) in the definition of morphism of sites.
	
	Let us now turn to condition (iii). 
	Given objects $(c_{1}, d_{1}, \alpha_{1}:d_{1}\to F(c_{1}))$ and $(c_{2}, d_{2}, \alpha_{2}:d_{2}\to F(c_{2}))$ and arrows $k_{1}:d\to \pi_{\cal D}((c_{1}, d_{1}, \alpha_{1}:d_{1}\to F(c_{1})))=d_{1}$ and $k_{2}:d\to \pi_{\cal D}((c_{2}, d_{2}, \alpha_{2}:d_{2}\to F(c_{2})))=d_{2}$ in $\cal D$, consider the arrows $g_{1}:=\alpha_{1}\circ k_{1}:d\to F(c_{1})$ and $g_{2}:=\alpha_{2}\circ k_{2}:d\to F(c_{2})$. Then, since $F$ is a morphism of sites, there is a $K$-covering family
	$$
	g'_i : d_i \longrightarrow d \, , \quad i \in I \, ,
	$$
	and a family of pairs of morphisms of ${\mathcal C}$
	$$
	f_1^i : c'_i \longrightarrow c_1 \, , \quad f_2^i : c'_i \to c_2 \, , \quad i \in I \,  ,
	$$
	and of morphisms of ${\mathcal D}$
	$$
	h_i : d_i \longrightarrow F(c'_i) \, , \quad i \in I \, ,
	$$
	making the following squares commutative:
	$$
	\xymatrix{
		d_i \ar[r]^{g'_i} \ar[d]_{h_i} & d \ar[d]^{g_1} \\
		F(c'_i) \ar[r]^{F(f_1^i)} &F(c_1)
	} 
	\qquad \qquad \qquad
	\xymatrix{
		d_i \ar[r]^{g'_i} \ar[d]_{h_i} & d \ar[d]^{g_2} \\
		F(c'_i) \ar[r]^{F(f_2^i)} &F(c_2)
	} 
	$$
	
	These yield objects $(c'_i, d_{i}, h_{i}:d_{i}\to F(c'_i))$ of the category $1_{\cal D}\downarrow F$ (for $i\in I$) and arrows $$(f_1^i, k_{1}\circ g'_i):(c'_i, d_{i}, h_{i}:d_{i}\to F(c'_i)) \to (c_{1}, d_{1}, \alpha_{1}:d_{1}\to F(c_{1}))$$ and $$(f_2^i, k_{2}\circ g'_i):(c'_i, d_{i}, h_{i}:d_{i}\to F(c'_i)) \to (c_{2}, d_{2}, \alpha_{1}:d_{2}\to F(c_{2}))$$ and morphisms of $\cal D$
	\[
	1_{d_{i}}:d_{i}\to \pi_{\cal D}((c'_i, d_{i}, h_{i}:d_{i}\to F(c'_i)))=d_{i}
	\]   
	(for each $i\in I$)
	making the following squares commutative:
	
	$$
	\xymatrix{
		d_i \ar[r]^{g'_i} \ar[d]_{1_{d_{i}}} & d \ar[d]^{k_1} \\
		d_{i} \ar[r]^{ k_{1}\circ g'_i} & d_{1}
	} 
	\qquad \qquad \qquad
	\xymatrix{
		d_i \ar[r]^{g'_i} \ar[d]_{1_{d_{i}}} & d \ar[d]^{k_2} \\
		d_{i} \ar[r]^{ k_{2}\circ g'_i} & d_{2}
	} 
	$$
	So $\pi_{\cal D}$ satisfies condition (iii) in the definition of morphism of sites. 
	
	Lastly, let us consider condition (iv). Let $(f_{1}, k_{1}), (f_{2}, k_{2}):(c_{1}, d_{1}, \alpha_{1}:d_{1}\to F(c_{1})) \rightrightarrows (c_{2}, d_{2}, \alpha_{2}:d_{2}\to F(c_{2}))$ be a pair of arrows in the category $(1_{\cal D}\downarrow F)$ and $z:d\to \pi_{\cal D}((c_{1}, d_{1}, \alpha_{1}:d_{1}\to F(c_{1})))=d_{1}$ be an arrow of $\cal D$ satisfying 
	\[
	k_{1}\circ z=k_{2}\circ z.
	\]
	Let $g$ be the arrow $\alpha_{1}\circ z$. We have $F(f_{1})\circ g=F(f_{2})\circ g$. Indeed, $F(f_{1})\circ g=F(f_{1})\circ \alpha_{1}\circ z=\alpha_{2}\circ k_{1}\circ z=\alpha_{2}\circ k_{2}\circ z=F(f_{2})\circ \alpha_{1}\circ z=F(f_{2})\circ g$. So the fact that $F$ satisfies condition (iv) in the definition of morphism of sites ensures that there exist a $K$-covering family
	$$
	g_i : d_i \longrightarrow d \, , \quad i \in I \, ,
	$$
	and a family of morphisms of ${\mathcal C}$
	$$
	w_i : c'_i \longrightarrow c_{1} \, , \quad i \in I \, ,
	$$
	satisfying
	$$
	f_1 \circ w_i = f_2 \circ w_i \, , \quad \forall \, i \in I \, ,
	$$
	and of morphisms of ${\mathcal D}$
	$$
	h_i : d_i \longrightarrow F(c'_i) \, , \quad i \in I \, ,
	$$
	making the following squares commutative:
	$$
	\xymatrix{
		d_i \ar[r]^{g_i} \ar[d]_{h_i} & d \ar[d]^{g} \\
		F(c'_i) \ar[r]^{F(w_i)} & F(c_{1})
	} 
	$$ 
	Then, for each $i\in I$, we have an object $(c'_i, d_{i}, h_{i}:d_i \to F(c'_i))$ of the category $(1_{\cal D}\downarrow F)$, and the commutativity of the above square ensures that $(w_{i}, z\circ g_{i})$ defines an arrow
	\[
	(w_{i}, z\circ g_{i}):(c'_i, d_{i}, h_{i}:d_i \to F(c'_i)) \to (c_{1}, d_{1}, \alpha_{1}:d_{1}\to F(c_1))
	\]
	in $(1_{\cal D}\downarrow F)$. These arrows satisfy the property 
	\[
	(f_{1}, k_{1})\circ (w_{i}, z\circ g_{i}) = (f_{2}, k_{2})\circ (w_{i}, z\circ g_{i})
	\]
	and the identity morphisms
	\[
	1_{d_{i}}:d_{i} \to \pi_{\cal D}((c'_i, d_{i}, h_{i}:d_i \to F(c'_i)))=d_{i}
	\]
	make the following squares commutative:
	$$
	\xymatrix{
		d_i \ar[r]^{g_i} \ar[d]_{1_{d_{i}}} & d \ar[d]^{z} \\
		d_{i} \ar[r]^{z\circ g_{i}} & d_{1}
	} 
	$$ 
	
	This shows that $\pi_{\cal D}$ satisfies condition (iv) in the definition of morphism of sites.
	
	The functor $\pi_{\cal D}$ is also a comorphism of sites $((1_{\cal D}\downarrow F), \tilde{K}) \to ({\cal D}, K)$. Indeed, given an object $(c, d, \alpha:d\to F(c))$ of $(1_{\cal D}\downarrow F)$ and a $K$-covering sieve $T$ on $\pi_{\cal D}((c, d, \alpha:d\to F(c)))=d$, the sieve generated by the family $\{(1_{c}, g):(c, \textup{dom}(g), \alpha\circ g: \textup{dom}(g) \to F(c)) \to (c, d, \alpha:d\to F(c)) \mid g\in T\}$ is $\tilde{K}$-covering and satisfies the required property since its image under $\pi_{\cal D}$ is precisely $T$.  
	
	Next, let us show that $\pi_{\cal D}$ is $\tilde{K}$-full and $K$-dense. This will imply, by Proposition \ref{propmorphismcomorphismequivalence}, that it induces an equivalence of toposes $\Sh((1_{\cal D}\downarrow F), \tilde{K}) \simeq \Sh({\cal D}, K)$ one half of which is the morphism $C_{\pi_{\cal D}}$ and whose other half is $\Sh(\pi_{\cal D})$.
	
	To show that $\pi_{\cal D}$ is $\tilde{K}$-full, let us consider an arrow $g$ in $\cal D$ from $\pi_{\cal D}((c_{1}, d_{1}, \alpha_{1}:d_{1}\to F(c_{1})))=d_{1}$ to $\pi_{\cal D}((c_{2}, d_{2}, \alpha_{2}:d_{2}\to F(c_{2})))=d_{2}$. By applying condition (iii) in the definition of morphism of sites for $F$ to the arrows $\alpha_{1}:d_{1}\to F(c_{1})$ and $\alpha_{2}\circ g:d_{1}\to F(c_{2})$ we obtain a $K$-covering family
	$$
	g'_i : d_i \longrightarrow d_1 \, , \quad i \in I \, ,
	$$
	and a family of pairs of morphisms of ${\mathcal C}$
	$$
	f_1^i : c'_i \longrightarrow c_1 \, , \quad f_2^i : c'_i \to c_2 \, , \quad i \in I \,  ,
	$$
	and of morphisms of ${\mathcal D}$
	$$
	h_i : d_i \longrightarrow F(c'_i) \, , \quad i \in I \, ,
	$$
	making the following squares commutative:
	$$
	\xymatrix{
		d_i \ar[r]^{g'_i} \ar[d]_{h_i} & d_1 \ar[d]^{\alpha_{1}} \\
		F(c'_i) \ar[r]^{F(f_1^i)} &F(c_1)
	} 
	\qquad \qquad \qquad
	\xymatrix{
		d_i \ar[r]^{g'_i} \ar[d]_{h_i} & d_1 \ar[d]^{\alpha_{2}\circ g} \\
		F(c'_i) \ar[r]^{F(f_2^i)} &F(c_2)
	} 
	$$
	
	Now, the commutativity of these squares ensures that we have, for each $i\in I$, arrows
	\[
	(f_1^i, g'_i):(c'_i, d_i, h_i:d_i \to F(c'_i)) \to (c_{1}, d_{1}, \alpha_{1}:d_{1}\to F(c_{1}))
	\]
	and
	\[
	(f_2^i, g\circ g'_i):(c'_i, d_i, h_i:d_i \to F(c'_i)) \to (c_{2}, d_{2}, \alpha_{2}:d_{2}\to F(c_{2}))
	\]
	arrows in the category $1_{\cal D}\downarrow F$. These arrows clearly satisfy the condition
	\[
	g\circ \pi_{\cal D}((f_1^i, g'_i))=\pi_{\cal D}((f_2^i, g\circ g'_i)),
	\]
	and the family $\{(f_1^i, g'_i) \mid i\in I\}$ is $\tilde{K}$-covering (since $\{ g'_i \mid i\in I \}$ is $K$-covering). This shows that $\pi_{\cal D}$ is $\tilde{K}$-full, as desired. 
	
	Let us now prove that $\pi_{\cal D}$ is $K$-dense. Given $d\in {\cal D}$, using the fact that $F$ satisfies condition (ii) in the definition of morphism of sites, we obtain a $K$-covering family 
	$$
	d_i \longrightarrow d \, , \quad i \in I \, ,
	$$
	by objects $d_i$ of ${\mathcal D}$ which have morphisms 
	$$
	h_{i}:d_i \longrightarrow F(c'_i)
	$$
	to the images under $F$ of objects $c'_i$ of ${\mathcal C}$. Since we can write $d_{i}=\pi_{\cal D}(c'_i, d_{i}, h_{i}:d_i \longrightarrow F(c'_i))$, it follows that $d$ can be $K$-covered by objects in the image of the functor $\pi_{\cal D}$, which proves that this functor is $K$-dense, as required.
	
	Let us now show that $\pi_{\cal C}$ is a comorphism of sites $((1_{\cal D}\downarrow F), \tilde{K}) \to ({\cal C}, J)$. Given a sieve $S$ on an object of the form $\pi_{\cal C}(c, d, \alpha:d\to F(c))=c$, since $F$ is cover-preserving, the family $\{F(f) \mid f\in S\}$ is $K$-covering on $F(c)$ and hence the pullback $\alpha^{\ast}(R)$ along $\alpha$ of the sieve $R$ generated by it is $K$-covering on $d$. So, for any $g\in \alpha^{\ast}(R)$, there is an element $i_{g}\in I$ and an arrow $\beta_{g}:\textup{dom}(g)\to F(\textup{dom}(f_{i_{g}}))$ such that $\alpha \circ g=F(\textup{dom}(f_{i_{g}}))\circ \beta_{g}$. So the sieve generated by the family of arrows 
	\[
	(f_{i_{g}}, g):(\textup{dom}(f_{i_{g}}), \textup{dom}(g), \beta_{g}: \textup{dom}(g) \to F(\textup{dom}(f_{i_{g}}))) \to (c, d, \alpha:d\to F(c))
	\] 
	satisfies the required property since its image under $\pi_{\cal C}$ is contained in $S$.  
	
	Since $\pi_{\cal C}$ is a comorphism of sites $((1_{\cal D}\downarrow F), \tilde{K}) \to ({\cal C}, J)$ and $\pi_{\cal C} \dashv i_{F}$, it thus follows from Proposition \ref{propadjointinduction}(iii) that $i_{F}$ is a morphism of sites $({\cal C}, J)\to ((1_{\cal D}\downarrow F), \tilde{K})$.

\end{proofs}

\begin{remark}\label{remcorrelationmorphisms}
	One can prove that the property of $F$ being a morphism of sites $({\cal C}, T') \to ({\cal D}, K)$ (where $T'$ is the trivial topology on $\cal C$) is not only a sufficient but also a necessary condition for $\pi_{\cal D}$ to be a morphism of sites $((1_{\cal C}\downarrow F), T) \to ({\cal D}, K)$. 
\end{remark}

\subsubsection{From comorphisms to morphisms of sites}

The following result represents a sort of \acc dual' of Theorem \ref{thmmorphismtocomorphism}, allowing to turn arbitrary comorphisms of sites into morphisms of sites.

Below, we shall abbreviate by $\hat{\cal D}$ the category of presheaves on a small category $\cal D$.

\begin{theorem}\label{thmcomorphismtomorphism}
	Let $F:({\cal D}, K)\to ({\cal C}, J)$ be a comorphism of small-generated sites. Let $\pi_{\cal C}':(F \downarrow 1_{\cal C})\to {\cal C}$ and  $\pi_{\cal D}':(F \downarrow 1_{\cal C})\to {\cal D}$ be the canonical projection functors and $j_{F}:{\cal D}\to (F \downarrow 1_{\cal C})$ the functor sending any object $d$ of $\cal D$ to the triplet $(d, F(d), 1_{F(d)})$. Let $\overline{K}$ be the Grothendieck topology on $(F \downarrow 1_{\cal C})$ whose covering families are those which are sent by $\pi_{\cal D}'$ to $K$-covering families. Then
	
	\begin{enumerate}[(i)]
		\item $j_{F}\dashv \pi_{\cal D}'$, $\pi_{\cal C}'\circ j_{F}=F$, $\pi_{\cal C}'$ is a comorphism of sites $(F \downarrow 1_{\cal C}, \overline{K})\to ({\cal C}, J)$ and $j_{F}$ is a (full and faithful) comorphism and dense morphism of sites $({\cal D}, K) \to (F \downarrow 1_{\cal C}, \overline{K})$; 
		
		\item $\pi_{\cal D}'$ is both a morphism and a comorphism of sites  $((F \downarrow 1_{\cal C}), \overline{K})\to ({\cal D}, K)$
		inducing equivalences
		\[
		C_{\pi_{D}'}:\Sh((F \downarrow 1_{\cal C}), \overline{K})\to \Sh({\cal D}, K)
		\]  
		and
		\[
		\Sh(\pi_{D}'):\Sh({\cal D}, K) \to \Sh((F \downarrow 1_{\cal C}), \overline{K})
		\] 
		which are quasi-inverse to each other and make the following triangle commute:
		\begin{equation*}
		\begin{tikzcd}
		\Sh((F \downarrow 1_{\cal C}), \overline{K}) \arrow[rr, "C_{\pi_{\cal D}'}\cong \Sh(j_{F})", yshift=1ex] \ar[rr, phantom, "\sim"{yshift=-0.25ex}]  \arrow[dr, "C_{\pi_{\cal C}'}"{below, xshift=-2ex}]   & &  \Sh({\cal D}, K) \ar[ll, "\Sh(\pi_{\cal D}')\cong C_{j_{F}}", yshift=-1ex] \ar[ld, "C_{F}", yshift=0.5ex, xshift=-1ex]            \\
		& \Sh({\cal C}, J)  & 
		\end{tikzcd}
		\end{equation*}	
		
		\item With $F$ we can associate a morphism of sites $m_{F}:({\cal C}, J)\to (\hat{\cal D}, \hat{K})$, where $m_{F}$ is the functor sending an object $c$ of $\cal C$ to the presheaf $\textup{Hom}_{\cal C}(F(-), c)$ and $\hat{K}$ is the extension of the Grothendieck topology $K$ along the Yoneda embedding ${\cal D}\rightarrow \hat{\cal D}$ (in the sense of section \ref{sec:presheaflifting}), which induces a geometric morphism $\Sh(m_{F})$ making the following triangle commute:  
		\begin{equation*}
		\begin{tikzcd}
		\Sh(\hat{D}, \hat{K}) \arrow[rr, "\Sh(y_{\cal D})", yshift=1ex] \ar[rr, phantom, "\sim"{yshift=-0.25ex}]  \arrow[dr, "\Sh(m_{F})"{below, xshift=-2ex}]   & &  \Sh({\cal D}, K) \ar[ll, "C_{y_{\cal D}}", yshift=-1ex] \ar[ld, "C_{F}", yshift=0.5ex, xshift=-1ex]            \\
		& \Sh({\cal C}, J)  & 
		\end{tikzcd}
		\end{equation*}	
	\end{enumerate}
\end{theorem}

\begin{proofs}
	(i) It is clear that $j_{F}\dashv \pi_{\cal D}'$ and that $\pi_{\cal C}'\circ j_{F}=F$. Since $\pi_{\cal D}'$ is cover-preserving, it follows from Proposition \ref{propadjointinduction}(i) that $j_{F}$ is a comorphism of sites $({\cal D}, K) \to ((F \downarrow 1_{\cal C}), \overline{K})$.	The fact that $\pi_{\cal C}'$ is a comorphism of sites $((F \downarrow 1_{\cal C}), \overline{K})\to ({\cal C}, J)$ follows from the fact that $F:({\cal D}, K)\to ({\cal C}, J)$ is. Indeed, given an object $(c,d, \alpha:F(d)\to c)$ of the category $(F \downarrow 1_{\cal C})$ and a $J$-covering sieve $S$ on $\pi_{\cal C}'((d, c, \alpha:F(d)\to c))=c$, $\alpha^{\ast}(S)$ is $J$-covering on $F(d)$ and hence contains the image under $F$ of a $K$-covering sieve $R$ on $d$. We thus have, for each $g\in R$, an arrow $f_{g}:F(\textup{dom}(g))\to c$ in $S$ such that $\alpha \circ F(g)=f_{g}$. The family $\{ (g, f_{g}):(\textup{dom}(g), F(\textup{dom}(g)), 1_{F(\textup{dom}(g))} ) \to (d, c, \alpha:F(d)\to c) \mid g\in R\}$ is therefore $\overline{K}$-covering and its image under $\pi_{\cal C}'$ is contained in $S$. It is clear that $j_{F}$ is full and faithful and cover-preserving, so to conclude that it is a dense morphism of sites, it suffices to show that it is $\overline{K}$-dense. But this follows from the following observation: for any object $(d,c, \alpha:F(d)\to c)$ of the category $(F \downarrow 1_{\cal C})$, the arrow $(1_{d}, \alpha):(d, F(d), 1_{F(d)}:F(d)\to F(d)) \to (d,c, \alpha:F(d)\to c)$ is $\overline{K}$-covering. 
	
	(ii) The fact that $\pi_{\cal D}'$ is a morphism of sites $((F \downarrow 1_{\cal C}), \overline{K})\to ({\cal D}, K)$ follows from the fact that its left adjoint $j_{F}$ is by (i) a comorphism of sites $({\cal D}, K) \to ((F \downarrow 1_{\cal C}), \overline{K})$ in light of Proposition \ref{propadjointinduction}(iii). Let us show that it is also a comorphism of sites. Given an object $(c,d, \alpha:F(d)\to c)$ of the category $(F \downarrow 1_{\cal C})$ and a $K$-covering sieve $T$ on $d$, the family $\{(\textup{dom}(g), 1_{c}, \alpha\circ F(g):F(\textup{dom}(g))\to c \mid g\in T\}$ is $\overline{K}$-covering and its image under $\pi_{\cal D}'$ is contained in $T$. By Proposition \ref{propmorphismcomorphismequivalence}, to show that $\pi_{\cal D}'$ induces an equivalence of toposes as in the statement of the theorem, it suffices to show that it is $K$-dense and $\overline{K}$-full. Since $\pi_{\cal D}'$ is surjective (for any $d\in {\cal D}$, $d=\pi_{\cal D}'((d, F(d), 1_{F(d)}))$), it is \emph{a fortiori} $K$-dense. To show that $\pi_{\cal D}'$ is $\overline{K}$-full, consider two objects $(d_{1}, c_{1}, \alpha_{1}:F(d_{1})\to c_{1})$ and $(d_{2}, c_{2}, \alpha_{2}:F(d_{2})\to c_{2})$ of the category $(F \downarrow 1_{\cal C})$ and an arrow $g:d_{1}\to d_{2}$ in $\cal D$. Then we have a $\overline{K}$-covering arrow $$(1_{d_{1}}, \alpha_{1}):(d_{1}, F(d_{1}), 1_{F(d_{1})}:F(d_{1})\to F(d_{1}))\to (d_{1}, c_{1}, \alpha_{1}:F(d_{1})\to c_{1})$$ (its image under $\pi_{\cal D}'$ is the identity arrow on $d_{1}$) and an arrow $$(g, \alpha_{2}\circ F(g)):(d_{1}, F(d_{1}), 1_{F(d_{1})}:F(d_{1})\to F(d_{1}))\to (d_{2}, c_{2}, \alpha_{2}:F(d_{2})\to c_{2})$$ in $(F \downarrow 1_{\cal C})$  satisfying the condition 
	\[
	\pi_{\cal D}'((g, \alpha_{2}\circ F(g)))=g\circ \pi_{\cal D}'((1_{d_{1}}, \alpha_{1})).
	\]   
	This shows that $\pi_{\cal D}'$ is $\overline{K}$-full, as desired.
	
	(iii) We know that the functor $m_{F}:{\cal C}\to [{\cal D}^{\textup{op}}, \Set]$ sending an object $c$ of $\cal C$ to the presheaf $\textup{Hom}_{\cal D}(F(-), c)$ is flat and that its composite with the associated sheaf functor $a_{K}:[{\cal D}^{\textup{op}}, \Set]\to \Sh({\cal D}, K)$ is $J$-continuous. So we can regard $m_{F}$ as a morphism of sites $({\cal C}, J)\to (\hat{\cal D}, \hat{K})$. The fact that the triangle in the statement of the theorem commutes follows from  Proposition \ref{propyonedamorcomor} and the fact that $a_{K}\circ m_{F}$ is the flat functor corresponding via Diaconescu's equivalence to the geometric morphism $C_{F}$.
\end{proofs}

\begin{remarks}\label{remcorrelationcomorphisms}
	\begin{enumerate}[(a)]
		\item The property that $F$ be a comorphism of sites $({\cal D}, K)\to ({\cal C}, J)$ is not only a necessary but also a sufficient condition for $\pi_{\cal C}'$ to be a comorphism of sites. Notice the analogy with Remark \ref{remcorrelationmorphisms}. 
		
		\item The canonical projection functor $\pi_{\cal C}':(F \downarrow 1_{\cal C}) \to {\cal C}$ is a split fibration, corresponding to the internal category in $[{\cal C}^{\textup{op}}, \Set]$ given by $c \mapsto (F(c)\downarrow {\cal D})$.
		
		\item If $F:{\cal D}\to {\cal C}$ has a right adjoint $G:{\cal C}\to {\cal D}$ then $G$ is a morphism of sites $({\cal C}, J)\to ({\cal D}, K)$ and $m_{F}\cong y_{\cal D}\circ G$; in other words, there is no need in this case to go to the presheaf topos $\hat{\cal D}$ in order to turn $F$ into a morphism of sites inducing the same morphism, since $C_{F}$ is isomorphic to $\Sh(G)$ (cf. Proposition \ref{propadjointinduction}). 
	\end{enumerate}	
\end{remarks}

We observe in passing that not only every comorphism of sites $({\cal D}, K)\to ({\cal C}, J)$ induces a morphism of sites $({\cal C}, J)\to (\hat{\cal D}, \hat{K})$, but the same is true for any comorphism of sites $({\cal D}, K)\to (\hat{\cal C}, \hat{J})$. Notice also that every comorphism of sites $({\cal D}, K)\to (\hat{\cal C}, \hat{J})$ can be extended along the Yoneda embedding $y_{\cal D}$ to a unique colimit-preserving comorphism of sites $(\hat{\cal D}, \hat{K})\to (\hat{\cal C}, \hat{J})$; this morphism is precisely the left Kan extension functor along $F$ which has as right adjoint the functor $F^{\ast}$, which is therefore a morphism of sites $(\hat{\cal C}, \hat{J}) \to (\hat{\cal D}, \hat{K})$; in fact, the composite of this functor with the Yoneda embedding is precisely the morphism of sites $m_{F}$ corresponding to $F$.  

\begin{theorem}\label{thmtoposadjunction}
	Let $({\cal C}, J)$ and $({\cal D}, K)$ be small-generated sites.
	
	\begin{enumerate}[(i)]
		\item Let $F:({\cal C}, J)\to ({\cal D}, K)$ be a morphism of sites, with corresponding comorphism of sites $c_{F}:((1_{\cal D}\downarrow F), \tilde{K}) \to ({\cal C}, J)$ as in Theorem \ref{thmmorphismtocomorphism}. Let $\pi_{\cal D}:((1_{\cal D} \downarrow F), \tilde{K})\to ({\cal D}, K)$ be the canonical projection functor, and let 
		\[
		w_{F}: (1_{\cal D}\downarrow F)\to (c_{F}\downarrow 1_{\cal D})
		\]
		be the functor $j_{c_{F}}$, sending an object $A$ of $(1_{\cal D}\downarrow F)$ to the object $(A, c_{F}(A), 1_{c_{F}(A)}:c_{F}(A)\to c_{F}(A))$. Then $w_{F}$ is both a (full and faithful) comorphism and a dense morphism of sites $((1_{\cal D}\downarrow F),  \tilde{K}) \to ((c_{F}\downarrow 1_{\cal D}), \overline{\tilde{K}})$ satisfying the relation $\pi_{\cal D}'''\circ w_{F}=\pi_{\cal D}$ and inducing an equivalence relating $F$ and $c_{F}$, which makes the following diagram commute (where $\pi_{\cal D}'''$ denotes the canonical projection functor $(c_{F}\downarrow 1_{\cal D}) \to {\cal D}$):
		\begin{equation*}
		\begin{tikzcd} 
		\Sh({\cal D}, K) \ar[d, "\Sh(\pi_{{\cal D}'''})", xshift=1ex] \arrow[rr, "", yshift=1ex] \ar[rr, phantom, "="{yshift=-0.25ex}] & & \Sh({\cal D}, K) \ar[ll, "", yshift=-1ex] \ar[d, "\Sh(\pi_{\cal D})", xshift=1ex]  \\
		\Sh((c_{F}\downarrow 1_{\cal D}), \overline{\tilde{K}})  \ar[u, "C_{\pi_{\cal D}'''}", xshift=-1ex] \ar[u, phantom, "\rotatebox{90}{$\sim$}"{xshift=0.25ex}]  \arrow[rr, "\Sh(w_{F})\cong C_{\pi_{(1_{\cal D}\downarrow F)}'}", yshift=1ex] \ar[rr, phantom, "\sim"{yshift=-0.25ex}]  \arrow[dr, "C_{\pi_{\cal C}'}"{below}]    & &  \Sh((1_{\cal D}\downarrow F), \tilde{K}) \ar[ll, "C_{w_{F}}\cong \Sh(\pi_{(1_{\cal D}\downarrow F)}') ", yshift=-1ex] \ar[u, "C_{\pi_{\cal D}}", xshift=-1ex ] \ar[u, phantom, "\rotatebox{90}{$\sim$}"{xshift=0.25ex}] \ar[ld, "C_{c_{F}}", yshift=0.5ex, xshift=-1ex]            \\
		& \Sh({\cal C}, J)  & 
		\end{tikzcd}
		\end{equation*}

		\item Let $G:({\cal D}, K) \to ({\cal C}, J)$ be a comorphism of sites, with corresponding morphism of sites $m_{G}:({\cal C}, J)\to (\hat{\cal D}, \hat{K})$ as in Theorem \ref{thmcomorphismtomorphism}. Let 
		\[
		z_{G}:(G\downarrow 1_{\cal C} )\to (1_{\hat{\cal D}}\downarrow m_{G})
		\]
		be the functor sending any object $(d,c, \alpha:G(d)\to c)$ of $(G\downarrow 1_{\cal C} )$ to the object $(y_{\cal D}(d), c, \overline{\alpha}:y_{\cal D}(d) \to m_{G}(c))$ of $(1_{\hat{\cal D}}\downarrow m_{G})$, where $\overline{\alpha}$ is the arrow corresponding to the element $\alpha$ of $m_{G}$ via the Yoneda Lemma. Then $z_{G}$ is both a (full and faithful) comorphism and a dense morphism of sites $((G\downarrow 1_{\cal C}), \overline{K} )\to ((1_{\hat{\cal D}}\downarrow m_{G}), \widetilde{\hat{K}})$ satisfying the relation $\pi_{\hat{\cal D}}\circ z_{G}=y_{\cal D}\circ \pi_{\cal D}'$ and inducing an equivalence relating $G$ and $m_{G}$, which makes the following diagram commute:
		\begin{equation*}
		\begin{tikzcd} 
		\Sh(\hat{\cal D}, \hat{K}) \ar[d, "\Sh(\pi_{\hat{\cal D}})", xshift=1ex] \arrow[rr, "\Sh(y_{\cal D})", yshift=1ex] \ar[rr, phantom, "\sim"{yshift=-0.25ex}] & & \Sh({\cal D}, K) \ar[ll, "C_{y_{\cal D}}", yshift=-1ex] \ar[d, "\Sh(\pi_{{\cal D}}')\cong C_{j_{G}}", xshift=1ex]  \\
		\Sh((1_{\hat{\cal D}}\downarrow m_{G}), \widetilde{\hat{K}}) \ar[u, "C_{\pi_{\hat{\cal D}}}", xshift=-1ex ] \ar[u, phantom, "\rotatebox{90}{$\sim$}"{xshift=0.25ex}] \arrow[rr, "\Sh(z_{G})", yshift=1ex] \ar[rr, phantom, "\sim"{yshift=-0.25ex}]  \arrow[dr, "C_{\pi_{\cal C}}\cong \Sh(i_{m_{G}})"{below, xshift=-5ex}]    & &  \Sh((G\downarrow 1_{\cal C}), \overline{K}) \ar[ll, "C_{z_{G}}", yshift=-1ex] \ar[u, "\Sh(j_{G})\cong C_{\pi_{{\cal D}'}}", xshift=-1ex ] \ar[u, phantom, "\rotatebox{90}{$\sim$}"{xshift=0.25ex}]  \ar[ld, "C_{\pi_{\cal C}'}", yshift=0.5ex, xshift=-1ex]            \\
		& \Sh({\cal C}, J)  & 
		\end{tikzcd}
		\end{equation*}	 
	\end{enumerate}   
	
\end{theorem}

\begin{proofs}
	(i) The fact that $w_{F}$ is both a morphism and a comorphism of sites	$((1_{\cal D}\downarrow F),  \tilde{K}) \to ((c_{F}\downarrow 1_{\cal D}), \overline{\tilde{K}})$ inducing equivalences which make the diagram in the statement of the theorem commute follows from the equality $\pi_{\cal D}'''\circ w_{F}=\pi_{\cal D}$ in light of Theorems \ref{thmcomorphismtomorphism} and \ref{thmmorphismtocomorphism}. The fact that $w_{F}=j_{c_{F}}$ is full and faithful and a dense morphism of sites follows from Theorem \ref{thmcomorphismtomorphism}. 
	
	(ii) The fact that $z_{G}$ is both a morphism and a comorphism of sites $((G\downarrow 1_{\cal C}), \overline{K})\to ((1_{\hat{\cal D}}\downarrow m_{G}), \widetilde{\hat{K}} )$ inducing equivalences making the diagram in the statement of the theorem commute follows from the equality $\pi_{\hat{\cal D}}\circ z_{G}=y_{\cal D}\circ \pi_{\cal D}'$ in light of Theorems \ref{thmcomorphismtomorphism} and \ref{thmmorphismtocomorphism}. It is clear that $z_{G}$ is full and faithful, so it remains to show that it is $\widetilde{\hat{K}}$-dense. But this easily follows from the fact that every presheaf on $\cal D$ is $\hat{K}$-covered by representables.	
\end{proofs} 

We shall call a functor which both a morphism and a comorphism of sites a \emph{bimorphism} of sites.

\begin{remark}
	Theorem \ref{thmtoposadjunction} shows that the relationship between a morphism $F$ (resp. comorphism $G$) of sites and the associated comorphism $c_{F}$ (resp. morphism $m_{F}$) of sites is captured by the equivalence
	\[
	\Sh((1_{\cal D}\downarrow F), \tilde{K})\simeq \Sh((c_{F}\downarrow 1_{\cal D}), \overline{\tilde{K}})
	\]
	(resp. 
	\[
	\Sh((G\downarrow 1_{\cal C}), \overline{K})\simeq \Sh((1_{\hat{\cal D}}\downarrow m_{G}), \widetilde{\hat{K}}))
	\]
	of toposes over $\Sh({\cal C}, J)$ induced by the bimorphism of sites $w_{F}$ (resp. $z_{G}$) over $\cal C$. Recall that two functors $H:{\cal A}\to {\cal B}$ and $K:{\cal B}\to {\cal A}$ are adjoint to each other ($K$ on the left and $H$ on the right) if and only if the comma categories $(K\downarrow 1_{\cal A})$ and $(1_{\cal B}\downarrow H)$ are equivalent over the product ${\cal A}\times {\cal B}$. Our theorem then tells us that $F$ and $c_{F}$ are not adjoint to each other in a concrete sense (that is, at the level of sites), since they are not defined between a pair of categories, nor the categories $(1_{\cal D}\downarrow F)$ and $(c_{F}\downarrow 1_{\cal D})$ (resp. the categories $(G\downarrow 1_{\cal C})$ and $(1_{\hat{\cal D}}\downarrow m_{G})$) are equivalent in general; nonetheless, they become \ac abstractly' adjoint in the world of toposes since toposes naturally attached to such categories are equivalent. This is a compelling illustration of the philosophy of toposes as \ac bridges'.   
\end{remark}

\subsubsection{The dual adjunction between morphisms and comorphisms of sites}

Recall that every morphism $F:({\cal C}, J)\to ({\cal D}, K)$ of small-generated sites induces a geometric morphism 
$\Sh(F):\Sh({\cal D}, K) \to \Sh({\cal C}, J)$ and that every comorphism of sites $G:({\cal D}, K) \to ({\cal C}, J)$ induces a geometric morphism $C_{G}:\Sh({\cal D}, K) \to \Sh({\cal C}, J)$. In fact, the constructions $\Sh$ and $C$ define functors from the category of small-generated sites and morphisms (resp. comorphisms) between them to the category of Grothendieck toposes.

Let us define the following categories:

\begin{definition}
	Let $({\cal C}, J)$ be a small-generated site. 
	\begin{enumerate}[(a)]
		\item The category ${\textup{\bf Mor}}_{({\cal C}, J)}$ has as objects the morphisms of sites from $({\cal C}, J)$ to a small generated site $({\cal D}, K)$ and as arrows $$(F:({\cal C}, J)\to ({\cal D}, K))\to (F':({\cal C}, J)\to ({\cal D}', K'))$$ between any two such morphisms the geometric morphisms $$f:\Sh({\cal D}', K')\to \Sh({\cal D}, K)$$ such that $\Sh(F)\circ f\cong \Sh(F')$:
		
		\begin{equation*}
		\begin{tikzcd} 
		\Sh({\cal C}, J) & \Sh({\cal D}, K) \ar[l, "\Sh(F)"{above}]  \\
		& \Sh({\cal D}', K') \ar[u, "f"] \ar[ul, "\Sh(F')"] 
		\end{tikzcd}
		\end{equation*}	 
			
		\item The category ${\textup{\bf Coh}}_{({\cal C}, J)}$ has as objects the comorphisms of sites from a small-generated site $({\cal D}, K)$ to $({\cal C}, J)$ and as arrows
		$$(U:({\cal D}, K)\to ({\cal C}, J))\to (U':({\cal D}', K') \to ({\cal C}, J))$$ between any two such comorphisms the geometric morphisms $$g:\Sh({\cal D}, K)\to \Sh({\cal D}', K')$$ such that $C_{U'}\circ g\cong C_{U}$:
		
		\begin{equation*}
		\begin{tikzcd} 
		\Sh({\cal D}, K) \ar[d, "g"] \ar[r, "C_{U}"{above}]  & \Sh({\cal C}, J)   \\
		\Sh({\cal D}', K') \ar[ur, "\quad C_{U'}"{below}] &   
		\end{tikzcd}
		\end{equation*}	
	\end{enumerate}
\end{definition}

Let us first show how the assignments $F\mapsto c_{F}$ and $G\mapsto m_{G}$ introduced in the last sections naturally define two functors
\[
C:	({\textup{\bf Mor}}_{({\cal C}, J)})^{\textup{op}} \to 	{\textup{\bf Coh}}_{({\cal C}, J)}
\] 
and
\[
M:{\textup{\bf Coh}}_{({\cal C}, J)} \to 	({\textup{\bf Mor}}_{({\cal C}, J)})^{\textup{op}}. 
\]

Given morphisms of sites $F:({\cal C}, J)\to ({\cal D}, K)$ and $F':({\cal C}, J)\to ({\cal D}', K')$ related by a geometric morphism $f:\Sh({\cal D}', K')\to \Sh({\cal D}, K)$ such that $\Sh(F)\circ f\cong \Sh(F')$, $f$ corresponds under the equivalences
\[
\Sh({\cal D}, K)\simeq \Sh((1_{\cal D}\downarrow F), \tilde{K})
\]
and
\[
\Sh({\cal D}', K')\simeq \Sh((1_{{\cal D}'}\downarrow F'), \tilde{K'})
\]
of Theorem \ref{thmmorphismtocomorphism} to a geometric morphism $C(f):\Sh((1_{{\cal D}'}\downarrow F'), \tilde{K'}) \to \Sh((1_{\cal D}\downarrow F), \tilde{K})$ such that $C_{c_{F}}\circ C(f)\cong C_{c_{F'}}$, which then defines an arrow $c_{F'}\to c_{F}$ in the category ${\textup{\bf Coh}}_{({\cal C}, J)}$. This defines the action of $C$ on arrows.

Notice that if the geometric morphism $f$ is of the form $\Sh(H)$ for some morphism of sites $H:({\cal D}, K)\to ({\cal D}', K')$ such that $F'\circ H\cong F$ then the functor $\overline{H}:(1_{{\cal D}}\downarrow F) \to (1_{{\cal D}'}\downarrow F')$ sending any object $(d, c, \alpha:d\to F(c))$ of $(1_{{\cal D}}\downarrow F)$ to the object $(H(d), c, H(\alpha):H(d)\to H(F(c))\cong F'(c))$ of $(1_{{\cal D}'}\downarrow F')$ is a morphism of sites $\overline{H}:((1_{{\cal D}}\downarrow F), \tilde{K}) \to ((1_{{\cal D}'}\downarrow F'), \tilde{K'})$ such that $\Sh(\overline{H})\cong f$, and $c_{F}\circ \overline{H}\cong c_{F'}$.
 
In the converse direction, given two comorphisms of sites $(U:({\cal D}, K)\to ({\cal C}, J))$ and $(U':({\cal D}', K') \to ({\cal C}, J))$ related by a geometric morphism $g:\Sh({\cal D}, K)\to \Sh({\cal D}', K')$ such that $C_{U'}\circ g\cong C_{U}$, $g$ corresponds, under the equivalences
\[
\Sh({\cal D}, K) \to \Sh(\hat{{\cal D}}, \hat{{K}})
\]
and
\[
\Sh({\cal D}', K') \to \Sh(\hat{{\cal D}'}, \hat{{K'}})
\]
of Theorem \ref{thmcomorphismtomorphism}, to a geometric morphism 
$M(g):\Sh((F \downarrow 1_{\cal C}), \overline{K}) \to \Sh((F' \downarrow 1_{\cal C}), \overline{K'})$ such that $\Sh(m_{F'})\circ M(g)\cong \Sh(m_{F})$. This defines the action of $M$ on arrows. 

Notice that if the geometric morphism $g$ is of the form $C_{H}$ for a comorphism of sites $H:({\cal D}, K)\to ({\cal D}', K')$ such that $U'\circ H\cong U$ then the functor $D_{H}:\hat{{\cal D}'}\to \hat{{\cal D}}$ is a morphism of sites $(\hat{{\cal D}'}, \hat{K'})\to (\hat{{\cal D}}, \hat{K})$ such that $g\cong \Sh(D_{H})$ (cf. Remark \ref{rempresheaflifting}) and $D_{H}\circ m_{U'}\cong m_{U}$.  

As shown by the following result, the functors $C$ and $M$ define a dual adjunction between the categories ${\textup{\bf Mor}}_{({\cal C}, J)}$ and ${\textup{\bf Coh}}_{({\cal C}, J)}$: 

\begin{theorem}\label{adjointequivalencemorphismscomorphisms}
	The functors 
	\[
	C:	({\textup{\bf Mor}}_{({\cal C}, J)})^{\textup{op}} \to 	{\textup{\bf Coh}}_{({\cal C}, J)}
	\] 
	and
	\[
	M:	{\textup{\bf Coh}}_{({\cal C}, J)} \to 	({\textup{\bf Mor}}_{({\cal C}, J)})^{\textup{op}} 
	\]  
	are ($2$-categorically) adjoint ($C$ on the right and $M$ on the left) and quasi-inverse to each other.
\end{theorem}	

\begin{proofs}
	Let us define the unit $\eta$ and counit $\epsilon$ of the dual adjunction.
	
	Given a morphism of sites $F:({\cal C}, J)\to ({\cal D}, K)$, the corresponding comorphism of sites $c_{F}:((1_{\cal D}\downarrow F), \tilde{K})\to ({\cal C}, J)$ is given by the canonical projection functor $(1_{\cal D}\downarrow F)\to {\cal C}$. The morphism of sites associated with $c_{F}$ is 
	\[
	m_{c_{F}}:({\cal C}, J)\to (\widehat{(1_{\cal D}\downarrow F)}, \widehat{\tilde{K}}), 
	\] 
	where $m_{c_{F}}(c)=\textup{Hom}_{\cal C}(c_{F}(-), c)$ (for any $c\in {\cal C}$).
	
	Let us define $\xi_{F}:{\cal D}\to \widehat{(1_{\cal D}\downarrow F)}$ as the functor sending any object $d\in {\cal D}$ to the presheaf $\textup{Hom}_{\cal D}(\pi_{\cal D}(-), d)$, where $\pi_{\cal D}$ is the canonical projection functor $(1_{\cal D}\downarrow F)\to {\cal D}$, and acting on the arrows in the obvious way. The functor $\xi_{F}$ yields a morphism of sites $({\cal D}, K)\to (\widehat{(1_{\cal D}\downarrow F)}, \widehat{\tilde{K}})$, since it is precisely the morphism of sites $m_{\pi_{\cal D}}$ associated as in Theorem \ref{thmcomorphismtomorphism}(iii) with the comorphism of sites $\pi_{\cal D}:((1_{\cal D}\downarrow F), \tilde{K})\to ({\cal D}, K)$ (note that $\pi_{\cal D}$ is indeed a comorphism of such sites by Theorem \ref{thmmorphismtocomorphism}(ii)). 
	
	Let us check that, whilst \emph{not} necessarily commuting at the level of morphisms of sites, the diagram 
	\begin{equation*}
	\begin{tikzcd}
	({\cal C}, J) \arrow[rr, "m_{c_{F}}"]   \arrow[drr, "F"]   & &  (\widehat{(1_{\cal D}\downarrow F)}, \hat{\tilde{K}})           \\
	&  & ({\cal D}, K) \ar[u, "\xi_{F}"] 
	\end{tikzcd}
	\end{equation*}	
	{\flushleft
	commutes after applying the functor $\Sh$ to it, i.e. $\Sh(F)\circ \Sh(\xi_{F})\cong \Sh(m_{c_{F}})$.}
	
	Since $C_{y_{{\cal D}'}}:\Sh((1_{\cal D}\downarrow F), \tilde{K})\to \Sh(\widehat{(1_{\cal D}\downarrow F)}, \widehat{\tilde{K}})$ is an equivalence (cf. Proposition \ref{propyonedamorcomor}), to show that $\Sh(F)\circ \Sh(\xi_{F})\cong \Sh(m_{c_{F}})$ it is equivalent to prove that $\Sh(F)\circ \Sh(\xi_{F})\circ C_{y_{(1_{\cal D}\downarrow F)}} \cong \Sh(m_{c_{F}})\circ C_{y_{(1_{\cal D}\downarrow F)}}$. Now, by Theorem \ref{thmcomorphismtomorphism}(iii), $\Sh(\xi_{F})\circ C_{y_{(1_{\cal D}\downarrow F)}}=\Sh(m_{\pi_{\cal D}})\circ C_{y_{(1_{\cal D}\downarrow F)}} \cong C_{\pi_{\cal D}}$ and $\Sh(m_{c_{F}})\circ C_{y_{(1_{\cal D}\downarrow F)}}\cong C_{c_{F}}$. But $\Sh(F)\circ C_{\pi_{\cal D}}\cong C_{c_{F}}$ by Theorem \ref{thmmorphismtocomorphism}(ii), as required.
	
	The morphism $\Sh(\xi_{F}):\Sh(\widehat{(1_{\cal D}\downarrow F)}, \widehat{\tilde{K}})\to \Sh({\cal D}, K)$ is actually an equivalence. Indeed, by Theorem \ref{thmcomorphismtomorphism}(iii), $\Sh(\xi_{F})=\Sh(m_{c_{\pi_{\cal D}}})$ is an equivalence if and only if $C_{c_{\pi_{\cal D}}}$ is, and by Theorem \ref{thmmorphismtocomorphism}(ii), $C_{c_{\pi_{\cal D}}}$ is an equivalence if and only if $\Sh(\pi_{\cal D})$ is, which is the case by Theorem \ref{thmmorphismtocomorphism}(ii). We set the component $\epsilon_{F}:m_{c_{F}} \to F$ at $F$ of the counit $\epsilon$ of our adjunction equal to the geometric morphism $\Sh(\xi_{F})$ (note that this is actually an arrow $F\to m_{c_{F}}$ in ${\textup{\bf Mor}}_{({\cal C}, J)}$, that is, an arrow $m_{c_{F}} \to F$ in $({\textup{\bf Mor}}_{({\cal C}, J)})^{\textup{op}}$).
	
	Let us now define the unit of our dual adjunction.
	
	Let $G$ be a comorphism of sites $({\cal D}, K)\to ({\cal C}, J)$. The comorphism of sites $c_{m_{G}}:((1_{\hat{D}}\downarrow m_{G}), \widetilde{\hat{K}})\to ({\cal C}, J)$ associated with the morphism of sites $m_{G}:({\cal C}, J)\to (\hat{D}, \hat{K})$ via Theorem \ref{thmmorphismtocomorphism}(iii)  is given by the canonical projection functor $(1_{\hat{D}}\downarrow m_{G}) \to {\cal C}$. We define the component $\eta_{G}:G\to c_{m_{G}}$ of the unit at $G$ as the geometric morphism
	\[
	C_{\chi_{G}}:\Sh({\cal D}, K)\to \Sh((1_{\hat{D}}\downarrow m_{G}), \widetilde{\hat{K}})
	\]
	induced by the comorphism of sites 
	\[
	\chi_{G}:({\cal D}, K)\to ((1_{\hat{D}}\downarrow m_{G}), \widetilde{\hat{K}})
	\] 
	given by the functor ${\cal D}\to (1_{\hat{D}}\downarrow m_{G})$ sending any object $d$ of $\cal D$ to the object $(y_{\cal D}(d), G(c), \alpha:y_{\cal D}(d)\to m_{G}(G(d))$, where $\alpha$ is the natural transformation corresponding to the identity on $G(d)$, regarded as an element of $m_{G}(G(d))(d)$, via the Yoneda lemma. Notice that $\chi_{G}$ is indeed a comorphism of sites $({\cal D}, K)\to ((1_{\hat{D}}\downarrow m_{G}), \widetilde{\hat{K}})$, since $\chi_{G}=z_{G}\circ j_{G}$, where $j_{G}$ and $z_{G}$ are the comorphisms of sites considered in Theorem \ref{thmtoposadjunction}(ii). Moreover, $C_{\chi_{G}}$ is an equivalence, since, by that theorem, both $C_{j_{G}}$ and $C_{z_{G}}$ are.
	
	The diagram
	\begin{equation*}
	\begin{tikzcd}
	({\cal D}, K) \arrow[rr, "G"] \ar[d, "\chi_{G}"]  & &  ({\cal C}, J)           \\ 
	((1_{\hat{D}}\downarrow m_{G}), \widetilde{\hat{K}}) \ar[urr, "c_{m_{G}}"{below, xshift=1.5ex}] &  &
	\end{tikzcd}
	\end{equation*}	
	{\flushleft
	of comorphisms of sites clearly commutes (without the need of applying the functor $C$ to it).}
	
	The verification that the unit and counit satisfy the triangular identities is an immediate consequence of their definition in light of the commutativity of the diagrams in Theorems \ref{thmmorphismtocomorphism} and \ref{thmcomorphismtomorphism}; the details are straightforward and left to the reader.	
\end{proofs}

The dual adjunction of Theorem \ref{adjointequivalencemorphismscomorphisms} shows that, from the point of view of toposes, morphisms and comorphisms of sites are dual to each other.

\subsubsection{The fibration of generalized elements of a functor}\label{sec:fibrgenelements}

Given a comorphism of sites $F:({\cal D}, K)\to ({\cal C}, J)$, the dual comma construction $(1_{\cal C}\downarrow F)$ also plays an important role, as shown by the following result:

\begin{theorem}\label{thmdualcommaconstruction}
	Let $F:({\cal D}, K)\to ({\cal C}, J)$ be a comorphism of small-generated sites, $i_{F}'$ the canonical functor ${\cal D}\to (1_{\cal C}\downarrow F)$ and $K^{i_{F}'}$ the Grothendieck topology on $(1_{\cal C}\downarrow F)$ coinduced by $K$ along $i_{F}'$ (in the sense of Proposition \ref{propimagetopology}). Let $\pi_{\cal C}^{F}$ and $\pi_{\cal D}^{F}$ be the canonical projections from $(1_{\cal C}\downarrow F)$ respectively to $\cal C$ and $\cal D$. Then 
	\begin{enumerate}[(i)]
		\item $\pi_{\cal D}^{F} \dashv i_{F}'$, $\pi_{\cal C}^{F}\circ i_{F}'=F$, $\pi_{\cal C}^{F}$ is a comorphism of sites $((1_{\cal C}\downarrow F), K^{i_{F}'}) \to ({\cal C}, J)$ and $\pi_{\cal D}^{F}$ is a comorphism of sites $((1_{\cal C}\downarrow F), K^{i_{F}'}) \to ({\cal D}, K)$;
		
		\item $i_{F}'$ is both a (full and faithful) comorphism of sites and a dense morphism of sites $({\cal D}, K) \to ((1_{\cal C}\downarrow F), K^{i_{F}'})$ inducing equivalences
		\[
		C_{i_{F}'}:\Sh({\cal D}, K) \to \Sh((1_{\cal C}\downarrow F), K^{i_{F}'})
		\]
		and
		\[
		\Sh(i_{F}'):\Sh((1_{\cal C}\downarrow F), K^{i_{F}'}) \to \Sh({\cal D}, K) 
		\]
		which are quasi-inverse to each other and make the following triangle commute:
		\begin{equation*}
		\begin{tikzcd}
		\Sh((1_{\cal C} \downarrow F), K^{i_{F}'}) \arrow[rr, "\Sh(i_{F}')\cong C_{\pi_{\cal D}^{F}}", yshift=1ex] \ar[rr, phantom, "\sim"{yshift=-0.25ex}]  \arrow[dr, "C_{{\pi_{\cal C}^{F}}}"{below, xshift=-2ex}]   & &  \Sh({\cal D}, K) \ar[ll, "C_{i_{F}'}", yshift=-1ex] \ar[ld, "C_{F}", yshift=0.5ex, xshift=-1ex]            \\
		& \Sh({\cal C}, J)  & 
		\end{tikzcd}
		\end{equation*}	
		
	\end{enumerate}
\end{theorem} 

\begin{proofs}
	(i) The fact that $\pi_{\cal D}^{F} \dashv i_{F}'$ is clear. By definition of coinduced topology $K^{i_{F}'}$, $i_{F}'$ is a comorphism of sites $({\cal D}, K) \to ((1_{\cal C}\downarrow F), K^{i_{F}'})$, which is also cover-preserving ($i_{F}'$ being full). So, by Proposition \ref{propadjointinduction}(i), its left adjoint $\pi_{\cal D}^{F}$ is a comorphism of sites  $((1_{\cal C}\downarrow F), K^{i_{F}'}) \to ({\cal D}, K)$ and $i_{F}'$ is a morphism of sites $({\cal D}, K) \to ((1_{\cal C}\downarrow F), K^{i_{F}'})$. It remains to show that $\pi_{\cal C}^{F}$ is a comorphism of sites $\Sh((1_{\cal C}\downarrow F), K^{i_{F}'}) \to ({\cal C}, J)$. Let $(c, d, \alpha:c\to F(d))$ be an object of $(1_{\cal C}\downarrow F)$, and let us suppose that $R$ is a $J$-covering sieve on the object $c=\pi_{\cal C}^{F}((c, d, \alpha:c\to F(d)))$. Let $\tilde{R}$ be the sieve in $(1_{\cal C}\downarrow F)$ on $(c, d, \alpha:c\to F(d))$ generated by by the arrows of the form $(f, 1_{d}):(\textup{dom}(f), d, \alpha\circ f:\textup{dom}(f) \to F(d))\to (c, d, \alpha:c\to F(d))$. Clearly, the image of $\tilde{R}$ under $\pi_{\cal C}^{F}$ is contained in $R$, so if we show that $\tilde{R}$ is $K^{i_{F}'}$-covering, we are done. Now, $\tilde{R}$ is $K^{i_{F}'}$-covering if and only if for any arrow $\xi=(\xi_{1}, \xi_{2})$ in $(1_{\cal C}\downarrow F)$ to $(c, d, \alpha:c\to F(d))$ from an object of the form $i_{F}'(d')=(F(d'),d', 1_{F(d')})$, $\xi^{\ast}(\tilde{R})$ contains the image under $i_{F}'$ of a $K$-covering sieve on $d'$. Notice that, since $\xi=(\xi_{1}, \xi_{2})$ is an arrow $(F(d'),d', 1_{F(d')})\to (c, d, \alpha:c\to F(d))$ in $(1_{\cal C}\downarrow F)$, we have $\alpha\circ \xi_{1}=F(\xi_{2})$.  
	Consider the sieve $\xi_{1}^{\ast}(R)$ on $F(d')$; it is $J$-covering as it is the pullback of a $J$-covering sieve. So, since $F$ is, by our hypothesis, a comorphism of sites $({\cal D}, K)\to ({\cal C}, J)$, there is a $K$-covering sieve $T$ on $d'$ such that $F(T)\subseteq \xi_{1}^{\ast}(R)$. Let us show that $i_{F}'(T)\subseteq \xi^{\ast}(\tilde{R})$. For any $t\in T$, $\xi\circ i_{F}'(t)=(\xi_{1}, \xi_{2})\circ (F(t), t)=(\xi_{1}\circ F(t), \xi_{2}\circ t):(F(\textup{dom}(t)),\textup{dom}(t), 1_{F(\textup{dom}(t))}) \to (c, d, \alpha:c\to F(d))$. But $\xi_{1}\circ F(t)\in R$ and hence the arrow $(\xi_{1}\circ F(t), 1_{d}):(F(\textup{dom}(t)), d, \alpha\circ \xi_{1}\circ F(t): F(\textup{dom}(t))\to F(d))\to (c, d, \alpha:c\to F(d))$ lies in $\tilde{R}$. But the arrow $(\xi_{1}\circ F(t), \xi_{2}\circ t)$ can be written as the composite of this arrow with the arrow $(1_{F(\textup{dom}(t))}, \xi_{2}\circ t):(F(\textup{dom}(t)), \textup{dom}(t), 1_{F(\textup{dom}(t))})\to (F(\textup{dom}(t)), d, \alpha\circ \xi_{1}\circ F(t):F(\textup{dom}(t))\to F(d))$, which is indeed an arrow in $(1_{\cal C}\downarrow F)$ since $F(\xi_{2}\circ t)=\alpha \circ \xi_{1}\circ F(t)$ (this follows from the previously remarked fact that $\alpha\circ \xi_{1}=F(\xi_{2})$). So $(\xi_{1}\circ F(t), \xi_{2}\circ t)$ also lies in $\tilde{R}$, as desired.   
	
	(ii) To show that $i_{F}'$ induces equivalences in the statement of the theorem, by Proposition \ref{propmorphismcomorphismequivalence} we are reduced, in light of (i), to prove that it is $K$-full and $K^{i_{F}'}$-dense. The functor $i_{F}'$ is actually full, so \emph{a fortiori} $K$-full, and it is $K^{i_{F}'}$-dense by Proposition \ref{propimagetopology}.  	
\end{proofs}

\begin{remark}
	The comorphisms of sites $F\circ \pi_{\cal D}^{F}$ and $\pi_{\cal C}^{F}$ are not in general isomorphic in spite of the fact that they induce isomorphic geometric morphisms; indeed, $C_{F}\circ C_{\pi_{\cal D}^{F}}\cong C_{\pi_{\cal C}^{F}}$ by the commutativity of the above diagram). This provides an example of a pair of comorphisms of sites which are non-trivially equivalent to each other (in the sense of Remark \ref{remdiffcomorphisms} and section \ref{sec:presheaflifting}).   
\end{remark}

Notice that the objects of the category $(1_{\cal C}\downarrow F)$ are the generalized elements of the functor $F$, and the canonical projection functor ${\pi_{\cal C}^{F}}:(1_{\cal C}\downarrow F) \to {\cal C}$ is a split fibration, since it is the result of applying the Grothendieck construction to the internal category $\mathbb F$ to $[{\cal C}^{\textup{op}}, \Set]$ given by $c\mapsto (c\downarrow F)$. This motivates the following definition:

\begin{definition}
	Let $F:{\cal D}\to {\cal C}$ be a functor. The \emph{fibration of generalized elements of $F$} is the canonical projection functor ${\pi_{\cal C}^{F}}:(1_{\cal C}\downarrow F) \to {\cal C}$. 
\end{definition}

\begin{remark}
	As shown in \cite{Kock}, the assignment $(F:{\cal C}\to {\cal D}) \mapsto ({\pi_{\cal C}^{F}}:(1_{\cal C}\downarrow F) \to {\cal C})$ yields an endofunctor $T$ on $\textbf{Cat}\slash {\cal C}$ which has the structure of a monad, whose unit $\eta_{T}$ is the functor $i'^{F}:{\cal D}\to (1_{\cal C}\downarrow F)$. If $F$ is a fibration then every cleavage for it defines an algebra $r:(1_{\cal C}\downarrow F) \to {\cal D}$ for this monad which is actually a morphism of fibrations; in particular, $r\circ \eta_{T}=1_{\cal D}$, whence $F$ is actually a retract of the split fibration $\pi_{\cal D}^{F}$ in the category $\textbf{Cat}\slash {\cal C}$ (but \emph{not} in the category of fibrations over $\cal C$ since $\eta_{T}$ is not in general a morphism of fibrations). Notice that if $F$ is a fibration then every object of $(1_{\cal C}\downarrow F)$ can be expressed as the cartesian image of an arrow in the image of the functor $\eta_{T}$. Indeed, given an object $(c,d, f:c\to F(d))$, take a cartesian arrow $g:d'\to d$ in $\cal D$ such that $F(g)\circ \alpha=f$ for some isomorphism $\alpha:c\to F(d')$; then the arrow $\eta_{T}(g)=(F(g), g)$ factors as $(f, 1_{d})\circ (\alpha, g)$, where $(f, 1_{d})$ is vertical and $(\alpha, g)$ is cartesian (with respect to the fibration $\pi_{\cal C}^{F}$).     
\end{remark}

The equivalence of Theorem \ref{thmdualcommaconstruction}(ii) between $C_{F}$ and $C_{{\pi_{\cal C}^{F}}}$ makes it possible to transfer many notions and constructions between $F$ and its fibrations of generalized elements. In particular, it is interesting to relate the Grothendieck topologies $M^{F}_{J}$ and $M^{\pi_{\cal C}^{F}}_{J}$ respectively on $\cal D$ and on $(1_{\cal C}\downarrow F)$ (see section \ref{sec:comorphismsofsites1} for the notation): 

\begin{proposition}\label{propthreepoints}
	Let $({\cal C}, J)$ be a small-generated site and $F:{\cal D}\to {\cal C}$ a functor. Then
	\begin{enumerate}[(i)]
		\item A sieve $S$ on an object $(c, d, \alpha:c\to F(d))$ of the category is $(1_{\cal C}\downarrow F)$ is $M^{\pi^{F}_{\cal C}}_{J}$-covering if and only if it contains a family of arrows of the form $\{(f_i, 1_{d}):(c_i, d, \alpha:c\to F(d)) \to (c, d, \alpha:c\to F(d))\}$ such that the family $\{f_i:c_i\to c \mid i\in I\}$ is $J$-covering.  
		
		\item $M^{F}_{J}=M_{M^{\pi_{\cal C}^{F}}_{J}}^{i_{F}'}$. 
		
		\item $({M^{F}_{J}})^{i_{F}'}=M^{\pi^{F}_{\cal C}}_{J}\vee R_{i_{F}'}$, where $R_{i_{F}'}={T_{\cal C}}^{i_{F}'}$ is the rigid topology associated with the inclusion $i_{F}'$ (here $T_{\cal C}$ denotes the trivial Grothendieck topology on $\cal C$).
	\end{enumerate}
\end{proposition}

\begin{proofs}
	(i) This is an immediate consequence of Theorem \ref{thmfibration1}.
	
	(ii) This follows from the ($2$-dimensional) pullback lemma, which implies that the left-hand square in the following diagram (where the unnamed vertical arrows are the canonical geometric inclusions) is a pullback (since the right-hand square and the outer rectangle are):
	$$
	\xymatrix{
		\Sh({\cal D}, M^{J}_{F}) \ar[d] \ar[r] & \Sh((1_{\cal C}\downarrow F), M^{\pi^{F}_{\cal C}}_{J}) \ar[r] \ar[d] & [{\cal D}^{\textup{op}}, \Set] \ar[d]^{C_{F}} \\	
		[{\cal D}^{\textup{op}}, \Set] \ar[r]^{C_{i_{F}'}\quad} & [(1_{\cal C}\downarrow F)^{\textup{op}}, \Set] \ar[r]^{\quad C_{\pi^{F}_{\cal C}}} & [{\cal C}^{\textup{op}}, \Set]	
	}
	$$
	(notice that the composite of the lower horizontal arrows is $C_{F}$).
	
	(iii) This follows from (ii) and Theorem \ref{thmcorrespondenceinclusion}. 
\end{proofs}

\begin{proposition}
	Let $({\cal C}, J)$ be a small-generated site and $F:{\cal D}\to {\cal C}$ a functor. Then, for any sieve $T$ on an object $d$ of $\cal D$, $T$ is $M_{J}^{F}$-covering if and only if it is sent by $i_{F}'$ to a $(M_{J}^{{\pi^{F}_{\cal C}}}\vee R_{i_{F}'})$-covering family.	
\end{proposition}

\begin{proofs}
	This follows from Corollary \ref{corexplicitdescriptionsmallesttopology}, applied to the morphism of sites $i_{F}':({\cal D}, T_{\cal D}) \to ((1_{\cal C}\downarrow F), ({M^{F}_{J}})^{i_{F}'})$, in light of Proposition \ref{propthreepoints}(iii). 		
\end{proofs}

\section{Weak morphisms and continuous comorphisms of sites}\label{sec:essentialcomorphisms}

In this section we first introduce the notion of \emph{weak morphism} of toposes (meaning simply a pair of adjoint functors between them, without the requirement that the left adjoint should preserve finite limits) and establish an equivalence result between the weak morphisms to a sheaf topos $\Sh({\cal C}, J)$ and the class of continuous functors ${\cal C}\to {\cal E}$, which notably includes that of $J$-continuous flat functors. Then we apply this result to obtain an equivalence between the essential geometric morphisms $\Sh({\cal C}, J)\to {\cal E}$ and the comorphisms of sites $({\cal C}, J)\to ({\cal E}, J^{\textup{can}}_{\cal E})$ which are $J$-continuous (recall that a geometric morphism $f:{\cal F}\to {\cal E}$ is said to be \emph{essential} if its inverse image $f^{\ast}:{\cal E}\to {\cal F}$ admits a left adjoint, usually denoted by $f_{!}:{\cal F}\to {\cal E}$, and which we shall call its \emph{essential image}). Then, after obtaining a number of equivalent characterizations for the notion of canonicity for comorphisms of sites, we discuss how the operation of extension of Grothendieck topologies on a category to the corresponding presheaf topos can be used for describing any essential geometric morphism $\Sh({\cal C}, J)\to \Sh({\cal D}, K)$ as induced by a comorphism of sites from $({\cal C}, J)$ to a site naturally attached to $({\cal D}, K)$ (namely, the site $([{\cal D}^{\textup{op}}, \Set], \hat{K})$, where $\hat{K}$ is the extension of $K$ from $\cal D$ to $[{\cal D}^{\textup{op}}, \Set]$).

\subsection{Weak morphisms of toposes}\label{sec:weakmorphisms}

It is sometimes useful to consider pairs of adjoint functors between toposes which do not yield geometric morphisms since the inverse image functor does not preserve finite limits. 

This motivates the following definition: 

\begin{definition}
	A \emph{weak morphism of toposes} $f:{\cal E}\to {\cal F}$ is a pair of adjoint functors $(f^{\ast} \dashv f_{\ast})$.
\end{definition}

As in the case of geometric morphism, we shall call the functor $f_{\ast}$ (resp. $f^{\ast}$) the \emph{direct image} (resp. the \emph{inverse image}) of the weak morphism $f$.

Clearly, the identity is a weak morphism and the composite of two weak morphisms is again a weak morphism. The category whose objects are the weak morphisms from $\cal E$ to $\cal F$ and whose arrows $f\to g$ are the natural transformations $f^{\ast}\to g^{\ast}$ will be denoted by $\textup{\bf Wmor}({\cal E}, {\cal F})$.

Weak morphisms of toposes allow us in particular to treat under the same roof both geometric morphisms and pairs of functors of the form $(f_{!}\dashv f^{\ast})$, where $f$ is an essential geometric morphism. In fact, some results and constructions valid for geometric morphisms naturally generalize to weak morphisms and then can be applied to derive results about essential geometric morphisms. An example of such a result is the equivalence between the property of a right adjoint functor to be faithful and the property that its left adjoint is locally surjective (cf. Corollary \ref{corinclusionequivalence} and Remark \ref{remsurjectionsinclusions}):

\begin{proposition}\label{propadjproperty}
	Let $\cal E$ and $\cal F$ categories and $L\dashv R$ a pair of adjoint functors $R:{\cal E}\to {\cal F}$ and $L:{\cal F}\to {\cal E}$. The the following conditions are equivalent:
	\begin{enumerate}[(i)]
		\item $R$ is faithful.
		
		\item For any $e\in {\cal E}$, there is $u\in {\cal F}$ and an epimorphism $L(u)\to e$.
	\end{enumerate}
	
	If $\cal E$ and $\cal F$ are Grothendieck toposes respectively with separating sets $\cal C$ and $\cal D$ then condition (ii) is equivalent to the following one: for any $c\in {\cal D}$, the family of arrows $L(d)\to c$ (where $d\in {\cal D}$) is epimorphic.
\end{proposition}

\begin{proofs}
	(i) $\imp$ (ii). Since $R$ is faithful, for any object $e$ of $\cal E$, the component $\epsilon_{e}:L(R(e))\to e$ of the counit of the adjunction $L\dashv R$ is an epimorphism. So condition (ii) is satisfied.
	
	(ii) $\imp$ (i). Any arrow $\alpha:L(u)\to e$, $\alpha$ factors through the component $\epsilon_{e}$ at $e$ of counit of the adjunction $L\dashv R$. So if $\alpha:L(u)\to e$ is an epimorphism, $\epsilon_{e}$ is also an epimorphism. Since this holds for every $e$, it implies that $R$ is faithful.
	
	It remains to prove the last assertion of the proposition. The fact that it suffices to require condition (ii) for every $e$ in $\cal C$ follows from the fact that every object $e$ of $\cal E$ admits an epimorphic family from objects in $\cal C$. On the other hand, since $\cal D$ is separating for $\cal F$, requiring, for a given object $c$, that there should be an epimorphism $L(u)\to c$, is equivalent to the condition that the family of all arrows $L(d)\to c$ (where $d\in {\cal D}$) is epimorphic (since $L$ preserves colimits, having a right adjoint, and hence epimorphic families).     	
\end{proofs}

\begin{remark}
	In general, for a given property $P$ of functors, in order to find a property $P_{l}$ such that whenever $R\dashv L$, $R$ satisfies $P$ if and only if, one can proceed to reformulate it in terms of the generalized elements of the objects in the image of $R$, and then rewrite the obtained condition in terms of $L$ by using the bijective correspondence between the generalized elements $b\to R(a)$ of $R(a)$ and the arrows $L(a)\to b$ provided by the adjunction between $R$ and $L$. Another possibility is to find a categorical property $I$ such that the category $(1\downarrow R)$ satisfies $I$ if and only if $R$ satisfies $P$, in which case one can obtain a property $P_{l}$ by unraveling what it means in terms of $L$ for the category $(L\downarrow 1)$ to satisfy $I$ (since having an adjunction $R\dashv L$ between two categories amounts precisely to having a categorical equivalence $(1\downarrow R)\simeq (L\downarrow 1)$ over the product of the two categories).
\end{remark}

The following result, which generalizes the classical one for geometric morphisms, characterizes the weak morphisms which factor through a geometric inclusion:

\begin{proposition}
	Let $i:{\cal F}\hookrightarrow {\cal E}$ be the geometric inclusion of a subtopos $\cal F$ of a Grothendieck topos $\cal E$ into $\cal E$, and let $f:{\cal G}\to {\cal E}$ be a weak morphism from a Grothendieck topos $\cal G$. Then the following conditions are equivalent:
	\begin{enumerate}[(i)]
		\item The weak morphism $f$ factors through $i$; 
		
		\item The direct image $f_{\ast}$ takes values in $\cal F$ (that is, factors through $i_{\ast}$);
		
		\item The inverse image $f^{\ast}$ factors (necessarily uniquely up to isomorphism) through $i^{\ast}$.
	\end{enumerate} 
\end{proposition}

\begin{proofs}
	The implications (i) $\imp$ (ii) and (i) $\imp$ (iii) are clear. The fact that, in (iii), if the factorization exists then it is unique (up to isomorphism) can be proved as follows. Since $i^{\ast}\circ i_{\ast}\cong 1_{\cal F}$ and $i^{\ast}$ preserves colimits, if $f^{\ast}$ factors as $r\circ i^{\ast}$ then the functor $r$ preserves colimits and hence has a right adjoint; let us call it $j$. Now, since $f^{\ast}\cong r\circ i^{\ast}$, passing to the right adjoints we obtain that $f_{\ast}\cong i_{\ast}\circ j$; in other words, $f_{\ast}$ factors (necessarily uniquely up to isomorphism) through $i_{\ast}$ and $j$ identifies with this factorization. So $r$ is uniquely determined as the left adjoint to this factorization. This argument actually proves that (iii) $\imp$ (i). Similarly, we can prove that (ii) $\imp$ (i).   
\end{proofs}

\begin{corollary}\label{corfact}
	Let $A:{\cal C}\to {\cal E}$ be a functor from an essentially small category $\cal C$ to a Grothendieck topos $\cal E$, and $J$ be a Grothendieck topology on $\cal C$. Then the following conditions are equivalent:
	\begin{enumerate}[(i)]
		\item The weak morphism $(L_{A} \dashv R_{A})$ factors through the canonical geometric inclusion $i:\Sh({\cal C}, J)\hookrightarrow [{\cal C}^{\textup{op}}, \Set]$; 
		
		\item The functor $R_{A}$ takes values in $\Sh({\cal C}, J)$;
		
		\item The functor $L_{A}$ factors (necessarily uniquely up to isomorphism) through the associated sheaf functor $a_{J}:[{\cal C}^{\textup{op}}, \Set]\to \Sh({\cal C}, J)$.
	\end{enumerate} 
\end{corollary}\qed

\begin{remark}\label{remflatnessandcanonicity}
	Interestingly, if $A$ is not flat the property of $J$-continuity, whilst a necessary condition for $R_{A}$ to take values in $\Sh({\cal C}, J)$, is no longer sufficient in general for this to hold: take for example the functor $A:{\cal C}\to \Set$ sending any object of $\cal C$ to the singleton $1_{\Set}$ and an arrow of $\cal C$ to the identity on $1_{\Set}$. The functor $R_{A}:\Set \to [{\cal C}^{\textup{op}}, \Set]$ thus sends any set $E$ to the constant presheaf on $\cal C$ with value $E$; now, for any Grothendieck topology $J$ on $\cal C$, these presheaves are all $J$-sheaves if and only if the site $({\cal C}, J)$ is locally connected in the sense of section C3.3 \cite{El} (that is, every $J$-covering sieve on an object $c$ is connected as a full subcategory of the slice category ${\cal C}\slash c$). So, if $({\cal C}, J)$ is not locally connected then the functor $R_{A}$ does not take values in $\Sh({\cal C}, J)$, in spite of $A$ being $J$-continuous.
\end{remark}

In the next section we shall investigate the functors between sites which induce weak morphisms between the associated toposes. 

\subsection{$(J, K)$-continuous functors}\label{subsec:continuousfunctors}

The results in the above section naturally lead us to consider the following notion: 

\begin{definition}\label{defcanonicalfunctor}
	\begin{enumerate}[(a)]
		\item Given a small-generated site $({\cal C}, J)$, we say that a functor $A:{\cal C}\to {\cal E}$ is \emph{$J$-continuous} if the hom functor $R_{A}:{\cal E}\to [{\cal C}^{\textup{op}}, \Set]$ takes values into $\Sh({\cal C}, J)$ (equivalently, by Corollary \ref{corfact}, if the functor $L_{A}: [{\cal C}^{\textup{op}}, \Set]\to {\cal E}$ factors through $a_{J}:[{\cal C}^{\textup{op}}, \Set] \to \Sh({\cal C}, J)$). 
		
		\item Given small-generated sites $({\cal C}, J)$ and $({\cal D}, K)$, a functor $A:{\cal C}\to {\cal D}$ is said to be $(J, K)$-continuous if $l'\circ A$ is $J$-continuous, where $l'$ is the canonical functor ${\cal D}\to \Sh({\cal D}, K)$.
	\end{enumerate}	
\end{definition}

As shown by Proposition \ref{propequivGroth} below, a functor $A:{\cal C}\to {\cal D}$ is $(J, K)$-continuous in the sense of Definition \ref{defcanonicalfunctor}(b) if and only if it is a continuous functor $({\cal C}, J)\to ({\cal D}, K)$ in the sense of section III.1 of \cite{grothendieck} (that is, satisfies condition (ii) of Proposition \ref{propequivGroth}), which explains our terminology.

\begin{proposition}\label{propequivGroth}
	Let $({\cal C}, J)$ and $({\cal D}, K)$ be small-generated sites and  $A:{\cal C} \to {\cal D}$ a functor. Then the following conditions are equivalent:
\begin{enumerate}[(i)]
	\item $A$ is $(J, K)$-continuous.
	
	\item The functor $$D_{A}:=(-\circ A^{\textup{op}}):[{\cal D}^{\textup{op}}, \Set]\to [{\cal C}^{\textup{op}}, \Set]$$ restricts to a functor $\Sh({\cal D}, K) \to \Sh({\cal C}, J)$.	
	\end{enumerate}
\end{proposition}

\begin{proofs}
By definition, $A$ is $(J, K)$-continuous if and only if the functor 
\[
R_{l'\circ A}:\Sh({\cal D}, K)\to [{\cal C}^{\textup{op}}, \Set]
\]
takes values in $\Sh({\cal C}, J)\subseteq [{\cal C}^{\textup{op}}, \Set]$. But we have
\[
R_{l'\circ A} \cong D_{A}\circ i_{K}.
\]  
Indeed, for any $E\in \Sh({\cal D}, K)$, 
\[
R_{l'\circ A}(E)= \textup{Hom}_{\Sh({\cal D}, K)}((l'\circ A)(-), E)\cong E\circ A^{\textup{op}}=(D_{A}\circ i_{K})(E)
\]
naturally in $E$. From this our thesis immediately follows.	
\end{proofs}

\begin{remark}\label{remhereditarycontinuity}
From the characterization of Proposition \ref{propequivGroth} it follows at once that any $(J, K)$-continuous functor $A:{\cal C}\to {\cal D}$ is $(J', K')$-continuous for every $J'\subseteq J$ and every $K'\supseteq K$. 	
\end{remark}

\begin{example}\label{examplemorphismsofdiscretefibrations}
	Let $P$ and $Q$ be presheaves on a category $\cal C$. Then, for any morphism $\alpha:P\to Q$, the comorphism of sites  ${\int \alpha}:({\int P}, J_{P})\to ({\int Q}, J_{Q})$ is $(J_{P}, J_{Q})$-continuous (where $J_{P}$ and $J_{Q}$ are the Grothendieck topologies defined in Remark \ref{remrelativetopology}). Indeed, this follows from Proposition \ref{propequivGroth} by observing that the functor $D_{\int \alpha}:=(-\circ {(\int \alpha)}^{\textup{op}}):[{(\int Q)}^{\textup{op}}, \Set]\to [{(\int P)}^{\textup{op}}, \Set]$ restricts to a functor $\Sh({\int Q}, J_{Q}) \to \Sh({\int P}, J_{P})$. In particular, all the comorphisms of the form $\pi_{P}:({\int P}, J_{P}) \to ({\cal C}, J)$ are $(J_{P}, J)$-continuous. Alternatively, the continuity of ${\int \alpha}:({\int P}, J_{P})\to ({\int Q}, J_{Q})$ follows from Corollary \ref{correstrictioncomorphism} below by observing that the geometric morphism $\Sh({\cal C}, J)\slash a_{J}(\alpha):\Sh({\cal C}, J)\slash a_{J}(P)\to \Sh({\cal C}, J)\slash a_{J}(Q)$ induced by the arrow $a_{J}(\alpha):a_{J}(P)\to a_{J}(Q)$ is essential and the following diagram, where $Z_{a_{J}(h)}$ is the functor $\Sh({\cal C}, J)\slash a_{J}(P)\to \Sh({\cal C}, J)\slash a_{J}(Q)$ given by composition with $a_{J}(\alpha)$ (that is, the essential image of $\Sh({\cal C}, J)\slash a_{J}(\alpha)$) and the vertical functors are the canonical ones modulo the equivalences of section \ref{sec:example}, commutes: 
	
	$$
	\xymatrix{
		{\int P} \ar[r]^{\int \alpha} \ar[d] & {\int Q} \ar[d] \\
		\Sh({\cal C}, J)\slash a_{J}(P) \ar[r]^{Z_{a_{J}(h)}} & \Sh({\cal C}, J)\slash a_{J}(Q)
	}
	$$ 
	
	See also Theorem \ref{thmlocallyconnected} below.
\end{example}

The following result shows in particular that, in the case of a flat functor $A:{\cal C}\to {\cal E}$ (but not in general), $A$ is $J$-continuous in the sense of Definition \ref{defcanonicalfunctor} if and only if it is $J$-continuous in the sense of section VII.7 \cite{MM}.   

\begin{proposition}\label{propcharJcanonical}
	Let $({\cal C}, J)$ and $({\cal D}, K)$ be small-generated sites and $\cal E$ a Grothendieck topos. Then 		
	\begin{enumerate}[(i)]
		\item A functor $A:{\cal C}\to {\cal E}$ is $J$-continuous if and only if for any $J$-covering sieve $S$ on an object $c$
		\[
		A(c)=\varinjlim_{f:d\to c \in S} A(d) 
		\]
		for each $J$-covering sieve $S$ on an object $c$ (where the colimit is indexed by the category $\int S$ of elements of $S$).
		
		\item A functor $A:{\cal C}\to {\cal D}$ is $(J, K)$-continuous if and only if for any $J$-covering sieve $S$ on an object $c$ the canonical cocone with vertex $A(c)$ on the diagram $\{A(\textup{dom}(f)) \mid f\in S\}$ indexed over $\int S$ is sent by $l'$ to a colimit in the topos $\Sh({\cal D}, K)$. 
		
		\item Every $J$-continuous functor $A:{\cal C}\to {\cal E}$ is $J$-continuous in the sense of Definition \ref{defcanonicalfunctor}(a) is $J$-continuous in the sense of section VII.7 of \cite{MM} (that is, sends $J$-covering families to epimorphic families), and the converse is true if $A$ is flat (but not in general, cf. Remark \ref{remflatnessandcanonicity}). More generally, every $(J, K)$-continuous functor $({\cal C}, J)\to ({\cal D}, K)$ is cover-preserving, and every morphism of sites $({\cal C}, J)\to ({\cal D}, K)$ is $(J, K)$-continuous.
		
		\item For any functor $A:{\cal C}\to {\cal E}$, there is a largest Grothendieck topology $J$ on the category $\cal C$ such that $A$ is $J$-continuous. This topology, denoted by $Z_{A}$, admits the following explicit description: for any sieve $S$ on an object $c$ of $\cal C$, $S$ is $Z_{A}$-covering if and only if for any arrow $g:d\to c$, $A(d)=\varinjlim_{h:e\to d \in g^{\ast}(S)} A(e)$. 	
	\end{enumerate}
\end{proposition}

\begin{proofs}
	(i) and (ii) are immediate, since stating that all the functors of the form $R_{A}(e)=\textup{Hom}_{\cal E}(A(-), e)$ (for $e\in {\cal E}$) satisfy the sheaf condition with respect to a sieve $S$ on an object $c$ amounts precisely to saying that $A(c)=\varinjlim_{f:d\to c \in S} A(d)$. 
	
	(iii) The equality $A(c)=\varinjlim_{f:d\to c \in S} A(d)$ implies in particular that the family of arrows $\{A(\textup{dom}(f)) \to A(c) \mid f\in S\}$ is epimorphic, i.e. that $A$ is $J$-continuous. If $A$ is flat then the geometric morphism $(L_{A}\dashv R_{A})$ factors through the canonical inclusion $\Sh({\cal C}, J)\hookrightarrow [{\cal C}^{\textup{op}}, \Set]$ if and only if $A$ is $J$-continuous (see, for instance, Corollary 4 at p. 394 of \cite{MM}). The second assertion, about morphisms of sites, follows from the first about flat functors by Remark \ref{remmorphismsofsites}(a).
	
	(iv) This follows from the fact that, for any collection $\cal P$ of presheaves on a given category $\cal C$, there is a largest Grothendieck topology for which all of them are sheaves. This topology can be explicitly characterized as the topology having as covering sieves precisely the sieves $S$ such that every presheaf $P$ in $\cal P$ satisfies the sheaf condition with respect to all the pullback sieves $g^{\ast}(S)$ (see, for instance, the proof of Proposition C2.1.9 \cite{El}). If $\cal P$ is the class of presheaves in the image of the functor $R_{A}$ then this condition amounts precisely to requiring that for any arrow $g:d\to c$, $A(d)=\varinjlim_{h:e\to d \in g^{\ast}(S)} A(e)$ (cf. the proof of (i)).
\end{proofs}

\begin{remarks}\label{remcontinuityprelim}
	\begin{enumerate}[(a)]		
		\item The condition for $A$ to be $J$-continuous can be interpreted as a sort of cofinality condition. Indeed, if $A$ is $J$-continuous then $A$ sends any $J$-covering sieve $S$ on an object $c$ of $\cal C$ to an epimorphic family and hence $A(c)$ is the colimit of the cocone under the diagram whose vertices are the objects of the form $A(d)$ where $d$ is the domain of an arrow $f:d\to c$ in $S$ and whose arrows are all the arrows in $\cal E$ over $A(c)$ between such objects. So the condition for $A$ to be $J$-continuous amounts precisely to the assertion that $A$ be $J$-continuous and that this colimit be equal to the colimit $\varinjlim_{f:d\to c \in S} A(d)$. See Remark \ref{remcharcanonicalfunctors}(a) below for more details.
		
		\item We shall obtain further characterizations for the property of a functor $A:{\cal C}\to {\cal D}$ (or more specifically a comorphism of sites $({\cal C}, J)\to ({\cal D}, K)$) to be $(J, K)$-continuous in Proposition \ref{propcharcanonicalfunctors} and Corollary \ref{corequivcharcanonicity}. See also Proposition III-1.2 \cite{grothendieck} for alternative characterizations. 
		
		\item The existence of the topology $Z_{A}$ as in point (iv) of Proposition \ref{propcharJcanonical} was already shown in \cite{grothendieck} (cf. Proposition III-3.2 therein), where this topology is called the Grothendieck topology induced by $A$, regarded as a functor to the site $({\cal E}, J^{\textup{can}}_{\cal E})$. For an explicit description of this topology, see Remark \ref{remcharcanonicalfunctors}(b) below. 
	\end{enumerate}
\end{remarks}

In general, the condition for a cocone under a diagram $D:{\cal I}\to {\cal D}$ with vertex $L$ to be sent by the canonical functor $l':{\cal D}\to \Sh({\cal D}, K)$ to a colimiting cocone can be expressed as the requirement that the canonical arrow $\textup{colim}(y_{\cal D}\circ D)\to y_{\cal D}(L)$ is sent by the associated sheaf functor $a_{K}$ to an isomorphism; these conditions are made explicit in Corollary \ref{corcolimitintopos}. In our case, the diagram to be considered is the functor 
\[
D^{A}_{S}:{\textstyle \int S} \to {\cal D}  
\]
sending any $(d, f)$ of $\int S$ to $A(d)$, together with the cocone $\xi_{A}$ with vertex $A(c)$ under it (whose legs are the arrows $A(f):A(d)=D^{A}_{S}((d, f))\to A(c)$ for any object $(d, f)$ of $\int S$). Applying this result in connection with Proposition \ref{propcharJcanonical}(i), we thus obtain the following explicit characterization of $(J, K)$-continuous functors:

\begin{proposition}\label{propcharcanonicalfunctors}
	Let $({\cal C}, J)$ and $({\cal D}, K)$ be small-generated sites. Then a functor $A:{\cal C}\to {\cal D}$ is $(J, K)$-continuous if and only if it is cover-preserving (i.e., sends $J$-covering families to $K$-covering ones) and for any $J$-covering sieve $S$ on an object $c$ and any commutative square of the form
	\[
	\xymatrix{
		d \ar[r] \ar[d]  & A(c') \ar[d]^{A(f)} \\
		A(c'') \ar[r]^{A(g)} & A(c),
	} 
	\]	
	where $f:c'\to c$ and $g:c''\to c$ are arbitrary arrows of $S$, there is a $K$-covering family $\{d_{i}\to d \mid i\in I\}$ such that for each $i\in I$, the composites $d_{i}\to A(c')$ and $d_{i}\to A(c'')$ belong to the same connected component of the category $(d_{i}\downarrow D^{A}_{S})$.  
\end{proposition}\qed 

\begin{remarks}\label{remcharcanonicalfunctors}
	\begin{enumerate}[(a)]
		\item The conditions of Proposition \ref{propcharcanonicalfunctors} are equivalent to the requirement that that the lift
		\[
		(D^{S}_{A})_{\xi_{A}}:{\textstyle \int S} \to {\cal D}\slash A(c)  
		\]
		of the diagram $D^{S}_{A}$ to ${\cal D}\slash A(c)$ induced by the cocone $\xi_{A}$ (in the sense of Corollary \ref{corcolimitintopos}) be $K_{A(c)}$-cofinal (where the notation is that of section \ref{sec:relativecofinality}). In fact, the property that $A$ be cover-preserving is equivalent to condition (i) of Corollary \ref{corcolimitintopos}, which reads as follows: for any $J$-covering sieve $S$ on an object $c$ of $\cal C$, for any arrow $y:d\to A(c)$ in $\cal D$ there are a $K$-covering family $\{g_i: d_i \to d \mid i\in I\}$ and for each $i\in I$ an arrow $f_{i}$ of $S$ and an arrow $y_i:d_{i}\to (A(\dom(f_i)))$ in $\cal D$ such that $y\circ g_{i}=A(f_{i})\circ y_{i}$.
		
		\item By Proposition \ref{propcharcanonicalfunctors}, the largest topology $Z_{A}$ which makes a functor $A:{\cal C}\to {\cal D}$ continuous to $({\cal D}, K)$ (as in point (iv) of Proposition \ref{propcharJcanonical}) consists precisely of the sieves $S$ such that the functor $A$ satisfies the conditions of Proposition \ref{propcharcanonicalfunctors} with respect to all the sieves $S'$ which are of the form $g^{\ast}(S)$ for some arrow $g$, in the sense that $A$ sends any such $S'$ to a $K$-covering family and for any commutative square of the form
		\[
		\xymatrix{
			d \ar[r] \ar[d]  & A(c') \ar[d]^{A(f)} \\
			A(c'') \ar[r]^{A(g)} & A(c),
		} 
		\]	
		where $f:c'\to c$ and $g:c''\to c$ are arbitrary arrows of $S'$, there is an epimorphic family $\{d_{i}\to d \mid i\in I\}$ in $\cal E$ such that for each $i\in I$, the composites $d_{i}\to A(c')$ and $d_{i}\to A(c'')$ belong to the same connected component of the category $(d_{i}\downarrow D^{A}_{S'})$.   
	\end{enumerate}

\end{remarks}

Recalling from section \ref{sec:smallestgrothendiecktopologycomorphism} the definition of the smallest Grothendieck topology making a functor a comorphism of sites, we deduce from Proposition \ref{propcharJcanonical}(iv) the following result:

\begin{proposition}
	Let $A:{\cal C}\to {\cal D}$ be a functor and $K$ a Grothendieck topology on $\cal D$. Then for any Grothendieck topology $J$ on $\cal C$ such that $M^{A}_{K}\subseteq J\subseteq Z_{A}$, the functor $A$ is a $(J, K)$-continuous comorphism of sites $({\cal C}, J)\to ({\cal D}, K)$. 
	We have $M^{A}_{K}\subseteq Z_{A}$ if and only if for any arrow $g:d\to c$ of $\cal C$ and any $K$-covering sieve $R$ on $A(c)$, 
	\[
	A(d)=\varinjlim_{h:e\to d \in g^{\ast}(S^{A}_{R})} A(e).
	\] 	   
\end{proposition}\qed 

The following example shows that a comorphism of sites $({\cal C}, J) \to  ({\cal D}, K)$ can be cover-preserving without being $(J, K)$-continuous.

\begin{example}
	Let $L$ be a locale. Then the comorphism of sites $(L, J^{\textup{can}_{L}})\to (\textbf{2}, J^{\textup{can}}_{\textbf{2}})$ sending $0$ to $0$ and every non-zero element of $L$ to $1$ is cover-preserving but not $(J^{\textup{can}}_{L}, J^{\textup{can}}_{\textbf{2}})$-continuous if $L$ is not locally connected (see the proof of Proposition \ref{propmisc}(iii) below).     
\end{example}

We shall denote by $[{\cal C}, {\cal E}]_{J}$ (resp. by $[{\cal C}, {\cal D}]_{J}$) the full subcategory of $[{\cal C}, {\cal E}]$ (resp. of $[{\cal C}, {\cal D}]$) on the $J$-continuous functors.

\begin{proposition}\label{propweakmorphismequivalence}
	Let $\cal C$ a locally small category and $\cal E$ a Grothendieck topos. 
	\begin{enumerate}[(i)]
		\item There is an equivalence
		\[
		\textup{\bf Wmor}({\cal E}, [{\cal C}^{\textup{op}}, \Set])\simeq [{\cal C}, {\cal E}]
		\]
		sending a weak morphism $f=(f^{\ast}\dashv f_{\ast})$ to the functor $f^{\ast}\circ y_{\cal C}$.
		
		\item For any Grothendieck topology $J$ on $\cal C$ making $({\cal C}, J)$ a small-generated site, the above equivalence restricts to an equivalence
		\[
		\textup{\bf Wmor}({\cal E}, \Sh({\cal C}, J))\simeq [{\cal C}, {\cal E}]_{J}
		\]
		sending a weak morphism $g=(g^{\ast}\dashv g_{\ast})$ to the functor $g^{\ast}\circ l$.
	\end{enumerate}
\end{proposition}

\begin{proofs}
	(i) One half of the equivalence is given by the assignment $A \mapsto (L_{A}\dashv R_{A})$; the other half sends a weak morphism $f$ to the functor $f^{\ast}\circ y_{\cal C}$. The two correspondences are quasi-inverses to each other since, on the one hand, for any functor $A:{\cal C}\to {\cal E}$, $L_{A}\circ y_{\cal C}\cong A$ and, on the other hand, for any weak morphism $f:{\cal E}\to [{\cal C}^{\textup{op}}, \Set]$, $f^{\ast}\cong L_{f^{\ast}\circ y_{\cal C}}$ since both functors preserve arbitrary colimits (as they have right adjoints) and their composites with the Yoneda embedding are isomorphic to the same functor (namely, $f^{\ast}\circ y_{\cal C}$).
	
	(ii) If $A$ is $J$-continuous then $R_{A}$ factors through the canonical inclusion $\Sh({\cal C}, J)\hookrightarrow [{\cal C}^{\textup{op}}, \Set]$. Regarded as a functor ${\cal E}\to \Sh({\cal C}, J)$, it preserves arbitrary limits (since $\Sh({\cal C}, J)$ is closed in $[{\cal C}^{\textup{op}}, \Set]$ under arbitrary limits), so it has a left adjoint and hence defines a weak morphism of toposes $g:{\cal E}\to \Sh({\cal C}, J)$ which, composed with the canonical geometric inclusion $i_{J}:\Sh({\cal C}, J)\hookrightarrow [{\cal C}^{\textup{op}}, \Set]$, gives the weak morphism induced by $A$; in particular, the inverse image $g^{\ast}:\Sh({\cal C}, J) \to {\cal E}$ of this morphism satisfies $g^{\ast}\circ a_{J}\cong L_{A}$. This also shows that the functor $A$ corresponding to $i_{J}\circ g$ is given by $g^{\ast}\circ l$; indeed, $A\cong L_{A}\circ y_{\cal C}\cong g^{\ast}\circ a_{J}\circ y_{\cal C}\cong g^{\ast}\circ l$. Conversely, given a weak morphism $f:{\cal E}\to [{\cal C}^{\textup{op}}, \Set]$ corresponding as in (i) to a functor $A:{\cal C}\to {\cal E}$ (so that $f^{\ast}\cong L_{A}$ and $f_{\ast}\cong R_{A}$), if $f$ factors through $\Sh({\cal C}, J)\hookrightarrow [{\cal C}^{\textup{op}}, \Set]$ then $f_{\ast}$ (that is, $R_{A}$) factors through the direct image of the geometric inclusion $\Sh({\cal C}, J)\hookrightarrow [{\cal C}^{\textup{op}}, \Set]$, that is, $A$ is $J$-continuous.	
\end{proofs}

Proposition \ref{propweakmorphismequivalence} motivates the following generalization of the notion of morphism of sites.

\begin{definition}
	Let $({\cal C}, J)$ and $({\cal D}, K)$ be small-generated sites. A functor $F:{\cal C}\to {\cal D}$ is said to be a \emph{weak morphism of sites} if it is $(J, K)$-continuous. 
\end{definition}

We can thus derive from Proposition \ref{propweakmorphismequivalence} the following result, whose first part was already established in Proposition III-1.2 \cite{grothendieck}:

\begin{proposition}
	Any weak morphism $F:({\cal C}, J)\to ({\cal D}, K)$ of small-generated sites induces a weak geometric morphism $\Sh(F):\Sh({\cal D}, K)\to \Sh({\cal C}, J)$ such that the following diagram commutes: 
	$$
	\xymatrix{
		{\mathcal C} \ar[r]^F \ar[d]_{l} &{\cal D} \ar[d]^{l'} \\
		\Sh({\cal C}, J) \ar[r]^{\Sh(F)^{\ast}} & \Sh({\cal D}, K)
	}
	$$ 
	
	Conversely, any weak geometric morphism $f=(f^{\ast}\dashv f_{\ast})$ such that $f^{\ast}\circ l$ factors through $l'$ is induced by a (necessarily unique, if $K$ is subcanonical) weak morphism of sites $({\cal C}, J)\to ({\cal D}, K)$.    
\end{proposition}\qed

\subsection{Continuous comorphisms of sites and essential geometric morphisms}\label{sec:essentialmorphismsandcomorphisms}

The following result provides an equivalence between continuous comorphisms of sites and essential geometric morphisms. 

\begin{theorem}\label{thmequivalencecomorphisms}
	Let $({\cal C}, J)$ be a small-generated site, $\cal E$ a Grothendieck topos. Let ${\textup{\bf Geom}}_{\textup{ess}}(\Sh({\cal C}, J), {\cal E})$ be the category of essential geometric morphisms, and $	{\textup{\bf Com}}_{\textup{cont}}(({\cal C}, J), ({\cal E}, J^{\textup{can}}_{\cal E}))$ the category of $J$-continuous comorphisms of sites $({\cal C}, J)\to ({\cal E}, J^{\textup{can}}_{\cal E})$. Then we have an equivalence
	\[
	{\textup{\bf Geom}}_{\textup{ess}}(\Sh({\cal C}, J), {\cal E}) \simeq 	{\textup{\bf Com}}_{\textup{cont}}(({\cal C}, J), ({\cal E}, J^{\textup{can}}_{\cal E}))
	\] 
	sending an essential geometric morphism $f=(f_{!}\dashv f^{\ast}\dashv f_{\ast})$ to the comorphism of sites $f_{!}\circ l$ and a $J$-continuous comorphism of sites $A$ to the geometric morphism $C_{A}$ induced by it. 
\end{theorem}

\begin{proofs}
	We have to prove that:
	\begin{enumerate}[(i)]
		\item For any $J$-continuous comorphism of sites $A:({\cal C}, J) \to ({\cal E}, J^{\textup{can}}_{\cal E})$, the geometric morphism $C_{A}$ it induces is essential.
		
		\item For any essential geometric morphism $\Sh({\cal C}, J) \to {\cal E}$, the functor $f_{!}\circ l$ is a $J$-continuous comorphism of sites $({\cal C}, J) \to ({\cal E}, J^{\textup{can}}_{\cal E})$.   
		
		\item Any essential geometric morphism $f:\Sh({\cal C}, J) \to {\cal E}$ is equal to $C_{f_{!}\circ l}$. 
		
		\item Any $J$-continuous comorphism of sites $A:({\cal C}, J) \to ({\cal E}, J^{\textup{can}}_{\cal E})$ satisfies $A\cong (C_{A})_{!}\circ l$.
	\end{enumerate}
	
	(i) By Proposition \ref{propweakmorphismequivalence}(ii), $A$ induces a pair of adjoint functors $L_{A} \dashv R_{A}$, where $R_{A}$ is the functor ${\cal E} \to \Sh({\cal C}, J)$ given by: $e\mapsto \textup{Hom}_{\cal E}(A(-), e)$. Now, $R_{A}$ corresponds to the inverse image of the geometric morphism $C_{A}$ induced by $A$ under the equivalence ${\cal E}\simeq \Sh({\cal E}, J^{\textup{can}}_{\cal E})$; indeed, by Theorem 4 at p. 412 of \cite{MM}, $C_{A}^{\ast}(P)\cong a_{J}(P\circ A^{\textup{op}})$ for any $P\in \Sh({\cal E}, J^{\textup{can}}_{\cal E})$, whence, by taking $P=\textup{Hom}_{\cal E}(-, e)$, $C_{A}^{\ast}(\textup{Hom}_{\cal E}(-, e))\cong R_{A}(e)$, as required. 
	
	(ii) Any essential geometric morphism of sites induces an adjunction $(f_{!}\dashv f^{\ast})$; so, by Proposition \ref{propweakmorphismequivalence}(ii), the functor $f_{!}\circ l$ is $J$-continuous. It remains to prove that $f_{!}\circ l$ is a $J$-continuous comorphism of sites $({\cal C}, J) \to ({\cal E}, J^{\textup{can}}_{\cal E})$. Since $f_{!}$ has a right adjoint, namely $f^{\ast}$, which preserves covers, it satisfies the covering lifting property; from this it immediately follows that $f_{!}\circ l$ also satisfies the covering lifting property.   
	
	(iii) By (ii), $f_{!}\circ l$ is a $J$-continuous comorphism of sites $({\cal C}, J) \to ({\cal E}, J^{\textup{can}}_{\cal E})$. So, by (i), the geometric morphism $C_{f_{!}\circ l}$ it induces is essential. To show that $f\cong C_{f_{!}\circ l}$, it thus suffices to show that the essential images of the two geometric morphisms are isomorphic. Now, we saw in the proof of (i) that the essential image $u$ of $C_{f_{!}\circ l}$ satisfies the relation $u\circ a_{J}\cong L_{f_{!}\circ l}$, while the essential image of $f$ is $f_{!}$. In order to prove that $f_{!}\cong L_{f_{!}\circ l}\circ a_{J}$, since both functors are colimit-preserving, it is equivalent to show that they give isomorphic functors when composed with $l$. But $u\circ l=u\circ a_{J}\circ y_{\cal C}\cong L_{f_{!}\circ l} \circ y_{\cal C}\cong f_{!}\circ l$, as required.  
	
	(iv) By Proposition \ref{propweakmorphismequivalence}(ii), $A:({\cal C}, J) \to ({\cal E}, J^{\textup{can}}_{\cal E})$ satisfies $A\cong (C_{A})_{!}\circ l$ if and only if $A$ and $(C_{A})_{!}\circ l$ induce the same weak morphisms ${\cal E}\to \Sh({\cal C}, J)$. Now, since $A$ and $(C_{A})_{!}\circ l$ are both $J$-continuous comorphisms of sites $({\cal C}, J) \to ({\cal E}, J^{\textup{can}}_{\cal E})$ (the latter by (i), which ensures that $C_{A}$ is essential, and then by (ii)), by the proof of (i) the inverse image of the weak morphism induced by $A$ is $(C_{A})_{!}$ and the inverse image of the weak morphism induced by $(C_{A})_{!}\circ l$ is $(C_{(C_{A})_{!}\circ l})_{!}$. But (iii) ensures that for any essential geometric morphism $f:\Sh({\cal C}, J) \to {\cal E}$, $f\cong C_{f_{!}\circ l}$, equivalently $f_{!}\cong (C_{f_{!}\circ l})_{!}$; applying this to $f=C_{A}$ thus yields our thesis.  
\end{proofs}

The following example shows that the condition for a comorphism of sites to be $J$-continuous is not necessary for the corresponding geometric morphism to be essential:

\begin{example}
	Let $({\cal C}, J)$ be any site and $\textbf{1}$ be the category with only one object and one morphism (the identity on it), endowed with the trivial Grothendieck topology $T$ on it. Then the unique functor $!_{\cal C}:{\cal C}\to \textbf{1}$ is clearly a comorphism of sites $({\cal C}, J)\to (\textbf{1}, T)$ and hence induces a geometric morphism $\Sh({\cal C}, J)\to \Sh(\textbf{1}, T)\simeq \Set$ which is necessarily the unique (up to isomorphism) geometric morphism $\Sh({\cal C}, J)\to \Set$. Now, the functor $R_{l'\circ !_{\cal C}}:\Set \to [{\cal C}^{\textup{op}}, \Set]$ clearly sends any set $E$ to the constant functor with value $E$. Now, all the constant functors are $J$-sheaves if and only if the site $({\cal C}, J)$ is locally connected (that is, every $J$-covering sieve on an object $c$ is connected as a full subcategory of the slice category ${\cal C}\slash c$). However, it is not necessary that $({\cal C}, J)$ be locally connected for the topos $\Sh({\cal C}, J)$ to be locally connected (in the sense that its unique geometric morphism to $\Set$ is essential); for instance, for any locally connected locale $L$, the canonical site $(L, J_{L}^{\textup{can}})$ is clearly not locally connected (the empty sieve covers the zero element of $L$) but the topos $\Sh(L)$ is locally connected (by Proposition C1.5.9 \cite{El}). We have thus found an example of a functor $!_{L}:L\to \textbf{1}$ which is cover-preserving but not canonical as a comorphism of sites $(L, J_{L}^{\textup{can}})\to (\textbf{1}, T)$, and which induces an essential geometric morphism. 
\end{example}

\begin{remark}\label{remdiffcomorphisms}
	Different comorphisms of sites $({\cal C}, J) \to ({\cal E}, J^{\textup{can}}_{\cal E})$ may induce the same geometric morphism. In fact, the equivalence relation $\sim$  given by $A\sim A'$ if and only if $C_{A}=C_{A'}$ can be explicitly characterized, by Theorem \ref{thmfunctionalrelations}, in terms of functional $J$-equivalence relations, as follows: $A\sim A'$ if and only if there are, for each $e\in {\cal E}$, $J$-functional relations between $\textup{Hom}_{\cal E}(A(-), e)$ and $\textup{Hom}_{\cal E}(A'(-), e)$ which are inverse to each other, naturally in $e\in {\cal E}$. In these terms, the proposition can be interpreted as saying that $J$-continuous comorphisms of sites are exactly the canonical representatives of equivalence classes of comorphisms of sites $({\cal C}, J)\to ({\cal E}, J^{\textup{can}}_{\cal E})$ inducing an essential geometric morphism $\Sh({\cal C}, J) \to {\cal E}$; more precisely, the $J$-continuous comorphism of sites representing the equivalence class of a comorphism of sites inducing an essential geometric morphism $f$ is given by $f_{!}\circ l$. In particular, a comorphism of sites $A:({\cal C}, J) \to ({\cal E}, J^{\textup{can}}_{\cal E})$ induces an essential geometric morphism if and only if there is a $J$-continuous comorphism of sites $A':({\cal C}, J) \to ({\cal E}, J^{\textup{can}}_{\cal E})$ such that $A\sim A'$. 
\end{remark}

\begin{corollary}\label{correstrictioncomorphism}
	Let $({\cal C}, J)$ and $({\cal D}, K)$ be small-generated sites. Then we have an equivalence between the essential geometric morphisms $f:\Sh({\cal C}, J)\to \Sh({\cal D}, K)$ such that $f_{!}\circ l:{\cal C}\to \Sh({\cal D}, K)$ factors through the canonical functor $l':{\cal D}\to \Sh({\cal D}, K)$ and the $(J, K)$-continuous comorphism of sites $({\cal C}, J) \to ({\cal D}, K)$, considered up to $K$-equivalence. 
\end{corollary}

\begin{remarks}
	
	\begin{enumerate}[(a)]
		\item If $A:({\cal C}, J)\to ({\cal D}, K)$ is both a morphism and a comorphism of sites then it is a $(J, K)$-continuous (by Proposition \ref{propcharJcanonical}(iii)) and hence by Corollary \ref{correstrictioncomorphism} the geometric morphism $C_{A}$ it induces is essential.
		
		\item The fact that a $(J, K)$-continuous comorphism of sites induces an essential geometric morphism was already proved as Proposition III-2.6 \cite{grothendieck}.
	\end{enumerate}

\end{remarks}

The following result provides several equivalent characterizations of the $(J, K)$-continuity property of a comorphism of sites $A:({\cal C}, J)\to ({\cal D}, K)$. Before stating it, we need to make a number of remarks.

The direct image of the geometric morphism $C_{A}$ induced by a comorphism of sites $A:({\cal C}, J)\to ({\cal D}, K)$ is the restriction to sheaves of the right Kan extension functor $\textup{Ran}_{A^{\textup{op}}}$ (see, for instance, Proposition 2.3.18 \cite{El}); in other words, the following diagram commutes:

\[
\xymatrix{
	[{\cal C}^{\textup{op}}, \Set]  \ar[r]^{\textup{Ran}_{A^{\textup{op}}}}  & [{\cal D}^{\textup{op}}, \Set]  \\
	\Sh({\cal C}, J) \ar@{^{(}->}[u]_{i_{J}}  \ar[r]^{(C_{A})_{\ast}} & \Sh({\cal D}, K) \ar@{^{(}->}[u]^{i_{K}}
} 
\] 
Recall that the functors $\textup{Ran}_{A^{\textup{op}}}$ and $\textup{Lan}_{A^{\textup{op}}}$ (namely, the right and left Kan extension functors along $A^{\textup{op}}$) are respectively right and left adjoints to the functor $D_{A}:=(-\circ A^{\textup{op}}):[{\cal D}^{\textup{op}}, \Set]\to [{\cal C}^{\textup{op}}, \Set]$.

By taking left adjoints to all the functors in the above square, we obtain that 
\[
C_{A}^{\ast} \circ a_{K}\cong a_{J}\circ D_{A}.
\]

Let us now suppose that the geometric morphism $C_{A}$ is essential, i.e. that $C_{A}^{\ast}$ has a left adjoint $(C_{A})_{!}$. Denoting by $\eta_{A}:1_{[{\cal C}^{\textup{op}}, \Set]}\to D_{A}\circ \textup{Lan}_{A^{\textup{op}}}$ the unit of the adjunction between $D_{A}$ and $\textup{Lan}_{A^{\textup{op}}}$, we have a morphism 
\[
((C_{A})_{!}\circ a_{J})\eta_{A}: (C_{A})_{!}\circ a_{J}\to (C_{A})_{!} \circ a_{J}\circ  D_{A}\circ \textup{Lan}_{A^{\textup{op}}}
\]
which, composed with the isomorphism
\[
(C_{A})_{!} \circ a_{J}\circ  D_{A}\circ \textup{Lan}_{A^{\textup{op}}}\cong (C_{A})_{!} \circ C_{A}^{\ast} \circ a_{K}\circ \textup{Lan}_{A^{\textup{op}}}
\]
induced by the isomorphism $C_{A}^{\ast} \circ a_{K}\cong a_{J}\circ D_{A}$, yields a morphism
\[
((C_{A})_{!}\circ a_{J})\eta_{A}:(C_{A})_{!} \circ a_{J} \to  (C_{A})_{!} \circ C_{A}^{\ast} \circ a_{K}\circ \textup{Lan}_{A^{\textup{op}}}.
\]
Composing this morphism with the morphism 
\[
\epsilon_{A} (a_{K}\circ \textup{Lan}_{A^{\textup{op}}}):(C_{A})_{!} \circ C_{A}^{\ast} \circ a_{K}\circ \textup{Lan}_{A^{\textup{op}}} \to a_{K}\circ \textup{Lan}_{A^{\textup{op}}} ,
\]
where   $\epsilon_{A}:(C_{A})_{!}\circ (C_{A})^{\ast} \to  1_{\Sh({\cal D}, K)}$ is the counit of the adjunction between $(C_{A})_{!}$ and $C_{A}^{\ast}$, thus yields a morphism
\[
(C_{A})_{!} \circ a_{J} \to a_{K}\circ \textup{Lan}_{A^{\textup{op}}},
\]
which we call $z_{A}$.

\begin{corollary}\label{corequivcharcanonicity}
	Let $A:({\cal C}, J)\to ({\cal D}, K)$ be a comorphism of sites. Then the following conditions are equivalent:
	\begin{enumerate}[(i)]
		\item $A$ is $(J, K)$-continuous.
		
		\item The left Kan extension functor $\textup{Lan}_{A^{\textup{op}}}:[{\cal C}^{\textup{op}}, \Set]\to [{\cal D}^{\textup{op}}, \Set]$ along $A^{\textup{op}}$ satisfies the property that $a_{K}\circ \textup{Lan}_{A^{\textup{op}}}$ factors (necessarily uniquely) through $a_{J}$.
		
		\item The geometric morphism $C_{A}$ induced by $A$ is essential and its essential image $(C_{A})_{!}$ makes the following diagram commute:
		\[
		\xymatrix{
			[{\cal C}^{\textup{op}}, \Set]  \ar[r]^{\textup{Lan}_{A^{\textup{op}}}} \ar[d]^{a_{J}}  & [{\cal D}^{\textup{op}}, \Set] \ar[d]^{a_{K}} \\
			\Sh({\cal C}, J) \ar[r]^{(C_{A})_{!}} & \Sh({\cal D}, K)
		} 
		\]		
	\end{enumerate}
	If $A$ induces an essential geometric morphism $C_{A}$ then there is a canonical morphism $(C_{A})_{!}\circ l \to l'\circ A$ (given by $z_{A}y_{\cal C}$, modulo the canonical isomorphism  $a_{K}\circ \textup{Lan}_{A^{\textup{op}}}\circ y_{\cal C}\cong l'\circ A$), and $A$ is $(J, K)$-continuous if and only if this morphism is an isomorphism, equivalently if and only if the canonical morphism 
	\[
	z_{A}:(C_{A})_{!} \circ a_{J} \to a_{K}\circ \textup{Lan}_{A^{\textup{op}}}
	\]
	defined above, is an isomorphism.
\end{corollary}

\begin{proofs}	
	(i) $\imp$ (ii). By Proposition \ref{propequivGroth}, $D_{A}$ restricts to a functor 
	$\Sh({\cal D}, K) \to \Sh({\cal C}, J)$. Then the following diagram commutes:
	\[
	\xymatrix{
		[{\cal D}^{\textup{op}}, \Set]  \ar[r]^{D_{A}}  & [{\cal C}^{\textup{op}}, \Set]  \\
		\Sh({\cal D}, K) \ar@{^{(}->}[u]_{i_{K}}  \ar[r]^{(C_{A})^{\ast}} & \Sh({\cal C}, J). \ar@{^{(}->}[u]^{i_{J}}
	} 
	\] 
	
	Indeed, $(C_{A})^{\ast}$ is given, by definition, by the functor $a_{J}\circ D_{A}\circ i_{K}$; so, if $D_{A}\circ i_{K}$ takes values in $\Sh({\cal C}, J)$ then $i_{J}\circ (C_{A})^{\ast}\cong D_{A}\circ i_{K}$.  
	
	Taking the left adjoints to all the functors in the above square (notice that $C_{A}$ is essential by Theorem \ref{thmequivalencecomorphisms}, since $A$ is $J$-continuous, equivalently satisfies condition (ii)) thus yields an isomorphism $a_{K}\circ \textup{Lan}_{A^{\textup{op}}} \cong (C_{A})_{!}\circ a_{J}$; in particular, $a_{K}\circ \textup{Lan}_{A^{\textup{op}}}$ factors through $a_{J}$, as required. 
	
	(ii) $\imp$ (iii) We first notice that, by the compatibility of adjoints with respects to composition, condition (ii) is equivalent to condition (ii) of Proposition \ref{propequivGroth} and hence, by that proposition, to condition (i). But we have already observed, in the proof that (i) $\imp$ (ii), that (i) actually implies (iii). 
	
	(iii) $\imp$ (i) By Theorem \ref{thmequivalencecomorphisms}, it suffices to show that $l'\circ A\cong (C_{A})_{!}\circ l$. But this immediately follows from the commutativity of the square in condition (iii), by composing both functors with $y_{\cal C}$. 
	
	It remains to prove the last part of the corollary. We preliminarily note that, since both functors $(C_{A})_{!} \circ a_{J}$ and $a_{K}\circ \textup{Lan}_{A^{\textup{op}}}$ preserve arbitrary colimits, $z_{A}$ is an isomorphism if and only if $z_{A}y_{\cal C}:(C_{A})_{!}\circ l\to l'\circ A$ is. If $z_{A}y_{\cal C}$ is an isomorphism then $A$ is $J$-continuous by Theorem \ref{thmequivalencecomorphisms}. Conversely, if $A$ is $J$-continuous then, as we have proved above, condition (iv) is satisfied, with the isomorphism $(C_{A})_{!} \circ a_{J} \cong a_{K}\circ \textup{Lan}_{A^{\textup{op}}}$ being precisely $z_{A}$.  
\end{proofs}

\begin{remark}
	From Corollaries \ref{correstrictioncomorphism} and \ref{corequivcharcanonicity} we may deduce the following criterion, alternative to that of Remark \ref{remdiffcomorphisms}, for a comorphism of sites $A:({\cal C}, J)\to ({\cal D}, K)$ to induce an essential geometric morphism: $C_{A}$ is essential if and only if there is a $J$-continuous comorphism of sites $A':({\cal C}, J)\to \Sh({\cal D}, K)$ and a natural transformation $\alpha:A'\to l'\circ A$ which is sent by $a_{K}$ to an isomorphism (notice that this latter condition admits an explicit reformulation obtained by applying Lemma \ref{lemmalift}(iii)). 
\end{remark}

The following result provides a natural relation between the sheafification operation and that of composition with the opposite of a functor inducing a comorphism of sites:

\begin{proposition}\label{propcompositionsheafificationcomorphism}
	Let $F:({\cal C}, J)\to ({\cal D}, K)$ be a comorphism of sites and $P$ a presheaf on $\cal D$. Then the canonical arrow
	\[
	P\circ F^{\textup{op}} \to a_{K}(P)\circ F^{\textup{op}}
	\]
	in $[{\cal C}^{\textup{op}}, \Set]$ is sent by the associated sheaf functor $a_{J}:[{\cal C}^{\textup{op}}, \Set] \to \Sh({\cal C}, J)$ to an isomorphism.
	
	If moreover $F$ is $(J, K)$-continuous (for instance, if $F$ is a morphism of sites $({\cal C}, J) \to ({\cal D}, K)$ -- cf. Proposition \ref{propcharJcanonical}(iii)) then the above arrow yields an isomorphism
	\[
	a_{J}(P\circ F^{\textup{op}})\cong a_{K}(P)\circ F^{\textup{op}}.
	\]
\end{proposition}

\begin{proofs}
	Let $(F^{\textup{op}})^{\ast}$ be the functor $D_{F}:=(-\circ F^{\textup{op}}):[{\cal D}^{\textup{op}}, \Set]  \to [{\cal C}^{\textup{op}}, \Set]$. Then the arrow in the statement of the proposition is precisely $\eta_{P}(F^{\textup{op}})^{\ast}$, where $\eta_{P}:P\to a_{K}(P)$ is the component at $P$ of the unit of the adjunction between $a_{K}$ and the inclusion functor of $K$-sheaves into presheaves. 
	
	As we observed in section \ref{sec:presheaflifting}, since $F$ is a comorphism of sites $({\cal C}, J) \to ({\cal D}, K)$, we have a commutative diagram
	\[
	\xymatrix{
		[{{\cal D}}^{\textup{op}}, \Set]  \ar[r]^{D_{F}} \ar[d]_{a_{K}} & [{\cal C}^{\textup{op}}, \Set]  \ar[d]^{a_{J}} \\
		\Sh({\cal D}, K)  \ar[r]^{C_{F}^{\ast}} & \Sh({\cal C}, J). 
	} 
	\]   
	From this it immediately follows that, since $\eta_{P}$ is sent by $a_{K}$ to an isomorphism, $a_{J}$ sends $\eta_{P}(F^{\textup{op}})^{\ast}$ to an isomorphism, as required. 
	
	Now, if $F$ is $(J, K)$-continuous then, by Proposition \ref{propequivGroth}, the presheaf $a_{K}(P)\circ F^{\textup{op}}$ is already a $J$-sheaf, whence the final assertion of the proposition follows.   	
\end{proofs}

Let us now consider localizations of geometric morphisms induced by a comorphism of sites.

For any comorphism of sites $F:({\cal C}, J)\to ({\cal D}, K)$ and any presheaf $Q$ on $\cal D$, we have a comorphism of sites 
\[
p^{F}_{Q}:(\textstyle \int (Q\circ F^{\textup{op}}), J_{Q\circ F^{\textup{op}}}) \to (\textstyle \int Q, J_{Q})
\]
sending every object $(c, x)\in \int (Q\circ F^{\textup{op}})$ to the object $(F(c), x)$ of $\int Q$ (and acting on arrows accordingly). 

\begin{proposition}\label{proplocalizationcanonicalcomorphism}
	Let $F:({\cal C}, J)\to ({\cal D}, K)$ be a comorphism of sites and $Q$ a presheaf on $\cal D$. The geometric morphism 
	\[
	C_{p^{F}_{Q}}:\Sh(\textstyle \int (Q\circ F^{\textup{op}}), J_{Q\circ F^{\textup{op}}}) \to \Sh(\textstyle \int Q, J_{Q})
	\]
	induced by $p^{F}_{Q}$ corresponds under the equivalences 
	\[
	\Sh({\textstyle \int Q\circ F^{\textup{op}}}, J_{Q\circ F^{\textup{op}}}) \simeq \Sh({\cal C}, J)\slash a_{J}(Q\circ F^{\textup{op}})
	\]
	and
	\[
	\Sh({\textstyle \int Q}, J_{Q}) \simeq \Sh({\cal D}, K) \slash a_{K}(Q)
	\]
	of section \ref{sec:example} to the geometric morphism
	\[
	C_{F}\slash a_{K}(Q):\Sh({\cal D}, K)\slash a_{K}(Q) \to \Sh({\cal C}, J)\slash C_{F}^{\ast}(a_{K}(Q))
	\]
	modulo the equivalence 
	\[
	\Sh({\cal C}, J)\slash C_{F}^{\ast}(a_{K}(Q))\cong \Sh({\cal C}, J)\slash a_{J}(Q\circ F^{\textup{op}})
	\]
	induced by the isomorphisms 
	\[
	C_{F}^{\ast}(a_{K}(Q))\cong a_{K}(Q)\circ F^{\textup{op}}\cong a_{J}(Q\circ F^{\textup{op}}),
	\] 
	the first of which results from the fact that $F$ is a continuous comorphism of sites $({\cal C}, J)\to ({\cal D}, K)$ and the second of which is that of Proposition \ref{propcompositionsheafificationcomorphism}.
\end{proposition}

\begin{proofs}
This is an immediate verification in light of the example of section \ref{sec:example}.
\end{proofs}

The following result shows that the property of a comorphism of sites to be continuous naturally behaves with respect to \ac relativization':

\begin{proposition}\label{propcancomorphismsstableunderpullback}
	Let $F:({\cal C}, J)\to ({\cal D}, K)$ be a $(J, K)$-continuous comorphism of sites. Then, for any presheaf $Q$ on $\cal D$, the comorphism of sites 
	\[
	p^{F}_{Q}:(\textstyle \int (Q\circ F^{\textup{op}}), J_{Q\circ F^{\textup{op}}}) \to (\textstyle \int Q, J_{Q})
	\]
	is $(J_{Q\circ F^{\textup{op}}}, J_{Q})$-continuous. 
\end{proposition}

\begin{proofs}
	Since $F$ is continuous, by Corollary \ref{correstrictioncomorphism} the geometric morphism $C_{F}:\Sh({\cal C}, J)\to \Sh({\cal D}, K)$ is essential. 
		  
	So, by Proposition \ref{proplocalizationcanonicalcomorphism}, the geometric morphism $C_{p^{F}_{Q}}$ induced by $p^{F}_{Q}$ is also essential, as it is (isomorphic to) a localisation of the essential geometric morphism $C_{F}$.
	
	By Corollary \ref{correstrictioncomorphism}, we can prove that $p^{F}_{Q}$ is continuous by showing that the diagram
		$$
	\xymatrix{
		\int (Q\circ F^{\textup{op}}) \ar[r]^{\quad \quad p^{F}_{Q}\,} \ar[d]_{l_{\int (Q\circ F^{\textup{op}})}} &{\int Q} \ar[d]^{l_{\int Q}} \\
		\Sh(\int (Q\circ F^{\textup{op}}), J_{Q\circ F^{\textup{op}}}) \ar[r]^{\quad\quad(C_{p^{F}_{Q}})_{!}} & \Sh(\int Q, J_{Q})
	}
	$$ 
	commutes.
	
	Now, consider the following cube:
		\begin{center}
		\begin{tikzcd}[row sep=large, column sep=small]
		& {\cal C} \arrow[dd, "l"{yshift=4ex}, rightarrow] \arrow[rr, "F"] &                                                                             & {\cal D} \arrow[dd, "l'"{yshift=2ex}, rightarrow] \\
		\int (Q\circ F^{\textup{op}}) \arrow[dd] \arrow[rr, "p^{F}_{Q}"{xshift=-4ex}, crossing over] \arrow[ru, "\pi_{Q\circ F^{\textup{op}}}"] &                                                                                             & {\int Q} \arrow[dd, crossing over]  \arrow[ru, "\pi_{Q}"] &                                                            \\
		& \Sh({\cal C}, J) \arrow[rr, "(C_{F})_{!}\quad\quad"{below, xshift=-3ex}]                                                &                                                                             & \Sh({\int Q}, J_{Q})                                 \\
		\Sh({\cal C}, J)\slash C_{F}^{\ast}(a_{K}(Q)) \arrow[ru, "V_{C_{F}^{\ast}(a_{K}(Q))}"'] \arrow[rr, "(C_{F}\slash a_{K}(Q))_{!}"']             &                                                                                             & \Sh({\cal D}, K)\slash a_{K}(Q) \arrow[ru, "V_{a_{K}(Q)}"']     &                                                           
		\end{tikzcd}
	\end{center} 
	where $V_{a_{K}(Q)}$ and $V_{C_{F}^{\ast}(a_{K}(Q))}$ are obvious forgetful functors and the vertical unnamed arrows are the canonical ones modulo the equivalences of Proposition \ref{proplocalizationcanonicalcomorphism}. 
	
	In light of the Proposition, we have to prove that the frontal face of the cube commutes. But this follows from the fact that the functor $V_{a_{K}(Q)}$ reflects isomorphisms, since the top face of the cube commutes by definition of $p^{F}_{Q}$, the lateral faces commute by Example \ref{examplemorphismsofdiscretefibrations}, the bottom face commutes since the square consisting of the right adjoints of the given functors clearly commutes, and the back face commutes by Corollary \ref{correstrictioncomorphism} as by our hypothesis $F$ is continuous.
\end{proofs}

The following proposition characterizes the morphisms of sites which induce an essential geometric morphism whose essential image restrict to representables:

\begin{proposition}
	Let $({\cal C},J)$ and $({\cal D}, K)$ be small-generated subcanonical sites and $g:\Sh({\cal C},J) \to \Sh({\cal D}, K)$ be a geometric morphism induced by a morphism of sites $G:({\cal D}, K) \to ({\cal C},J)$. Then the following conditions are equivalent:
	\begin{enumerate}[(i)]
		\item $g$ is essential and its essential image restricts to the representables; 
		
		\item $G:{\cal D}\to {\cal C}$ has a left adjoint $F:{\cal C}\to {\cal D}$ which is $(J, K)$-continuous. 
	\end{enumerate} 
\end{proposition}

\begin{proofs}
	(i) $\imp$ (ii) Suppose that $g=\Sh(G)$ is essential. Let us denote by $F$ the restriction of $g_{!}$ to the representables. We know from Corollary \ref{correstrictioncomorphism} that $F$ is $(J, K)$-continuous.  It thus remains to prove that $F\dashv G$. But this follows from the fact that $g_{!}\dashv g^{\ast}$, since $F$ is the restriction of $g_{!}$ to the representables and $G$ is the restriction of $g^{\ast}$ to the representables.
	
	(ii) $\imp$ (i) By Proposition \ref{propadjointinduction}(iv), $C_{F}= \Sh(G)$. Since $F$ is $(J, K)$-continuous then, by Corollary \ref{correstrictioncomorphism}, $C_{F}=\Sh(G)$ is essential and its essential image restricts to the representables, as desired. 
\end{proofs}

\begin{remark}
	Proposition \ref{propmisc}(iii) shows that, if the essential image of $g$ does not restrict to representables, its restriction ${\cal C}\to \Sh({\cal D}, K)$ to $\cal C$ does not necessarily take the same values as $l'\circ F$, if a left adjoint $F$ to $G$ exists (in fact, this happens precisely when $F$ is not $(J, K)$-continuous).  
\end{remark}

Let us now prove a useful result concerning cofinality conditions (in the sense of Proposition \ref{procofinality}) in the context of a continuous comorphism of sites whose essential image of the associated geometric morphism is conservative.

\begin{proposition}\label{propcofinalcontinuous}
	Let $p:({\cal C}, J)\to ({\cal D}, K)$ be a continuous comorphism of small-generated sites such that $(C_{p})_{!}$ is conservative, and $F:{\cal A}\to {\cal C}$ and $F:{\cal A}'\to {\cal C}$ two functors to $\cal C$ related by a functor $\xi:{\cal A}\to {\cal A}'$ and a natural transformation $\alpha:F\to F'\circ \xi$. Then $(\xi, \alpha)$ satisfies the conditions of Proposition \ref{procofinality} if and only if $(\xi, p\alpha)$ does.
\end{proposition}

\begin{proofs}
The canonical arrow
\[
a_{J}(\tilde{\alpha}): \textup{colim}_{\Sh({\cal C}, J)}(l\circ F) \to \textup{colim}_{\Sh({\cal C}, J)}(l\circ F) 
\]
in $\Sh({\cal C}, J)$ induced by $(\xi, \alpha)$ correponds, under the isomorphisms
\[
(C_{p})_{!}(\textup{colim}_{\Sh({\cal C}, J)}(l\circ F))\cong \textup{colim}_{\Sh({\cal D}, K)}(C_{p})_{!}\circ l\circ F)\cong \textup{colim}_{\Sh({\cal D}, K)}(l'\circ p\circ F)
\]
and 
\[
(C_{p})_{!}(\textup{colim}_{\Sh({\cal C}, J)}(l\circ F'))\cong \textup{colim}_{\Sh({\cal D}, K)}(C_{p})_{!}\circ l\circ F')\cong \textup{colim}_{\Sh({\cal D}, K)}(l'\circ p\circ F')
\]
resulting from the fact that $(C_{p})_{!}$ preserves colimits and that $(C_{p})_{!}\circ l\cong l'\circ p$ (cf. Corollary \ref{correstrictioncomorphism}), to the canonical arrow
\[
a_{K}(\tilde{p\alpha}): \textup{colim}_{\Sh({\cal D}, K)}(l'\circ p\circ F) \to \textup{colim}_{\Sh({\cal D}, K)}(l'\circ p\circ F) 
\]
in $\Sh({\cal D}, K)$ induced by $(\xi, p\alpha)$. The thesis thus follows from the conservativity of $(C_{p})_{!}$. 
\end{proofs}

\begin{remark}\label{remlocalhomeocofinality}
	Any a fibration $p$ of the form $\pi_{Q}:{\int Q}\to {\cal D}$, where $Q$ is a presheaf on $\cal D$ satisfies the hypotheses of the proposition (that is, $(C_{p})_{!}$ is conservative). Indeed, the geometric morphism $C_{p}$ is isomorphic (over $\Sh({\cal D}, K)$) to the local homeomorphism $\Sh({\cal D}, K)\slash a_{K}(Q)$, and hence its essential image corresponds to the forgetful functor $\Sh({\cal D}, K)\slash a_{K}(Q)\to \Sh({\cal D}, K)$, which is clearly conservative.  
\end{remark}

\subsubsection{Generalized compactness conditions}\label{subsec:generalizedcompactnessconditions}

In relation to Theorem \ref{thmequivalencecomorphisms}, it is worth remarking that, as shown by the following proposition, the essential images of essential geometric morphisms preserve a wide class of generalized compactness properties. Hence, in the presence of subcanonical sites $({\cal C}, J)$ and $({\cal D}, K)$ such that the objects of the form $l(c)$ (resp. $l'(d)$) can be characterized as the objects of the topos $\Sh({\cal C}, J)$ (resp. $\Sh({\cal D}, K)$) which satisfy such compactness properties, they restrict along $l$ and $l'$, thus yielding ($(J, K)$-continuous) comorphisms of sites $({\cal C}, J)\to ({\cal D}, K)$ inducing these morphisms. This is relevant in connection with the method, introduced in \cite{OC10}, for constructing dualities or equivalences by means of functorializing topos-theoretic \ac bridges'.

We say that a property $P$ of presieves is \emph{stable} if whenever $f:{\cal F}\to {\cal E}$ is a weak  morphism of toposes, if a presieve $S$ satisfies $P$ then $f^{\ast}(S)$ satisfies $P$. Examples of stable properties include to be finite, to be a singleton, to contain the identity arrow, etc. 

We say that a family $R$ of arrows to a given object $e$ of a topos $\cal E$ \emph{refines} a family $S$ of arrows to $e$ in $\cal E$ if every arrow in $R$ factors through an arrow in $S$.  Given a property $P$ of presieves, we say that an object $e$ of a topos is \emph{$P$-compact} if every epimorphic family of arrows to $e$ admits an epimorphic refinement satisfying property $P$.

\begin{proposition}\label{proppreservationcompactness}
	Let $P$ be a stable property of presieves and let $f$ be an essential geometric morphism. Then $f_{!}$ sends $P$-compact objects to $P$-compact objects.
\end{proposition}

\begin{proofs}
	Let $f:{\cal F}\to {\cal E}$ be an essential geometric morphism and $a$ a $P$-compact object of $\cal F$. We want to prove that $f_{!}(a)$ is $P$-compact. Let $\{e_{i}:u_{i}\to f_{!}(a) \mid i\in I\}$ be an epimorphic family to $f_{!}(a)$. Then the family $\{f^{\ast}(e_{i}):f^{\ast}(u_{i})\to f^{\ast}(f_{!}(a)) \mid i\in I\}$ is also epimorphic, and so is its pullback along the unit arrow $\eta_{a}:a\to f^{\ast}(f_{!}(a))$. In other words, given the pullbacks
	\[
	\xymatrix{
		c_{i}  \ar[r]^{\xi_{i}} \ar[d]_{s_{i}} & a \ar[d]^{\eta_{a}} \\
		f^{\ast}(u_{i}) \ar[r]^{f^{\ast}(e_{i})} & f^{\ast}(f_{!}(a))
	} 
	\]
	the family $\{\xi_{i}:c_{i}\to a\}$ is epimorphic. Since $a$ is $P$-compact, there is an epimorphic family of arrows $\{\chi_{k}:b_{k}\to a \mid k\in K\}$ satisfying $P$ and, for each $k$ an index $i_{k}$ and a factorization $\xi_{i_{k}}\circ z_{k}$ of $\chi_{k}$ through $\xi_{i_{k}}$. Therefore the image $\{f_{!}(\chi_{k}):f_{!}(b_{k})\to f_{!}(a) \mid k\in K\}$ under $f_{!}$ of this family is epimorphic (since $f_{!}$ preserves colimits) and satisfies $P$. It remains to show that this family refines $\{e_{i}:u_{i}\to f_{!}(a) \mid i\in I\}$. For each $k\in K$, let $w_{k}$ be the arrow $s_{i_{k}}\circ z_{k}:b_{k}\to f^{\ast}(u_{i_{k}})$, and $\overline{w_{k}}:f_{!}(b_{k})\to u_{i_{k}}$ be its transpose along the adjunction $f_{!}\dashv f^{\ast}$. Let us show that, for each $k\in K$, $f_{!}(\chi_{k})=e_{i_{k}}\circ \overline{w_{k}}$; this will establish our thesis. By applying $f^{\ast}$ to both sides and composing the obtained arrow with the unit $\eta_{b_{k}}:b_{k}\to f^{\ast}(f_{!}(b_{k}))$, our equality rewrites as $f^{\ast}(f_{!}(\chi_{k}))\circ \eta_{b_{k}} =  f^{\ast}(e_{i_{k}}\circ \overline{w_{k}})\circ \eta_{b_{k}}$. But $f^{\ast}(f_{!}(\chi_{k}))\circ \eta_{b_{k}}=\eta_{a}\circ \chi_{k}$, while $f^{\ast}(e_{i_{k}}\circ \overline{w_{k}})\circ \eta_{b_{k}}=f^{\ast}(e_{i_{k}})\circ f^{\ast}(\overline{w_{k}})\circ \eta_{b_{k}}=f^{\ast}(e_{i_{k}})\circ w_{k}=f^{\ast}(e_{i_{k}}) \circ  s_{i_{k}}\circ z_{k}=\eta_{a}\circ \xi_{i_{k}}\circ z_{k}=\eta_{a}\circ \chi_{k}$.      
\end{proofs}

While we are on the topic of $P$-compact objects, let us record a result about them which will be useful too in connection with the method for generating dualities through \ac bridges' introduced in \cite{OC10}.

Recall that a sieve $S$ on an object $c$ of a category $\cal C$ is said to be \emph{effective-epimorphic} if the canonical cone with vertex $c$ under
the (large) diagram $D_{S}$ given by the canonical projection ${\int S}\to {\cal C}$ is a colimit. We shall say that a family $T$ of arrows to $c$ is effective-epimorphic if the sieve $S_{T}$ generated by it is effective-epimorphic. For a family $T$ of arrows to $c$, we shall denote the diagram $D_{S_{T}}$ associated with the sieve $S_{T}$ generated by $T$ simply by $D_{T}$. A sieve in $\cal C$ is said to be \emph{universally effective-epimorphic} if its pullback along arbitrary arrows in $\cal C$ is effective-epimorphic. The canonical topology on a category $\cal C$, defined as the largest one for which all the representable functors are sheaves, can be characterized as the topology whose covering sieves are precisely the universally effective-epimorphic ones (see, for instance, pp. 542-4 of \cite{El}).  

\begin{lemma}\label{lemmainducedcanonicaltopology}
	Let $\cal C$ be a separating set for a Grothendieck topos $\cal E$ and  $S$ a sieve in $\cal C$ (regarded as a full subcategory of $\cal E$) such that the colimit in $\cal E$ of the diagram $D_{\overline{S}}$, where $\overline{S}$ is the sieve generated by $S$ in $\cal E$, lies in $\cal C$. Then the following conditions are equivalent:
	\begin{enumerate}[(i)]
		\item $S$ is $J_{\cal E}^{\textup{can}}|_{\cal C}$-covering. 
		
		\item $S$ is effective-epimorphic in $\cal C$.
		
		\item $S$ is universally effective-epimorphic in $\cal C$.
	\end{enumerate} 
\end{lemma}

\begin{proofs}
	(i) $\imp$ (iii) This follows from the fact that $J_{\cal E}^{\textup{can}}|_{\cal C}$ is subcanonical and hence contained in the canonical topology on $\cal C$.
	
	(iii) $\imp$ (ii) This is obvious.
	
	(ii) $\imp$ (i) If $S$ is effective-epimorphic in $\cal C$ then $c$ is the colimit of $D_{S}$ in $\cal C$. Now, since the colimit of $D_{\overline{S}}$ in $\cal E$ lies in $\cal C$ by our hypothesis, it necessarily yields a colimit of $D_{S}$ in $\cal C$ and hence it coincides with $c$. Indeed, since $\cal C$ is separating for $\cal E$, any cone over $D_{S}$ in $\cal C$ can be uniquely extended to a cone over $D_{\overline{S}}$ in $\cal E$ (with the same vertex). So $\overline{S}$ is epimorphic in $\cal E$, that is, $S\in J_{\cal E}^{\textup{can}}|_{\cal C}(c)$.
\end{proofs}

We shall use the notation $\ast$ to denote the operation of multicomposition of families of arrows, defined as
\[
{\cal F}\ast \{{\cal G}_{i} \mid i\in I \}=\{f_i\circ g^{i}_{j} \mid i\in I, j\in J_{i}\}
\]
for any family of arrows ${\cal F}:=\{f_i:\textup{dom}(f_i) \to c \mid i\in I\}$ and families of arrows ${\cal G}_{i}:=\{g^{i}_{j}:\textup{dom}(g^{i}_{j}) \to \textup{dom}(f_i) \mid j\in J_{i}\}$ for each $i\in I$.

Before stating the following proposition, we need to introduce some terminology.

We shall say that a property $P$ of families of arrows in a topos with common codomain is \emph{fractal} if it is stable under arbitrary multicompositions, and that it is \emph{stable under pullback} if for any family $\{f_{i}:\textup{dom}(f_i)\to u \mid i\in I\}$ satisfying $P$ and any arrow $g:v\to u$, the family $\{g^{\ast}(f_i):\textup{dom}(g^{\ast}(f_i)) \to v \mid i\in I \}$ also satisfies $P$. We shall say that $P$ is \emph{hereditary} if whenever we have families $\{g_i:a_i \to c \mid i\in I\}$ and $\{f_i:b_i \to c \mid i\in I\}$ and monomorphisms $\{m_i:a_i \mono b_i \mid i\in I\}$ such that $f_i\circ m_i=g_i$ for every $i \in I$, if $\{f_i:b_i \to c \mid i\in I\}$ satisfies $P$ then $\{g_i:a_i \to c \mid i\in I\}$ also satisfies $P$; notice that this is a form of stability under \ac pointwise' refinement. We shall say that $P$ \emph{relativizes to a separating set} $\cal C$ of a Grothendieck topos $\cal E$ if there is a property $P|_{\cal C}$ of families of arrows in $\cal C$ (regarded as a full subcategory of $\cal E$) with common codomain such that such a family $\cal F$ satisfies $P|_{\cal C}$ if and only if $\cal F$, regarded as a family of arrows in $\cal E$, satisfies $P$. We say that $P$ \emph{relativizes to separating sets} if it relativizes to every separating set of a Grothendieck topos.

\begin{proposition}
	Let $P$ be a property of presieves in a topos, $\cal E$ a Grothendieck topos and ${\cal E}^{c}_{P}$ the full subcategory of $\cal E$ on the $P$-compact objects. Suppose that 
	\begin{enumerate}[(1)]
		\item  ${\cal E}^{c}_{P}$ is separating for $\cal E$,
		\item  $P$ is fractal and relativizes to separating sets (or at least to the separating set ${\cal E}^{c}_{P}$ of $\cal E$),
		\item for any effective-epimorphic family ${\cal R}=\{r_i:\dom(r_i) \to c \mid i\in I\}$ in ${\cal E}^{c}_{P}$ satisfying $P|_{{\cal E}^{c}_{P}}$, the family of canonical colimit arrows $\{\lambda_i:\dom(r_i) \to e \mid i\in I\}$, where $e$ is the colimit of the diagram $D_{{\cal R}}$ in $\cal E$, also satisfies $P$.
	\end{enumerate} Then the Grothendieck topology $J^{c}_{P}$ on ${\cal E}^{c}_{P}$ induced by the canonical topology on $\cal E$ can be intrinsically characterized (in terms of ${\cal E}^{c}_{P}$) as follows: the $J^{c}_{P}$-covering sieves on an object $c$ of ${\cal E}^{c}_{P}$ are precisely those which contain a family of arrows to $c$ which satisfies $P|_{{\cal E}^{c}_{P}}$ and is effective-epimorphic (equivalently, universally effective-epimorphic). 
	
	Conditions (2) and (3) can be replaced by the following condition: 
	\begin{enumerate}
		\item[(4)] The category ${\cal E}^{c}_{P}$ is closed in $\cal E$ under colimits of diagrams of the form $D_{\cal R}$, where $\cal R$ is an effective-epimorphic family of arrows in ${\cal E}^{c}_{P}$ satisfying $P|_{{\cal E}^{c}_{P}}$.
	\end{enumerate} 
	which in fact implies (3), and is implied by (3) assuming (1) and (2).
	
	Condition (3) is satisfied if $P$ is stable under pullback and hereditary. 
\end{proposition} 

\begin{proofs}
	Since the topology $J^{c}_{P}$ is subcanonical, it is contained in the canonical topology on ${\cal E}^{c}_{P}$. In particular, all its covering sieves are universally effective-epimorphic. To show that every $J^{c}_{P}$-covering sieve contains a family of arrows to $c$ which satisfies $P$ and which is universally effective-epimorphic, it suffices to prove that every $J^{c}_{P}$-covering sieve contains a $J^{c}_{P}$-covering family which satisfies property $P$. But this follows at once from the $P$-compactness of $c$. 
	
	Conversely, we want to show that if a sieve $S$ on an object $c$ of ${\cal E}^{c}_{P}$ contains a family of arrows to $c$ which satisfies $P$ and which is effective-epimorphic then $S$ is $J^{c}_{P}$-covering. This will follow at once from Lemma \ref{lemmainducedcanonicaltopology} once we have proved that if $R$ is a sieve generated by a family of arrows $\{r_{i} \mid i\in I\}$ in ${\cal E}^{c}_{P}$ satisfying $P$ then the colimit $e$ in $\cal E$ of the diagram $D_{\overline{R}}$ lies in $\cal C$. We have to show that $e$ is $P$-compact. The universal property of the colimit yields colimit arrows $\{\lambda_{i}:\textup{dom}(r_i)\to E \mid i\in I\}$ and an arrow $\alpha:e\to c$ such that $\alpha \circ \lambda_{i}=r_{i}$ for each $i\in I$.
	By condition (3), the family $\{\lambda_{i}:\textup{dom}(r_i)\to e \mid i\in I\}$ satisfies property $P$. 
	
	Let ${\cal F}$ be an epimorphic family on $e$, which can suppose to be a sieve without loss of generality. We want to prove that it admits an epimorphic refinement satisfying $P$. For each $i\in I$, the family $\lambda_{i}^{\ast}({\cal F})$ is epimorphic on $\dom(r_i)$. Since $\dom(r_i)$ is $P$-compact, $\lambda_{i}^{\ast}({\cal F})$ is refined by an epimorphic family of arrows ${\cal F}_{i}$ to $\textup{dom}(r_{i})$ satisfying $P$. So the multicomposite family $\{\lambda_i \mid i\in I\} \ast \{{\cal F}_{i} \mid i\in I\}$ satisfies $P$ (since $P$ is fractal). Since this family is also epimorphic (notice that the family $\{\lambda_i \mid i\in I\}$ is epimorphic by the universal property of the colimit) and clearly refines $\cal F$, our thesis is proved.  
	
	By the above arguments, if conditions (1) and (2) are satisfied, then condition (3) implies condition (4). Conversely, (4) implies (3). Indeed, if ${\cal R}=\{r_i:\dom(r_i) \to c \mid i\in I\}$ is such a family then the family of canonical colimit arrows $\{\lambda_i:\dom(r_i) \to e \mid i\in I\}$ corresponds to $\cal R$ under the identification $e\cong c$ and hence \emph{a fortiori} satisfies $P$ since $\cal R$ does. 
	
	On the other hand, conditions (1) and (4) alone ensure, by Lemma \ref{lemmainducedcanonicaltopology}, that if a sieve $S$ on an object $c$ of ${\cal E}^{c}_{P}$ contains a family of arrows to $c$ which satisfies $P$ and which is effective-epimorphic then $S$ is $J^{c}_{P}$-covering, so they allow to conclude the thesis of the proposition.  
	
	Finally, let us show that if $P$ is stable under pullback and hereditary then condition (3) is satisfied.
	
	Consider the pullbacks $s_{i}$ of the arrows $r_{i}$ along $\alpha:e\to c$. By the universal property of the pullback, there is a unique arrow $m_{i}:\dom(r_i)\to u_i$ such that $s_i\circ m_i=\lambda_i$ and $t_i\circ m_i=1_{\dom(r_i)}$:
	
	\[
	\begin{tikzcd}
	\dom(r_i) \ar[ddr, "\lambda_{i}"', bend right] \ar[dr, "m_{i}"description]\ar[drr, "1_{\dom(r_i)}", bend left]&&\\
	&u_{i} \ar[d, "s_{i}"] \ar[r, "t_{i}"] & \textup{dom}(r_i) \ar[d, "r_i"] \\
	& e \ar[r, "\alpha"']& c \ar[ul, "\lrcorner"{xshift=-1.5ex, yshift=0.5ex, near end}, phantom]
	\end{tikzcd}
	\]
	
	It follows in particular that $m_i$ is a (split) monomorphism (for each $i\in I$).
	
	Since $P$ is stable under pullback and $\{r_{i} \mid i\in I\}$ satisfies $P$, the family ${\cal S}:=\{s_i \mid i\in I\}$ also satisfies $P$. Then, the fact that $P$ is hereditary implies that the family $\{\lambda_i \mid i\in I\}$ also satisfies $P$. 
\end{proofs}

\begin{remark}
	The proposition generalizes to any full dense subcategory $\cal C$ of ${\cal E}^{c}_{P}$ (to which $P$ relativizes) which is closed in $\cal E$ under colimits of diagrams of the form $D_{\cal R}$, where $\cal R$ is an effective-epimorphic family of arrows of $\cal C$ satisfying $P|_{\cal C}$, giving an intrinsic description of the induced topology $J^{\textup{can}}|_{\cal C}$ as the topology on $\cal C$ whose covering sieves are precisely those which contain a family of arrows to $c$ which satisfies $P|_{\cal C}$ and is effective-epimorphic (equivalently, universally effective-epimorphic).    
\end{remark}

\begin{examples}
	\begin{enumerate}[(a)]
		\item The property of a presieve to be finite (resp. a  singleton, to consist of pairwise disjoint arrows) is fractal, stable under pullback and hereditary. The first two properties also relativize to arbitrary separating sets.
		
		\item The property $P$ of being a singleton family only consisting of the identity arrow is fractal, stable under pullback, relativizes to separating sets but is not hereditary. Still, it clearly satisfies condition (4) (the given colimit is trivial). Notice that $P$-compactness coincides in this case with irreducibility. 
		
		\item Many important classes of toposes admit separating sets consisting of $P$-compact objects (for some property $P$ of families of arrows with common codomain). For instance, presheaf toposes are generated by $P$-compact objects where $P$ is the property \ac \ac to contain the identity'' (that is, by irreducible objects), regular toposes are generated by $P$-compact objects where $P$ is the property \ac\ac to be a singleton'' (that is, by supercompact objects) and coherent toposes are generated by $P$-compact objects where $P$ is the property \ac\ac to be finite'' (that is, by compact objects). Also, as shown in \cite{MorganRogers}, toposes of continuous actions of a topological monoid on discrete sets are generated by supercompact objects. 
	\end{enumerate}	
\end{examples}

\subsubsection{Continuous locale homomorphisms}

In this section we shall address the problem of characterizing the locale homomorphisms $L\to L'$ whose corresponding geometric morphism $\Sh(L)\to \Sh(L')$ is essential.

Recall that a frame is a complete partially ordered sets in which arbitrary joins distribute over finite meets. The category $\textbf{Loc}$ of \emph{locales} is the opposite of the category $\textbf{Frm}$ of frames and frame homomorphisms between them; so a locale is just a frame, but considered as an object of the category $\textbf{Loc}$. Notice that the canonical topology $J^{\textup{can}}_{L}$ on a frame $L$ has as covering families on a given element $l\in L$ the families $\{l_i \to l \mid i\in I\}$ such that $l=\bigvee_{i\in I}l_i$. A frame $L$ can be recovered from the corresponding topos $\Sh(L):=\Sh(L, J^{\textup{can}}_{L})$ as the frame of its subterminal objects. Since every frame homomorphism $L'\to L$ defines a morphism of sites $(L', J^{\textup{can}}_{L'}) \to (L, J^{\textup{can}}_{L})$, every locale homomorphism $u:L\to L'$ induces a geometric morphism $\Sh(u):\Sh(L)\to \Sh(L')$, from which one can recover the associated frame homomorphism $L'\to L$ as the restriction of $\Sh(u)^{\ast}$ to the subterminal objects. In fact, in this way the category $\textbf{Loc}$ fully embeds as a subcategory of the category of Grothendieck toposes and geometric morphisms between them; the toposes of the form $\Sh(L)$ for a locale $L$ are called \emph{localic}. We shall say that a locale homomorphism $f:L\to L'$ (or the associated frame homomorphism $f:L'\to L$) is \emph{complete} if $f:L'\to L$ preserves arbitrary meets; notice that this condition is equivalent, by the Special Adjoint Functor Theorem, to it having a left adjoint $L\to L'$; note that such a left adjoint need not be a frame homomorphism since it might not preserve finite meets.

First, we notice that the essential image of an essential geometric morphism does not in general restrict to subterminal objects. Take for example a locally connected but not connected topos $\cal E$; then the unique geometric morphism from $\cal E$ to $\Set$ is essential and its essential image sends an object of $\cal E$ to the set of its connected components - if applied to $1$, it yields the set of connected components of $1$, which is not a (subobject of the) singleton since $\cal E$ is by our hypothesis non-connected. Interestingly, not even essential geometric morphisms induced by functors between posets satisfy in general the property that their essential image sends subterminal objects to subterminal objects (in spite of the fact that it sends the subterminal \ac generators' to subterminal objects). Take for instance the unique functor $!_{\cal P}:{\cal P}\to \textbf{1}$. The topos $[{\cal P}, \Set]$ is locally connected (being a presheaf topos), and hence the essential image of the unique (up to isomorphism) geometric morphism from $[{\cal P}, \Set]$ to $\Set$ sends any object to the set of its connected components; now, the terminal object of $[{\cal P}, \Set]$ is not necessarily indecomposable (take, for instance, ${\cal P}$ equal to the discrete preorder with two elements), and hence it is not sent by the essential image to a subterminal object. 

We notice that if $A:({\cal C}, J)\to ({\cal D}, K)$ is a cover-preserving comorphism of sites then $A$ might be not $J$-continuous but its localic reflection satisfies the analogue of this property for locales. Indeed, while the functor $D_{A}$ might not restrict to all sheaves, it does restrict to subterminal objects, and sends subterminal sheaves to subterminal sheaves, that is, $J$-ideals to $K$-ideals.

This motivates the consideration, in the context of locales, of a weaker form of the notion of essential geometric morphism, namely the property of a morphism consisting in the fact that the restriction of the inverse image functor to subterminals admits a left adjoint (of frames):

\begin{definition}
	A geometric morphism $f:{\cal F}\to {\cal E}$ is said to be \emph{weakly essential} if the frame homomorphism $f^{\ast}|_{\Sub_{\cal E}(1_{\cal E})}:\Sub_{\cal E}(1_{\cal E}) \to \Sub_{\cal F}(1_{\cal F})$ admits a left adjoint.
\end{definition}

One might wonder whether any complete frame homomorphism $g:L\to L'$ induces an essential geometric morphism $\Sh(L')\to \Sh(L)$. The answer to this question is negative, as shown by the following proposition. 

\begin{proposition}\label{propmisc}
	\begin{enumerate}[(i)]
		\item Every essential geometric morphism $f:{\cal F}\to {\cal E}$ is weakly essential; in fact, more generally, all the frame homomorphisms of the form $f^{\ast}|_{\Sub_{\cal E}(A)}:\Sub_{\cal E}(A) \to \Sub_{\cal F}(f^{\ast}(A))$, where $A$ is an object of $A$, are complete (equivalently, they admit a left adjoint).
		
		\item For any locales $L$ and $L'$, the weakly essential geometric morphisms $\Sh(L)\to \Sh(L')$ correspond precisely to the complete frame homomorphisms $L'\to L$. 
		
		\item Let $u:L'\to L$ be a complete locale homomorphism whose corresponding geometric morphism $\Sh(u):\Sh(L')\to \Sh(L)$ is essential. Then the left adjoint to the frame homomorphism associated with $u$ does not in general take the same values as the restriction of the essential image of $\Sh(u)$ to the subterminal objects of $\Sh(L)$. 
		
		\item There are weakly essential morphisms between localic toposes which are not essential; in other words, the property of a frame homomorphism being complete is a necessary but not sufficient condition for the associated geometric morphism to be essential. 
		
		\item Any weakly essential morphism between toposes which can be represented as the toposes of presheaves on some posets is essential.  
	\end{enumerate}
\end{proposition} 

\begin{proofs}
	(i) It is easy to see that a left adjoint to $$f^{\ast}|_{\textup{Sub}_{\cal E}(A)}:\textup{Sub}_{\cal E}(A)\to \textup{Sub}_{\cal F}(f^{\ast}(A))$$ is given by $m\in \textup{Sub}_{\cal E}(A) \mapsto \textup{Im}(\epsilon_{A}\circ f_{!}(m))$, where $\epsilon_{A}:f_{!}(f^{\ast}(A))\to A$ is the component at $A$ of the counit of the adjunction $f_{!}\dashv f^{\ast}$.  
	
	(ii) This is clear by definition of weakly essential geometric morphism.
	
	(iii) Take $L=\textbf{2}=\{0,1\}$, $L'$ arbitrary and $u$ equal to the unique frame homomorphism $L \to L'$; then $u$ is clearly complete, and hence it has a left adjoint $L'\to 2$. By Proposition \ref{propadjointinduction}(i), this preorder map is therefore a comorphism of sites $(L', J^{\textup{can}}_{L'})\to (\textbf{2}, J^{\textup{can}}_{\textbf{2}})$, which induces the unique geometric morphism $\Sh(L')\to \Sh(\textbf{2}, J^{\textup{can}}_{\textbf{2}})\simeq \Set$. This comorphism sends $0$ to $0$ and any non-zero element of $L'$ to $1$. Now, if $L'$ is locally connected then this morphism is essential but if $L'$ is not connected then the essential image of this morphism sends the terminal object to the set of its connected components rather than to an element of $\textbf{2}$. 
	
	(iv) Take $L=\textbf{2}=\{0,1\}$ and $L'$ to be any not locally connected locale and $u$ equal to the unique frame homomorphism $L \to L'$. Then $u$ is complete but the geometric morphism $\Sh(u):\Sh(L')\to \Sh(\textbf{2})\simeq \Set$ which it induces is not essential (cf. point (iii) above).	
	
	(v) Let $f:{\cal F}\to {\cal E}$ be a weakly essential geometric morphism and ${\cal F}\simeq [{\cal Q}^{\textup{op}}, \Set]$, ${\cal E}\simeq [{\cal P}^{\textup{op}}, \Set]$, where $\cal P$ and $\cal Q$ are posets. Then ${\cal F}\simeq \Sh(\textup{Id}(\cal Q))$, ${\cal E}\simeq \Sh(\textup{Id}(\cal P))$, where $\textup{Id}(\cal P)$ and $\textup{Id}(\cal Q)$ are respectively the frame of ideals of $\cal P$ and of $\cal Q$. Indeed, these are exactly the frames of subterminal objects of $[{\cal P}^{\textup{op}}, \Set]$ and $[{\cal Q}^{\textup{op}}, \Set]$. Therefore the geometric morphism $f$ is induced by a frame homorphism $n:\textup{Id}(\cal P)\to \textup{Id}(\cal Q)$, which is given by the restriction of its inverse image $f^{\ast}$ to subterminals. Since $f$ is by our hypothesis weakly essential, it thus follows $n$ admits a left adjoint. This implies that the left adjoint to $n$ restricts to a poset homomorphism $u:{\cal Q}\to {\cal P}$ (cf. the proof of Proposition \ref{proppreservationcompactness} below) and hence that $n$ can be identified with the frame homomorphism $\textup{Id}(u):\textup{Id}({\cal P}) \to \textup{Id}({\cal Q})$. So $f$ can be identified with the geometric morphism induced by $u$, since they have isomorphic inverse image functors; in particular, $f$ is essential, as required. 
\end{proofs}

\subsection{Presheaf lifting of Grothendieck topologies}\label{sec:presheaflifting}

We have seen that, given two small-generated sites $({\cal C}, J)$ and $({\cal D}, K)$, not every essential geometric morphism $f:\Sh({\cal C}, J)\to \Sh({\cal D}, K)$ is induced by a comorphism of sites $({\cal C}, J)\to ({\cal D}, K)$; however, as shown by Corollary \ref{correstrictioncomorphism}, the only condition for this is that the composite of $f$ with the canonical functor $l:{\cal C}\to \Sh({\cal C}, J)$ factor through the canonical functor $l':{\cal D}\to \Sh({\cal D}, K)$. Now, as shown below, we can replace $({\cal D}, K)$ with a larger site giving rise to the same topos so that this condition is always satisfied and so every essential geometric morphism $\Sh({\cal C}, J)\to \Sh({\cal D}, K)$ can be seen as being induced by a comorphism from $({\cal C}, J)$ to this larger site. This involves the well-known operation of extension of a Grothendieck topology on a category to a Grothendieck topology on the corresponding presheaf topos.    

Recall from \cite{grothendieck} (Expos\'e II.5) that, for any small-generated site $({\cal C}, J)$, there is a Grothendieck topology $\hat{J}$ on the topos $[{\cal C}^{\textup{op}}, \Set]$ such that the site $([{\cal C}^{\textup{op}}, \Set], \hat{J})$ is small-generated and we have (by the Comparison Lemma) an equivalence
\[
\Sh({\cal C}, J) \simeq \Sh([{\cal C}^{\textup{op}}, \Set], \hat{J}).
\]
The Grothendieck topology $\hat{J}$ can be characterized as the topology coinduced by $J$ along the Yoneda embedding $y_{\cal C}:{\cal C}\to [{\cal C}^{\textup{op}}, \Set]$ (in the sense of Proposition \ref{propimagetopology}). Concretely, it has as covering sieves those sent by $a_{J}$ to epimorphic families.

Denoting by $l':[{\cal C}^{\textup{op}}, \Set]\to \Sh([{\cal C}^{\textup{op}}, \Set], \hat{J})$ the canonical functor, it sends any object of the form $l'(P)$, where $P\in [{\cal C}^{\textup{op}}, \Set]$, to $a_{J}(P)\in \Sh({\cal C}, J)$; in particular, the other half of the equivalence sends any $Q\in \Sh({\cal C}, J)$ to $l'(Q)$. Indeed, this follows at once from the general fact that for any (co)morphism of sites $i:({\cal D}, J_{\cal D})\to ({\cal C}, J)$ satisfying the hypotheses of the Comparison Lemma and inducing an equivalence $C_{i}:\Sh({\cal D}, J_{\cal D})\to \Sh({\cal C}, J)$, the diagram
\[
\xymatrix{
	\Sh({\cal C}, J)  \ar[r]^{C_{i}^{\ast}}_{\sim}  & \Sh({\cal D}, J_{\cal D}) \\
	[{\cal C}^{\textup{op}}, \Set] \ar[u]^{a_{J}}  \ar[r]^{-\circ i^{\textup{op}}} & [{\cal D}^{\textup{op}}, \Set] \ar[u]^{a_{J_{\cal D}}} 
} 
\] 
commutes, observing that $(-\circ y_{\cal C})\circ y_{[{\cal C}^{\textup{op}}, \Set]}\cong 1_{[{\cal C}^{\textup{op}}, \Set]}$.

The Yoneda embedding $y_{\cal C}:{\cal C}\to [{\cal C}^{\textup{op}}, \Set]$ and the associated sheaf functor $a_{J}$ respectively
define comorphisms of sites   
\[
y_{\cal C}:({\cal C}, J)\to ([{\cal C}^{\textup{op}}, \Set], \hat{J})
\]
and
\[
a_{J}:([{\cal C}^{\textup{op}}, \Set], \hat{J}) \to (\Sh({\cal C}, J), J^{\textup{can}}_{\Sh({\cal C}, J)})
\]
(cf. also Proposition \ref{propyonedamorcomor} below).

We shall call \emph{$J$-equivalence} the relation identifying two presheaves $P$ and $Q$ on $\cal C$ which have isomorphic $J$-sheafifications. By Theorem \ref{thmfunctionalrelations}, the relation of $J$-equivalence can be explicitly characterized in terms of $J$-functional relations: $P$ and $Q$ are $J$-equivalent if and only if there exist $J$-functional relations $R$ from $P$ to $Q$ and $R'$ from $Q$ to $P$ such that $R\circ R'$ is the identical relation on $Q$ and $R'\circ R$ is the identical relation on $P$.

As we had anticipated, we can profitably use this construction for representing essential morphisms to a topos $\Sh({\cal D}, K)$ of sheaves on a site $({\cal D}, K)$ as induced by comorphisms of sites towards a topos of sheaves on a site through which they factor. Indeed, by Theorem \ref{thmequivalencecomorphisms}, any essential geometric morphism $f:\Sh({\cal C}, J)\to \Sh({\cal D}, K)$ induces a comorphism of sites $({\cal C}, J)\to ([{\cal D}^{\textup{op}}, \Set], \hat{K})$ (given by $f_{!}\circ l$, regarded as taking values in $[{\cal D}^{\textup{op}}, \Set]$), since the composite of the latter with the canonical functor $[{\cal D}^{\textup{op}}, \Set]\to \Sh([{\cal D}^{\textup{op}}, \Set], \hat{K})$ is isomorphic to the composite of $f_{!}$ with the equivalence $\Sh({\cal D}, K) \simeq \Sh([{\cal D}^{\textup{op}}, \Set], \hat{K})$. Clearly, again by the above remark (this time applied in the converse direction), two comorphisms of sites $A, A':({\cal C}, J)\to ([{\cal D}^{\textup{op}}, \Set], \hat{K})$ induce the same geometric morphism $\Sh({\cal C}, J)\to \Sh({\cal D}, K) \simeq \Sh([{\cal D}^{\textup{op}}, \Set], \hat{K})$ if and only if $A(c)$ and $A'(c)$ are $K$-equivalent for any $c$, naturally in $c\in {\cal C}$. In particular, if $A=y_{\cal D}\circ F$ and $A'=y_{\cal D}\circ F'$ then $A$ and $A'$ induce the same geometric morphism if and only if, for any $c\in {\cal C}$, there is an isomorphism between the images of $F(c)$ and $F'(c)$ under the canonical functor $l':{\cal D}\to \Sh({\cal D}, K)$, natural in $c$; if $K$ is subcanonical then $A$ and $A'$ induce the same geometric morphism if and only if $F$ and $F'$ are isomorphic.

Let us now describe the functorial behaviour of the passage from a site $({\cal C}, J)$ to $([{\cal C}^{\textup{op}}, \Set], \hat{J})$.

\begin{proposition}\label{proppresheaflifting}
Any comorphism of sites $f:({\cal C}, J)\to ({\cal C}', J')$ extends to a comorphism of sites $\hat{f}:([{\cal C}^{\textup{op}}, \Set], \hat{J}) \to ([{\cal C}'^{\textup{op}}, \Set], \hat{J'})$ given by the left Kan extension along $f^{\textup{op}}$, and the following diagram commutes:
\[
\xymatrix{
	\Sh({\cal C}, J)  \ar[r] \ar[d]_{C_{f}} & \Sh([{\cal C}^{\textup{op}}, \Set], \hat{J}) \ar[d]^{C_{\hat{f}}} \\
	\Sh({\cal C}', J')  \ar[r] & \Sh([{{\cal C}'}^{\textup{op}}, \Set], \hat{J'})
} 
\]  	
\end{proposition}
\begin{proofs}
	The commutativity of the above diagram will follow from the functoriality of the geometric morphisms induced by comorphism of sites, once we have proved that $\hat{f}:([{\cal C}^{\textup{op}}, \Set], \hat{J}) \to ([{\cal C}'^{\textup{op}}, \Set], \hat{J'})$ is a comorphism of sites. Since $\hat{f}$ has a right adjoint, namely the inverse image functor $E(f)^{\ast}=D_{f}:=(-\circ f^{\textup{op}})$ of the essential geometric morphism $E(f):[{{\cal C}'}^{\textup{op}}, \Set] \to [{\cal C}^{\textup{op}}, \Set]$ induced by $f$, it is equivalent to show that this latter functor sends $\hat{J'}$-covering families to $\hat{J}$-covering families (cf. Proposition \ref{propadjointinduction}(i)). For this, we observe that, since $f$ is a comorphism of sites $({\cal C}, J)\to ({\cal C}', J')$, the direct image of the geometric morphism $C_{f}$ induced by it is given by the restriction to sheaves of the direct image of $E(f)$; passing to the left adjoints, we obtain the commutative diagram
	\[
	\xymatrix{
		[{{\cal C}'}^{\textup{op}}, \Set]  \ar[r]^{E(f)^{\ast}} \ar[d]_{a_{J'}} & [{\cal C}^{\textup{op}}, \Set]  \ar[d]^{a_{J}} \\
		\Sh({\cal C}', J')  \ar[r]^{C_{f}^{\ast}} & \Sh({\cal C}, J), 
	} 
	\] 
	from which it immediately follows that $E(f)^{\ast}$ sends $\hat{J'}$-covering families to $\hat{J}$-covering families.  
\end{proofs}

\begin{remark}\label{rempresheaflifting}
	As we observed in the proof of Proposition \ref{proppresheaflifting}, if $f:({\cal C}, J)\to ({\cal C}', J')$ is a comorphism of sites then the functor $D_{f}=-\circ f^{\textup{op}}:\hat{{\cal C}'}\to \hat{\cal C}$ is a morphism of sites $(\hat{{\cal C}'}, \hat{J'})\to (\hat{\cal C}, \hat{J})$ (since it clearly preserves finite limits and is cover-preserving), and the geometric morphism $\Sh(D_{f})$ can be identified with $C_{\hat{f}}$ (cf. Proposition \ref{propadjointinduction}).  
\end{remark}

We shall abbreviate by $\hat{\cal C}$ the category $[{\cal C}^{\textup{op}}, \Set]$ of presheaves on a category $\cal C$. 

As shown by the following result, the Yoneda embedding $y_{\cal C}:{\cal C}\to [{\cal C}^{\textup{op}}, \Set]$ can be used not only for turning comorphisms of sites towards a site with underlying category $\cal C$ to a site with underlying category $[{\cal C}^{\textup{op}}, \Set]$ but also for doing the same for morphisms of sites:

\begin{proposition}\label{propyonedamorcomor}
	Let $({\cal C}, J)$ be a small-generated site. Then the Yoneda embedding $y_{\cal C}:{\cal C}\to [{\cal C}^{\textup{op}}, \Set]$ is both a morphism and a comorphism of sites $({\cal C}, J)\to (\hat{\cal C}, \hat{J})$. It is moreover $\hat{J}$-dense (and $J$-full), so (cf. Proposition \ref{propmorphismcomorphismequivalence}) it induces an equivalence of toposes 	
	\[
	\Sh({\cal C}, J) \simeq \Sh(\hat{\cal C}, \hat{J})
	\]
	one half of which is $C_{y_{\cal C}}$ and whose other half is $\Sh(y_{\cal C})$.
\end{proposition}

\begin{proofs}
	Since a sieve $S$ on an object $c$ is $J$-covering if and only if it is sent by the associated sheaf functor $a_{J}$ to an isomorphism, $y_{\cal C}$ is a comorphism of sites $({\cal C}, J) \to (\hat{\cal C}, \hat{J})$. To see that it is also a morphism of sites  $({\cal C}, J) \to (\hat{\cal C}, \hat{J})$, we have to show that $l_{\hat{\cal C}}\circ y_{\cal C}$ is a $J$-continuous flat functor ${\cal C} \to \Sh(\hat{\cal C}, \hat{J})$. But the Comparison Lemma yields an equivalence
	\[
	\Sh({\cal C}, J) \simeq \Sh(\hat{\cal C}, \hat{J}),
	\]    
	under which this functor corresponds precisely to the functor $a_{J}\circ y_{\cal C}$, which is clearly flat and $J$-continuous. 
	
	It remains to show that $y_{\cal C}$ is $\hat{J}$-dense and $J$-full. But $y_{\cal C}$ is full and hence \emph{a fortiori} $J$-full, while the $\hat{J}$-denseness follows from the fact that every presheaf is a colimit of representables, whence it is $\hat{J}$-covered by them in $\hat{\cal C}$. 
\end{proofs}

\subsection{Fibrations as continuous functors}

In this section we prove the continuity of fibrations, regarded as comorphisms of sites as in section \ref{subsubsecfibrationscomorphisms}. We refer to section \ref{subsec:continuousfunctors} for the notation.

\begin{theorem}\label{thmfibration2}
	Let $p:{\cal C}\to {\cal D}$ be a fibration. Then, for any Grothendieck topology $K$ on $\cal D$, the comorphism of sites $p:(M^{p}_{K}, {\cal C})\to ({\cal D}, K)$ is continuous.
\end{theorem}

\begin{proofs}
	To show that $p:(M^{p}_{K}, {\cal C})\to ({\cal D}, K)$ is continuous, we shall apply the criterion of Proposition \ref{propcharcanonicalfunctors}. First, it is clear from the characterization of the topology $M^{p}_{K}$ given by Theorem \ref{thmfibration1} that $p$ is cover-preserving. So it remains to prove that, given a $M^{p}_{K}$-covering sieve $S$ on an object $c$ and a commutative diagram of the form 
	\[
	\xymatrix{
		d \ar[r]^{w} \ar[d]^{z}  & p(c') \ar[d]^{p(f)} \\
		p(c'') \ar[r]^{p(g)} & p(c),
	} 
	\]	
	where $f:c'\to c$ and $g:c''\to c$ are two arrows of $S$, there is a $K$-covering family $\{d_{i}\to d \mid i\in I\}$ such that for each $i\in I$, the composites $d_{i}\to p(c')$ and $d_{i}\to p(c'')$ belong to the same connected component of the category $(d_{i}\downarrow D^{p}_{S})$. Our $K$-covering family on $d$ will be given by the sieve $\chi^{\ast}(\langle p(P_{S})\rangle )$, where $P_{S}$ is the presieve consisting of the cartesian arrows in $S$ and $\chi$ is the arrow $p(f)\circ w=p(g)\circ z$ witnessing the commutativity of the above square. Notice that, by Theorem \ref{thmfibration1}, this sieve is indeed $K$-covering. We have to check that, for any $\alpha\in \chi^{\ast}(\langle p(P_{S})\rangle )$, the arrows $w\circ \alpha:\dom(\alpha)\to p(c')$ and $z\circ \alpha:\dom(\alpha)\to p(c'')$, regarded as objects of the category $(\dom(\alpha)\downarrow D^{p}_{S})$, belong to the same connected component. Since $\alpha\in \chi^{\ast}(\langle p(P_{S})\rangle )$, there is a cartesian arrow $\xi$ in $S$ and an arrow $\gamma$ such that $\chi\circ \alpha =p(\xi)\circ \gamma$. The arrow $\xi$ will actually serve as a means for connecting our two arrows; indeed, we will show that, on the one hand, $\gamma:\dom(\alpha)\to p(\dom(\xi))$ and $w\circ \alpha:\dom(\alpha)\to p(c')$ belong to the same connected component of $(\dom(\alpha)\downarrow D^{p}_{S})$ and that, on the other hand, $\gamma:\dom(\alpha)\to p(\dom(\xi))$ and $z\circ \alpha:\dom(\alpha)\to p(c'')$ belong to the same connected component of $(\dom(\alpha)\downarrow D^{p}_{S})$, from which it follows by transitivity that $w\circ \alpha:\dom(\alpha)\to p(c')$ and $z\circ \alpha:\dom(\alpha)\to p(c'')$ belong to the same connected component, as desired. 
	
	To prove that $\gamma:\dom(\alpha)\to p(\dom(\xi))$ and $w\circ \alpha:\dom(\alpha)\to p(c')$ belong to the same connected component of $(\dom(\alpha)\downarrow D^{p}_{S})$, we observe that, since $p$ is a fibration, there are an object $\tilde{c'}$ of $\cal C$, an isomorphism $\tau:\dom(\alpha)\to p(\tilde{c'})$ and an arrow $h:\tilde{c'}\to c'$ such that $p(h)\circ \tau=w\circ \alpha$. Next, we notice that, since $\chi=p(f)\circ w$ and $\chi\circ \alpha =p(\xi)\circ \gamma$, we have $p(f)\circ w\circ \alpha=p(\xi)\circ \gamma$ and hence $p(f)\circ p(h)\circ \tau=p(\xi)\circ \gamma$, equivalently, $p(f\circ h)=p(\xi)\circ \gamma \circ \tau^{-1}$. From the fact that $\xi$ is cartesian, it thus follows that there is an arrow $\psi:\tilde{c'}\to \dom(\xi)$ such that $p(\psi)=\gamma\circ \tau^{-1}$ and $\xi\circ \psi=f\circ h$. We thus see from the following diagram that $\gamma$ and $w\circ \alpha$ belong to the same component of $(\dom(\alpha)\downarrow D^{p}_{S})$ (where the vertical unnamed arrow is $p(f\circ h)=p(\xi\circ \psi)$):
	
	\begin{equation*}
	\begin{tikzcd} 
	& \dom(\alpha) \ar[d, "\tau"] \ar[ld, "\gamma"{above}] \ar[rd, "w\circ \alpha"] &  \\
	p(\dom(\xi)) \ar[dr, "p(\xi)\quad"{below}]   & p(\tilde{c'}) \ar[d] \ar[l, "p(\psi)"{above}] \ar[r, "p(h)"] &  p(c')  \ar[ld, "p(f)"]            \\
	& p(c)  & 
	\end{tikzcd}
	\end{equation*}	 	
	
	The proof that $\gamma:\dom(\alpha)\to p(\dom(\xi))$ and $z\circ \alpha:\dom(\alpha)\to p(c'')$ also belong to the same connected component of $(\dom(\alpha)\downarrow D^{p}_{S})$ is completely analogous. 
\end{proofs}

The following results enlighten a context where the continuity of a comorphism of sites can be deduced, under appropriate hypotheses, from that of a comorphism of sites of which it represents a \ac lifting' along a fibration.

\begin{lemma}\label{lemmacartesianrelmono}
	Let $C:{\cal D}\to {\cal D}'$ be a fibration.
	\begin{enumerate}[(i)]
		\item Let $f:d\to d'$ be a cartesian arrow. Then for any parallel arrows $h, k:d''\to d$, if $C(f\circ h)=C(f\circ k)$ then $h=k$.
		
		\item Every arrow of $\cal D$ can be functorially (and uniquely up to unique isomorphism) factored as a vertical arrow followed by a cartesian arrow.
		
		\item Let $A:{\cal C}\to {\cal D}$ be a functor and $S$ a sieve in $\cal C$ which supports a functorial factorization $g=g_{c}\circ g_{v}$ for any arrow $g$ in $S$, such that $g_{c}$ lies in $S$, $A(g_{c})$ is a cartesian arrow and $A(g_{v})$ is a vertical arrow. [Note that if $B:{\cal C}\to {\cal D'}$ is a fibration and $A$ is a morphism of fibrations from $B$ to $C$, then every sieve $S$ containing the cartesian images of all its arrows admits such a functorial factorization of arrows.] Then any two objects $(f:c' \to c, \alpha:d\to A(c'))$ and $(g:c''\to c, \beta:d\to A(c''))$ of the category $(d\downarrow D^{A}_{S})$ belong to the same connected component if and only if the objects $(f:c' \to c, C(\alpha):C(d)\to C(A(c')))$ and $(g:c''\to c, C(\beta):C(d)\to C(A(c'')))$ of the category $(C(d)\downarrow D^{C\circ A}_{S})$ belong to the same connected component.
	\end{enumerate}	 
\end{lemma}

\begin{proofs}
	(i) This is an immediate consequence of the uniqueness of the filling arrow in the definition of cartesian arrow.
	
	(ii) It is well-known that, since $C$ is a fibration, every arrow of $\cal D$ can be (uniquely up to unique isomorphism) factored as a vertical arrow followed by a cartesian arrow. It remains to show that this factorization is functorial. But this follows immediately from the universal property of cartesian arrows. 
	
	(iii) The \ac only if' direction is trivial, so we only have to show the converse one. So, let us suppose that $(f:c' \to c, C(\alpha):C(d)\to C(A(c')))$ and $(g:c''\to c, C(\beta):C(d)\to C(A(c'')))$ belong to the same connected component of the category $(C(d)\downarrow D^{C\circ A}_{S})$. Given a zig-zag connecting these two objects, we can obtain, by the functoriality of the factorizations of arrows in $S$, a zig-zag connecting the objects $(f_{c}:\dom(f_{c})\to c, C(A(f_{v})\circ \alpha):C(d)\to C(A(\dom(f_{c}))))$ and $(g_{c}:\dom(g_{c})\to c, C(A(g_{v})\circ \beta):C(d)\to C(A(\dom(g_{c}))))$, where the intermediate objects have the form $(h:\dom(h)\to c, u:C(d)\to C(A(\dom(h))))$ for a cartesian arrow $h$ in $S$. Now, it is easy to see, by means of an inductive argument starting from one or the other of the two extremes of the zig-zag connecting $(f_{c}:\dom(f_{c})\to c, C(A(f_{v})\circ \alpha):C(d)\to C(A(\dom(f_{c}))))$ and $(g_{c}:\dom(g_{c})\to c, C(A(g_{v})\circ \beta):C(d)\to C(A(\dom(g_{c}))))$, by using the cartesianness of the arrows $h$, that each arrow $u$ is actually the image under $C$ of an arrow $\gamma:d\to A(\dom(h))$ which makes the relevant triangle commute in the category $\cal D$. Indeed, suppose given two objects $(r:\dom(r)\to c, u:C(d')\to C(A(\dom(r))))$ and $(s:\dom(s)\to c, v:C(d')\to C(A(\dom(s))))$ of the category $(C(d')\downarrow D^{C\circ A}_{S})$, where the arrows $A(r)$ and $A(s)$ are cartesian, related by an arrow $t:\dom(s)\to \dom(r)$ in $\cal C$ such that $r\circ t=s$, in the sense the following diagram commutes:
	
	\begin{equation*}
	\begin{tikzcd} 
	& C(A(\dom(r)))   \ar[rd, "C(A(r))"] &  \\
	C(d') \ar[ur, "u"] \ar[dr, "v\quad"{below}]   &   &  C(A(c))              \\
	& C(A(\dom(s))) \ar[uu, "\quad\quad\quad C(A(t))"{below}]  \ar[ur, "\quad\quad\quad C(A(s))"{below}] & 
	\end{tikzcd}
	\end{equation*}	 
	
	If $v=C(\rho)$ for some arrow $\rho:d'\to A(\dom(s))$ then $u=C(A(t))\circ v=C(A(t)\circ \rho)$ and the following diagram commutes:
	
	\begin{equation*}
	\begin{tikzcd} 
	& A(\dom(r))   \ar[rd, "A(r)"] &  \\
	d' \ar[ur, "A(t)\circ \rho"] \ar[dr, "\rho\quad"{below}]   &   &  A(c)              \\
	& A(\dom(s)) \ar[uu, "\quad\quad A(t)"{below}]  \ar[ur, "\quad\quad A(s)"{below}] & 
	\end{tikzcd}
	\end{equation*} 
	
	If $u=C(\delta)$ for some arrow $\delta:d'\to A(\dom(r))$ then, by the cartesiannes of the arrow $A(s)$ and part (i) of the lemma (applied to the cartesian arrow $A(r)$), there is a unique arrow $\rho:d'\to A(\dom(s))$ such that $A(t)\circ \rho=\delta$. The following diagram thus commutes:
	\begin{equation*}
	\begin{tikzcd} 
	& A(\dom(r))   \ar[rd, "A(r)"] &  \\
	d' \ar[ur, "\delta"] \ar[dr, "\rho\quad"{below}]   &   &  A(c)              \\
	& A(\dom(s)) \ar[uu, "\quad\quad A(t)"{below}]  \ar[ur, "\quad\quad A(s)"{below}] & 
	\end{tikzcd}
	\end{equation*} 
	
	So in both cases we have a lift of the connection between the objects $(r:\dom(r)\to c, u:C(d')\to C(A(\dom(r))))$ and $(s:\dom(s)\to c, v:C(d')\to C(A(\dom(s))))$ from the category $D^{C\circ A}_{S}$ to the category $D^{A}_{S}$.
	
	So our zig-zag can be lifted to a zig-zag connecting the objects $(f_{c}:\dom(f_{c})\to c, A(f_{v})\circ \alpha:d\to A(\dom(f_{c})))$ and $(g_{c}:\dom(g_{c})\to c, A(g_{v})\circ \beta:d\to A(\dom(g_{c})))$, and hence between the objects $(f:\dom(f_{c})\to c, \alpha:d\to A(c'))$ and $(g:c''\to c, \beta:d\to A(c''))$, since $(f_{c}:\dom(f_{c})\to c, A(f_{v})\circ \alpha:d\to A(\dom(f_{c})))$ (resp. $(g_{c}:\dom(g_{c})\to c, A(g_{v})\circ \beta:d\to A(\dom(g_{c})))$) is in the same connected component as $(f:\dom(f_{c})\to c, \alpha:d\to A(c'))$ (resp. $(g:c''\to c, \beta:d\to A(c''))$).   
	
	The last statement of (iii) follows from (ii), since every morphism of fibrations sends cartesian arrows to cartesian arrows and vertical arrows to vertical arrows. 
\end{proofs}

\begin{proposition}\label{propliftcontinuity}
	Let $A:({\cal C}, J)\to ({\cal D}, K)$, $B:({\cal C}, J) \to ({\cal E}, L)$ and $C:({\cal D}, K)\to ({\cal E}, L)$ be comorphisms of sites related by a natural isomorphism $C\circ A\cong B$. If $B$ is $(J, L)$-continuous and $C$ is a fibration then $A$ satisfies the conditions of Proposition \ref{propcharcanonicalfunctors} with respect to every $J$-covering sieve $S$ which satisfies the property that there is a functorial factorization $g=g_{c}\circ g_{v}$ for any arrow $g$ in $S$ such that $g_{c}$ lies in $S$, $A(g_{c})$ is a cartesian arrow and $A(g_{v})$ is a vertical arrow (whence the factorization $A(g)=A(g_{c})\circ A(g_{v})$ yields \ac the' decomposition of $A(g)$ in $\cal D$ as a vertical arrow followed by a cartesian arrow); that is, for any such sieve $S$ on an object $c$, $A(S)$ is $K$-covering and for any commutative square of the form
	\[
	\xymatrix{
		d \ar[r] \ar[d]  & A(c') \ar[d]^{A(f)} \\
		A(c'') \ar[r]^{A(g)} & A(c),
	} 
	\]	
	where $f:c'\to c$ and $g:c''\to c$ are arbitrary arrows of $S$, there is a $K$-covering family $\{d_{i}\to d \mid i\in I\}$ such that for each $i\in I$, the composites $d_{i}\to A(c')$ and $d_{i}\to A(c'')$ belong to the same connected component of the category $(d_{i}\downarrow D^{A}_{S})$.	
\end{proposition}

\begin{proofs}
	The fact that $S$ admits a functorial statement as in the statement of the proposition implies that $S$ which is generated by a presieve $P$ of arrows $f$ such that $A(f)$ is cartesian (indeed, one can take $P$ to consist precisely of the arrows $g_{c}$ for $g\in S$). Let us first prove that $A(S)$ is $K$-covering.  From the existence of an isomorphism $C\circ A\cong B$ and the fact that $B$ is cover-preserving, it follows that the family $\{C(A(f))\mid f\in S\}$ is $L$-covering on $C(A(c))$. Now, since $C$ is a comorphism of sites, the sieve $S^{C}_{\langle \{C(A(f))\mid f\in S\}\rangle }$ is $K$-covering (we refer to section \ref{sec:smallestgrothendiecktopologycomorphism} for the notation). But by Lemma \ref{lemmafibrations}(iii), since the presieve $\{A(f) \mid f\in P\}$ entirely consists of cartesian arrows, $S^{C}_{\langle \{C(A(f))\mid f\in S\}\rangle }=\langle \{A(f) \mid f\in P\}\rangle $; so $\langle \{A(f) \mid f\in P\}\rangle $ is $K$-covering, as required.
	
	Let us now suppose that we have a commutative square of the form
		\[
	\xymatrix{
		d \ar[r]^{\chi} \ar[d]^{\xi}  & A(c') \ar[d]^{A(f)} \\
		A(c'') \ar[r]^{A(g)} & A(c),
	} 
	\]
	where $f,g\in S$.
	
	Consider the image of this square under the functor $C$: using the isomorphism $C\circ A\cong B$, we obtain a commutative square
			\[
	\xymatrix{
		C(d) \ar[r] \ar[d]  & B(c') \ar[d]^{B(f)} \\
		B(c'') \ar[r]^{B(g)} & B(c).
	} 
	\]
	The $(J, L)$-continuity of $B$ thus ensures the existence of a $L$-covering family $\{r_{i}:e_{i}\to C(d) \mid i\in I\}$ such that the composite arrows $e_{i}\to B(c')$ and $e_{i}\to B(c'')$ belong to the same connected component of the category $(e_{i} \downarrow D^{B}_{S})$. Now, since $C$ is a comorphism of sites, there is a $K$-covering family $\{g_{k}:d_{z}\to d \mid z\in Z\}$ whose image under $C$ is contained in the sieve generated by the family $\{e_{i}\to C(d) \mid i\in I\}$. So for each $z\in Z$ there are an element $i_{z}\in I$ and an arrow $h_{z}:d_{z}\to e_{i_{z}}$ such that $r_{i_{z}}\circ h_{z}=C(g_{z})$
	
	We thus have a commutative square
	\[
	\xymatrix{
		C(d_{z}) \ar[rr]^{C(\chi \circ g_{z})} \ar[d]^{C(\xi \circ g_{z})}  & & C(A(c')) \ar[d]^{C(A(f))} \\
		C(A(c'')) \ar[rr]^{C(A(g))} & & C(A(c)),
	} 
	\]
	where the objects $C(\chi\circ g_{z})$ and $C(\xi \circ g_{z})$ belong to the same connected component of the category $(C(d_{z})\downarrow D^{C\circ A}_{S})$. Let us show that the objects $\chi\circ g_{z}$ and $\xi\circ g_{z}$ belong to the same connected component of the category $(d_{z}\downarrow D^{A}_{S})$. This will clearly imply our thesis. But this follows from the fact that, more generally, for any arrows $\alpha:d'\to A(c')$ and $\beta:d'\to A(c'')$ in $\cal D$ such that $C(\alpha)$ and $C(\beta)$ belong to the same connected component of the category $(C(d')\downarrow D^{C\circ A}_{S})$, $\alpha$ and $\beta$ belong to the same connected component of the category $(d'\downarrow D^{C\circ A}_{S})$, which is an immediate consequence of Lemma \ref{lemmacartesianrelmono}(iii). 	     
\end{proofs}

\begin{corollary}
	Let $F:{\cal C}\to {\cal E}$ and $G:{\cal D}\to {\cal E}$ be fibrations and $H$ a morphism of fibrations $F\to G$, that is, a functor $H:{\cal C}\to {\cal D}$ which sends cartesian arrows to cartesian arrows and such that there is an isomorphism $G\circ H\cong F$. Then $H$ is $(M^{F}_{L}, M^{G}_{L})$-continuous.   	
\end{corollary}

\begin{proofs}
	Lemma \ref{lemmacartesianrelmono}(iii) ensures that $H$ satisfies the hypothesis of Proposition \ref{propliftcontinuity} with respect to   every sieve which contains all the cartesian images of its arrows (notice that a sieve has this property if and only if it is generated by a presieve entirely consisting of cartesian arrows). Now, not every $M^{F}_{L}$-covering sieve has this property, but $M^{F}_{L}$ is generated by sieves with this property. So, proving that $H$ is $(M^{F}_{L}, M^{G}_{L})$-continuous amounts precisely (by Proposition \ref{propcharJcanonical}(iv)) to showing that $H$ satisfies the continuity condition with respect to all the pullbacks of the sieves with the above property. But by Lemma \ref{lemmafibrations}(vii) this class of sieves is closed under pullback, whence our thesis follows.   
\end{proofs}

\begin{corollary}\label{corliftalongdiscretefibration}
	Let $p:{\cal D}\to {\cal D}'$ be a discrete fibration, $K$ a Grothendieck topology on ${\cal D}'$ and $F:({\cal C}, J)\to ({\cal D}, M^{p}_{K})$ a comorphism of sites. Then $F\circ p$ is continuous if and only if $F$ is.
\end{corollary}

\begin{proofs}
	Since $p$ is a discrete fibration, every arrow of $\cal D$ is cartesian and hence, by Proposition \ref{propliftcontinuity}, if $F\circ p$ is continuous then $F$ trivially satisfies the hypotheses of Proposition \ref{propcharcanonicalfunctors} with respect to every ($J$-covering) sieve on $\cal C$, that is, $F$ is continuous. Conversely, if $F$ is continuous then $F\circ p$ also is, since $p$ is continuous being a fibration.  
\end{proofs}

The following other corollary of Proposition \ref{propliftcontinuity} shows that the property of continuity of a comorphism of sites under relativisation. It should be regarded as a natural companion to Proposition \ref{propcancomorphismsstableunderpullback}.

\begin{corollary}
	Let $F:({\cal C}, J)\to ({\cal D}, K)$ be a continuous comorphism of sites. Then, for any object $c$ of $\cal C$, the functor
	\[
	F\slash c: ({\cal C}\slash c, J_{c}) \to ({\cal D}\slash F(c), K_{F(c)})
	\]
	is a continuous comorphism of sites $({\cal C}\slash c, J_{c}) \to ({\cal D}\slash F(c), K_{F(c)})$.
\end{corollary}

\begin{proofs}
	Since the canonical projection ${\cal D}\slash F(c)\to {\cal D}$ is a discrete fibration, all the arrows of ${\cal D}\slash F(c)$ are cartesian. Therefore every sieve $S$ in $\cal C$ trivially satisfies the property in Proposition \ref{propliftcontinuity}, namely that there is a functorial factorization $g=g_{c}\circ g_{v}$ for any arrow $g$ in $S$ such that $g_{c}$ lies in $S$, $A(g_{c})$ is a cartesian arrow and $A(g_{v})$ is a vertical arrow for \emph{any} functor $A:{\cal C}\to {\cal D}\slash F(c)$. Our thesis thus follows from the proposition by taking $A$ to be the functor $F\slash c$, and $B$ to be the composite of $F$ and the projection functor ${\cal C}\slash c \to {\cal C}$, and $C$ to be the projection functor ${\cal D}\to {\cal D}\slash F(c)$; note that $B$ is continuous as it is the composite of two continuous comorphisms of sites.  
\end{proofs}

The following technical result about morphisms of fibrations will be instrumental in the next section, for proving that such morphisms induce locally connected morphisms of toposes:

\begin{proposition}\label{proppropertymorphismsoffibrations}
Let $F:{\cal C}\to {\cal E}$ and $G:{\cal D}\to {\cal E}$ be fibrations and $H$ a morphism of fibrations $F\to G$. Then $H$ satisfies the following property: for any object $c$ of $\cal C$ and arrow $g:d\to H(c)$ in $\cal D$, there is a cartesian (with respect to $F$) arrow $f:c'\to c$ and a vertical arrow $v:d'\to H(c')$ such that $g=H(f)\circ v$ (note that $H(f)$ is cartesian with respect to $G$ since $f$ is cartesian and $H$ is a morphism of fibrations). In other words, the cartesian image of an arrow towards an object in the image of $H$ also lies in the image of $H$. 	
\end{proposition}

\begin{proofs}
Let $\alpha$ be an isomorphism $F\to G\circ H$. Given an arrow $g:d\to H(c)$ in $\cal D$, consider its image 
\[
G(g):G(d)\to G(H(c))
\] 	
under $G$. Then its composite with $\alpha(c)^{-1}:G(H(c))\to F(c)$ yields an arrow $G(d)\to F(c)$ in $\cal E$. The fact that $F$ is a fibration implies that there is an isomorphism $\gamma:G(d)\to F(c')$ and a cartesian arrow $f:c'\to c$ such that $\alpha(c)^{-1}\circ G(g)=F(f)\circ \gamma$. 

Now consider the arrow 
\[
\gamma^{-1}:F(c')\to G(d).
\]
Since $G$ is a fibration, there are a cartesian arrow $k:d'\to d$ and an isomorphism $\beta:G(d')\to F(c')$ such that $\gamma^{-1}\circ \beta=G(k)$. Since $k$ is cartesian and $G(k)$ is an isomorphism, it follows that $k$ is an isomorphism. On the other hand, since $f$ is cartesian (with respect to $F$) and $H$ is a morphism of fibrations, the arrow $H(f):H(c')\to H(c)$ is cartesian (with respect to $G$). Now, we have 
\[
G(H(f))\circ (\alpha(c')\circ \beta)=G(g\circ k).
\] 
Indeed, $G(H(f))\circ \alpha(c')\circ \beta=\alpha(c)\circ F(f)\circ \beta$ (by the naturality of $\alpha$), and $\alpha(c)\circ F(f)\circ \beta=G(g)\circ \gamma^{-1}\circ \beta=G(g)\circ G(k)=G(g\circ k)$.

So the cartesianness of $H(f)$ implies that there is a unique arrow $u:d'\to H(c')$ such that $H(f)\circ u=g\circ k$, equivalently $g=H(f)\circ (u\circ k^{-1})$. To conclude our proof, it remains to observe that the arrow $u\circ k^{-1}$ is vertical; but this follows from the fact that $G(u\circ k^{-1})=G(u)\circ G(k^{-1})=\alpha(c')\circ \beta  \circ \beta^{-1} \circ \gamma=\alpha(c')\circ \gamma$, which is an isomorphism.  
\end{proofs}

\begin{remark}
	Given a fibration $p:{\cal C}\to {\cal D}$, we denote by ${\cal C}^{\textup{cart}}_{p}$ the subcategory of $\cal C$ whose objects are the objects of $\cal C$ and whose arrows are the arrows of $\cal C$ which are cartesian with respect to the fibration $p$ (note that this is a well-defined category since, for any fibration, identical arrows are cartesian and the composite of two cartesian arrows is cartesian). Then a morphism of fibrations $H$ as in the statement of Proposition \ref{proppropertymorphismsoffibrations} yields by restriction a discrete fibration ${\cal C}^{\textup{cart}}_{F} \to {\cal D}^{\textup{cart}}_{G}$, by the property of $H$ stated in Proposition \ref{proppropertymorphismsoffibrations}.
	
	Notice that, conversely, any fibration $p:{\cal D}\to {\cal C}$ can be seen as a morphism of fibrations from $p$ to the identical fibration on $\cal C$ (since all the arrows of $\cal C$ are cartesian with respect to the latter).
\end{remark}

\subsection{Locally connected morphisms}\label{sec:locallyconnectedmorphisms}

The class of locally connected morphisms consistutes an important class of essential geometric morphisms. In light of the results of section \ref{sec:essentialmorphismsandcomorphisms}, it is interesting to investigate them from the point of view of comorphisms of sites. Recall (see, for instance, section C3.3 of \cite{El}) that a geometric morphism $f:{\cal F}\to {\cal E}$ is said to be \emph{locally connected} if $f^{\ast}$ has an $\cal E$-indexed left adjoint. By Proposition C3.3.1 \cite{El}, a geometric morphism $f:{\cal F}\to {\cal E}$ if locally connected if and only if, for any arrow $h:A\to B$ in $\cal E$, the square
\[
\xymatrix{
	{\cal E}\slash A \ar[r]^{\Pi_{h}} \ar[d]^{(f\slash A)^{\ast}}  & {\cal E}\slash B \ar[d]^{(f\slash B)^{\ast}} \\
	{\cal F}\slash f^{\ast}(A) \ar[r]^{\Pi_{f^{\ast}(h)}} & {\cal F}\slash f^{\ast}(B),
} 
\]
where $\Pi_{h}$ and $\Pi_{f^{\ast}(h)}$ are respectively the direct images of the geometric morphisms ${\cal E}\slash A \to {\cal E}\slash B$ and  ${\cal F}\slash f^{\ast}(A) \to {\cal F}\slash f^{\ast}(B)$ induced respectively by $h$ and by $f^{\ast}(h)$, commutes.

Let us preliminarily notice that if $f:{\cal F}\to {\cal E}$ is essential then $f\slash E: {\cal F}\slash f^{\ast}(E)\to {\cal E}\slash E$ is also essential for every object $E$ of $\cal E$. Indeed, if $f_{!}$ is left adjoint to $f^{\ast}$ then the functor ${\cal F}\slash f^{\ast}(E)\to {\cal E}\slash E$ sending an object $[y:Y\to f^{\ast}(E)]$ of ${\cal F}\slash f^{\ast}(E)$ to the object $[\epsilon_{E}\circ f_{!}(y):f_{!}(Y)\to E]$, where $\epsilon_{E}:f_{!}(f^{\ast}(E))\to E$ is the component at $E$ of the counit $\epsilon:f_{!}\circ f^{\ast} \to 1$ of the adjunction between $f_{!}$ and $f^{\ast}$, and acting on arrows in the obvious way, is left adjoint to $(f\slash E)^{\ast}$.  

In light of this observation, we can rephrase the above condition in terms of the left adjoints, as follows: for any arrow $h:A\to B$ in $\cal E$, the square
\[
\xymatrix{
	{\cal F}\slash f^{\ast}(B) \ar[rr]^{(f\slash B)_{!}} \ar[d]^{(f^{\ast}(h))^{\ast}} & &  {\cal E}\slash B    \ar[d]^{h^{\ast}} \\
	{\cal F}\slash f^{\ast}(A) \ar[rr]^{(f\slash A)_{!}} & & {\cal E}\slash A
} 
\]
commutes. We shall refer to this condition as $\Gamma^{f}_{h}$.

\begin{remark}\label{remlocallyconnectedcondition}
Recall that every \emph{local homeomorphism} (that is, geometric morphism of the form $U^{\cal E}_{E}:{\cal E}\slash E\to {\cal E}$ for some object $E$ of a topos $\cal E$) is locally connected. Notice that, for any arrow $u:E'\to E$ in a topos $\cal E$, the geometric morphism $U^{\cal E}_{u}:{\cal E}\slash E'\to {\cal E}\slash E$ induced by it (whose inverse image is given by the pullback functor along $u$) is a local homeomorphism (since it can be identified with the functor $U^{{\cal E}\slash B}_{[u:E'\to E]}$). Let us denote by $\Sigma^{\cal E}_{u}$ the essential image of $U^{\cal E}_{u}$, that is, the functor ${\cal E}\slash E'\to {\cal E}\slash E$ acting as composition with $u$.  Given an arrow $h:[a:A\to E]\to [b:B\to E]$ in ${\cal E}\slash E$, that is an arrow $h:A\to B$ in $\cal E$ such that $b\circ h=a$, condition $\Gamma^{U^{\cal E}_{u}}_{h}$ rewrites as follows: given pullback squares
$$
\xymatrix {
	A' \ar[d]^{a'} \ar[r]^{u_{A}} &  A \ar[d]^{a} & \textup{ and } & B' \ar[d]^{b'} \ar[r]^{u_{B}} &  B \ar[d]^{b} \\ 
	E' \ar[r]^{u}  & E & &
	E' \ar[r]^{u} & E,}
$$
the diagram	
\[
\xymatrix{
	{\cal E}\slash B' \ar[rr]^{\Sigma^{\cal E}_{u_{B}}} \ar[d]^{(u^{\ast}(h))^{\ast}} & &  {\cal E}\slash B    \ar[d]^{h^{\ast}} \\
	{\cal E}\slash A' \ar[rr]^{\Sigma^{\cal E}_{u_{A}}} & & {\cal E}\slash A
} 
\]
commutes. In the case $b=1_{E}$, the arrow $h$ is simply an arrow $a:A\to E$ and this square specializes to the following one:
\[
\xymatrix{
	{\cal E}\slash E' \ar[rr]^{\Sigma^{\cal E}_{u}} \ar[d]^{(u^{\ast}(a))^{\ast}} & &  {\cal E}\slash E    \ar[d]^{a^{\ast}} \\
	{\cal E}\slash A' \ar[rr]^{\Sigma^{\cal E}_{u_{A}}} & & {\cal E}\slash A.
} 
\]
\end{remark}

Note that the commutativity of this square ultimately results from the pullback lemma.

\begin{proposition}\label{propsimplificationlocallyconnected}
	Let  $f:{\cal F}\to {\cal E}$ be an essential geometric morphism between Grothendieck toposes and $h:A\to B$ an arrow of $\cal E$.
	\begin{enumerate}[(i)]
		\item For any epimorphic family $\{f_{i}:A_{i}\to A \mid i\in I\}$ in $\cal E$, $\Gamma^{f}_{h}$ holds if and only if $\Gamma^{f}_{h\circ f_{i}}$ holds for every $i\in I$.
		
		\item For any epimorphic family $\{g_{j}:B_{j}\to B \mid j\in J\}$ in $\cal E$, $\Gamma^{f}_{h}$ holds if and only if $\Gamma^{f}_{g_{j}^{\ast}(h)}$ holds for every $j\in J$.
	\end{enumerate} 
\end{proposition}	

\begin{proofs}
	(i) From the fact that the canonical indexing of a Grothendieck topos over itself is a stack with respect to the canonical topology it follows that the square
	\[
	\xymatrix{
	{\cal F}\slash f^{\ast}(B) \ar[rr]^{(f\slash B)_{!}} \ar[d]^{(f^{\ast}(h))^{\ast}} & &  {\cal E}\slash B    \ar[d]^{h^{\ast}} \\
	{\cal F}\slash f^{\ast}(A) \ar[rr]^{(f\slash A)_{!}} & & {\cal E}\slash A
	} 
	\]
	commutes if and only if each of the outer rectangles
	\[
	\xymatrix{
	{\cal F}\slash f^{\ast}(B) \ar[rr]^{(f\slash B)_{!}} \ar[d]^{(f^{\ast}(h))^{\ast}} & &  {\cal E}\slash B    \ar[d]^{h^{\ast}} \\
	{\cal F}\slash f^{\ast}(A) \ar[rr]^{(f\slash A)_{!}} \ar[d]^{(f^{\ast}(f_{i}))^{\ast}} & & {\cal E}\slash A \ar[d]^{f_{i}^{\ast}}\\
	{\cal F}\slash f^{\ast}(A_{i}) \ar[rr]^{(f\slash A_{i})_{!}} & & {\cal E}\slash A_{i}
	} 
	\]
	does. But this rectangle is precisely the square corresponding to the arrow $h\circ f_{i}$. This shows that the commutativity of the square with respect to an arrow $h$ is equivalent to that of all the squares with respect to the arrows of the form $h\circ f_{i}$ for $i\in I$.
	
	(ii) Let us preliminarily notice that, for any epimorphic family $\{c_{k}:C_{k}\to C \mid k\in K\}$ in a Grothendieck topos $\cal E$, the family of functors $\{\Sigma^{\cal E}_{c_{k}}:{\cal E}\slash C_{k} \to {\cal E}\slash C\}$ satisfies the property that for any two colimit-preserving functors $H, K:{\cal E}\slash C \to {\cal F}$, there are isomorphisms $H\circ \Sigma^{\cal E}_{c_{k}}\cong K\circ \Sigma^{\cal E}_{c_{k}}$ satisfying the obvious coherence relations if and only if $H\cong K$. Indeed, for any object $[u:U\to C]$ of ${\cal E}\slash C$, consider the pullback squares
		\[
	\xymatrix{
		U_{k} \ar[rr]^{\xi_{k}} \ar[d]^{u_{k}} & &  U    \ar[d]^{u} \\
		C_{k} \ar[rr]^{c_{k}} & & C.
	} 
	\]    
	
	  Since the family $\{\xi_{k}:U_{k} \to U \mid k\in K\}$ is epimorphic, denoting by $S_{u}$ the sieve generated by the family $\{\xi_{k}:U_{k} \to U \mid k\in K\}$, we can represent $[u:U\to C]$ as the colimit in ${\cal E}\slash C$ of the diagram $D_{u}:{\int S}_{u}\to {\cal E}\slash C$ which sends any object $(Z, z:Z\to U)$ of ${\int S_{u}}$ to the object $[u\circ z:Z\to C]$ of ${\cal E}\slash C$. So, for any colimit-preserving functors $H$ and $K$ as above, $H([u])\cong K([u])$ if $H\circ D_{u}\cong K\circ D_{u}$. Now, we have such isomorphisms $H\circ D_{u}\cong K\circ D_{u}$ by the coherence relations existing by our hypotheses between the isomorphisms between the functors $H\circ \Sigma^{\cal E}_{c_{k}}$ and $K\circ \Sigma^{\cal E}_{c_{k}}$. Moreover, these isomorphisms $H([u])\cong K([u])$ (note that any object in the image of $D_{u}$ belongs to the image of some $\Sigma^{\cal E}_{c_{k}}$) are clearly natural in $[u:U\to C]\in {\cal E}\slash C$, which establishes our claim. 
	
	Now, given, for each $j\in J$, the pullback square
	\[
	\xymatrix{
		A_{j} \ar[rr]^{h_{j}=g_{j}^{\ast}(h)} \ar[d]^{h^{\ast}(g_{j})} & &  B_{j}    \ar[d]^{g_{j}} \\
		A \ar[rr]^{h} & & B,
	} 
	\]
	consider the following cube:  
	
	\begin{center}
		\begin{tikzcd}[row sep=large, column sep=small]
		& {\cal F}\slash f^{\ast}(B) \arrow[dd, "(f^{\ast}(h))^{\ast}"{yshift=4ex}, rightarrow] \arrow[rr, "(f\slash B)_{!}"] &                                                                             & {\cal E}\slash B \arrow[dd, "h^{\ast}"{yshift=2ex}, rightarrow] \\
		{\cal F}\slash f^{\ast}(B_{j}) \arrow[dd, "(f^{\ast}(h_{j}))^{\ast}"{yshift=2.5ex, xshift=0ex}, rightarrow] \arrow[rr, "(f\slash B_{j})_{!}"{xshift=-4ex}, crossing over] \arrow[ru, "\Sigma^{\cal F}_{f^{\ast}(g_{j})}"] &                                                                                             & {\cal E}\slash B_{j} \arrow[dd, "h_{j}^{\ast}"{yshift=4ex}, rightarrow, crossing over]  \arrow[ru, "\Sigma^{\cal E}_{g_{j}}"] &                                                            \\
		& {\cal F}\slash f^{\ast}(A) \arrow[rr, "(f\slash A)_{!}"{below, xshift=-3ex}]                                                &                                                                             & {\cal E}\slash A                                 \\
		{\cal F}\slash f^{\ast}(A_{j}) \arrow[ru, "\Sigma^{\cal F}_{f^{\ast}(h^{\ast}(g_{j}))}"'] \arrow[rr, "(f\slash A_{j})_{!}"']             &                                                                                             & {\cal E}\slash A_{j} \arrow[ru, "\Sigma^{\cal E}_{h^{\ast}(g_{j})}"']     &                                                           
		\end{tikzcd}
	\end{center} 

We want to show that the back square is commutative (that is, $\Gamma^{f}_{h}$ holds) if and only if the front square is for each $j\in J$ (that is, if and only if $\Gamma^{f}_{g_{j}^{\ast}(h)}$ holds for each $j\in J$). This follows from the above preliminary remark and the fact that all the other squares in the cube commute. More specifically, the lower and upper squares commute since the squares obtained by taking the right adjoints of all functors commute, while the left-hand and right-hand lateral faces commute since they are commutative squares associated with local homeomorphisms of toposes (cf. Remark \ref{remlocallyconnectedcondition}).
\end{proofs}

\begin{corollary}\label{corsimplifiedlocalconnectedness}
	Let $({\cal C}, J)$ be a small-generated site and $f:{\cal F}\to \Sh({\cal C}, J)$ an essential geometric morphism. Then $f$ is locally connected if and only if $\Gamma^{f}_{l(g)}$ holds for every arrow $g$ in $\cal C$ (where $l$ is the canonical functor ${\cal C}\to \Sh({\cal C}, J)$). 
\end{corollary}

\begin{proofs}
	Let $h:Q\to Q'$ be an arrow in $\Sh({\cal C}, J)$. The object $Q'$ of $\Sh({\cal C}, J)$ can be represented as $\textup{colim}_{\Sh({\cal C}, J)}(l\circ \pi_{Q'})$; we shall denote by $\xi_{(c, x)}:l(c)\to Q'$ the colimit arrow corresponding to the object $(c, x)$ of $\int Q'$. Then, by Proposition \ref{propsimplificationlocallyconnected}(ii), $\Gamma^{f}_{h}$ holds if and only if $\Gamma^{f}_{\xi_{(c, x)}^{\ast}(h)}$ holds for every object $(c, x)$ of $\int Q'$. Now, every arrow $k:Q''\to l(c)$ in $\Sh({\cal C}, J)$ can be locally represented in terms of arrows of the form $l(g)$ for $g$ an arrow of $\cal C$ with codomain $c$, since $Q''$ can be covered by objects of the form $l(c')$ for $c'\in {\cal C}$ and then one can apply Proposition \ref{propexplicit}(i). So Proposition \ref{propsimplificationlocallyconnected}(i) ensures that $f$ is locally connected (that is, $\Gamma^{f}_{h}$ holds for every $h$) if and only if $\Gamma^{f}_{l(g)}$ holds for every arrow $g$ in $\cal C$, as required. 
\end{proofs}

We shall now proceed to identify necessary and sufficient conditions on a $(J, K)$-continuous comorphism of sites $({\cal C}, J)\to ({\cal D}, K)$ to induce a locally connected geometric morphism $\Sh({\cal C}, J)\to \Sh({\cal D}, K)$.

Suppose that $f$ is the geometric morphism $C_{F}:\Sh({\cal C}, J)\to \Sh({\cal D}, K)$ induced by a comorphism of sites $F:({\cal C}, J)\to ({\cal D}, K)$. Given an arrow $h:Q\to Q'$ in $[{\cal D}^{\textup{op}}, \Set]$, we have to investigate under which conditions the square 
\[
\xymatrix{
	\Sh({\cal C}, J)\slash (C_{F})^{\ast}(a_{K}(Q')) \ar[rr]^{\quad(C_{F}\slash a_{K}(Q'))_{!}} \ar[d]^{(C_{F}^{\ast}(a_{K}(h)))^{\ast}} & &  \Sh({\cal D}, K)\slash a_{K}(Q')    \ar[d]^{a_{K}(h)^{\ast}} \\
	\Sh({\cal C}, J)\slash C_{F}^{\ast}(a_{K}(Q)) \ar[rr]^{\quad(C_{F}\slash a_{K}(Q))_{!}} & & \Sh({\cal D}, K) \slash a_{K}(Q)
} 
\]
relative to the arrow $a_{K}(h)$ in $\Sh({\cal D}, K)$ commutes.

Notice that if $F$ is $(J, K)$-continuous then by Theorem 2.2 \cite{dependent} and Proposition \ref{proplocalizationcanonicalcomorphism}, this square is isomorphic to the following one, where $p^{F}_{Q}$ (resp. $p^{F}_{Q'}$) is the $(J_{Q\circ F^{\textup{op}}}, J_{Q})$-continuous (resp.  $(J_{Q'\circ F^{\textup{op}}}, J_{Q'})$-continuous) comorphism of sites $(\int (Q\circ F^{\textup{op}}), J_{Q\circ F^{\textup{op}}}) \to (\int Q, J_{Q})$ (resp. $(\int (Q'\circ F^{\textup{op}}), J_{Q'\circ F^{\textup{op}}}) \to (\int Q', J_{Q'})$) of Proposition \ref{propcancomorphismsstableunderpullback}:
\[
\xymatrix{
	\Sh(\int (Q'\circ F^{\textup{op}}), J_{Q'\circ F^{\textup{op}}}) \ar[rr]^{\quad\quad(C_{p^{F}_{Q'}})_{!}} \ar[d]^{(C_{\int hF^{\textup{op}}})^{\ast}} & &  \Sh(\int Q', J_{Q'}) \ar[d]^{(C_{\int h})^{\ast}} \\
	\Sh(\int (Q\circ F^{\textup{op}}), J_{Q\circ F^{\textup{op}}}) \ar[rr]^{\quad\quad(C_{p^{F}_{Q}})_{!}} & & \Sh(\int Q, J_{Q})
} 
\] 

The following result provides an explicit characterization of the continuous comorphisms of sites $F$ whose associated geometric morphism $C_{F}$ is locally connected. In order to state it, we need to introduce the following constructions.

For any pair $(c, x)$ consisting of an object $c$ of $\cal C$ and an element $x$ of $Q'(F(c))$, and any arrow $h:Q\to Q'$ in $[{\cal D}^{\textup{op}}, \Set]$, we define two categories ${\cal A}^{y_{\cal D}(h)}_{(c, x)}$ and ${\cal B}^{y_{\cal D}(h)}_{(c, x)}$ as follows:

\begin{enumerate}[(1)]
	\item The objects of ${\cal A}^{y_{\cal D}(h)}_{(c, x)}$ are the triplets $(c', y, f)$ where $c'$ is an object of $\cal C$, $y$ is an element of $Q(F(c'))$ and $f:c\to c$ is an arrow of $\cal C$ such that $h(F(c'))(y)=Q'(F(f))(x)$.
	
	The arrows $(c_{1}, y_{1}, f_{1})\to (c_{2}, y_{2}, f_{2})$ in ${\cal A}^{y_{\cal D}(h)}_{(c, x)}$ are the arrows $t:c_{1} \to c_{2}$ in $\cal C$ such that $f_{2}\circ t=f_{1}$ and $Q(F(t))(y_{2})=y_{1}$. 
	
	\item The objects of ${\cal B}^{y_{\cal D}(h)}_{(c, x)}$ are the triplets $(d, z, g)$ where $d$ is an object of $\cal D$, $z$ is an element of $Q(z)$ and $g:d\to F(c)$ is an arrow of $\cal D$ such that $h(d)(z)=Q'(g)(x)$.
	
	The arrows $(d_{1}, z_{1}, g_{1})\to (d_{2}, z_{2}, g_{2})$ in ${\cal B}^{y_{\cal D}(h)}_{(c, x)}$ are the arrows $s:d_{1} \to d_{2}$ in $\cal D$ such that $g_{2}\circ s=g_{1}$ and $Q(s)(z_{2})=z_{1}$.
\end{enumerate}

The categories ${\cal A}^{y_{\cal D}(h)}_{(c, x)}$ and ${\cal B}^{y_{\cal D}(h)}_{(c, x)}$ are related to each other by a functor $$\xi^{y_{\cal D}(h)}_{(c, x)}:{\cal A}^{y_{\cal D}(h)}_{(c, x)} \to {\cal B}^{y_{\cal D}(h)}_{(c, x)}$$ sending any object $(c', y, f)$ of ${\cal A}^{y_{\cal D}(h)}_{(c, x)}$ to the object $(F(c'), y, F(f))$ of ${\cal B}^{y_{\cal D}(h)}_{(c, x)}$ and any arrow $s:(c_{1}, y_{1}, f_{1})\to (c_{2}, y_{2}, f_{2})$ in ${\cal A}^{y_{\cal D}(h)}_{(c, x)}$ to the arrow $F(s):(F(c_{1}), y_{1}, F(f_{1}))\to (F(c_{2}), y_{2}, F(f_{2}))$ of ${\cal B}^{y_{\cal D}(h)}_{(c, x)}$. 

As we shall see in the proof of Theorem \ref{thmlocallyconnected}, the categories ${\cal A}^{y_{\cal D}(h)}_{(c, x)}$ and ${\cal B}^{y_{\cal D}(h)}_{(c, x)}$ are actually fibered over $\cal D$.
We shall denote by $a^{y_{\cal D}(h)}_{(c, x)}:{\cal A}^{y_{\cal D}(h)}_{(c, x)} \to {\cal D}$ the functor sending any object  $(c', y, f)$ of ${\cal A}^{y_{\cal D}(h)}_{(c, x)}$ to the object $F(c')$ of $\cal D$ and any arrow $s:(c_{1}, y_{1}, f_{1})\to (c_{2}, y_{2}, f_{2})$ in ${\cal A}^{y_{\cal D}(h)}_{(c, x)}$ to the arrow $F(s):F(c_1)\to F(c_2)$ of $\cal D$, and by $b^{y_{\cal D}(h)}_{(c, x)}:{\cal B}^{y_{\cal D}(h)}_{(c, x)} \to {\cal D}$ the canonical projection functor. In fact, $\xi^{y_{\cal D}(h)}_{(c, x)}$ is a morphism of fibrations $a^{y_{\cal D}(h)}_{(c, x)} \to b^{y_{\cal D}(h)}_{(c, x)}$.  

Let us now consider the particular case in which the arrow $h:Q\to Q'$ is of the form $y_{\cal D}(h):y_{\cal D}(d_{0})\to y_{\cal D}(d_{1})$, where $h:d_{0}\to d_{1}$ is an arrow in $\cal D$.  

\begin{enumerate}[(1)]
\item The category ${\cal A}^{y_{\cal D}(h)}_{(c, x)}$, where $c$ is an object of $\cal D$ and $x$ is an arrow $F(c)\to d_{1}$ in $\cal D$, has as objects the triplets $(c', y, f)$ where $c'$ is an object of $\cal C$, $y$ is an arrow $F(c')\to d_{0}$ in $\cal D$ and $f:c'\to c$ is an arrow of $\cal C$ such that $x\circ F(f)=h\circ y$, and as arrows $(c_{1}, y_{1}, f_{1})\to (c_{2}, y_{2}, f_{2})$ the arrows $t:c_{1} \to c_{2}$ in $\cal C$ such that $f_{2}\circ t=f_{1}$ and $y_{2} \circ F(t)=y_{1}$.	
	
\item The category ${\cal B}^{y_{\cal D}(h)}_{(c, x)}$ has as objects the triplets $(d, z, g)$ where $d$ is an object of $\cal D$, $z$ is an arrow $d\to d_{0}$ in $\cal D$ and $g:d\to F(c)$ is an arrow of $\cal D$ such that $x\circ g=h\circ z$, and as arrows $(d_{1}, z_{1}, g_{1})\to (d_{2}, z_{2}, g_{2})$ in ${\cal B}^{y_{\cal D}(h)}_{(c, x)}$ the arrows $s:d_{1} \to d_{2}$ in $\cal D$ such that $g_{2}\circ s=g_{1}$ and $z_{2} \circ s=z_{1}$.

\item The functor 
\[
\xi^{y_{\cal D}(h)}_{(c, x)}:a^{y_{\cal D}(h)}_{(c, x)} \to b^{y_{\cal D}(h)}_{(c, x)}  
\]
sends any object $(c', y, f)$ of ${\cal A}^{y_{\cal D}(h)}_{(c, x)}$ to the object $(F(c'), y, F(f))$ of ${\cal B}^{y_{\cal D}(h)}_{(c, x)}$ and any arrow $s:(c_{1}, y_{1}, f_{1})\to (c_{2}, y_{2}, f_{2})$ in ${\cal A}^{y_{\cal D}(h)}_{(c, x)}$ to the arrow 
\[
F(s):(F(c_{1}), y_{1}, F(f_{1}))\to (F(c_{2}), y_{2}, F(f_{2}))
\] 
of ${\cal B}^{y_{\cal D}(h)}_{(c, x)}$. 
\end{enumerate}

\begin{theorem}\label{thmlocallyconnected}
	Let $F:({\cal C}, J)\to ({\cal D}, K)$ be a continuous comorphism of small-generated sites. Then the following conditions are equivalent:
	\begin{enumerate}[(i)]
		\item The geometric morphism $C_{F}:\Sh({\cal C}, J)\to \Sh({\cal D}, K)$ induced by $F$ is locally connected.
		
		\item For any arrow $h:Q\to Q'$ in $[{\cal D}^{\textup{op}}, \Set]$, the morphism of fibrations 
		\[
		\xi^{y_{\cal D}(h)}_{(c, x)}: a^{y_{\cal D}(h)}_{(c, x)}\to b^{y_{\cal D}(h)}_{(c, x)}
		\]
		to $\cal D$ satisfies (together with the identical natural transformation  $a^{y_{\cal D}(h)}_{(c, x)} \to b^{y_{\cal D}(h)}_{(c, x)}\circ \xi^{y_{\cal D}(h)}_{(c, x)}$) the \ac cofinality' conditions of Proposition \ref{procofinality}, that is (cf. Remark \ref{rempropcofinality}(b)), using the notation of the proposition, the following conditions: 
		\begin{enumerate}[(i)]
			\item For any object $(d, z, g)$ of the category ${\cal B}^{y_{\cal D}(h)}_{(c, x)}$ there are a $K$-covering family $\{g_{i}:d_{i}\to d \mid i\in I\}$ and for each $i\in I$ an object $(c_{i}, y_{i}, f_{i})$ of the category ${\cal A}^{y_{\cal D}(h)}_{(c, x)}$ and an arrow $s_{i}:d_{i}\to F(c_{i})= a^{y_{\cal D}(h)}_{(c, x)}((c_{i}, y_{i}, f_{i}))=b^{y_{\cal D}(h)}_{(c, x)}(\xi^{y_{\cal D}(h)}_{(c, x)}((c_{i}, y_{i}, f_{i})))=b^{y_{\cal D}(h)}_{(c, x)}((F(c_{i}), y_{i}, F(f_{i})))$ such that $(g_{i}:d_{i}\to d=b^{y_{\cal D}(h)}_{(c, x)}((d, z, g)),s_{i}:d_{i}\to F(c_{i})=b^{y_{\cal D}(h)}_{(c, x)}((F(c_{i}), y_{i}, F(f_{i}))))$ belongs to $R'_{d_{i}}$. 
			
			\item For any object $d$ of $\cal D$ and any arrows $\alpha:d\to a^{y_{\cal D}(h)}_{(c, x)}((a, y, f))=F(a)$ and $\beta:d\to a^{y_{\cal D}(h)}_{(c, x)}((b, y', f'))=F(b)$ in $\cal D$ such that $(\alpha:d\to b^{y_{\cal D}(h)}_{(c, x)}((F(a), y, F(f)))=F(a), (\beta:d\to b^{y_{\cal D}(h)}_{(c, x)}((F(b), y', F(f')))=F(b)))$ belongs to $R'_{d}$, there is a $K$-covering family $\{g_i:d_i\to d\mid i\in I\}$ such that $((\alpha\circ g_{i}:d_{i}\to a^{y_{\cal D}(h)}_{(c, x)}((a, y, f))=F(a), (\beta\circ g_i:d_{i}\to b^{y_{\cal D}(h)}_{(c, x)}((b, y', f'))=F(b)))$ belongs to $R_{d}$.
		\end{enumerate}
		
		\item For any arrow $h:d_{0}\to d_{1}$ in $\cal D$, object $c$ of $\cal C$ and arrow $x:F(c)\to d_{1}$, the following conditions hold:
		\begin{enumerate}[(a)]
			\item For any object $(d, z, g)$ of the category ${\cal B}^{y_{\cal D}(h)}_{(c, x)}$ there is a $K$-covering family $\{g_{i}:d_{i}\to d \mid i\in I\}$ and for each $i\in I$ an object $(c_{i}, y_{i}, f_{i})$ of the category ${\cal A}^{y_{\cal D}(h)}_{(c, x)}$ and an arrow $$s_{i}:d_{i}\to F(c_{i})= a^{y_{\cal D}(h)}_{(c, x)}((c_{i}, y_{i}, f_{i}))$$ such that $$g_{i}:d_{i}\to d=b^{y_{\cal D}(h)}_{(c, x)}((d, z, g))$$ and $$s_{i}:d_{i}\to F(c_{i})=b^{y_{\cal D}(h)}_{(c, x)}((F(c_{i}), y_{i}, F(f_{i}))))$$ belong to the same connected component of the category $(d_{i}\downarrow b^{y_{\cal D}(h)}_{(c, x)})$. 
			
			\item For any object $d$ of $\cal D$ and any arrows $\alpha:d\to a^{y_{\cal D}(h)}_{(c, x)}((a, y, f))=F(a)$ and $\beta:d\to a^{y_{\cal D}(h)}_{(c, x)}((b, y', f'))=F(b)$ in $\cal D$ such that $$\alpha:d\to b^{y_{\cal D}(h)}_{(c, x)}((F(a), y, F(f)))=F(a)$$ and $$\beta:d\to b^{y_{\cal D}(h)}_{(c, x)}((F(b), y', F(f')))=F(b)$$ belong to the same connected component of the category $(d\downarrow b^{y_{\cal D}(h)}_{(c, x)})$, there is a $K$-covering family $\{g_i:d_i\to d\mid i\in I\}$ such that $$\alpha\circ g_{i}:d_{i}\to a^{y_{\cal D}(h)}_{(c, x)}((a, y, f))=F(a)$$ and $$\beta\circ g_i:d_{i}\to b^{y_{\cal D}(h)}_{(c, x)}((b, y', f'))=F(b)$$ belong to the same connected component of the category $(d\downarrow a^{y_{\cal D}(h)}_{(c, x)})$.
		\end{enumerate}
		
	\end{enumerate}
\end{theorem}

\begin{proofs}
	(i) $\Leftrightarrow$ (ii) As we observed at the beginning of this section, $C_{F}$ is locally connected if and only if for every arrow $h:Q\to Q'$ in $[{\cal D}^{\textup{op}}, \Set]$, the square
	\[
	\xymatrix{
		\Sh(\int (Q'\circ F^{\textup{op}}), J_{Q'\circ F^{\textup{op}}}) \ar[rr]^{\quad\quad(C_{p^{F}_{Q'}})_{!}} \ar[d]^{(C_{\int hF^{\textup{op}}})^{\ast}} & &  \Sh(\int Q', J_{Q'}) \ar[d]^{(C_{\int h})^{\ast}} \\
		\Sh(\int (Q\circ F^{\textup{op}}), J_{Q\circ F^{\textup{op}}}) \ar[rr]^{\quad\quad(C_{p^{F}_{Q}})_{!}} & & \Sh(\int Q, J_{Q})
	} 
	\] 
	commutes, that is, if and only if the canonical arrow
	\[
	(C_{p^{F}_{Q}})_{!} \circ (C_{\int hF^{\textup{op}}})^{\ast} \to (C_{\int h})^{\ast} \circ (C_{p^{F}_{Q'}})_{!}
	\]
	is an isomorphism.
	
	Since all the functors appearing in the above square are colimit-preserving, this arrow is an isomorphism if and only if for every object $(c, x)$ of the category $\int (Q'\circ F^{\textup{op}})$, the canonical arrow
	\[
	((C_{p^{F}_{Q}})_{!} \circ (C_{\int hF^{\textup{op}}})^{\ast})(l_{\int (Q'\circ F^{\textup{op}})}((c, x))) \to (C_{\int h})^{\ast} \circ (C_{p^{F}_{Q'}})_{!}(l_{\int (Q'\circ F^{\textup{op}})}((c, x)))
	\]
	is an isomorphism.
	
	Now, since $\int hF^{\textup{op}}$ is a comorphism of sites (cf. section \ref{sec:comorphismsofsites1}), we have
	\[
	(C_{\int hF^{\textup{op}}})^{\ast}(l_{\int (Q'\circ F^{\textup{op}})}((c, x)))\cong a_{J_{Q\circ F^{\textup{op}}}}(y_{\int (Q'\circ F^{\textup{op}})}((c, x)) \circ (hF^{\textup{op}})^{\textup{op}}).
	\]
	So, by Corollary \ref{corequivcharcanonicity}(iii),
	\[
	((C_{p^{F}_{Q}})_{!} \circ (C_{\int hF^{\textup{op}}})^{\ast})(l_{\int (Q'\circ F^{\textup{op}})}((c, x)))\cong a_{J_{Q}}(\textup{Lan}_{(p^{F}_{Q})^{\textup{op}}}(y_{\int (Q'\circ F^{\textup{op}})}((c, x)) \circ (hF^{\textup{op}})^{\textup{op}}))
	\]
	(where $\textup{Lan}$ denotes, as usual, the relevant left Kan extension functor).
	
	On the other hand, since $F$ is continuous by our hypothesis, $p^{F}_{Q'}$ is also a continuous comorphism of sites $(\int (Q'\circ F^{\textup{op}}), J_{Q'\circ F^{\textup{op}}}) \to (\int Q', J_{Q'})$  (by Proposition \ref{propcancomorphismsstableunderpullback}) and hence by Corollary \ref{correstrictioncomorphism} 	
	\[
	(C_{p^{F}_{Q'}})_{!}(l_{\int (Q'\circ F^{\textup{op}})}((c, x)))\cong (l_{\int Q'}\circ p^{F}_{Q'})((c, x))=a_{J_{Q'}}(y_{\int Q'}((F(c), x))). 
	\]
	So, since $\int h$ is a comorphism of sites (cf. section \ref{sec:comorphismsofsites1}),
	\[
	(C_{\int h})^{\ast}((C_{p^{F}_{Q'}})_{!}(l_{\int (Q'\circ F^{\textup{op}})}((c, x))))\cong a_{J_{Q}}(y_{\int Q'}((F(c), x)) \circ (\textstyle \int h)^{\textup{op}}).
	\]
	
	Next, we note that the discrete fibration associated with the presheaf $\textup{Lan}_{(p^{F}_{Q})^{\textup{op}}}(y_{\int (Q'\circ F^{\textup{op}})}((c, x)) \circ (hF^{\textup{op}})^{\textup{op}})$ on $\int Q$ is precisely the canonical functor ${\cal A}^{y_{\cal D}(h)}_{(c, x)}\to {\int Q}$ (sending an object $(c', y, f)$ of ${\cal A}^{y_{\cal D}(h)}_{(c, x)}$ to the object $(F(c'), y)$ of $\int Q$), while the discrete fibration associated with the presheaf $y_{\int Q'}((F(c), x)) \circ (\int h)^{\textup{op}}$ on $\int Q$ is the canonical projection functor ${\cal B}^{y_{\cal D}(h)}_{(c, x)} \to {\int Q}$. Moreover, the canonical arrow 
	\[
	((C_{p^{F}_{Q}})_{!} \circ (C_{\int hF^{\textup{op}}})^{\ast})(l_{\int (Q'\circ F^{\textup{op}})}((c, x))) \to (C_{\int h})^{\ast} \circ (C_{p^{F}_{Q'}})_{!}(l_{\int (Q'\circ F^{\textup{op}})}((c, x)))
	\]
	is induced by the functor $\xi^{y_{\cal D}(h)}_{(c, x)}:{\cal A}^{y_{\cal D}(h)}_{(c, x)} \to {\cal B}^{y_{\cal D}(h)}_{(c, x)}$ introduced above. So this arrow is an isomorphism if and only if the pair $(\xi^{y_{\cal D}(h)}_{(c, x)}, 1_{[{\cal A}^{y_{\cal D}(h)}_{(c, x)}], {\int Q}})$ satisfies the conditions of Proposition \ref{procofinality}. But by Proposition \ref{propcofinalcontinuous}, this holds if and only if $(\xi^{y_{\cal D}(h)}_{(c, x)}, 1_{[{\cal A}^{y_{\cal D}(h)}_{(c, x)}, {\cal D}]})$ satisfies them, where $\xi^{y_{\cal D}(h)}_{(c, x)}$ is regarded here as a morphism of fibrations $a^{y_{\cal D}(h)}_{(c, x)}\to b^{y_{\cal D}(h)}_{(c, x)}$ to $\cal D$. This concludes the proof of our thesis.
	
	(ii) $\Leftrightarrow$ (iii). Condition (iii) is the particular case of (ii) for the arrows of the form $y_{\cal D}(h)$ where $h$ is an arrow in $\cal D$. The fact that it is equivalent to the general condition follows from Corollary \ref{corsimplifiedlocalconnectedness}.
\end{proofs}

It is instructive to apply Theorem \ref{thmlocallyconnected} in the context of essential geometric morphisms between presheaf toposes induced by a functor. The problem of obtaining necessary and sufficient conditions for such a morphism to be locally connected, to which the next corollary provides a complete solution, was judged \ac\ac hard'' by Johnstone in \cite{El} (see page 653 therein), who only provided partial sufficient conditions.

\begin{corollary}\label{corlocallyconnectedpresheaf}
	Let $F:{\cal C}\to {\cal D}$ be a functor. Then the geometric morphism
	\[
	E(F):[{\cal C}^{\textup{op}}, \Set]\to [{\cal D}^{\textup{op}}, \Set]
	\]
	induced by $F$ is locally connected if and only if for any arrow $h:d_{0}\to d_{1}$ in $\cal D$, object $c$ of $\cal C$ and arrow $x:F(c)\to d_{1}$ in $\cal D$, the following conditions hold:
	\begin{enumerate}[(a)]
		\item For any object $(d, z, g)$ of the category ${\cal B}^{y_{\cal D}(h)}_{(c, x)}$ there is an object $(c', y, f)$ of the category ${\cal A}^{y_{\cal D}(h)}_{(c, x)}$ and an arrow $$s:d\to F(c')= a^{y_{\cal D}(h)}_{(c, x)}((c', y, f))$$ such that $$1_{d}:d\to d=b^{y_{\cal D}(h)}_{(c, x)}((d, z, g))$$ and $$s:d\to F(c')=b^{y_{\cal D}(h)}_{(c, x)}((F(c'), y, F(f))))$$ belong to the same connected component of the category $(d\downarrow b^{y_{\cal D}(h)}_{(c, x)})$.  
		
		\item For any object $d$ of $\cal D$ and any arrows $\alpha:d\to a^{y_{\cal D}(h)}_{(c, x)}((a, y, f))=F(a)$ and $\beta:d\to a^{y_{\cal D}(h)}_{(c, x)}((b, y', f'))=F(b)$ in $\cal D$ such that $$\alpha:d\to b^{y_{\cal D}(h)}_{(c, x)}((F(a), y, F(f)))=F(a)$$ and $$\beta:d\to b^{y_{\cal D}(h)}_{(c, x)}((F(b), y', F(f')))=F(b)$$ belong to the same connected component of the category $(d\downarrow b^{y_{\cal D}(h)}_{(c, x)})$, $$\alpha:d\to a^{y_{\cal D}(h)}_{(c, x)}((a, y, f))=F(a)$$ and $$\beta:d\to b^{y_{\cal D}(h)}_{(c, x)}((b, y', f'))=F(b)$$ belong to the same connected component of the category $(d\downarrow a^{y_{\cal D}(h)}_{(c, x)})$.
	\end{enumerate}
\end{corollary}

\begin{remark}\label{remconditionalocallyconnected}
	Condition (a) of Corollary \ref{corlocallyconnectedpresheaf} is satisfied if for any object $(d, z, g)$ of the category ${\cal B}^{y_{\cal D}(h)}_{(c, x)}$ there is an object $(c', y, f)$ of the category ${\cal A}^{y_{\cal D}(h)}_{(c, x)}$ and an arrow $$s:d\to F(c')= a^{y_{\cal D}(h)}_{(c, x)}((c', y, f))$$ such that $z=y\circ s$ and $g=F(f)\circ s$. Note that this condition is always satisfied if $F$ is a fibration (in which case $s$ can be taken to be an isomorphism).
\end{remark}

\begin{corollary}\label{corfibrationslocallyconnected}
	Let $F:{\cal C}\to {\cal D}$ be a fibration and $K$ a Grothendieck topology on $\cal D$. Then the geometric morphism $C_{F}:\Sh({\cal C}, M^{F}_{K})\to \Sh({\cal D}, K)$ induced by $F$, regarded as a (continuous) comorphism of sites $({\cal C}, M^{F}_{K})\to ({\cal D}, K)$, is locally connected.
\end{corollary}

\begin{proofs}
	As observed in Remark \ref{remconditionalocallyconnected}, condition (a) of Corollary \ref{corlocallyconnectedpresheaf} is satisfied by $F$, so it only remains to show that also condition (b) holds.
	
	To establish this, we show that any zig-zag connecting two objects $\alpha:d\to a^{y_{\cal D}(h)}_{(c, x)}((a, y, f))=b^{y_{\cal D}(h)}_{(c, x)}((F(a), y, F(f)))=F(a)$ and $\beta:d\to a^{y_{\cal D}(h)}_{(c, x)}((b, y', f'))=b^{y_{\cal D}(h)}_{(c, x)}((F(b), y', F(f')))=F(b)$ in the category $(d\downarrow b^{y_{\cal D}(h)}_{(c, x)})$ can be \ac lifted' to a zig-zag connecting them in the category $(d\downarrow a^{y_{\cal D}(h)}_{(c, x)})$. Note that we can suppose without loss of generality the arrows $f$ and $f'$ to be cartesian (at the cost of possibly replacing them with their cartesian images, which does not affect the connected component at which the given object belongs).
	
	Given an arrow $\delta:d\to a^{y_{\cal D}(h)}_{(c, x)}((a, y, f))=b^{y_{\cal D}(h)}_{(c, x)}((F(a), y, F(f)))=F(a)$ in $\cal D$ and an arrow $\gamma:d \to b^{y_{\cal D}(h)}_{(c, x)}((d', z, g))$ connected to it in the category $(d\downarrow b^{y_{\cal D}(h)}_{(c, x)})$, there are two basic cases:
	\begin{enumerate}[(1)]
		\item There is an arrow $s:(F(a), y, F(f))\to (d', z, g)$ in the category ${\cal B}^{y_{\cal D}(h)}_{(c, x)}$ such that $s\circ \delta=\gamma$; 
		
		\item There is an arrow $t:(d', z, g) \to (F(a), y, F(f))$ in the category ${\cal B}^{y_{\cal D}(h)}_{(c, x)}$ such that $t\circ \gamma=\delta$. 
	\end{enumerate} 

	Let us show that in each of this cases we can lift this connection to a connection in the category $(d\downarrow a^{y_{\cal D}(h)}_{(c, x)})$.
	\begin{enumerate}
		\item[(1)] Note preliminarily that, since $s:(F(a), y, F(f))\to (d', z, g)$ is an arrow in the category ${\cal B}^{y_{\cal D}(h)}_{(c, x)}$, we have $g\circ s=F(f)$ and $z\circ s=y$. Since $F$ is a fibration, there are an isomorphism $\tau:F(a')\to d'$ and a cartesian arrow $k:a'\to a$ such that $g\circ \tau=F(k)$. Now, since $k$ is cartesian and $F(k)\circ (\tau^{-1}\circ s)=F(f)$, there is a unique arrow $h:a\to a'$ in $\cal C$ such that $F(h)=\tau^{-1}\circ s$ and $k\circ h=f$. Therefore $h$ is an arrow $(a, y, f)\to (a', z\circ \tau, k)$ in the category ${\cal A}^{y_{\cal D}(h)}_{(c, x)}$ which connects $\delta:d\to a^{y_{\cal D}(h)}_{(c, x)}((a, y, f))=F(a)$ with $F(h)\circ \delta:d\to a^{y_{\cal D}(h)}_{(c, x)}((a', z\circ \tau, k))$ and which \ac lifts' our given connection in the sense that $F(h)\circ \delta:d\to a^{y_{\cal D}(h)}_{(c, x)}((a', z\circ \tau, k))=b^{y_{\cal D}(h)}_{(c, x)}((F(a'), z\circ \tau, F(k)))$ is isomorphic in the category $(d\downarrow b^{y_{\cal D}(h)}_{(c, x)})$ to the object $\gamma:d \to b^{y_{\cal D}(h)}_{(c, x)}((d', z, g))$. 
		
		\item[(2)] Note preliminarily that, since $t:(d', z, g)\to (F(a), y, F(f))$ is an arrow in the category ${\cal B}^{y_{\cal D}(h)}_{(c, x)}$, we have $y\circ t=z$ and $F(f)\circ t=g$. Since $F$ is a fibration, there are an object $a''$ of $\cal C$, an isomorphism $\chi:F(a'')\to d'$ and a cartesian arrow $u$ such that $t\circ \chi=F(u)$. Then $t$ is an arrow $(a'', y\circ F(u), f\circ u)\to (a, y, f)$ in the category ${\cal A}^{y_{\cal D}(h)}_{(c, x)}$ which connects $\delta:d\to a^{y_{\cal D}(h)}_{(c, x)}((a, y, f))=F(a)$ with $\chi^{-1}\circ \gamma:d\to a^{y_{\cal D}(h)}_{(c, x)}((a'', y\circ F(u), f\circ u))$ and \ac lifts' our given connection in the sense that $\chi^{-1}\circ \gamma:d\to a^{y_{\cal D}(h)}_{(c, x)}((a'', y\circ F(u), f\circ u))=a^{y_{\cal D}(h)}_{(c, x)}((F(a''), y\circ F(u), F(f\circ u)))$ is isomorphic in the category $(d\downarrow b^{y_{\cal D}(h)}_{(c, x)})$ to the object $\gamma:d \to b^{y_{\cal D}(h)}_{(c, x)}((d', z, g))$.    
	\end{enumerate}

	Note that all the \ac lifts' constructed in this way have automatically the property that the arrow in the third component is cartesian (both $k$ and $u$ are taken to be cartesian). Now, by inductively applying the argument starting from one of the extremes of the zig-zag connecting the two given objects, we are eventually led to an arrow in the category ${\cal B}^{y_{\cal D}(h)}_{(c, x)}$ under $d$ between the object at the other extreme of the zig-sag and an object in the image of the canonical functor $(d\downarrow a^{y_{\cal D}(h)}_{(c, x)}) \to (d\downarrow b^{y_{\cal D}(h)}_{(c, x)})$, whose arrow in the third component is cartesian. Now, similarly to the proof of Lemma \ref{lemmacartesianrelmono}(iii), the cartesianness of the arrows in the third components of these two objects ensures that the arrow in ${\cal B}^{y_{\cal D}(h)}_{(c, x)}$ actually lifts uniquely to an arrow in ${\cal A}^{y_{\cal D}(h)}_{(c, x)}$ satisfying the required commutativity conditions. This completes the lift of the zig-zag and hence our proof.
\end{proofs}

\subsection{The relative comprehensive factorization}\label{sec:relativecomprehensivefactorization}

Recall that a geometric morphism $f:{\cal F}\to {\cal E}$ is said to be \emph{connected} if $f^{\ast}$ is full and faithful. We shall say that an essential geometric morphism $f:{\cal F}\to {\cal E}$ is \emph{terminally connected} if $f_{!}(1_{\cal F})$ is the terminal object of $\cal E$. 

As observed in the proof of Lemma C3.3.3 \cite{El}, if $f$ is connected then the counit of the adjunction $(f_{!}\dashv f^{\ast})$ is an isomorphism and hence $f_{!}(1)\cong f_{!}(f^{\ast}(1))\cong 1$, that is, $f$ is terminally connected. By that lemma, if $f$ is locally connected, then $f$ is connected if and only if it is terminally connected. 

We will need the following result (where, given a topos $\cal E$ and an object $A$ of it, we shall denote by $U_{A}$ the canonical local homeomorphism ${\cal E}\slash A \to {\cal E}$, whose inverse image $U_{A}^{\ast}$ sends any object $X$ of $\cal E$ to the object $[A\times X\to X]$ of ${\cal E}\slash A$ given by the canonical projection arrow):

\begin{lemma}\label{lemmaslicetoposclassification}
	Let $A$ be an object of a topos $\cal E$. Then, for any topos $\cal F$, the geometric morphisms $g:{\cal F}\to {\cal E}\slash A$ correspond bijectively to the pairs $(f, \alpha)$, where $f$ is a geometric morphism ${\cal F}\to {\cal E}$ and $\alpha$ is an arrow $1_{\cal F}\to f^{\ast}(A)$ in $\cal F$ (naturally in $A$ and $\cal F$). More specifically:	
	\begin{itemize}
		\item given a geometric morphism $g:{\cal F}\to {\cal E}\slash A$, the pair associated with it is $(U_{A}\circ g, \alpha_{g})$, where $\alpha_{g}:1_{\cal F}\to g^{\ast}(U_{A}^{\ast}(A))$ is the arrow $g^{\ast}(\Delta_{A})$, where $\Delta_{A}:1_{{\cal E}\slash A}\to U_{A}^{\ast}(A)$ is the arrow given by the diagonal $A\to A\times A$ on the object $A$; 
		
		\item given a pair $(f, \alpha)$ consisting of a geometric morphism $f:{\cal F}\to {\cal E}$ and an arrow $\alpha:1_{\cal F} \to f^{\ast}(A)$, the geometric morphism $f_{\alpha}:{\cal F}\to {\cal E}\slash A$ corresponding to the pair $(f, \alpha)$ is defined by setting, for any object $[x:X\to A]$ of ${\cal E}\slash A$, $f_{\alpha}^{\ast}([x])\mono f^{\ast}(X)$ equal to the equalizer of the arrows $\alpha\circ !_{f^{\ast}(X)}$ and $f^{\ast}(x)$ (where $!_{f^{\ast}(X)}$ is the unique arrow $f^{\ast}(X)\to 1_{\cal F}$).
	\end{itemize} 		 
\end{lemma}

\begin{proofs}
	The proof is a straightforward consequence of the following easy-to-prove fact: any object $[x:X\to A]$ of ${\cal E}\slash A$ can be canonically identified, via the monomorphism $[x]\mono U_{A}^{\ast}(X)$ given by the arrow $<x, 1_{X}>:X\to A\times X$, with the equalizer of the arrows $U_{A}^{\ast}(x), \Delta_{A}\circ !_{U_{A}^{\ast}(X)}:U_{A}^{\ast}(X) \doublerightarrow{}{}   U_{A}^{\ast}(A)$, where $\Delta_{A}:1_{{\cal E}\slash A}\to U_{A}^{\ast}(A)$ is the arrow given by the diagonal $A\to A\times A$ on $A$ and $!_{U_{A}^{\ast}(X)}$ is the unique arrow $U_{A}^{\ast}(X) \to 1_{{\cal E}\slash A}$. The details are straightforward and left to the reader.
\end{proofs}

\begin{remark}\label{remlemmaslice}
	For any geometric morphisms $g:{\cal F}\to {\cal E}\slash A$, $g':{\cal F}'\to {\cal E}\slash A$ and $h:{\cal F}'\to {\cal F}$ such that $g\circ h\cong g'$, we have $\alpha_{g'}=h^{\ast}(\alpha_{g})$. In particular, for any composable geometric morphisms $g$ and $h$, we have
	\[
	\alpha_{g\circ h}=h^{\ast}(\alpha_{g}).
	\] 
\end{remark}

The following result is a generalization and a strengthening of Lemma C3.3.4 \cite{El}: 

\begin{proposition}\label{proporthogonality}
Terminally connected morphisms are orthogonal to local homeomorphisms
in the $2$-category of Grothendieck toposes; that is, for any commutative square
\begin{equation*}
\begin{tikzcd}
{\cal F} \arrow[rr, "f"]   \arrow[d, "m"]   & &  {\cal E} \arrow[d, "n"] \arrow[dll, dashed, "k"{above}]  \\
{\cal G} \ar[rr, "g"] & & {\cal H},   
\end{tikzcd}
\end{equation*}	
where $f$ is terminally connected and $g$ is a local homeomorphism, there exists a morphism $k:{\cal E}\to {\cal G}$ (unique up
to unique $2$-isomorphism) making both triangles commute. 	
\end{proposition}

\begin{proofs}
Since $g$ is a local homeomorphism, we can suppose it to be of the form $U_{A}:{\cal H}\slash A \to {\cal H}$ for some object $A$ of $\cal H$. By Lemma \ref{lemmaslicetoposclassification}, the geometric morphism $m$ corresponds to an arrow $\alpha:1_{\cal F}\to (g\circ m)^{\ast}(A)\cong (n\circ f)^{\ast}(A)=f^{\ast}(n^{\ast}(A))$. Again by the Lemma, giving a geometric morphism $k:{\cal E}\to {\cal G}$ such that $g\circ k\cong n$ corresponds to giving an arrow $\beta:1_{\cal E}\to n^{\ast}(A)$ in $\cal E$. By Remark \ref{remlemmaslice}, the commutativity condition $k\circ f\cong m$ can be formulated in terms of pairs corresponding to these morphisms via Lemma \ref{lemmaslicetoposclassification} as the condition $f^{\ast}(\beta)=\alpha$. We thus have to show that there is exactly one arrow $\beta:1_{\cal E}\to n^{\ast}(A)$ such that $f^{\ast}(\beta)=\alpha$. If we apply the functor $f_{!}$ to the equality $f^{\ast}(\beta)=\alpha$, we obtain the equality $f_{!}(f^{\ast}(\beta))=f_{!}(\alpha)$. But the naturality of the counit $\epsilon$ of the adjunction between $f^{\ast}$ and $f_{!}$ implies the commutativity of the diagram
  \begin{equation*}
  \begin{tikzcd}
  f_{!}(f^{\ast}(1_{\cal E})) \arrow[rr, "f_{!}(f^{\ast}(\beta))"]   \arrow[d, "\epsilon_{1_{\cal E}}"]   & &  f_{!}(f^{\ast}(n^{\ast}(A))  \arrow[d, "\epsilon_{n^{\ast}(A)}"]  \\
  1_{\cal E} \ar[rr, "\beta"] & & n^{\ast}(A),   
  \end{tikzcd}
  \end{equation*} 
which forces, since $f_{!}(f^{\ast}(1_{\cal E}))=1_{\cal E}$ and $\epsilon_{1_{\cal E}}=1_{1_{\cal E}}$ by the terminal connectedness of $f$, $\beta$ to be equal to the arrow $\epsilon_{n^{\ast}(A)}\circ f_{!}(f^{\ast}(\beta))$, in other words to the transpose $1_{\cal E}=f_{!}(1_{\cal E})\to n^{\ast}(A)$ of the arrow $\alpha$.
\end{proofs}

The following result generalizes Proposition C3.3.5 \cite{El}.

\begin{proposition}\label{propfactorizationdiscrete}
Any essential geometric morphism can be factored, uniquely up to 
equivalence, as a terminally connected morphism followed by a local
homeomorphism. In particular, an essential geometric morphism $f:{\cal F}\to {\cal E}$ factors as a terminally connected morphism $f':{\cal F}\to {\cal E}\slash f_{!}(1)$ followed by the canonical local homeomorphism $U_{f_{!}(1)}:{\cal E}\slash f_{!}(1) \to {\cal E}$. 
\end{proposition}

\begin{proofs}
The uniqueness of the factorization follows from Proposition \ref{proporthogonality}, so it remains to prove its existence. 
If $A=f_{!}(1)$ and $\alpha$ is the component at $1$ of the unit $1\to f^{\ast}\circ f_{!}$ of the adjunction between $f_{!}$ and $f^{\ast}$ then the geometric morphism $f_{\alpha}$ associated with the pair $(f, \alpha)$ as in Lemma \ref{lemmaslicetoposclassification} is essential and its essential image $(f_{\alpha})_{!}$ sends any object $F$ of $\cal F$ to the object $[f_{!}(!_{F}):f_{!}(F)\to f_{!}(1)]$ of ${\cal E}\slash f_{!}(1)$ (where $!_{F}$ is the unique arrow $F\to 1$). Indeed, it is immediate to see that, under the identification of $f_{\alpha}^{\ast}([x])$ with a subobject of $f^{\ast}(X)$, the natural bijection
\[
\textup{Hom}_{\cal F}(F, f^{\ast}(X))\cong \textup{Hom}_{\cal E}(f_{!}(F), X) 
\]
restricts to a bijection 
\[
\textup{Hom}_{\cal F}(F, f_{\alpha}^{\ast}([x]))\cong \textup{Hom}_{{\cal E}\slash f_{!}(1)}(f_{!}(F), x)
\]
(for any object $x:X\to f_{!}(1)$ of ${\cal E}\slash f_{!}(1)$).  

In particular, $(f_{\alpha})_{!}$ sends $1$ to the terminal object of ${\cal E}\slash f_{!}(1)$, that is $f_{\alpha}$ is terminally connected.  
Moreover, composing $f_{\alpha}:{\cal F}\to {\cal E}\slash f_{!}(1)$ with the local homeomorphism ${\cal E}\slash f_{!}(1) \to {\cal E}$ yields precisely $f$, as required.	
\end{proofs}

\begin{remark}
	The composite of any local homeomorphism with a terminally connected morphism is always an essential geometric morphism.
\end{remark}

The following result enlightens a link between relative cofinality of a continuous comorphism of sites and terminal connectedness of the associated geometric morphism: 

\begin{proposition}\label{proplinkterminallyconnected}
	Let $F:({\cal C}, J)\to ({\cal D}, K)$ be a continuous comorphism of small-generated sites. Then the geometric morphism $C_{F}:\Sh({\cal C}, J)\to \Sh({\cal D}, K)$ induced by $F$ is terminally connected (resp. connected and locally connected) if and only if $F$ is $K$-cofinal (resp. $F$ is $K$-cofinal and satisfies the conditions of Theorem \ref{thmlocallyconnected}).
\end{proposition}

\begin{proofs}
Since $F$ is continuous, we have  $$(C_{F})_{!}(1)=(C_{F})_{!}(\textup{colim}_{\Sh({\cal C}, J)}(l))\cong \textup{colim}(l'\circ F),$$ whence our thesis follows from Corollary \ref{corcofinalcharacterization}.	
\end{proofs}

Proposition \ref{propfactorizationdiscrete} will allow us to interpret in terms of a topos-theoretic invariant the comprehensive factorization of a functor, as introduced in \cite{RossWalters}. 

Recall that discrete fibrations on a category $\cal D$ correpond precisely to presheaves on $\cal D$. 

\begin{definition}
	We say that a discrete fibration $p:{\cal E}\to {\cal D}$ is \emph{$K$-glueing} if the presheaf corresponding to it is a $K$-sheaf, equivalently, if for any $K$-covering sieve $T$ on an object $d$ of $\cal D$, the map
	\[
	\textup{ \bf Fib}_{\cal D}({\cal D}\slash d, p) \to \textup{ \bf Fib}_{\cal D}({\textstyle \int T}, p)
	\] 
	is a bijection (where $\textup{ \bf Fib}_{\cal D}(-,-)$ denotes the set of morphisms of fibrations over $\cal D$). 
\end{definition}

Given a category $\cal C$, we define $\textbf{Cat}\slash\slash {\cal C}$ to be the category whose objects are the functors with codomain $\cal C$ and whose arrows $(F:{\cal A}\to {\cal C}) \to (F':{\cal B} \to {\cal C})$ and the pairs $(G, \alpha)$, where $G$ is a functor ${\cal A}\to {\cal B}$ and $\alpha$ is an isomorphism $F \to F'\circ G$ (with the obvious notion of composition between them).  

The following result clarifies the identification between sheaves and the associated discrete fibrations. 

\begin{proposition}\label{propadjunctiondiscrete}
Let $({\cal C}, J)$ be a small-generated site. Then	the functors 
\[
U:\Sh({\cal C}, J)\to \textup{\bf Cat}\slash\slash {\cal C}
\]
and
\[
V:\textup{\bf Cat}\slash\slash {\cal C} \to \Sh({\cal C}, J)
\]
which send respectively a $J$-sheaf $P$ on $\cal C$ to the associated discrete fibration $$\pi_{P}:{\textstyle \int P}\to {\cal C}$$ and a functor $F:{\cal A}\to {\cal C}$ to the $K$-sheaf $$F_{K}:=\textup{colim}_{\Sh({\cal C}, J)}(l_{\cal C}\circ F)$$ on $\cal D$ (and act on arrows in the obvious way) are adjoint to each other ($U$ on the right and $V$ on the left), and this adjunction restricts to an equivalence between $\Sh({\cal C}, J)$ and the full subcategory of $\textup{\bf Cat} \slash\slash {\cal C}$ on the discrete $J$-glueing fibrations.  
\end{proposition}

\begin{proofs}
	Straightforward and left to the reader.
\end{proofs}

Given a functor $F:{\cal C}\to {\cal D}$, we can associate with it, as in Proposition \ref{propadjunctiondiscrete}, the $K$-sheaf 
$F_{K}=\textup{colim}_{\Sh({\cal D}, K)}(l'\circ F).$

Note that there is a canonical functor $\xi^{F}_{K}:{\cal C}\to {\int F_{K}}$ such that $\pi_{F_{K}}\circ \xi^{F}_{K}=F$ (which assigns to each object $c$ of $\cal C$ the element of $F_{K}$ corresponding to the colimit arrow $l'(F(c))\to F_{K}$ indexed by the object $c$). As we shall see below, the functor $\xi^{F}_{K}$ is $M^{F_{K}}_{K}$-cofinal. This motivates the following 

\begin{definition}
Let $F:{\cal C}\to {\cal D}$ be a functor and $K$ a Grothendieck topology  on $\cal D$. The \emph{$K$-comprehensive factorization of $F$} is given by the composite $F_{K}\circ \xi^{F}_{K}$:
\begin{equation*}
\begin{tikzcd}
{\cal C} \arrow[rr, "F"]   \arrow[dr, "\xi^{F}_{K}"]   & &  {\cal D}           \\
&  \int F_{K} \ar[ur, "\pi_{F_{K}}"] &  
\end{tikzcd}
\end{equation*}	
\end{definition}

\begin{proposition}\label{proporthogonality}
	The following orthogonality property of cofinal functors with respect of glueing fibrations (along continuous functors) holds: for any commutative square
	\begin{equation*}
	\begin{tikzcd}
	{\cal A} \arrow[rr, "F"]   \arrow[d, "u"]   & &  {\cal C} \arrow[d, "v"] \arrow[dll, dashed, "G"{above}]  \\
	{\cal B} \ar[rr, "p"] & & {\cal C}'   
	\end{tikzcd}
	\end{equation*}	
	where $\cal C$ and ${\cal C}'$ are respectively endowed with Grothendieck topologies $J$ and $J'$ such that $v$ is $(J, J')$-continuous, $F$ is $J$-cofinal and $p$ is a $J'$-glueing discrete fibration, there is a unique (up to isomorphism) functor $G:{\cal C}\to {\cal B}$ which makes both triangles commute (up to isomorphism). 
\end{proposition} 
 
\begin{proofs}
Since $p$ is $J'$-glueing, by Proposition \ref{propadjunctiondiscrete} giving a functor $G:{\cal C}\to {\cal B}$ such that $p\circ G\cong v$ amounts precisely to giving	an arrow $$\textup{colim}_{\Sh({\cal C}', J')}(l_{{\cal C}'}\circ v) \to P_{p},$$ where $P_{p}$ is the $J'$-sheaf corresponding to $p$ as in the proposition. 

By Proposition \ref{propcontinuitycofinality},
\[
\textup{colim}_{\Sh({\cal C}', J')}(l_{{\cal C}'}\circ v)\cong \textup{colim}_{\Sh({\cal C}', J')}(l_{{\cal C}'}\circ v\circ F)=\textup{colim}_{\Sh({\cal C}', J')}(l_{{\cal C}'}\circ p\circ u),
\] 
and we clearly have an arrow 
\[
\textup{colim}_{\Sh({\cal C}', J')}(l_{{\cal C}'}\circ p\circ u)\to \textup{colim}_{\Sh({\cal C}', J')}(l_{{\cal C}'}\circ p)=P_{p}
\]
induced by $u$, so by composition we obtain an arrow
\[
\textup{colim}_{\Sh({\cal C}', J')}(l_{{\cal C}'}\circ v) \to P_{p}.
\]

It is immediate to see, by exploiting the adjunction of  Proposition \ref{propadjunctiondiscrete}, that the functor $G$ corresponding to it satisfies, besides the condition $p\circ G\cong v$, the relation $G\circ F\cong u$.  

It remains to show the uniqueness (up to isomorphism) of a functor $G$ making the two triangles commute. This can be deduced from the ($2$-)universal property of the comma construction, noticing that for any  discrete fibration $p:{\cal B}\to {\cal C}'$ associated with a presheaf $P_{p}$ on ${\cal C}'$, the following diagram, where $\textbf{1}$ is the category with one object and one morphism and $\ast$ is the functor sending this object to the singleton $1_{\Set}$, is a comma square in $\textbf{Cat}$:
\begin{equation*}
\begin{tikzcd}
{\cal B}^{\textup{op}} \arrow[rr, "p^{\textup{op}}"]   \arrow[d]   & &  {{\cal C}'}^{\textup{op}} \arrow[d, "P_{p}"] \\
\textbf{1} \ar[rr, "\ast"] & & \Set.   
\end{tikzcd}
\end{equation*}
 
Indeed, $p^{\textup{op}}:{\cal B}^{\textup{op}}\to {{\cal C}'}^{\textup{op}}$ is isomorphic to ${\pi_{P_{p}}}^{\textup{op}}: {\int P_{p}}^{\textup{op}}\to {{\cal C}'}^{\textup{op}}$, which is in turn isomorphic to the opposite of the canonical projection functor $(\ast \downarrow P_{p})\to {\cal C}'$ (cf. the proof of Theorem 4 \cite{RossWalters}). From this universal property and the fact that $p$ (equivalently, $p^{\textup{op}}$) reflects isomorphisms, it thus follows that for any category ${\cal Z}$ and any functors $Z_{1}, Z_{2}:{\cal Z} \to {\cal B}^{\textup{op}}$, if there is an isomorphism $\chi:p^{\textup{op}}\circ Z_{1}\cong p^{\textup{op}}\circ Z_{2}$ such that $P_{p}\chi$ is the identity natural transformation on the functor $\ast \circ !_{\cal Z}$, where $!_{\cal Z}$ is the unique functor ${\cal Z}\to \textbf{1}$, then $Z_{1}$ and $Z_{2}$ are naturally isomorphic.   
\end{proofs} 
 
\begin{remark}
	If $J$ and $J'$ are the trivial Grothendieck topologies then the orthogonality property of Proposition \ref{proporthogonality} specializes to the well-known orthogonality property of cofinal functors with respect to discrete fibrations, first established (in the dual form of initial functors and discrete $0$-fibrations) in \cite{RossWalters}.
\end{remark} 
\begin{proposition}
Let	$F:{\cal C}\to {\cal D}$ be a functor and $K$ a Grothendieck topology  on $\cal D$. 
	\begin{enumerate}[(i)]
		\item The $K$-comprehensive factorization of $F$ is characterized by being the unique (up to equivalence) factorization of $F$ as a  $M^{p}_{K}$-cofinal functor ${\cal C}\to {\cal E}$ followed by a $K$-glueing fibration $p:{\cal E}\to {\cal D}$.
		
		\item If $F$ is a continuous comorphism of sites $({\cal C}, J)\to ({\cal D}, K)$ then $\xi^{F}_{K}:({\cal C}, J)\to ({\int F_{K}}, M^{\pi_{F_{K}}}_{K})$ and $F_{K}:({\int F_{K}}, M^{\pi_{F_{K}}}_{K}) \to ({\cal D}, K)$ are continuous comorphism of sites and $C_{F}\cong C_{F_{K}}\circ C_{\xi^{F}_{K}}$ is the terminally connected-local homeomorphism factorization of the geometric morphism $C_{F}:\Sh({\cal C}, J)\to \Sh({\cal D}, K)$ (cf. Proposition \ref{propfactorizationdiscrete}). 
	\end{enumerate}
\end{proposition} 

\begin{proofs}
	(i) Recall that we have
	\[
	F=F_{K}\circ \xi^{F}_{K}.
	\]
	First, let us show that the functor $\xi^{F}_{K}:{\cal C}\to \int F_{K}$ is $M^{F_{K}}_{K}$-cofinal. By Corollary \ref{corcofinalcharacterization}, $\xi^{F}_{K}$ is $M^{F_{K}}_{K}$-cofinal if and only if the unique arrow
	\[
	\textup{colim}_{\Sh({\cal C}, J)}(l_{\int F_{K}}\circ \xi^{F}_{K}) \to 1_{\Sh({\int F_{K}}, M^{F_{K}}_{K})}
	\] 
	is an isomorphism. Now, since $(C_{\pi_{F_{K}}})_{!}$ reflects isomorphisms (as it is isomorphic to the essential image of the geometric morphism $U^{\Sh({\cal D}, K)}_{F_{K}}$, namely the forgetful functor $\Sh({\cal D}, K)\slash F_{K}\to \Sh({\cal D}, K)$) and the square
	$$
	\xymatrix{
		{\int F_{K}} \ar[rr]^{\quad\pi_{F_{K}}} \ar[d]_{l_{{\int F_{K}}}} & & {\cal D} \ar[d]^{l_{\cal D}} \\
		\Sh({{\int F_{K}}}, M^{F_{K}}_{K}) \ar[rr]^{\quad(C_{\pi_{F_{K}}})_{!}} & & \Sh({\cal D}, K)
	}
	$$ 
	commutes (since $\pi_{F_{K}}$ is a continuous comorphism of sites $({{\int F_{K}}}, M^{F_{K}}_{K}) \to ({\cal D}, K)$ -- cf. Corollary \ref{correstrictioncomorphism}), the above arrow is an isomorphism if and only if its image 
	\[
	(C_{\pi_{F_{K}}})_{!}(\textup{colim}_{\Sh({\cal C}, J)}(l_{\int F_{K}}\circ \xi^{F}_{K})) \to (C_{\pi_{F_{K}}})_{!}(1_{\Sh({\int F_{K}}, M^{F_{K}}_{K})}) 
	\]
	under $(C_{\pi_{F_{K}}})_{!}$ is an isomorphism in $\Sh({\cal D}, K)$. But this arrow is (isomorphic to) the canonical one
	\[
	\textup{colim}_{\Sh({\cal D}, K)}(l_{\cal D}\circ F) \to F_{k},
	\]
	and this is an isomorphism by definition of $F_{k}$. This argument, rewritten for a general factorization of $F$ as a $K$-glueing fibration $p:{\cal E}\to {\cal D}$ composed with a $M^{p}_{K}$-cofinal functor ${\cal C}\to {\cal E}$, shows that $p$ must necessarily be isomorphic to $F_{k}$. The uniqueness of a $M^{F_{K}}_{K}$-cofinal factorization of $F$ through $\pi_{F_{K}}$ follows from Proposition \ref{proporthogonality}.
	
	(ii) First, we notice that $F_{K}:({\int F_{K}}, M^{\pi_{F_{K}}}_{K}) \to ({\cal D}, K)$ is a continuous comorphism of sites since $F:({\cal C}, J)\to ({\cal D}, K)$ is, by Proposition \ref{corliftalongdiscretefibration}. To conclude our thesis, it thus suffices to show that $F_{K}$ induces a terminally connected geometric morphism and that $\pi_{F_{K}}$ induces a local homeomorphism. But the first claim follows from Proposition \ref{proplinkterminallyconnected}, while the second follows from the equivalence of section \ref{sec:example}.
\end{proofs}

\begin{remark}
	In the case $K$ is the trivial topology on $\cal D$ then $M^{p}_{K}$ is also the trivial topology (for any fibration $p$), so the $K$-comprehensive factorization of $F$ specializes to the opposite of the usual comprehensive factorization (in the sense of \cite{RossWalters}). 
\end{remark}

\section{Denseness conditions and equivalences of toposes}\label{sec:denseness}

\begin{definition}
	A morphism of sites $F:({\cal C}, J)\to ({\cal D}, K)$ is said to be \emph{dense} if it satisfies the following properties:
	\begin{enumerate}[(i)]	
		\item $P$ is a $J$-covering family in $\cal C$ if and only if $F(P)$ is a $K$-covering family in $\cal D$;
		\item $F$ is \emph{$K$-dense} in the sense that for any object $d$ of $\cal D$ there exists a $K$-covering family of arrows $d_{i}\to d$ whose domains $d_{i}$ are in the image of $F$;
		\item for every $c_{1}, c_{2}\in {\cal C}$ and any arrow $g:F(c_{1})\to F(c_{2})$ in $\cal D$, there exist a $J$-covering family of arrows $f_{i}:c'_{i}\to c_{1}$ and a family of arrows $k_{i}:c'_{i}\to c_{2}$ such that $g\circ F(f_{i})=F(k_{i})$ for all $i$.
	\end{enumerate}
\end{definition}

\begin{remark}\label{remmorphismsofsitesfaithfulness}
	This definition is adapted from section 11 of \cite{Shulman}, where a functor $F:{\cal C}\to {\cal D}$ is defined to be a dense morphism of sites if it satisfes the three conditions above plus the following one (which we shall call $J$-faithfulness in section \ref{sec:localconditions}): 
	\begin{enumerate}[(i)]
	\item[(iv)] for any arrows $f_{1}, f_{2}:c_{1}\to c_{2}$ in $\cal C$ such that $F(f_1)=F(f_2)$ there exists a $J$-covering family of arrows $k_{i}:c'_{i}\to c_{1}$ such that $f_1 \circ k_{i}=f_{2}\circ k_{i}$ for all $i$. 
    \end{enumerate}
	In fact, Theorem 11.2 \cite{Shulman}(b) proves that every such functor $F$ is a morphism of sites $({\cal C}, J)\to ({\cal D}, K)$. On the other hand, if $F$ is already a morphism of sites, condition (iv) is unnecessary as it is implied by condition (i) (cf. condition (iv) in Definition \ref{defmorphismsites}).  
	
\end{remark}

By Theorem 11.8 \cite{Shulman}, if $F$ is a dense morphism of sites then the associated geometric morphism $\Sh(F):\Sh({\cal D}, K)\to \Sh({\cal C}, J)$ is an equivalence. The following proposition shows that, conversely, if $F$ is a morphism of sites such that $({\cal D}, K)$ is subcanonical and $\Sh(F)$ is an equivalence then $F$ is a dense morphism of sites.

\begin{proposition}\label{propdensitysubcanonical}
	Let $F:({\cal C}, J)\to ({\cal D}, K)$ be a morphism of sites. Suppose that $K$ is subcanonical. Then, if the geometric morphism $\Sh(F):\Sh({\cal D}, K)\to \Sh({\cal C}, J)$ is an equivalence, $F$ is a dense morphism of sites $({\cal C}, J)\to ({\cal D}, K)$. 
\end{proposition}

\begin{proofs}
	Let us check that all the conditions in the definition of a dense morphism of sites are satisfied.
	
	(i) This condition follows from the commutativity of the diagram
	$$
	\xymatrix{
		{\mathcal C} \ar[r]^F \ar[d]_{l} &{\cal D} \ar[d]^{l'} \\
		\Sh({\cal C}, J) \ar[r]^{\Sh(F)^{\ast}} & \Sh({\cal D}, K)
	}
	$$
	by exploiting the fact that a sieve $P$ (resp. $Q$) in $\cal C$ (resp. in $\cal D$) is $J$-covering (resp. $K$-covering) if and only if it is sent by $l$ (resp. by $l'$) to an epimorphic family, since for any arrow $f$ in $\cal C$, $l'(F(f))\cong \Sh(F)^{\ast}(l(f))$ and $\Sh(F)^{\ast}$ is an equivalence.
	
	(ii) By the commutativity of the above square, the objects of the form $l'(F(c))\simeq \Sh(F)^{\ast}(l(c))$ form a separating set for the topos $\Sh({\cal D}, K)$. Therefore, for any $d\in {\cal D}$, the family of arrows from objects of the form $l'(F(c))$ (for $c\in {\cal C}$) to $l'(d)$ is epimorphic. Since, by the subcanonicity of $K$, all these arrows are of the form $l'(f)$ for $f$ an arrow in $\cal D$, this means that the family of arrows from objects of the form $F(c)$ to $d$ is $K$-covering in $\cal D$, as required.  
	
	(iii) By the commutativity of the above square, we have $l'(g):l'(F(x))\to l'(F(y))$ is the image under $\Sh(F)^{\ast}$ of an arrow $\xi:l(x)\to l(y)$ in $\Sh({\cal C}, J)$. By Proposition \ref{propexplicit}(i), there exists a $J$-covering family of arrows $f_{i}:x_{i}\to x$ and a family of arrows $g_{i}:x_{i}\to y$ such that $\xi \circ l(f_{i})=l(g_{i})$ for all $i$. Applying $\Sh(F)^{\ast}$ yields $l'(g)\circ  l'(F(f_{i}))=l'(F(g_{i}))$ for all $i$, whence, since $l'$ is faithful ($K$ being subcanonical), $g\circ  F(f_{i})= F(g_{i})$ for all $i\in I$, as required.   
\end{proofs}

Recall that if $\cal C$ is the full subcategory of a Grothendieck topos $\cal E$ on a family of objects which is separating for it then $\cal E$ is equivalent to the topos $\Sh({\cal C}, J^{\textup{can}}_{\cal E}|_{\cal C})$, where $J^{\textup{can}}_{\cal E}|_{\cal C}$ is the Grothendieck topology on $\cal C$ induced by the canonical topology $J^{\textup{can}}_{\cal E}$ on $\cal E$: the $J^{\textup{can}}_{\cal E}|_{\cal C}$-covering sieves are those which generate $J^{\textup{can}}_{\cal E}$-covering sieves in $\cal E$. Given a small-generated site $({\cal C}, J)$, we shall denote by $C_{J}^{\cal C}$ the Grothendieck topology $J^{\textup{can}}_{\Sh({\cal C}, J)}|_{a_{J}({\cal C})}$; so we have an equivalence
\[
\Sh({\cal C}, J)\simeq \Sh(a_{J}({\cal C}), C_{J}^{\cal C}).
\] 
This equivalence is clearly induced by the morphism of sites 
\[
l:({\cal C}, J)\to (a_{J}({\cal C}), C_{J}^{\cal C}).
\]

\subsection{Weakly dense morphisms of sites}

In light of Proposition \ref{propdensitysubcanonical}, it is natural to give the following definition.

\begin{definition}\label{defweaklydensemorphism}
	A morphism of sites $F:({\cal C}, J)\to ({\cal D}, K)$ is \emph{weakly dense} if the morphism of sites $l'\circ F:({\cal C}, J)\to (a_{K}({\cal D}), C_{K}^{\cal D})$ is dense.
\end{definition}

Notice that, if $K$ is subcanonical then any weakly dense morphism of sites
$({\cal C}, J)\to ({\cal D}, K)$ is dense, but the converse does not hold in general (see Example \ref{exampleweaklydensenotdense} below).

The following proposition gives an explicit characterization of weakly dense morphisms of sites (for conciseness we give the constructively stronger formulation relying on the axiom of choice, but the reader who prefers the constructive phrasing in terms of more general families satisfying the conditions in Proposition \ref{propexplicit}(i)-(ii) can obtain it by replacing all the families of arrows indexed by a covering sieve occurring in the statement with these more general families).
 
\begin{proposition}\label{propweaklydensemorphism}
	Let $F:({\cal C}, J)\to ({\cal D}, K)$ be a morphism of sites. Then $F$ is a weakly dense morphism of sites if and only if it satisfies the following conditions:
	\begin{enumerate}[(i)]
	\item $P$ is a $J$-covering family in $\cal C$ if and only if $F(P)$ is a $K$-covering family in $\cal D$;
	
    \item for any object $d$ of $\cal D$ there exist a family $\{S_{i} \mid i\in I\}$ of $K$-covering sieves on objects of the form $F(c_{i})$ (where $c_{i}$ is an object of $\cal C$) and for each $f\in S_{i}$ an arrow $g_{f}:\textup{dom}(f)\to d$ such that $g_{f\circ z}\equiv_{K} g_{f}\circ z$ whenever $z$ is composable with $f$, such that the family of arrows $g_{f}$ (for $f\in S_{i}$ for some $i$) is $K$-covering;  
   
    \item for any objects $x, y$ of ${\cal C}$ and any family of arrows $g_{h}:\textup{dom}(h)\to F(y)$ indexed by the arrows of a $K$-covering sieve $U$ on $F(x)$ such that $g_{h\circ k}\equiv_{K} g_{h}\circ k$ for every arrow $k$ composable with $h$, there exist a $J$-covering family of arrows $\{f_{i}:x_{i}\to x \mid  i\in I\}$ and arrows $k_{i}:x_{i}\to y$ (for each $i\in I$) such that for every arrows $w$ and $z$ such that $F(f_{i})\circ w=h\circ z$, we have $g_{h}\circ z\equiv_{K} F(k_{i})\circ w$ (for every $h\in U$ and $i\in I$).
	\end{enumerate}
\end{proposition}

\begin{proofs}
Conditions (i) for $F$ is clearly equivalent to condition (i) in the definition of a dense morphism of sites for the functor $l'\circ F$.

Let us now reformulate condition (ii) in the definition of a dense morphism of sites for the functor $l'\circ F$, namely the condition that for every object $d$ of $\cal D$ there exists an epimorphic family of arrows  $\{l'(F(c_{i}))\to l'(d) \mid i\in I\}$ whose domains are objects in the image of the functor $l'\circ F$. Our thesis thus immediately follows from the description of such arrows provided by Proposition \ref{propexplicit}. 

Let us show that we can reformulate condition (iii) in the definition of a dense morphism of sites for the functor $l'\circ F$ as property (iii) of the proposition. For any arrow $\xi:l'(F(x))\to l'(F(y))$ in $\Sh({\cal D}, K)$, by Proposition \ref{propexplicit} there exists a $K$-covering sieve $\{h:\textup{dom}(h)\to F(x) \mid h\in U \}$ and for each $h\in U$ an arrow $g_{h}:\textup{dom}(h)\to F(y)$ such that $\xi \circ l'(h)=l'(g_{h})$ for each $h$ and $g_{h\circ k}\equiv_{K} g_{h}\circ k$ for any arrow $k$ composable with $h$. The existence of a $J$-covering family of arrows $f_{i}:x_{i}\to x$ and a family of arrows $k_{i}:x_{i}\to y$ such that $\xi\circ l'(F(f_{i}))=l'(k_{i})$ for each $i$ can be reformulated in terms of the arrows $h$ and $g_{h}$, as follows. Let us consider, for each $i\in I$ and $h\in U$, the following pullback square:
$$
 \xymatrix{
 	P_{i, h} \ar[rd]^{e_{i, h}} \ar[r]^{\pi^{1}_{i,h}\,\,\,\,\,\,\,} \ar[d]_{\pi^{2}_{i,h}} & l'(F(x_{i})) \ar[d]^{l'(F(f_{i}))} \\
 	l'(\textup{dom}(h)) \ar[r]^{l'(h)} & l'(F(x))
 } 
$$      	

It is easy to see (cf. Proposition \ref{proplifttosite}) that the collection ${\cal F}_{i,h}$ of arrows $\chi:l'(c)\to P_{i, h}$ such that both $\pi^{1}_{i,h}\circ \chi$ and $\pi^{2}_{i,h}\circ \chi$ are images $l'(w):l'(c)\to l'(F(x_{i}))$ and $l'(z):l'(c)\to l'(\textup{dom}(h))$ under $l'$ of arrows $w$ and $z$ in $\cal D$ (where $c$ varies among the objects of $\cal C$) is epimorphic. At the cost of composing with another covering family on $l'(c)$, we can suppose that $F(f_{i})\circ w=h\circ z$ without loss of generality. Therefore, since for any $i\in I$ the collection of arrows $\{\pi^{1}_{i,h}\circ \chi \mid h\in U, \chi\in {\cal F}_{i,h}\}$ is epimorphic, $\xi\circ l'(F(f_{i}))=l'(k_{i})$ if and only if for any $h\in U$ and arrows $w$ and $z$ such that $F(f_{i})\circ w=h\circ z$, $\xi\circ e_{i,h}\circ \chi=l'(k_{i})\circ \pi^{1}_{i,h}\circ \chi$. But $\pi^{1}_{i,h}\circ \chi=l'(w)$ and $\xi\circ e_{i,h}\circ \chi=\xi \circ l'(h)\circ \pi^{2}_{i,h}\circ \chi=l'(g_{h})\circ l'(z)$. So the condition $\xi\circ l'(F(f_{i}))=l'(k_{i})$ is equivalent to the condition $g_{h}\circ z\equiv_{K} k_{i}\circ w$ for every arrow $h\in U$ and arrows $w$ and $z$ such that $F(f_{i})\circ w=h\circ z$.     

\end{proofs}

\begin{remark}
		Condition (iv) in the definition of dense morphism of sites (cf. Remark \ref{remmorphismsofsitesfaithfulness}) for the functor $l'\circ F$ admits the following reformulation in terms of $F$:
		
		\begin{enumerate}[(i)]
			\item[(iv)] For any arrows $f_{1}, f_{2}:c_{1}\to c_{2}$ in $\cal C$ such that $F(f_1)\equiv_{K} F(f_2)$ there exists a $J$-covering family of arrows $k_{i}:c'_{i}\to c_{1}$ such that $f_1 \circ k_{i}=f_{2}\circ k_{i}$ for all $i$.
		\end{enumerate}
	
		If a functor $F$ satisfies the conditions of Proposition \ref{propweaklydensemorphism} plus condition (iv) then $l'\circ F$ is a dense morphism of sites and hence $F:({\cal C}, J)\to ({\cal D}, K)$ is a weakly dense morphism of sites (cf. Remark \ref{remmorphismsofsites}(a)).  	
\end{remark}	

Summarizing, we have the following result:
\begin{theorem}\label{thmweakdensityequivalence}
	Let $F:({\cal C}, J)\to ({\cal D}, K)$ be a morphism of sites. Then the following conditions are equivalent:
	\begin{enumerate}[(i)]
		\item The geometric morphism $\Sh(F):\Sh({\cal D}, K)\to \Sh({\cal C}, J)$ is an equivalence.
		
		\item $l'\circ F:({\cal C}, J)\to (a_{K}({\cal D}), C_{K}^{\cal D})$ is a dense morphism of sites.
		
		\item $F$ is a weakly dense morphism of sites $({\cal C}, J)\to ({\cal D}, K)$ (i.e. it satisfies the conditions of Proposition \ref{propweaklydensemorphism}).
	\end{enumerate}
\end{theorem}

\begin{proofs}
Since $\Sh(l'\circ F)$ is the composite of $\Sh(F)$ with the canonical equivalence $\Sh(a_{K}({\cal D}), C_{K}^{\cal D})\simeq \Sh({\cal D}, K)$ and the site $(a_{K}({\cal D}), C_{K}^{\cal D})$ is subcanonical, the direction (ii) $\Rightarrow$ (i) follows from Theorem 11.8 \cite{Shulman}. The equivalence between (ii) and (iii) follows by definition of a  weakly dense morphism of sites (in light of Proposition \ref{propweaklydensemorphism}). Lastly, the implication (i) $\Rightarrow$ (ii) follows from Proposition \ref{propdensitysubcanonical}. 	
\end{proofs}

\begin{remark}
	Theorem \ref{thmweakdensityequivalence} constitutes a vast generalization of Grothendieck's Comparison Lemma (Theorem 4.1 from section III of \cite{grothendieck}).
\end{remark}

\begin{example}\label{exampleweaklydensenotdense}
Let us discuss an example of a weakly dense morphism of sites which is not dense.
Let $\bf{2}$ be the preorder category with two distinct objects $0$ and $1$ and just one arrow $0\to 1$ apart from the identities. Notice that $\bf{2}$ is cartesian, since it is a meet-semilattice, and the functor $F:\bf{2}\to \bf{2}$ sending $0$ to $1$, $1$ to $1$ and the arrow $0\to 1$ to the identity arrow on $1$ is cartesian as it is a meet-semilattice homomorphism. Let us equip $\bf{2}$ with the atomic topology $J_{\textup{at}}$, whose covering sieves are the maximal ones plus the sieve on $1$ consisting of the arrow $0\to 1$. Since it is cartesian and cover-preserving, the functor $F$ is a morphism of sites $({\mathbf{2}}, J_{\textup{at}})\to ({\mathbf{2}}, J_{\textup{at}})$. Now, by the Comparison Lemma, the topos $\Sh({\mathbf{2}}, J_{\textup{at}})$ is equivalent to the topos of sheaves on the one-object full subcategory $\{0\}$ of $\bf{2}$ with respect to the induced topology on it, namely the maximal one; so we have $\Sh({\mathbf{2}}, J_{\textup{at}})\simeq \Set$. Therefore the geometric morphism $\Sh(F):\Sh({\mathbf{2}}, J_{\textup{at}}) \to \Sh({\mathbf{2}}, J_{\textup{at}})$ is necessarily (isomorphic) to the identity morphism (since the only geometric morphism $\Set \to \Set$ is the identity one); in particular, it is an equivalence. Therefore, by Theorem \ref{thmweakdensityequivalence}, $F$ is a weakly dense morphism of sites. However, $F$ is \emph{not} dense; for instance, it does not satisfy condition (ii) in the definition of a dense morphism of sites. Indeed, there is no $J_{\textup{at}}$-covering family of arrows to the object $0$ whose domain is in the image of the functor $F$. On the other hand, it satisfies condition (ii) in the characterization of weakly dense morphisms of sites provided by Proposition \ref{propweaklydensemorphism}. Indeed, if $d=1$ then the arrow $k=1_{1}$ satisfies the condition (by taking $h=1_{1}$, $x=0$ or $x=1$, $S=\{0\to 1\}$ and $g_{0\to 1}=0\to 1$), while if $d=0$ then the arrow $k=1_{0}$ satisfies the condition (by taking $h=0\to 1$, $x=0$ or $x=1$, $S=\{0\to 1\}$ and $g_{0\to 1}=1_{0}$).
\end{example}

\begin{remark}
	In section \ref{sec:equivalencetoposes2} we shall obtain an alternative criterion, given by Corollary \ref{corquivalencesites}, for a morphism of sites to induce an equivalence of toposes. 
\end{remark}

\subsection{Two criteria for equivalence}

Given a flat functor $F:{\cal C}\to {\cal E}$ defined on an essentially small category $\cal C$ with values in a Grothendieck topos $\cal E$, if $F$ is $J$-continuous for a Grothendieck topology $J$ on $\cal C$ then we know by Diaconescu's equivalence that it induces a geometric morphism $f:{\cal E}\to \Sh({\cal C}, J)$. The following result provides a necessary and sufficient condition, phrased entirely in terms of $F$, for $f$ to be an equivalence:

\begin{corollary}\label{flatequivalence}
	Let $({\cal C}, J)$ be a small-generated site, $\cal E$ a Grothendieck topos and $F:{\cal C}\to {\cal E}$ a $J$-continuous flat functor. Then the geometric morphism $f:{\cal E}\to \Sh({\cal C}, J)$ induced by $F$ is an equivalence if and only if $F$ satisfies the following conditions:
		\begin{enumerate}[(i)]	
			\item If the image under $F$ of a sieve $S$ in $\cal C$ is epimorphic in $\cal E$ then $S$ is $J$-covering;
			\item the family of objects of the form $F(c)$ for $c\in {\cal C}$ is separating for $\cal E$;
			\item for every $x, y\in {\cal C}$ and any arrow $g:F(x)\to F(y)$ in $\cal E$, there exist a $J$-covering family of arrows $f_{i}:x_{i}\to x$ and a family of arrows $g_{i}:x_{i}\to y$ such that $g\circ F(f_{i})=F(g_{i})$ for all $i$.			
		\end{enumerate}
\end{corollary}

\begin{proofs}
	A $J$-continuous flat  functor $F:{\cal C}\to {\cal E}$ is clearly the same thing as a morphism of sites $F:({\cal C}, J)\to ({\cal E}, J_{\cal E}^{\textup{can}})$, and the geometric morphism $f:{\cal E}\to \Sh({\cal C}, J)$ corresponding to $F$ via Diaconescu's equivalence is precisely $\Sh(F)$. Our thesis thus follows from Theorem \ref{thmweakdensityequivalence}, since the conditions of the corollary are precisely those for $F:({\cal C}, J)\to ({\cal E}, J_{\cal E}^{\textup{can}})$ to be dense (notice that the site $({\cal E}, J_{\cal E}^{\textup{can}})$ is subcanonical).
\end{proofs}

\begin{remark}
	In section \ref{sec:equivalencetoposes2}, we shall encounter an alternative criterion, provided by Corollary \ref{corflatequivalence2}, for a continuous flat functor to induce an equivalence of toposes. 
\end{remark}

We can deduce from Theorem \ref{thmweakdensityequivalence} a criterion for two (essentially small) sites to give rise to equivalent toposes of sheaves on them.

\begin{theorem}
	Let $({\cal C}, J)$ and $({\cal D}, K)$ be two small-generated sites. Then the following conditions are equivalent:
	
	\begin{enumerate}[(i)]
		\item The toposes $\Sh({\cal C}, J)$ and $\Sh({\cal D}, K)$ are equivalent.
		
		\item There exist a category (resp. an essentially small category, if $\cal C$ and $\cal D$ are essentially small) $\cal A$, a Grothendieck topology $Z$ on $\cal A$ (which can be supposed subcanonical) and two functors $H:{\cal C}\to {\cal A}$ and $K:{\cal D}\to {\cal A}$ satisfying the following conditions:
		\begin{enumerate}[(i)]	
		\item $P$ is a $J$-covering family in $\cal C$ if and only if $H(P)$ is a $Z$-covering family in $\cal A$;
		\item $Q$ is a $K$-covering family in $\cal D$ if and only if $K(Q)$ is a $Z$-covering family in $\cal A$;
		\item for any object $a$ of $\cal A$ there exists a $Z$-covering sieve whose arrows factor both through an arrow whose domain is in the image of $H$ and through an arrow whose domain is in the image of $K$;
		\item for every $x, y\in {\cal C}$ (resp. $x',y'\in {\cal D}$) and any arrow $g:H(x)\to H(y)$ (resp. $g':K(x')\to K(y')$) in $\cal A$, there exist a $J$-covering family of arrows $f_{i}:x_{i}\to x$ (resp. a $K$-covering family of arrows $f_{j}':x_{j}'\to x'$) and a family of arrows $g_{i}:x_{i}\to y$ (resp. a family of arrows $g'_{j}:x_{j}'\to y'$) such that $g\circ H(f_{i})=H(g_{i})$ for all $i$ (resp. $g'\circ K(f_{j}')=K(g'_{j})$ for all $j$);
		\item for any arrows $h, k:x\to y$ (resp. $h',k':x'\to y'$) in $\cal C$ (resp. in $\cal D$) such that $H(h)=H(k)$ (resp. $K(h')=K'(k')$) there exists a $J$-covering (resp. $K$-covering) family of arrows $f_{i}:x_{i}\to x$ (resp. $f'_{j}:x'_{j}\to x'$) such that $h\circ f_{i}=k\circ f_{i}$ for all $i$ (resp. $h'\circ f'_{j}=k'\circ f'_{j}$ for all $j$). 
		\end{enumerate}		 

	\end{enumerate} 

\end{theorem}

\begin{proofs}
The given conditions are precisely those for the functor $H$ and $K$ to respectively define dense morphisms of sites $({\cal C}, J)\to ({\cal A}, Z)$ and $({\cal D}, K)\to ({\cal A}, Z)$. Such morphisms induce by Theorem 11.8 \cite{Shulman} (cf. also Theorem \ref{thmweakdensityequivalence}) equivalences of toposes $\Sh({\cal C}, J)\simeq \Sh({\cal A}, Z)$ and $\Sh({\cal D}, K)\to \Sh({\cal A}, Z)$, whence an equivalence $\Sh({\cal C}, J)\simeq \Sh({\cal D}, K)$. Conversely, if $\Sh({\cal C}, J)\simeq \Sh({\cal D}, K)$ then, by taking $\cal A$ to be the full subcategory of this topos on the objects that are either coming from the site $({\cal C}, J)$ or from the site $({\cal D}, K)$ with the Grothendieck topology $Z$ induced on it by the canonical topology on the topos, we obtain by the Comparison Lemma equivalences $\Sh({\cal C}, J)\simeq \Sh({\cal A}, Z)$ and $\Sh({\cal D}, K)\to \Sh({\cal A}, Z)$, and hence by Theorem \ref{thmweakdensityequivalence} the canonical functors ${\cal C}\to {\cal A}$ and ${\cal D}\to {\cal A}$ are respectively dense morphisms of sites $({\cal C}, J)\to ({\cal A}, Z)$ and $({\cal D}, K)\to ({\cal A}, Z)$. 
\end{proofs}

\subsection{Local faithfulness, local fullness and local surjectivity}\label{sec:localconditions}

In this section we shall introduce some notions which are naturally related to the denseness conditions considered above.

\begin{definition}\label{deflocalproperties}
	Let $F:{\cal C}\to {\cal D}$ be a functor and $J$ (resp. $K$) a Grothendieck topology on $\cal C$ (resp. on $\cal D$). Then $F$ is said to be 
	\begin{enumerate}[(a)]
		\item \emph{$(J, K)$-faithful} (resp. $J$-faithful) if whenever $F(h)\equiv_{K} F(k)$ (resp. $F(h)=F(k)$), $h\equiv_{J} k$;
		
		\item \emph{$(J, K)$-full} (resp. \emph{$J$-full}) if for every $x, y\in {\cal C}$ and any arrow $g:F(x)\to F(y)$ in $\cal D$, there exist a $J$-covering family of arrows $f_{i}:x_{i}\to x$ and arrows $g_{i}:x_{i}\to y$ (for each $i\in I$) such that $g\circ F(f_{i})\equiv_{K} F(g_{i})$ (resp. $g\circ F(f_{i})= F(g_{i})$) for all $i$;
		
		\item \emph{$K$-dense} if for every $d\in {\cal D}$, there exists a $K$-covering family of arrows whose domains are in the image of $F$. 
	\end{enumerate} 
\end{definition}

\begin{remarks}\label{remcoveringliftingstrictequality}
	\begin{enumerate}[(a)]
		\item If $K$ is the canonical topology on $\cal D$ then the relation $\equiv_{K}$ reduces to equality and $(J, K)$-faithfulness (resp. $(J, K)$-fullness) reduces to $J$-faithfulness (resp. $J$-fullness). 
		
		\item If $F$ satisfies the covering-lifting property (that is, the property that for any $c\in {\cal C}$ and any $K$-covering sieve $S$ on $F(c)$ there is a $J$-covering sieve $R$ on $c$ such that $F(R)\subseteq S$) then the local equality $\equiv_{K}$ in conditions (a) and (b) of Definition \ref{deflocalproperties} can be replaced by strict equality, whence $(J, K)$-faithfulness (resp. $(J, K)$-fullness) coincides with $J$-faithfulness (resp. $J$-fullness) for $F$.  
		
		\item Any dense morphism of sites $({\cal C}, J)\to ({\cal D}, K)$ is $J$-faithful, $J$-full and $K$-dense.
		
		\item Inverse images functors of geometric inclusions of toposes (in particular, associated sheaf functors) satisfy a form of local fullness (cf. Corollary \ref{corinclusionequivalence} below).
	\end{enumerate}	
\end{remarks}

\begin{proposition}\label{proplocalfullnessandfaithfulness}
	Let $F:{\cal C}\to {\cal D}$ be a functor and $J$ a Grothendieck topology on $\cal C$. Then $F$ is $J$-full and $J$-faithful if and only if for any $c\in {\cal C}$, the canonical arrow $l(c)\to a_{J}(\textup{Hom}_{\cal D}(F(-), F(c)))$ is an isomorphism in $\Sh({\cal C}, J)$.
\end{proposition}

\begin{proofs}
	The canonical arrow $l(c)\to l(\textup{Hom}_{\cal D}(F(-), F(c)))$ is the image under the associated sheaf functor $a_{J}:[{\cal C}^{\textup{op}}, \Set]\to \Sh({\cal C}, J)$ of the canonical arrow $y_{\cal C}(c)\to (\textup{Hom}_{\cal D}(F(-), F(c))$. Our thesis thus immediately follows from Lemma \ref{lemmalift}. 	
\end{proofs}

\begin{proposition}\label{propcovlifting}
	Let $F:{\cal C}\to {\cal D}$ be a $J$-full, $K$-dense, cover-reflecting and cover-preserving functor (for instance, a dense morphism of sites $F:({\cal C}, J)\to ({\cal D}, K)$). Then $F:({\cal C}, J)\to ({\cal D}, K)$ has the covering-lifting property. 
\end{proposition}

\begin{proofs}
	Let $c$ be an object of $\cal C$ and $S$ be a $K$-covering sieve on $F(c)$. Then, for each $g\in S$, by $K$-denseness there exists a $K$-covering sieve $U_{g}$ on $\textup{dom}(g)$ generated by a family of arrows $\{\xi^{g}_{i} \mid i\in I \}$ whose domains are of the form $F(c^{g}_{i})$, where $c^{g}_{i}$ are objects of $\cal C$. Now, by $(J, K)$-fullness, there is a $J$-covering family of arrows $h^{(g, i)}_{k}$ to $c^{g}_{i}$ and for each $h^{(g, i)}_{k}$ an arrow $\chi_{h^{(g, i)}_{k}}$ from $\textup{dom}(h^{(g, i)}_{k})$ to $c$ such that $g\circ \xi^{g}_{i} \circ F(h^{(g, i)}_{k})= F(\chi_{h^{(g, i)}_{k}})$. So the arrows $\chi_{h^{(g, i)}_{k}}$ generate a sieve whose image under $F$ factors through $g$ and hence belongs to $S$. This sieve is $J$-covering since $F$ is both cover-preserving and cover-reflecting; indeed, the sieve generated by the arrows $g\circ \xi^{g}_{i} \circ F(h^{(g, i)}_{k})$ is $K$-covering since $S$ is $K$-covering, the family of arrows $\{\xi^{g}_{i} \mid i\in I \}$ is $K$-covering and the family of arrows $h^{(g, i)}_{k}$ is $J$-covering (whence its image under $F$ is $K$-covering). So $F$ satisfies the covering-lifting property.     	
\end{proofs}

\begin{corollary}
	Let $F:{\cal C}\to {\cal D}$ be a $J$-full, $K$-dense, cover-reflecting and cover-preserving functor (for instance, a dense morphism of sites $F:({\cal C}, J)\to ({\cal D}, K)$). Then, for any sieve $R$ on an object $F(c)$ in $\cal D$, $R\in K(F(c))$ if and only if $\{f:\dom(f)\to c \mid F(f)\in R\}\in J(c)$.  
\end{corollary}

\begin{proofs}
If 	$\{f:\dom(f)\to c\ mid F(f)\in R\}\in J(c)$ then, since $F$ is cover-preserving, the family $F(\{f:\dom(f)\to c\ mid F(f)\in R\})\subseteq R$ is $K$-covering and hence $R$ is also $K$-covering. Conversely, if $R\in K(F(c))$ then, since $F$ satisfies the covering-lifting property by Proposition \ref{propcovlifting}, there is a $J$-covering sieve $S$ on $c$ such that $F(S)\subseteq R$. Therefore $S\subseteq \{f:\dom(f)\to c\ mid F(f)\in R\}$, and since $S$ is $J$-covering, it follows that $\{f:\dom(f)\to c\ mid F(f)\in R\}$ is also $J$-covering.  
\end{proofs}

\begin{remark}
	For any small-generated site $({\cal C}, J)$, the corollary notably applies to the canonical functor $l:{\cal C}\to \Sh({\cal C}, J)$, which is a dense morphism of sites $({\cal C}, J)\to (\Sh({\cal C}, J), J^{\textup{can}}_{\Sh({\cal C}, J)})$, yielding the following (well-known) result: for any sieve $R$ on an object $l(c)$ in $\Sh({\cal C}, J)$, $R$ is covering for the canonical topology on $\Sh({\cal C}, J)$ (that is, it yields an epimorphic family) if and only if $\{f:d\to c \mid l(f)\in R\}\in J(c)$.		
\end{remark}

Given a morphism of sites $F:({\cal C}, J)\to ({\cal D}, K)$, we have an induced functor $a_{F}:a_{J}({\cal C})\to a_{K}({\cal D})$ given by the restriction of the inverse image of the induced geometric morphism $\Sh(F)$, which is a morphism of sites if $a_{J}({\cal C})$ and $a_{K}({\cal D})$ are endowed with the Grothendieck topologies $C^{\cal C}_{J}$ and $C^{\cal D}_{K}$ respectively induced by the canonical ones on $\Sh({\cal C}, J)$ and $\Sh({\cal D}, K)$.

\begin{proposition}\label{proplocalproperties}
	Let $F:({\cal C}, J)\to ({\cal D}, K)$ be a morphism of sites. Then
	
	\begin{enumerate}[(i)]
		\item $F$ is $(J, K)$-faithful if and only if $a_{F}$ is faithful;
		
		\item If $F$ is $(J, K)$-faithful and satisfies the covering-lifting property then $F$ is $(J, K)$-full if and only if $a_{F}$ is full;
		
		\item If $F$ satisfies the covering-lifting property then $F$ is $K$-dense if and only if $a_{F}$ is $C^{\cal D}_{K}$-dense. 
	\end{enumerate}
\end{proposition}

\begin{proofs}
	(i) If $a_{F}$ is faithful then for any parallel arrows $h, k$ in $\cal C$, the condition $l(h)=l(k)$ (that is, $h\equiv_{J} k$) is equivalent to the condition $l'(F(h))=a_{F}(l(h))=a_{F}(l(k))=l'(F(k))$ (that is, $F(h)\equiv_{K} F(k)$). 
	
	Conversely, let us now suppose that $F$ is $(J, K)$-faithful, and prove that $a_{F}$ is faithful. If $u, v:l(c)\to l(c')$ are arrows in $\Sh({\cal C}, J)$ then there are $J$-covering sieves $S$ and $T$ on $c$ such that for each $f\in S$ there is an arrow $\xi_{f}:\textup{dom}(f)\to c'$ such that $u\circ l(f)=l(\xi_{f})$ for each $f$, and for each $g\in T$ there is an arrow $\chi_{g}:\textup{dom}(g)\to c'$ such that $v\circ l(g)=l(\chi_{g})$. Notice that the sieve $S\cap T$ is $J$-covering as it is the intersection of two $J$-covering sieves. Suppose that $a_{F}(u)=a_{F}(v)$. Then $a_{F}(u)\circ a_{F}(l(h))=a_{F}(v)\circ a_{F}(l(h))$ for each $h\in S\cap T$. Therefore $l'(F(\xi_{h}))= a_{F}(l(\xi_{h}))=a_{F}(l(\chi_{h}))=l'(F(\chi_{h}))$ for each $h\in S\cap T$, that is $F(\xi_{h}) \equiv_{K} F(\chi_{h})$. Since $F$ is $(J, K)$-faithful, it follows that $\xi_{h}\equiv_{J} \chi_{h}$ for each $h\in S\cap T$. But this means that $u\circ l(h)=l(\xi_{h})=l(\chi_{h})=v\circ l(h)$ for each $h\in S\cap T$. Since $S\cap T$ is $J$-covering, this implies that $u=v$, as required.        
	
	(ii) Suppose that $a_{F}$ is full. Then for every $x, y\in {\cal C}$ and any arrow $g:F(x)\to F(y)$ in $\cal D$, there is an arrow $\xi:l(x)\to l(y)$ such that $a_{F}(\xi)=l'(g)$. By Proposition \ref{propexplicit} there are a $J$-covering family of arrows $f_{i}:x_{i}\to x$ and arrows $g_{i}:x_{i}\to y$ (for each $i\in I$) such that $\xi \circ l(f_{i}) = l(g_{i})$. Applying $a_{F}$, we thus obtain that $l'(g)\circ l'(F(f_{i}))=l'(F(g_{i}))$, equivalently $g\circ F(f_{i})\equiv_{K} F(g_{i})$. So $F$ is $(J, K)$-full. In fact, this argument shows that the fullness condition applied to arrows of the form $l'(g)$ for some arrow $g$ in $\cal D$ is precisely equivalent to the  $(J, K)$-fullness condition. Conversely, let us suppose that $F$ is $(J, K)$-full and satisfies the covering-lifting property. Given an arrow $\xi:l'(F(c))\to l'(F(c'))$, by Proposition \ref{propexplicit} there exist a $K$-covering sieve $S$ on $F(c)$ and for each arrow $g$ in $S$ an arrow $t_{g}:\textup{dom}(g)\to F(c')$ such that $\xi\circ l'(g)=l'(t_{g})$. Since $F$ satisfies the covering-lifting property, there is a $J$-covering sieve $R$ on $c$ such that $F(R)\subseteq S$. For any $r\in R$, consider the arrow $t_{F(r)}:F(\textup{dom}(r))\to F(c')$. Since $F$ is  $(J, K)$-full, there exist a $J$-covering family of arrows $\{ f^{r}_{j}:y_{j}\to \textup{dom}(r) \mid j\in H_{r}\}$ and, for each $j\in H_{r}$, an arrow $k^{r}_{j}:y_{j}\to c'$ such that $t_{F(r)}\circ F(f^{r}_{j})\equiv_{K} F(k^{r}_{j})$. We want to show that the arrows $k^{r}_{j}:y_{j}\to c'$ indexed by the $J$-covering family of arrows $r\circ f^{r}_{j}$ to $c$ define an arrow $\chi:l(c)\to l(c')$ such that $a_{F}(\chi)=\xi$. For this, by Proposition \ref{propexplicit} we have to verify that if $m$ and $n$ are arrows such that $r\circ f^{r}_{j}\circ m=r'\circ f^{r'}_{j'}\circ n$, then $k^{r}_{j}\circ m \equiv_{J} k^{r'}_{j'}\circ n$. Now, since $F$ is $(J, K)$-faithful, the latter condition is equivalent to $F(k^{r}_{j}\circ m) \equiv_{K} F(k^{r'}_{j'}\circ n)$. But $l'(F(k^{r}_{j}\circ m))=l'(F(k^{r}_{j}))\circ l'(F(m))=l'(t_{F(r)})\circ l'(F(f^{r}_{j}))\circ l'(F(m))=\xi \circ l'(F(r))\circ l'(F(f^{r}_{j}))\circ l'(F(m))=\xi \circ l'(F(r\circ f^{r}_{j}\circ m))=\xi \circ l'(F(r'\circ f^{r'}_{j'}\circ n))=\xi \circ l'(F(r'))\circ l'(F(f^{r'}_{j'}))\circ l'(F(n))=l'(t_{F(r')})\circ l'(F(f^{r'}_{j'}))\circ l'(F(n))=l'(F(k^{r'}_{j'}))\circ l'(F(n))=l'(F(k^{r'}_{j'}\circ n))$, as required. The fact that $a_{F}(\chi)=\xi$ follows immediately from the definition of $\chi$, by using the fact that $F$ is cover-preserving.   
	
	(iii) If $F$ is $K$-dense then, clearly, $a_{F}$ is $C^{\cal D}_{K}$-dense. Suppose instead that $a_{F}$ is $C^{\cal D}_{K}$-dense. Given an object $d$ of $\cal D$, we want to show that there is a $K$-covering family of arrows whose domains are in the image of $F$. Consider the object $l'(d)$. Then there is an epimorphic family $\cal F$ in $\Sh({\cal D}, K)$ of arrows $\xi$ to $l'(d)$ whose domains are of the form $l'(F(c))$ for some $c\in {\cal C}$. For any such $\xi:l'(F(c_{\xi}))\to l'(d)$, by Proposition \ref{propexplicit} there exists a $K$-covering sieve $S_{\xi}$ on $F(c_{\xi})$ and for each arrow $g$ in $S_{\xi}$ an arrow $t_{g}:\textup{dom}(g)\to d$ such that $\xi\circ l'(g)=l'(t_{g})$. Since $F$ satisfies the covering lifting property, there is a $J$-covering sieve $R_{\xi}$ on $c$ such that $F(R_{\xi})\subseteq S_{\xi}$. Now, the collection of arrows  $\{t_{F(r)}:F(\textup{dom}(r))\to d \mid r\in R_{\xi}, \xi\in {\cal F}\}$ is $K$-covering. Indeed, by using the fact that $F$ is cover-preserving, one immediately sees that the image of this family under $l'$ is epimorphic in $\Sh({\cal D}, K)$. This completes our proof. 
\end{proofs}

\begin{remark}
Proposition \ref{proplocalproperties} implies that every morphism of sites $F:({\cal C}, J)\to ({\cal D}, K)$ whose underlying functor has the covering-lifting property, is $J$-faithful, $J$-full and $K$-dense induces an equivalence of toposes. In fact, as shown by the following result, these morphisms are precisely the dense morphisms of sites $({\cal C}, J)\to ({\cal D}, K)$.	
\end{remark}

\begin{corollary}\label{corequivalencecovliftproperty}
	Let $F$ be a morphism of sites $({\cal C}, J)\to ({\cal D}, K)$. Then the following conditions are equivalent:
	\begin{enumerate}[(i)]
		\item $F$ is a weakly dense morphism of sites which has the covering-lifting property;
		
		\item $F$ is dense (that is, $K$-dense, $J$-full and cover-reflecting).
	\end{enumerate}
	
\end{corollary}

\begin{proofs}
	(i) $\imp$ (ii)  By definition, any weakly dense morphism of sites $F:({\cal C}, J)\to ({\cal D}, K)$ is cover-reflecting and satisfies the property that $a_{F}$ is $C^{\cal D}_{K}$-dense; so by Proposition \ref{proplocalproperties}(iii) $F$ is $K$-dense. The fact that $F$ is $J$-full follows from the proof of Proposition \ref{proplocalproperties}(ii) in light of Theorem \ref{thmweakdensityequivalence} and Remark \ref{remcoveringliftingstrictequality}(b).
	
	(ii) $\imp$ (i) This follows from Propositions \ref{propcovlifting} and \ref{proplocalproperties}.   
\end{proofs}

\begin{remarks}\label{remweaklydense}	
	\begin{enumerate}[(a)]
		\item The morphism of sites $F:({\mathbf{2}}, J_{\textup{at}})\to ({\mathbf{2}}, J_{\textup{at}})$ of Example \ref{exampleweaklydensenotdense} is weakly dense but does not satisfy the covering-lifting property (the sieve on $1=F(0)=F(1)$ generated by the arrow $0\to 1$ cannot be lifted neither to a $J_{\textup{at}}$-covering sieve on $0$ nor on $1$ since the image under $F$ of any sieve is the maximal one).
		
		\item Every weakly dense morphism of sites $F:({\cal C}, J)\to ({\cal D}, T)$ is $(J, T)$-faithful if $T$ is a subcanononical topology on $\cal D$ (in which cases the local equality $\equiv_{T}$ reduces to strict equality). On the other hand, every weakly dense morphism of sites $F:({\cal C}, J)\to ({\cal D}, K)$ satisfying the covering lifting property is $(J, K)$-faithful (see Proposition \ref{propweaklydensemorphism} and Remark \ref{remmorphismsofsitesfaithfulness}).  
	\end{enumerate}
	
\end{remarks}

\section{Site characterizations of some properties of geometric morphisms}\label{sec:propgeometricmorphisms}

Corollary \ref{flatequivalence} gives a characterization of the $J$-continuous flat  functors ${\cal C}\to {\cal E}$ whose corresponding geometric morphism $f:{\cal E}\to \Sh({\cal C}, J)$ is an equivalence. In this section we shall characterize the $J$-continuous flat  functors whose associated geometric morphism $f$ is a geometric inclusion (resp. a surjection, localic, hyperconnected).  

\subsection{Surjections and inclusions}\label{secsurjectioninclusion}

\begin{proposition}\label{propcoverreflecting}
	Let $({\cal C}, J)$ be a small-generated site, $\cal E$ a Grothendieck topos and $F:{\cal C}\to {\cal E}$ a $J$-continuous flat  functor. Then the geometric morphism $f:{\cal E}\to \Sh({\cal C}, J)$ induced by $F$ is a surjection (that is, $f^{\ast}$ is faithful) if and only if $F$ is \emph{cover-reflecting} (in the sense that the $J$-covering families in $\cal C$ are precisely those which are sent by $F$ to epimorphic families in $\cal E$).
\end{proposition}

\begin{proofs}
Let us preliminarily notice that $f^{\ast}$ is faithful if and only if for every monomorphism $m$, if $f^{\ast}(m)$ is an isomorphism then $m$ is an isomorphism. Indeed, one direction follows by considering the equalizer of the given pair of arrows, while the other follows from the balancedness of $\Sh({\cal C}, J)$ (which allows to check that $m$ is an isomorphism by showing that it is both a monomorphism and an epimorphism).
	
Suppose that $f^{\ast}$ is faithful. Given a sieve $S$ in $\cal C$ on an object $c$, consider the corresponding monomorphism $m_{S}:S\mono \Hom_{\cal C}(-, c)$. The sieve $S$ is $J$-covering if and only if $a_{J}(m_{S})$ is an isomorphism. But this is equivalent, since $f^{\ast}$ is faithful, to $f^{\ast}(a_{J}(m_{S}))$ being an isomorphism, that is, to the family $\{F(f) \mid f\in S \}$ being epimorphic.

Conversely, let us suppose that $F$ is cover-reflecting. Given a monomorphism $m:P\to Q$ in $\Sh({\cal C}, J)$, we have to prove that if $f^{\ast}(m)$ is an isomorphism then $m$ is an isomorphism. By considering the pullbacks $S_{\alpha}\mono l(c)$ of $m$ along all the arrows $\alpha:l(c)\to Q$ where $c$ is an object of $\cal C$, we obtain that, since $F$ is cover-reflecting, the monomorphisms $S_{\alpha}\mono l(c)$ are all isomorphisms, whence $m$ is an isomorphism as well (the family of arrows $\alpha$ being epimorphic and $\cal E$ being balanced).
\end{proofs}

Let us apply this result in the context of interpretations between geometric theories. Recall from section 2.1.3 of \cite{OCBook} that an interpretation of a geometric theory $\mathbb T$ into a geometric theory $\mathbb S$ is a geometric functor $I:{\cal C}_{\mathbb T}\to {\cal C}_{\mathbb S}$ between their geometric syntactic categories. Notice that $I$ acts on contexts $\vec{x}$ for geometric formulae over the signature of $\mathbb T$ sending a context $\vec{x}$ to the context $I(\vec{x})$ appearing in the formula $I(\{\vec{x}. \top\})=\{I(\vec{x}). \xi\}$; since $I$ preserves monomorphisms, we have $I(\{\vec{x}. \phi\})=\{I(\vec{x}).\chi\}$ for some geometric formula $\chi$ in the context $I(\vec{x})$; we shall write $I(\phi)$ for $\chi$ (when the context for $\phi$ can be unambigously inferred). Therefore, with any sequent $\sigma\equiv (\phi \vdash_{\cal x} \psi)$ over the signature of $\mathbb T$ can be associated a sequent $I(\sigma)\equiv (I(\phi) \vdash_{I(\vec{x})} I(\psi))$ over the signature of $\mathbb S$ which is provable in $\mathbb S$ if the former sequent is provable in $\mathbb T$.  

\begin{corollary}\label{corinterpretationsfaithful}
\begin{enumerate}[(i)]
	\item An intepretation $I:{\cal C}_{\mathbb T}\to {\cal C}_{\mathbb S}$ between geometric theories $\mathbb T$ and $\mathbb S$ induces a geometric surjection $\Set[{\mathbb S}]\to \Set[{\mathbb T}]$ between their classifying toposes if and only if, with the above notation, for every sequent $\sigma$ over the signature of $\mathbb T$, if the sequent $I(\sigma)$ is provable in $\mathbb S$ then $\sigma$ is provable in $\mathbb T$. 
	
	\item In particular, if $I$ is the canonical interpretation given by an expansion ${\mathbb T}'$ of ${\mathbb T}$ (in the sense of section 7.1 of \cite{OCBook}) then $I$ induces a geometric surjection if and only if ${\mathbb T}'$ is a \emph{conservative expansion} of $\mathbb T$, that is, every geometric sequent over the signature of $\mathbb T$ which is provable in ${\mathbb T}'$ is provable in ${\mathbb T}$.
\end{enumerate}
\end{corollary}

\begin{proofs}
It suffices to apply Proposition \ref{propcoverreflecting}, observing that the $J_{\mathbb T}$-continuous flat functor given by the composite $y_{\mathbb S}\circ I$ is cover-reflecting if and only if for every subobject $[\psi]:\{\vec{x}. \psi\} \mono \{\vec{x}. \phi \}$ in ${\mathbb T}$, if $I([\psi])$ is the identity subobject (equivalently, the sequent $(I(\phi) \vdash_{I(\vec{x})} I(\psi))$ is provable in $\mathbb S$) then $[\psi]$ is the identity subobject (equivalently, the sequent $(\phi \vdash_{\vec{x}} \psi)$ is provable in $\mathbb T$). 	
\end{proofs}

Let us recall from Theorem A4.2.10 \cite{El} that every geometric morphism can be factored, uniquely up to commuting equivalence, as a surjection followed by an inclusion. Let us describe this factorization in terms of sites.

In the proof of the following theorem, we exploit the possibility of regarding the inverse image functor $f^{\ast}:{\cal E}\to {\cal F}$ of a geometric morphism $f:{\cal F}\to {\cal E}$ as a morphism of small-generated sites $({\cal E}, J_{\cal E}^{\textup{can}})\to ({\cal F}, J_{\cal F}^{\textup{can}})$.

\begin{theorem}\label{Morphismsmaalgenerated}
	Let $F:({\cal C}, J)\to ({\cal D}, K)$ be a morphism of small-generated sites. Then:
	
	\begin{enumerate}[(i)]
		\item The geometric morphism $\Sh(F):\Sh({\cal D}, K) \to \Sh({\cal C}, J)$ induced by $F$ is a surjection if and only if $F$ is \emph{cover-reflecting} (that is, if the image of a family of arrows with a fixed codomain is $K$-covering then the family is $J$-covering).
		
		\item The surjection-inclusion factorization of the geometric morphism $\Sh(F):\Sh({\cal D}, K) \to \Sh({\cal C}, J)$ induced by $F$ can be identified with the factorization $\Sh(i_{J_{F}})\circ \Sh(F_{r})$, 
		where $J_{F}$ is the Grothendieck topology on $\cal C$ whose covering sieves are exactly those whose image under $F$ are $K$-covering families, $i_{J_{F}}:({\cal C}, J)\to ({\cal C}, J_{F})$ is the morphism of sites given by the canonical inclusion functor and $F_{r}: ({\cal C}, J_{F}) \to ({\cal D}, K)$ is the morphism of sites given by $F$. 
		
		\item The geometric morphism $\Sh(F):\Sh({\cal D}, K) \to \Sh({\cal C}, J)$ induced by $F$ is an inclusion if and only if $F_{r}: ({\cal C}, J_{F}) \to ({\cal D}, K)$ is a weakly dense morphism of sites (equivalently, if $K$ is subcanonical, a dense morphism of sites); in particular, if $K$ is subcanonical then $\Sh(F)$ is an inclusion if and only if the following conditions are satisfied:
		
		\begin{enumerate}[(i)]	
			\item for any object $d$ of $\cal D$ there exists a $K$-covering family of arrows $d_{i}\to d$ whose domains $d_{i}$ are in the image of $F$;
			\item for every $x, y\in {\cal C}$ and any arrow $g:F(x)\to F(y)$ in $\cal D$, there exist a $J_{F}$-covering family of arrows $f_{i}:x_{i}\to x$ and a family of arrows $g_{i}:x_{i}\to y$ such that $g\circ F(f_{i})=F(g_{i})$ for all $i$.
		\end{enumerate}
	\end{enumerate}
\end{theorem}

\begin{proofs}
(i) This follows as an immediate consequence of Proposition \ref{propcoverreflecting} by observing that a family of arrows in $\cal D$ with common codomain is sent by the canonical functor ${\cal D}\to \Sh({\cal D}, K)$ to an epimorphic family if and only if it is $K$-covering. 	

(ii) The morphism $\Sh({i_{J_{F}}})$ is the canonical inclusion of $\Sh({\cal C}, J_{F})$ into $\Sh({\cal C}, J)$. On the other hand, by definition of $J_{F}$, the morphism of sites $F_{r}$ is cover-reflecting and hence, by point (i), induces a surjective geometric morphism. So $\Sh(F_{r})$ and $\Sh({i_{J_{F}}})$ give a factorization of $\Sh(F)$ as a surjection followed by an inclusion.

(iii) A geometric morphism is an inclusion if and only if the surjection part of its surjection-inclusion factorization is an equivalence; our thesis thus follows from Theorem \ref{thmweakdensityequivalence}.
\end{proofs}

Theorem \ref{Morphismsmaalgenerated} can be notably applied to a $J$-continuous flat  functor $F:{\cal C} \to {\cal E}$, regarded as a morphism of sites $({\cal C}, J)\to ({\cal E}, J_{\cal E}^{\textup{can}})$, giving the following result: 

\begin{corollary}\label{corinclusionflatfunctor}
	A $J$-continuous flat  functor $F:{\cal C} \to {\cal E}$ induces a geometric inclusion if and only if it satisfies the following conditions:  
	\begin{enumerate}[(i)]	
		\item every object of $\cal E$ can be covered by objects which are in the image of $F$;
		\item for every $x, y\in {\cal C}$ and any arrow $g:F(x)\to F(y)$ in $\cal E$, there exist a family of arrows $f_{i}:x_{i}\to x$ which is sent by $F$ to an epimorphic family and a family of arrows $g_{i}:x_{i}\to y$ such that $g\circ F(f_{i})=F(g_{i})$ for all $i$.
	\end{enumerate}  
\end{corollary}

\subsection{Induced topologies}\label{sec:inducedtopology}

The following result, which is a corollary of Theorem \ref{Morphismsmaalgenerated}(ii), shows that the natural context for defining `induced Grothendieck topologies' is that provided by flat functors (equivalently, geometric morphisms) or more generally morphisms of sites.

\begin{proposition}\label{reminducedtopology}
	Let $f:{\cal E}\to [{\cal C}^{\textup{op}}, \Set]$ be a geometric morphism (equivalently, a flat functor $F:{\cal C}\to {\cal E}$). Then there exists a Grothendieck topology $J_{f}$ (resp. $J_{F}$) on $\cal C$, called the Grothendieck topology \emph{induced} by $f$ (resp. $F$) whose covering sieves are precisely the sieves which are sent by $f^{\ast}$ (resp. by $F$) to epimorphic families in $\cal E$. This applies in particular to a morphism of small-generated sites $G:({\cal C}, J)\to ({\cal D}, K)$, yielding a Grothendieck topology $J_{G}$ on $\cal C$ whose covering sieves are exactly those whose image under $G$ are $K$-covering families (cf. Theorem \ref{Morphismsmaalgenerated}(ii)).	  
\end{proposition}\qed

\begin{remarks}\label{rem1}
	\begin{enumerate}[(a)]
		\item If $F:({\cal C}, J)\to ({\cal D}, K)$ is a cover-reflecting morphism of sites (for instance, a weakly dense morphism of sites) then $J$ is a Grothendieck topology (since it coincides with $J_{F}$).
		
		\item The classical definition of Grothendieck topology induced on a full $K$-dense subcategory $\cal C$ of a category $\cal D$ is subsumed by this general definition of induced topology; indeed, it is readily seen that the embedding of $\cal C$ into $\cal D$ defines a morphism of sites $({\cal C}, T)\to ({\cal D}, K)$, where $T$ is the trivial Grothendieck topology on $\cal C$.
		
		\item The topology $J_{f}$ (resp. $J_{F}$) on $\cal C$ induced by  by $f$ (resp. $F$) in the sense of Proposition \ref{reminducedtopology} coincides with the Grothendieck topology induced by $f$ (resp. $F$), regarded as a functor to the canonical site $({\cal E}, J^{\textup{can}}_{\cal E})$, in the sense of section III.3 of \cite{grothendieck}. In fact, more generally, for any morphism of small-generated sites $G:({\cal C}, J)\to ({\cal D}, K)$, the Grothendieck topology induced by $u$ in the sense of  \cite{grothendieck} (that is, the largest topology which makes $u$, regarded as a functor to the site $({\cal D}, K)$, continuous) coincides with the topology $J_{G}$ in our Proposition \ref{reminducedtopology} (cf. Corollary III-3.3 \cite{grothendieck}).      
		
		\item We can characterize the toposes which can be represented as subtoposes of a given presheaf topos as follows: a topos $\cal E$ is equivalent to a subtopos of a presheaf topos $[{\cal C}^{\textup{op}}, \Set]$ if and only if there is a flat (equivalently, filtering) functor $F:{\cal C}\to {\cal E}$ satisfying the conditions of Corollary \ref{corinclusionflatfunctor}, in which case $\cal E$ is equivalent to $\Sh({\cal C}, J_{F})$, where $J_{F}$ is the Grothendieck topology induced by $F$ (in the sense of Proposition \ref{reminducedtopology}). Notice that this generalizes Giraud's result that every Grothendieck topos with a separating set of objects  $\cal C$  can be represented as $\Sh({\cal C}, J^{\textup{can}}_{\cal E}|_{\cal C})$, where $\cal C$ is regarded as a full subcategory of $\cal E$; indeed, for any separating set $\cal C$ of objects of a Grothendieck topos $\cal E$, the inclusion functor of $\cal C$ into $\cal E$ is flat.   
	\end{enumerate}
	
\end{remarks}

The following result shows in particular that if $F:({\cal C}, J)\to ({\cal D}, K)$ is a dense morphism of sites then the Grothendieck topology $K$ can be recovered from $J$. Notice that, by Corollary \ref{corequivalencecovliftproperty}, any flat functor $F:{\cal C}\to {\cal E}$ on a small-generated category $\cal C$ with values in a Grothendieck topos $\cal E$ which induces an equivalence ${\cal E} \simeq \Sh({\cal C}, J_{F})$ (equivalently, satisfies the conditions of Corollary \ref{flatequivalence}, cf. also Corollary \ref{corflatequivalenceinducedtopology} below) is a dense morphism of sites $({\cal C}, J_{F})\to ({\cal E}, J_{\cal E}^{\textup{can}})$. 

\begin{proposition}\label{proprecoverytopology}
	Let $F:({\cal C}, J)\to ({\cal D}, K)$ be a $K$-dense functor with the covering-lifting property (for instance, a dense morphism of sites $({\cal C}, J)\to ({\cal D}, K)$ -- cf. Corollary \ref{corequivalencecovliftproperty}). Then the Grothendieck topology $K$ can be recovered from $J$, as follows: for any sieve $T$ on an object $d$ of $\cal D$, $T\in K(d)$ if and only if for any object $c\in {\cal C}$ and any arrow $\xi:F(c)\to d$ in $\cal D$, there exists a $J$-covering sieve $R$ on $c$ such that $F(R)\subseteq \xi^{\ast}(T)$.  
\end{proposition}

\begin{proofs}
 By the transitivity and stability axioms for Grothendieck topologies, $T$ is $K$-covering if and only if for any object $c\in {\cal C}$ and any arrow $\xi:F(c)\to d$ in $\cal D$, $\xi^{\ast}(T)\in K(\textup{dom}(\xi))$. But, $F$ having the covering-lifting property, $\xi^{\ast}(T)\in K(\textup{dom}(\xi))$ if and only if there exists a $J$-covering sieve $R$ on $c$ such that $F(R)\subseteq \xi^{\ast}(T)$. 
\end{proofs}

\begin{remark}
	Every dense morphism of sites $F:({\cal C}, J)\to ({\cal D}, K)$ is cover-reflecting, so $J$ can be recovered from $K$ as the collection of sieves which are sent by $F$ to $K$-covering families.
\end{remark}

The notion of induced topology can be very profitably applied for establishing equivalences of toposes. Indeed, if we have a flat functor $F:{\cal C}\to {\cal E}$, we dispose of easily applicable criteria for $F$ to induce an equivalence of toposes 
\[
{\cal E}\simeq \Sh({\cal C}, J_{F}),
\]
as shown by the following result (which is an immediate consequence of Corollary \ref{flatequivalence}):

\begin{corollary}\label{corflatequivalenceinducedtopology}
	Let ${\cal C}$ be an essentially small category, $\cal E$ a Grothendieck topos and $F:{\cal C}\to {\cal E}$ a flat functor. Then $F$ induces an equivalence 
	\[
	{\cal E}\simeq \Sh({\cal C}, J_{F}),
	\]
	if and only if $F$ is $J_{F}$-full and the objects of the form $F(c)$ for $c\in {\cal C}$ form a separating set for the topos $\cal E$
\end{corollary}\qed

The following proposition describes a general setup for building equivalences of toposes. 

\begin{proposition}\label{propgeneratingequivalence}
	Let $\cal E$ be a Grothendieck topos, represented as the category $\Sh({\cal D}, K)$ of sheaves on a small-generated site $({\cal D}, K)$. Let $F:{\cal C}\to {\cal D}$ be a functor and $J_{F}$ the collection of sieves in $\cal C$ whose image under $F$ is $K$-covering. If $F$ is a weakly dense morphism of sites $({\cal C}, J_{F})\to ({\cal D}, K)$ then $J_{F}$ is a Grothendieck topology on $\cal C$, and we have an equivalence of toposes 
	\[
	\Sh({\cal D}, K) \simeq \Sh({\cal C}, J_{F}).
	\]
	
	The following conditions are equivalent to $F$ being a dense (or weakly dense, if $K$ is subcanonical) morphism of sites $({\cal C}, J_{F})\to ({\cal D}, K)$:
	\begin{enumerate}[(i)]
		\item[(i)'] For any object $d$ of $\cal D$ there exists a $K$-covering family of arrows $d_{i}\to d$ whose domains $d_{i}$ are in the image of $F$.
		
		\item[(ii)'] For every $c_{1}, c_{2}\in {\cal C}$ and any arrow $g:F(c_{1})\to F(c_{2})$ in $\cal D$, there exist a family of arrows $f_{i}:c'_{i}\to c_{1}$ which is sent by $F$ to a $K$-covering family and a family of arrows $k_{i}:c'_{i}\to c_{2}$ such that $g\circ F(f_{i})=F(k_{i})$ for all $i$. 
		
		\item[(iii)'] For any arrows $f_{1}, f_{2}:c_{1}\to c_{2}$ in $\cal C$ such that $F(f_1)=F(f_2)$ there exists a family of arrows $k_{i}:c'_{i}\to c_{1}$ which is sent by $F$ to a $K$-covering family such that $f_1 \circ k_{i}=f_{2}\circ k_{i}$ for all $i$.	
	\end{enumerate}
\end{proposition} 

\begin{proofs}
	The first part of the proposition follows from Theorem \ref{thmweakdensityequivalence} and Remark \ref{rem1}(a).
	
	The second part of the proposition follows from the fact that every dense morphism is weakly dense and that the two notions are equivalent when $K$ is subcanonical.
\end{proofs}

\subsection{Coinduced topologies}

The following result, established in the proof of Lemma C2.3.19(ii) \cite{El}, represents a kind of converse to Proposition \ref{proprecoverytopology}. It involves a notion of \emph{image} of a Grothendieck topology under a functor.

\begin{proposition}[cf. Lemma C2.3.19(ii) \cite{El}]\label{propimagetopology}
	Let $F:{\cal C}\to {\cal D}$ be a functor and $J$ a Grothendieck topology on $\cal C$. Then there is a Grothendieck topology $J^{F}$ on $\cal D$, called the \emph{image of $J$ along $F$} (or the Grothendieck topology \emph{coinduced} by $J$ along $F$), such that for any sieve $T$ on an object $d$ of $\cal D$, $T\in J^{F}(d)$ if and only if for any object $c\in {\cal C}$ and any arrow $\xi:F(c)\to d$ in $\cal D$, there exists a $J$-covering sieve $R$ on $c$ such that $F(R)\subseteq \xi^{\ast}(T)$. Then the functor $F:({\cal C}, J)\to ({\cal D}, J^{F})$ is $J^{F}$-dense and has the covering-lifting property.
	
	The Grothendieck topology $J^{F}$ enjoys the following universal property: the subtopos $\Sh({\cal D}, J^{F})\hookrightarrow [{\cal D}^{\textup{op}}, \Set]$ is the image (in the sense of surjection-inclusion factorization) of the composite of the induced geometric morphism $[F^{\textup{op}}, \Set]:[{\cal C}^{\textup{op}}, \Set]\to [{\cal D}^{\textup{op}}, \Set]$ with the subtopos inclusion $\Sh({\cal C}, J)\hookrightarrow [{\cal C}^{\textup{op}}, \Set]$.
\end{proposition}
	
\begin{proofs}
First, let us show that $J^{F}$ is indeed a Grothendieck topology on $\cal C$. The maximality and pullback-stability axiom trivially hold, so it remains to prove the transitivity one. It is clear that any sieve containing a $J^{F}$-covering one is also $J^{F}$-covering. Let us denote, given a sieve $S$ on an object $d$ and, for each arrow $f\in S$, a sieve $S_{f}$ on $\textup{dom}(f)$, by $S\ast \{S_{f} \mid f\in S\}$ the sieve on $d$ consisting of the arrows of the form $f\circ h$ where $f\in S$ and $h\in S_{f}$. We have to prove that, if $S\in J^{F}(d)$ and, for each $f\in S$, $S_{f}\in J^{F}(\textup{dom}(f))$, then $S\ast \{S_{f} \mid f\in S\}\in J^{F}(d)$. Given an arrow $\xi:F(c)\to d$, $F(R) \subseteq \xi^{\ast}(S)$ for some $R\in J(c)$. For any $r\in R$, $\xi \circ F(r)\in S$, so, since $S_{\xi \circ F(r)}\in J^{F}(F(\textup{dom}(r)))$, there is a sieve $R_{r}\in J(\textup{dom}(r))$ such that $F(R_{r})\subseteq S_{\xi \circ F(r)}$. So the $J$-covering sieve $R\ast \{R_{r}\mid r\in R \}$ satisfies the condition $F(R\ast \{R_{r}\mid r\in R \})\subseteq \xi^{\ast}(S\ast \{S_{f} \mid f\in S\})$.

The fact that $F:({\cal C}, J)\to ({\cal D}, J^{F})$ is $J^{F}$-dense and has the covering-lifting property is obvious.

The last part of the proposition follows at once from Proposition \ref{propsurjectioncomorphismofsites} below. Indeed, the construction of the geometric morphism induced by a comorphism of sites is functorial, the geometric morphism $[G^{\textup{op}}, \Set]:[{\cal A}^{\textup{op}}, \Set] \to [{\cal B}^{\textup{op}}, \Set]$ induced by a functor $G:{\cal A}\to {\cal B}$ is precisely $L_{G}$, regarding $G$ as a comorphism of sites $({\cal A}, T_{\cal A})\to ({\cal B}, T_{\cal B})$ where $T_{\cal A}$ and $T_{\cal B}$ are the trivial Grothendieck topologies respectively on $\cal A$ and $\cal B$, and any canonical geometric inclusion $\Sh({\cal A}, Z) \hookrightarrow [{\cal A}^{\textup{op}}, \Set]$ is induced by the comorphism of sites $1_{\cal A}:({\cal A}, Z) \to ({\cal A}, T)$ given by the identity functor $1_{\cal A}$ on $\cal A$. 
\end{proofs}

\begin{remark}\label{remcovlift}
If $F:({\cal C}, J)\to ({\cal D}, K)$ satisfies the covering-lifting property then $K\subseteq J^{F}$. This follows immediately from the fact that $K$ satisfies the pullback-stability property. In fact, in Lemma C2.3.19(ii) \cite{El}, $J^{F}$ is characterized as the largest Grothendieck topology on $\cal C$ for which $F$ satisfies the covering-lifting property.    	
\end{remark}

\begin{lemma}\label{lemmafullfaithfulcoverreflecting}
Let $F:{\cal C}\to {\cal D}$ be a $J$-full and $J$-faithful functor, for a Grothendieck topology $J$ on $\cal C$. Then:
\begin{enumerate}[(i)]
	\item for any sieve $R$ on an object $c$ of $\cal C$ and any arrow $t$ with codomain $c$, if $F(t)$ belongs to the sieve generated by $F(R)$ then there exists a $J$-covering sieve $U$ on $\textup{dom}(t)$ such that $t\circ u\in R$ for each $u\in U$;  
	
	\item the functor $F:({\cal C}, J)\to ({\cal D}, J^{F})$ is cover-reflecting.
\end{enumerate}
\end{lemma}

\begin{proofs}
(i) Suppose that $F(t)=F(r)\circ \xi$ for some $r\in R$, where $t:a\to c$, $r:b\to c$ and $\xi:F(a)\to F(b)$. Since $F$ is $J$-full, there exist a $J$-covering sieve $S$ on $a$ and for each $s\in S$ an arrow $t_{s}:\textup{dom}(s)\to b$ such that $\xi \circ F(s)=F(t_{s})$. Since $F(t)=F(r)\circ \xi$, we have $F(t\circ s)=F(t)\circ F(s)=F(r)\circ \xi \circ F(s)= F(r)\circ F(t_{s})=F(r\circ t_{s})$ for each $s\in S$. The fact that $F$ is $J$-faithful then implies that $t\circ s\equiv_{J} r\circ t_{s}$, that is, there is a $J$-covering sieve $V_{s}$ on $\textup{dom}(s)$ such that for each $v\in V_{s}$, $t\circ s\circ v=r\circ t_{s}\circ v$. Therefore the $J$-covering sieve $U:=S\ast \{V_{s} \mid s\in S\}$ on $a$ satisfies the property that for each $u\in U$, $t\circ u\in R$. 

(ii) Let $R$ be a sieve on an object $c$ of $\cal C$. If $F(R)$ generates a $J^{F}$-covering sieve then there is a $J$-covering sieve $T$ on $c$ such that the sieve generated by $F(T)$ is contained in the sieve generated by $F(R)$. So, for each $t\in T$, by point (i) of the Lemma, there is a $J$-covering sieve $U_{t}$ on $\textup{dom}(t)$ such that $t\circ u\in R$ for each $u\in U_{t}$. So the sieve $R$ is $J$-covering, since it contains the $J$-covering sieve $T\ast \{U_{t} \mid t\in T\}$.  	
\end{proofs}

\begin{theorem}\label{corliftinggeneration}
Let $F:({\cal C}, J)\to ({\cal D}, K)$ be a $K$-dense, cover-preserving and cover-reflecting functor with the covering-lifting property (for instance, a dense morphism of sites $({\cal C}, J)\to ({\cal D}, K)$ - cf. Corollary \ref{corequivalencecovliftproperty}). Suppose that $K$ is generated by a collection $K'$ of families of arrows which is stable under pullback and such that every $K'$-covering family on an object in the image of $F$ contains the image $F(R)$ under $F$ of some $J$-covering family $R$. If $F$ is $J'$-full and $J'$-faithful, where $J'$ is the Grothendieck topology generated by the sieves $R$ such that $F(R)\in K'$, then $J=J'$.	
\end{theorem}

\begin{proofs}
Let $J'$ be the Grothendieck topology on $\cal C$ generated by the sieves $R$ such that $F(R)$ generates a $K$-covering sieve. We have to prove that $J'=J$. Consider the topology $J'^{F}$. We clearly have $K'\subseteq J'^{F}$, whence $K\subseteq J'^{F}$. On the other hand, $K=J^{F}$ by Proposition \ref{proprecoverytopology}.  Now, the fact that $F$ is cover-reflecting implies that $J'\subseteq J$; so $J'^{F}\subseteq J^{F}=K$, and hence $J'^{F}=K$. Now, since both $F:({\cal C}, J)\to ({\cal D}, K)$ and $F:({\cal C}, J')\to ({\cal D}, J'^{F})$ are cover-reflecting (the former by our hypothesis and the latter by Lemma \ref{lemmafullfaithfulcoverreflecting}), $J'^{F}=K$ implies that $J=J'$.
\end{proofs}

\begin{remark}
	If $F$ is a weakly dense morphism of sites $({\cal C}, J)\to ({\cal D}, K)$, the condition that $F$ be $J'$-full and $J'$-faithful is also necessary for $J$ to be equal to $J'$ (see Proposition \ref{propmorphismcomorphismequivalence}).  
\end{remark}

Theorem \ref{corliftinggeneration} has the following immediate corollary, which allows to obtain, from a set of generators of a Grothendieck topology, a set of generators for a topology induced by it on a full dense subcategory.

\begin{corollary}\label{corgensep}
	Let $({\cal D}, K)$ be a small-generated site, $\cal C$ a full $K$-dense subcategory of $\cal D$ and $J$ a Grothendieck topology on $\cal C$. Suppose that $K$ is generated by a collection $K'$ of families of arrows which is stable under pullback and such that every $K'$-covering family on an object in $\cal C$ is refined by a $J$-covering family. Then $J$ is generated by the families in $K'$ whose arrows lie all in $\cal C$. 	
\end{corollary}

\subsection{Intrinsic characterization of geometric inclusions}\label{sec:intrinsiccharinclusions}

The first part of the following result characterizes the property of a geometric morphism being an inclusion entirely in terms of its inverse image. It can be useful, for instance, in cases where one disposes of a description of the inverse image of the morphism which is simpler or more tractable than that of its direct image.	

\begin{corollary}\label{corinclusionequivalence}
	Let $f:{\cal F}\to {\cal E}$ be a geometric morphism. Then 
	
	\begin{enumerate}[(i)]
	\item $f$ is an inclusion if and only if $f^{\ast}$ satisfies the following conditions:

		\begin{enumerate}[(i)]	
			\item \textup{$f^{\ast}$ is locally surjective}, that is every object of $\cal F$ can be covered by objects in the image of $f^{\ast}$;
			
			\item \textup{$f^{\ast}$ is locally full}, that is for every $x, y\in {\cal E}$ and any arrow $g:f^{\ast}(x)\to f^{\ast}(y)$ in $\cal F$, there exists a family of arrows $s_{i}:x_{i}\to x$ in $\cal E$ which is sent by $f^{\ast}$ to an epimorphic family and a family of arrows $g_{i}:x_{i}\to y$ such that $g\circ f^{\ast}(s_{i})=f^{\ast}(g_{i})$ for all $i$.
		\end{enumerate}
	
	\item $f$ is an equivalence if and only if $f^{\ast}$ is faithful, locally full and locally surjective.

	\end{enumerate}
\end{corollary}

\begin{proofs}
	(i) The above conditions amount precisely to the property of the morphism of sites $({\cal E}, {J_{\cal E}^{\textup{can}}}_{f^{\ast}})\to ({\cal F}, J_{\cal F}^{\textup{can}})$ to be dense, so the thesis follows from Theorem \ref{Morphismsmaalgenerated}. 
	
	(ii) This follows from (i), given that a geometric morphism is an equivalence if and only if it is a surjection and an inclusion.
\end{proofs}

\begin{remark}\label{remsurjectionsinclusions}
	Since a Grothendieck topos has all coproducts, all the covering families arising in the formulations of the properties in part (i) of Corollary \ref{corinclusionequivalence} can be supposed, without loss of generality, to consist of a single element. More precisely, the condition for $f$ to be an inclusion is equivalent to the conjunction of the following ones:  
		\begin{enumerate}[(i)]	
			\item every object of $\cal F$ is a quotient of an object in the image of $f^{\ast}$;
			
			\item for every $x, y\in {\cal E}$ and any arrow $g:f^{\ast}(x)\to f^{\ast}(y)$ in $\cal F$, there exist an arrow $s:x'\to x$ in $\cal E$ which is sent by $f^{\ast}$ to an epimorphism and an arrow $g':x'\to y$ such that $g\circ f^{\ast}(s)=f^{\ast}(g')$.
		\end{enumerate}
		
		So, $f$ is an equivalence if and only if $f^{\ast}$ is faithful and satisfies the above conditions.

\end{remark}

\begin{example}\label{exgeometricinclusion}
	Let $f:{\cal F}\to {\cal E}$ be a geometric morphism between Grothendieck toposes and $A$ an object of $\cal E$. Then we have a geometric morphism $f_{A}:{\cal F}\slash f^{\ast}(A)\to {\cal E}\slash A$ such that the diagram
	$$
	\xymatrix{
		{\cal F}\slash f^{\ast}(A)  \ar[r] \ar[d]_{f_{A}} & {\cal F} \ar[d]^{f} \\
		{\cal E}\slash A \ar[r] & {\cal E}
	} 
	$$  
	where the horizontal morphisms are the canonical ones, commutes. The inverse image of $f_{A}$ is the functor ${\cal E}\slash A \to {\cal F}\slash f^{\ast}(A)$ sending an object $u:B\to A$ of ${\cal E}\slash A$ to $f^{\ast}(u):f^{\ast}(B)\to f^{\ast}(A)$ (and acting accordingly on arrows). Let us show, by applying Corollary \ref{corinclusionequivalence}(i), that if $f$ is an inclusion then $f_{A}$ is also an inclusion. First, let us show that ${f_{A}}^{\ast}$ is locally surjective. Given an object $u:C\to f^{\ast}(A)$ of ${\cal F}\slash f^{\ast}(A)$, since $f^{\ast}$ is locally surjective there is an object $B$ of $\cal E$ and an epimorphism $q:f^{\ast}(B)\to C$. Consider the composite arrow $u\circ q:f^{\ast}(B)\to f^{\ast}(A)$. Since $f^{\ast}$ satisfies condition (ii) of Corollary \ref{corinclusionequivalence}(i), there is an arrow $s:B'\to B$ in $\cal E$ and an arrow $g:B'\to A$ in $\cal E$ such that $f^{\ast}(s)$ is an epimorphism and $u\circ q\circ f^{\ast}(s)=f^{\ast}(g)$ (cf. Remark \ref{remsurjectionsinclusions}). Therefore the arrow $q\circ f^{\ast}(s)$ is an epimorphism from ${f_{A}}^{\ast}(g)$ to $u$ in ${\cal F}\slash f^{\ast}(A)$. The fact that ${f_{A}}^{\ast}$ inherits local fullness from $f^{\ast}$ is obvious. 
\end{example}

\subsection{Hyperconnected and localic morphisms}

Recall that a geometric morphism $f:{\cal F}\to {\cal E}$ is said to be
\emph{hyperconnected} if $f^{\ast}$ is full and faithful and its image is closed under subobjects in $\cal F$, and is said to be \emph{localic} if every object of $\cal F$ is a subquotient (that is, a quotient of a subobject) of an object of the form $f^{\ast}(A)$ for $A\in {\cal E}$. By Theorem A4.6.5 \cite{El}, every geometric morphism can be factored, uniquely up to commuting equivalence, as the composite of a hyperconnected morphism followed by a localic one.  

In this section we shall characterize the property of a geometric morphism to be hyperconnected (or localic) in terms of sites, and also provide a natural site-level description of the hyperconnected-localic factorization.

\subsubsection{Hyperconnected morphisms}

The following lemma will be instrumental in the sequel.

\begin{lemma}\label{lempb}
	Let $q:X\to X'$ be the coequalizer of a pair of arrows $r_{1}, r_{2}:R\to X$  in a Grothendieck topos $\cal E$ and $f:Y\to X'$ a monomorphism. Then, the arrows $t_{1}, t_{2}:R\times_{X'} Y \to Z$ which make the diagram
	
		$$
	\xymatrix{
		R\times_{X'} Y  \ar@<-.5ex>[d]_{t_{1}} \ar@<.5ex>[d]^{t_{2}} \ar[rr] & & R \ar@<-.5ex>[d]_{r_{1}} \ar@<.5ex>[d]^{r_{2}}   \\
		Z \ar[d]^{k} \ar[rr]^{h} & & X \ar[d]^{q} \\
		Y \ar[rr]^{f} & & X'		
	}
	$$ 
	commutative, where the bottom square is a pullback, satisfy the property that $\langle t_{1}, t_{2}\rangle :R\times_{X'} Y \to Z\times Z$ is (isomorphic to) the pullback of $\langle r_{1}, r_{2} \rangle  $ along $h\times h$.
	
	Moreover, $k$ is the coequalizer of the arrows $t_{1}$ and $t_{2}$.  
\end{lemma}

\begin{proofs}
	By using the internal language, we shall argue as if $\cal E$ were the topos $\Set$ of sets.
	
	Let $U\to Z\times Z$ be the pullback of $\langle r_{1}, r_{2} \rangle  $ along $h\times h$. We have 
	\[
	U=\{ ((z_{1}, z_{2}), \xi) \mid  (z_{1}, z_{2})\in Z, \, \xi\in R, \, h(z_{1})=r_{1}(\xi), \, h(z_{2})=r_{2}(\xi) \}.
	\]
	
	We want to prove that $U\to Z\times Z$ is isomorphic to  $\langle t_{1}, t_{2}\rangle\,:R\times_{X'} Y \to Z\times Z$ as arrows over $Z\times Z$. 
	
	We have $$R\times_{X'} Y=\{(\xi, y) \mid \xi \in R, \, y\in Y, \, q(r_{1}(\xi))=q(r_{2}(\xi))=f(y) \}.$$
	Now, the arrow $t_{1}:R\times_{X'} Y \to Z$ (resp. $t_{2}:R\times_{X'} Y \to Z$) sends an element $(\xi, y)\in R\times_{X'} Y$ to the unique element $z_{1}\in Z$ such that $k(z_{1})=y$ and $h(z_{1})=r_{1}(\xi)$ (resp. to the unique element $z_{2}\in Z$ such that $k(z_{2})=y$ and $h(z_{2})=r_{2}(\xi)$). Let us define arrows $R\times_{X'} Y \to U$ and $U \to R\times_{X'} Y$ over $Z\times Z$ that are inverse to each other. In one direction, we send an element $(\xi, y)\in R\times_{X'} Y$ to the element $(t_{1}(\xi, y), t_{2}(\xi, y), \xi)$. This is well-defined by definition of the arrows $t_{1}$ and $t_{2}$. In the converse direction, we observe that, given an element $((z_{1}, z_{2}), \xi) $ of $U$, we have $k(z_{1})=k(z_{2})$. Indeed, since $f$ is monic, this condition is equivalent to $f(k(z_{1}))=f(k(z_{2}))$; but $f(k(z_{1}))=q(h(z_{1}))=q(r_{1}(\xi))=q(r_{2}(\xi))=q(h(z_{2}))=f(k(z_{2}))$. We can therefore define an arrow $U\to R\times_{X'} Y$ which sends an element $((z_{1}, z_{2}), \xi) $ of $U$ to the element $(\xi, k(z_{1}))$. This is well-defined since $q(r_{1}(\xi))=q(h(z_{1}))=f(k(z_{1}))$ and $q(r_{2}(\xi))=q(h(z_{2}))=f(k(z_{2}))=f(k(z_{1}))$ (as $k(z_{1})=k(z_{2})$). It is now straightforward to check that the two arrows just defined are inverse to each other and compatible with the structure arrows to $Z\times Z$. This completes the proof of the first part of the lemma. The fact that $k$ is the coequalizer of the arrows $t_{1}$ and $t_{2}$ follows from the fact that in a Grothendieck topos colimits are stable under pullback.
\end{proofs}

\begin{remark}
	The proof of the lemma shows that the general version of it (for $f$ not necessarily monic) reads as follows: $\langle t_{1}, t_{2}\rangle$ is isomorphic over $Z\times Z$ to the composite with the canonical embedding $z_{Y}:Z\times_{Y} Z \hookrightarrow Z\times Z$ of the pullback of $\langle r_{1}, r_{2} \rangle  $ along $(h\times h)\circ z_{Y}$. 
\end{remark}

\begin{proposition}\label{prophypersubobjects}
	Let ${\cal C}$ be a separating set for a Grothendieck topos $\cal E$ and $f:{\cal F}\to {\cal E}$ a geometric morphism. Then the image of $f^{\ast}$ is closed under subobjects if and only if for any $c\in {\cal C}$, every subobject of $f^{\ast}(c)$ is in the image of $f^{\ast}$.
\end{proposition}

\begin{proofs}
	Let $m:B\mono f^{\ast}(A)$ be a monomorphism. Let us represent $A$ as a quotient $q:\coprod_{c\in {\cal C}_{A}}c\to A$ of a coproduct of objects coming from $\cal C$. Then, since the image of $f^{\ast}|_{\cal C}$ is closed under subobjects and coproducts are stable under pullbacks, the pullback of $m$ along $f^{\ast}(q)$ is of the form $f^{\ast}(z)$ for some subobject $z:Z\mono \coprod_{c\in {\cal C}_{A}}c$ in $\cal E$. Now, if $q$ is the coequalizer of its kernel pair $\langle r_{1},r_{2} \rangle  $ in $\cal E$ then $f^{\ast}(q)$ is the coequalizer of the pair $f^{\ast}(r_{1}), f^{\ast}(r_{2})$. Let $\langle t_{1}, t_{2}\rangle:T\to Z$ be the pullback in $\cal E$ of $\langle r_{1}, r_{2}\rangle  $ along $z$. Then, since $f^{\ast}$ preserves pullbacks, $\langle f^{\ast}(t_{1}), f^{\ast}(t_{2})\rangle:f^{\ast}(T)\to f^{\ast}(Z)$ is the pullback of $\langle f^{\ast}(r_{1}), f^{\ast}(r_{2})\rangle  $ along $f^{\ast}(z)$.  
	
	 By Lemma \ref{lempb}, it follows that $m$ is isomorphic to the canonical arrow $U \to f^{\ast}(A)$ from the codomain $U$ of the coequalizer $f^{\ast}(Z)\to U$ of the arrows $f^{\ast}(t_{1}), f^{\ast}(t_{2})$. But, since $f^{\ast}$ preserves coequalizers, this arrow is the image under $f^{\ast}$ of the canonical arrow $\xi:Z_{T} \to A$ from the codomain $Z_{T}$ of the coequalizer $Z\to Z_{T}$ of the two arrows $t_{1}, t_{2}$ to $A$, that is, $m\cong \textup{Im}(f^{\ast}(\xi))\cong f^{\ast}(\textup{Im}(\xi))$. 
\end{proofs}

\begin{proposition}\label{prophyperconnected}
	Let $({\cal C}, J)$ be a small-generated site, $\cal E$ a Grothendieck topos and $F:{\cal C}\to {\cal E}$ a $J$-continuous flat  functor. Then the geometric morphism $f:{\cal E}\to \Sh({\cal C}, J)$ induced by $F$ is hyperconnected if and only if $F$ is cover-reflecting and for every subobject $A\mono F(c)$ in $\cal E$ there exists a ($J$-closed) sieve $R$ on $c$ such that $A$ is the union of the images of the arrows $F(f)$ for $f\in R$. 
\end{proposition}

\begin{proofs}
	It is clear that if $f^{\ast}$ is faithful and its image is closed under subobjects then $f^{\ast}$ is also full (see the proof of the implication  (vi) $\imp$ (i) at p. 229 of \cite{El}). So it follows from Propositions \ref{propcoverreflecting} and  \ref{prophypersubobjects} that $f$ is hyperconnected if and only if $F$ is cover-reflecting and, for every object $c$ of $\cal C$, all the subobjects of $f^{\ast}(l(c))$ are images under $f^{\ast}$ of subobjects of $l(c)$ in $\Sh({\cal C}, J)$. Now, since the subobjects of $l(c)$ in $\Sh({\cal C}, J)$ are all images under the associated sheaf functor $a_{J}:[{\cal C}^{\textup{op}}, \Set] \to \Sh({\cal C}, J) $ of subobjects of the form $R \mono \Hom_{\cal C}(-,c)$ (which can be supposed $J$-closed without loss of generality), and this subobject is the union of the images of the arrows of the form $y_{\cal C}(f)$ for $f\in R$ (where $y_{\cal C}$ is the Yoneda embedding), we can reformulate the latter condition as the requirement that for any object $c$ of $\cal C$ and any subobject $A\mono F(c)$ in $\cal E$ there should exist a ($J$-closed) sieve $R$ on $c$ such that $A$ is the union of the images of the arrows $F(f)$ for $f\in R$.  
\end{proofs}

Let us apply this result in the context of interpretations of geometric theories (see section \ref{secsurjectioninclusion} for the notation).

\begin{corollary}\label{corsyntactichyperconnected}
An interpretation $I:{\cal C}_{\mathbb T}\to {\cal C}_{\mathbb S}$ between geometric theories $\mathbb T$ and $\mathbb S$ induces an hyperconnected geometric morphism $\Set[{\mathbb S}]\to \Set[{\mathbb T}]$ if and only if for every sequent $\sigma$ over the signature of $\mathbb T$, if the sequent $I(\sigma)$ is provable in $\mathbb S$ then $\sigma$ is provable in $\mathbb T$, and for any sorts $A_{1}, \ldots, A_{n}$ of the signature of $\mathbb T$, every geometric formula $\psi$ over the signature of $\mathbb S$ in the context $I((x_{1}^{A_{1}}, \ldots, x_{n}^{A_{n}}))$ is $\mathbb T$-provably equivalent (in its context) to a formula $I(\phi)$, where $\phi$ is a geometric formula over the signature of $\mathbb T$ in the context $(x_{1}^{A_{1}}, \ldots, x_{n}^{A_{n}})$. 	

In particular, if $I$ is the canonical interpretation given by an expansion ${\mathbb T}'$ of ${\mathbb T}$ (in the sense of section 7.1 of \cite{OCBook}) then the geometric morphism induced by $I$ is hyperconnected if and only if every sequent over the signature of $\mathbb T$ which is provable in ${\mathbb T}'$ is provable in $\mathbb T$, and for any sorts $A_{1}, \ldots, A_{n}$ of the signature of $\mathbb T$, every geometric formula $\psi$ in the context $(x_{1}^{A_{1}}, \ldots, x_{n}^{A_{n}})$ over the signature of ${\mathbb T}'$ is $\mathbb T$-provably equivalent (in its context) to a geometric formula over the signature of $\mathbb T$.
\end{corollary}

\begin{proofs}
	Our thesis follows immediately from Proposition \ref{prophyperconnected} in light of Corollary \ref{corinterpretationsfaithful}. 
\end{proofs}

Let us now apply Proposition \ref{prophyperconnected} in the context of morphisms of sites.

\begin{proposition}\label{prophyperconnectedsites}
		Let $F:({\cal C}, J)\to ({\cal D}, K)$ be a morphism of small-generated sites. Then the geometric morphism $\Sh(F):\Sh({\cal D}, K) \to \Sh({\cal C}, J)$ induced by $F$ is hyperconnected if and only if $F$ is cover-reflecting and \emph{closed-sieve-lifting}, in the sense that for every object $c$ of $\cal C$ and any $K$-closed sieve $S$ on $F(c)$ there exists a ($J$-closed) sieve $R$ on $c$ such that $S$ coincides with the ($K$-closed) sieve on $F(c)$ generated by the arrows $F(f)$ for $f\in R$.   
\end{proposition}

\begin{proofs}
	The thesis follows as an immediate consequence of Proposition \ref{prophyperconnected} in light of Proposition \ref{propclosedsubobjects}. Indeed, the condition of Proposition \ref{prophyperconnected} amounts precisely to the requirement that for every subobject $A\mono l'(F(c))$, where $l'$ is the canonical functor ${\cal D}\to \Sh({\cal D}, K)$, there should be a ($J$-closed) sieve $R$ on $c$ such that  $A$ is the union of the images of the arrows $l'(F(f))$ for $f\in R$. Now, the subobjects of $l'(F(c))$ in $\Sh({\cal D}, K)$ are the images under the associated sheaf functor $a_{K}:[{\cal D}^{\textup{op}}, \Set]\to \Sh({\cal D}, K)$ of $K$-closed sieves $S$ on $F(c)$, and the union of the images of the arrows $l'(F(f))$ for $f\in R$ is the image under $a_{K}$ of the union of the images of the arrows $y_{\cal D}(F(f))$ for $f\in R$; our claim thus follows from the fact that, by Proposition \ref{propclosedsubobjects}, two subobjects of a given object in $[{\cal D}^{\textup{op}}, \Set]$ have isomorphic images under $a_{K}$ if and only if their $K$-closures coincide.   
\end{proofs}

\subsubsection{Essential surjectivity}

\begin{proposition}\label{propessentiallysurjective}
	Let $({\cal C}, J)$ be a small-generated site, $\cal E$ a Grothendieck topos and $F:{\cal C}\to {\cal E}$ a $J$-continuous flat functor inducing a geometric morphism $f:{\cal E}\to \Sh({\cal C}, J)$. If the image of $f^{\ast}$ is closed under subobjects then $f^{\ast}$ is essentially surjective if and only if the objects of the form $F(c)$ for $c\in {\cal C}$ form a separating set for the topos $\cal E$.  
\end{proposition}

\begin{proofs}
	If the objects of the form $F(c)$ for $c\in {\cal C}$ form a separating set of objects for $\cal E$ then for every $A\in {\cal E}$, we have an epimorphism $q:f^{\ast}(\coprod_{i\in I}l(c_{i}))\cong \coprod_{i\in I}f^{\ast}(l(c_{i}))\epi A$. Since the image of $f^{\ast}$ is closed under subobjects, the kernel pair of $q$ is of the form $f^{\ast}(R)$, where $R$ is a relation on $\coprod_{i\in I}l(c_{i})$. Then $q$ is isomorphic to the image under $f^{\ast}$ of the coequalizer of $R$; in particular, $A$ is isomorphic to an object in the image of $f^{\ast}$, as desired.
	
	Conversely, if $f^{\ast}$ is essentially surjective then for every object $A$ of $\cal E$, there exists an object $B$ of $\Sh({\cal C}, J)$ such that $A\cong f^{\ast}(B)$. But $B$ can be covered in $\Sh({\cal C}, J)$ by objects of the form $l(c)$ for $\cal C$, whence $A$ can be covered in $\cal E$ by objects of the form $F(c)$ for $c\in {\cal C}$, as required.  
\end{proofs}

Let us apply Proposition \ref{propessentiallysurjective} in the context of interpretations of theories.

\begin{corollary}
		Let $I:{\cal C}_{\mathbb T}\to {\cal C}_{\mathbb S}$ be an interpretation of a geometric theory $\mathbb T$ into a geometric theory ${\mathbb S}$. Then the geometric morphism $\Set[{\mathbb S}]\to \Set[{\mathbb T}]$ induced by $I$ satisfies the property that its inverse image is essentially surjective and its image is closed under subobjects if and only if for every sort $B$ of the signature of $\mathbb S$, the object $\{x^{B}. \top\}$ is $J_{\mathbb S}$-covered by objects in the image of $I$ and the image of $I$ is closed under subobjects (that is, for any sorts $A_{1}, \ldots, A_{n}$ of the signature of $\mathbb T$, every geometric formula $\psi$ over the signature of $\mathbb S$ in the context $I((x_{1}^{A_{1}}, \ldots, x_{n}^{A_{n}}))$ is $\mathbb T$-provably equivalent (in its context) to a formula $I(\phi)$, where $\phi$ is a geometric formula over the signature of $\mathbb T$ in the context $(x_{1}^{A_{1}}, \ldots, x_{n}^{A_{n}})$). 	    
\end{corollary}

\begin{proofs}
	By Propositions \ref{prophypersubobjects} and \ref{propessentiallysurjective}, since the geometric syntactic category ${\cal C}_{\mathbb S}$ is closed under subobjects in the topos $\Sh({\cal C}_{\mathbb S}, J_{\mathbb S})$, it suffices to show, assuming that the image of $I$ is closed under subobjects, that the objects in the image of $I$ form a separating set for the topos $\Set[{\mathbb S}]$ if and only if every object of the form $\{x^{B}. \top\}$ (for $B$ a sort of the signature of $\mathbb S$) is $J_{\mathbb S}$-covered by objects in the image of $I$. But this follows from the fact that every object of ${\cal C}_{\mathbb S}$ is a subobject of a finite product of objects of the form $\{x^{B}. \top\}$, using the fact that $I$ preserves finite products and its image is closed under subobjects.  
\end{proofs}

Let us now apply Proposition \ref{propessentiallysurjective} in the context of morphisms of sites.

\begin{proposition}\label{propessentialsurjectivitysites}
	Let $F:({\cal C}, J)\to ({\cal D}, K)$ be a morphism of small-generated sites. Then the geometric morphism $\Sh(F):\Sh({\cal D}, K) \to \Sh({\cal C}, J)$ induced by $F$ satisfies the property that $\Sh(F)^{\ast}$ is essentially surjective and its image is closed under subobjects if and only if $F$ is closed-sieve-lifting (in the sense of Proposition \ref{prophyperconnectedsites}) and satisfies condition (ii) of Proposition \ref{propweaklydensemorphism} (equivalently, if $K$ is subcanonical, every object of $\cal D$ admits a $K$-covering sieve generated by arrows whose domain lies in the image of $F$).
\end{proposition}

\begin{proofs}
	The equivalence between the property of the image of $\Sh(F)^{\ast}$ to be closed under subobjects and the property of $F$ to be closed-sieve-lifting follows from the proof of Proposition \ref{prophyperconnectedsites}, while the equivalence between the property of $\Sh(F)^{\ast}$ to be essentially surjective and condition (ii) of Proposition \ref{propweaklydensemorphism} follows from the proof of Proposition \ref{propweaklydensemorphism}.
\end{proofs}

\subsubsection{Localic morphisms}

The following proposition shows the natural behavior of the property of a geometric morphism being localic in terms of separating sets for toposes.

\begin{proposition}\label{proplocalic}
	Let ${\cal C}$ (resp. $\cal D$) be a separating set for a Grothendieck topos $\cal E$ (resp. $\cal F$) and $f:{\cal F}\to {\cal E}$ a geometric morphism.
	Then: 
	\begin{enumerate}[(i)]
		\item The full subcategory $\cal G$ of $\cal F$ on the objects which are subquotients of objects of the form $f^{\ast}(A)$ for $A\in {\cal E}$ coincides with the full subcategory of $\cal F$ on the objects which can be covered by objects which are domains of subobjects of objects of the form $f^{\ast}(A)$ for $A$ in $\cal C$. 
		
		\item $f$ is localic if and only if every object of $\cal D$ can be covered by objects which are domains of subobjects of objects of the form $f^{\ast}(A)$ for $A$ in $\cal C$. 
	\end{enumerate}	
\end{proposition}

\begin{proofs}
	(i) Let $B\in {\cal G}$. Then $B$ is a quotient $d\to B$ of a subobject $d\mono f^{\ast}(A)$. But $A$ is a quotient $q:\coprod_{i\in I}c_{i}\to A$ of a coproduct of objects $c_{i}$ in $\cal C$, so, by considering the pullback
	\[
	\xymatrix{
		d' \ar[d] \ar[r] & f^{\ast}(\coprod_{i\in I}c_{i})\cong \coprod_{i\in I}f^{\ast}(c_i) \ar[d]^{q}  \\
		d \ar[r] & f^{\ast}(A) } 
	\] 
	we see that $B$ is a quotient of $d'$, which is a coproduct of subobjects of the objects $f^{\ast}(c_i)$. Therefore $B$ can be covered by the domains of these subobjects, as required.
	
	Conversely, let us suppose that $B$ is an object of $\cal F$ that can be covered by objects $u_{i}$ (for $i\in I$) which are domains of subobjects $u_{i} \mono f^{\ast}(c_{i})$ of objects of the form $f^{\ast}(c_{i})$ for $c_{i}$ in $\cal C$. Then $B$ is a quotient of the coproduct of the $u_{i}$ (for $i\in I$), which is a subobject of the coproduct of the $f^{\ast}(c_{i})$; but this coproduct is in the image of $f^{\ast}$ as $f^{\ast}$ preserves coproducts. 
	
	(ii) By definition of a localic morphism, $f$ is localic if and only if ${\cal G}={\cal F}$. So, by point (i), if $f$ is localic then every object of $\cal D$ can be covered by objects which are domains of subobjects of objects of the form $f^{\ast}(A)$ for $A$ in $\cal C$. Conversely, again by (i), we have to prove that if every object of $\cal D$ can be covered by objects which are domains of subobjects of objects of the form $f^{\ast}(A)$ for $A$ in $\cal C$ then every object of $\cal F$ can be covered by objects which are domains of subobjects of objects of the form $f^{\ast}(A)$ for $A$ in $\cal C$; but this immediately follows from the fact that $\cal D$ is separating for $\cal F$.  
\end{proofs}

\begin{remark}\label{remseparatingsetlocalic}
	We know from section A4.6 of \cite{El} that if $\cal E$ and $\cal F$ are Grothendieck toposes then the category $\cal G$ defined in point (i) of Proposition \ref{proplocalic} is a Grothendieck topos; so this part of the proposition actually asserts that the collection of objects which are subobjects of an object of the form $f^{\ast}(A)$ for $A\in {\cal C}$ is a separating set for $\cal G$. 
\end{remark}

The following result is an immediate corollary of Proposition \ref{proplocalic}.

\begin{proposition}
	Let $({\cal C}, J)$ be a small-generated site, $\cal E$ a Grothendieck topos and $F:{\cal C}\to {\cal E}$ a $J$-continuous flat functor inducing a geometric morphism $f:{\cal E}\to \Sh({\cal C}, J)$. Then $f$ is localic if and only if the subobjects of objects of the form $F(c)$ for $c\in {\cal C}$ form a separating set for the topos $\cal E$.
\end{proposition}\qed

Let us now apply Proposition \ref{proplocalic} in the general context of morphisms of sites.

\begin{proposition}\label{proplocalicsites}
	Let $F:({\cal C}, J)\to ({\cal D}, K)$ be a morphism of small-generated sites. Then the geometric morphism $\Sh(F):\Sh({\cal D}, K) \to \Sh({\cal C}, J)$ induced by $F$ is localic if and only if for any object $d$ of $\cal D$ there exist a family $\{S_{i} \mid i\in I\}$ of sieves on objects of the form $F(c_{i})$ (where $c_{i}$ is an object of $\cal C$) and for each $f\in S_{i}$ an arrow $g_{f}:\textup{dom}(f)\to d$ such that $g_{f\circ z}\equiv_{K} g_{f}\circ z$ whenever $z$ is composable with $f$, such that the family of arrows $g_{f}$ (for $f\in S_{i}$ for some $i$) is $K$-covering.  
\end{proposition}

\begin{proofs}
Let us suppose that $\Sh(F)$ is localic. Let $l'$ be the canonical functor ${\cal D}\to \Sh({\cal D}, K)$. By Proposition \ref{proplocalic}, for every $d\in {\cal D}$ there exist objects $c_{i}$ of $\cal C$ (for $i\in I$), subobjects $A_{i}\mono l'(F(c_i))$ and a jointly epimorphic family of arrows $A_{i}\to l'(d)$ (for $i\in I$). Now, we can write $A_{i}=a_{K}(S_{i})$, where $S_i$ is a sieve on $F(c_i)$. For each $f\in S_{i}$, composing the arrow $a_{K}(S_{i})\to l'(d)$ with the canonical arrow $l'(\textup{dom}(f))\to a_{K}(S_{i})$ thus yields an arrow $l'(\textup{dom}(f))\to l'(d)$ which we can suppose without loss of generality of the form $l'(g_{f})$ for an arrow $g_{f}:\textup{dom}(f)\to d$ (at the cost of replacing $S_{i}$, by using Proposition \ref{propexplicit}, which a smaller sieve with the same image under $a_{K}$); clearly, for any $z$ composable with $f$, $l'(g_{f\circ z})=l'(g_{f}\circ z)$, equivalently $g_{f\circ z}\equiv_{K} g_{f}\circ z$. The property of our arrows $A_{i}\to l'(d)$ to be jointly epimorphic can thus be reformulated as the property of the arrows $g_{f}$ to be sent by $l'$ to an epimorphic family, equivalently to be $K$-covering.

Conversely, any family of arrows $g_{f}:\textup{dom}(f)\to d$ indexed by a sieve $S_{i}$ such that $g_{f\circ z}\equiv_{K} g_{f}\circ z$ whenever $z$ is composable with $f$ induces an arrow $a_K(S_i)\to l'(d)$. So if we have a family of sieves $S_{i}$ on objects of the form $F(c_i)$ and one such family of arrows for each sieve, the resulting arrows $a_K(S_i)\to l'(d)$ will be jointly epimorphic if (and only if) the family of arrows $g_{f}$ (for $f\in S_{i}$ for some $i$) is $K$-covering.    	
\end{proofs}

The following result characterizes the interpretations between geometric theories which induce a localic geometric morphism.

\begin{proposition}
	Let $I:{\cal C}_{\mathbb T}\to {\cal C}_{\mathbb S}$ be an interpretation of a geometric theory $\mathbb T$ into a geometric theory ${\mathbb S}$. Then the geometric morphism $\Set[{\mathbb S}]\to \Set[{\mathbb T}]$ induced by $I$ is localic if and only if every object of ${\cal C}_{\mathbb S}$ can be $J_{\mathbb S}$-covered by formulae over the signature of $\mathbb S$ in a context of the form $I(\vec{x})$ which $\mathbb T$-provably entail a formula of the form $I(\phi)$, where $\phi$ is a geometric formula in the context $\vec{x}$ over the signature of $\mathbb T$. 
	
	In particular, the canonical interpretation of $\mathbb T$ into an expansion of it over a signature which does not contain any new sorts with respect to the signature of $\mathbb T$ induces a localic geometric morphism.  
\end{proposition}

\begin{proofs}
	This is an immediate consequence of Proposition \ref{proplocalic}, by using the fact that the subobjects in $\Set[{\mathbb S}]\simeq \Sh({\cal C}_{\mathbb S}, J_{\mathbb S})$ of an object of the form $I(\{\vec{x}. \phi\})$ for a geometric formula-in-context $\phi(\vec{x})$ over the signature of $\mathbb T$ can be identified with the ($\mathbb S$-provable equivalence classes of) geometric formulae in the context $I(\vec{x})$ over the signature of $\mathbb S$ which $\mathbb S$-provably entail $I(\phi)$.
\end{proofs}

\subsubsection{The hyperconnected-localic factorization}\label{sec:hyperconnfact}

In this section we shall provide a site-level description of the hyperconnected-localic factorization (in the sense of \cite{Johnstonefactorization} - cf. also section A4.6 of \cite{El}) of a geometric morphism between Grothendieck toposes.

Let us refer to Corollary \ref{corhyperconnectedsieves} for the notation. Let $C_{J}^{s}$ be the Grothendieck topology induced on ${\cal C}_{J}^{s}$ by the canonical topology on the topos $\Sh({\cal C}, J)$.

We have a canonical functor ${\cal C}\to {\cal C}_{J}^{s}$ sending an object $c$ of $\cal C$ to the pair $(c, M_{c})$, where $M_{c}$ is the maximal sieve on $c$, which yields a morphism of sites $i_{J}^{s}:({\cal C}, J)\to ({\cal C}_{J}^{s}, C_{J}^{s})$. This morphism induces an equivalence $\Sh(i_{J}^{s}):\Sh({\cal C}_{J}^{s}, C_{J}^{s}) \simeq \Sh({\cal C}, J)$.  

\begin{theorem}\label{thmhyperconnectedlocalicfactorization}
	Let $F:({\cal C}, J)\to ({\cal D}, K)$ be a morphism of small-generated sites. Then the hyperconnected-localic factorization of the geometric morphism $\Sh(F)\circ \Sh(i_{K}^{s}):\Sh({\cal D}_{K}^{s}, C_{K}^{s})\simeq \Sh({\cal D}, K) \to \Sh({\cal C}, J)$ induced by $F$ can be identified with the factorization $\Sh(F^{s})\circ \Sh(i_{F}^{s})$, where
	${\cal D}_{F}^{s}$ is the full subcategory of ${\cal D}_{K}^{s}$ on the objects of the form $(F(c), S)$ for an object $c$ of $\cal C$ and a $K$-closed sieve $S$ on $F(c)$, $C_{K}^{s}|_{{\cal D}_{F}^{s}}$ is the Grothendieck topology induced by $C_{K}^{s}$ on ${\cal D}_{F}^{s}$, $i_{F}^{s}:({\cal D}_{F}^{s}, C_{K}^{s}|_{{\cal D}_{F}^{s}})\to ({\cal D}_{K}^{s}, C_{K}^{s})$ is the morphism of sites given by the canonical inclusion functor ${\cal D}_{F}^{s}\hookrightarrow {\cal D}_{K}^{s}$, and $F_{s}$ is the morphism of sites $({\cal C}, J)\to ({\cal D}_{F}^{s}, C_{K}^{s}|_{{\cal D}_{F}^{s}})$ given by the functor $F$.
\end{theorem}

\begin{proofs}
		By Remark \ref{remseparatingsetlocalic}, the topos $\cal G$ arising in the hyperconnected-localic factorization $p'\circ p$ of the geometric morphism $\Sh(F)$ admits as separating set of objects the family of domains of subobjects of objects of the form $\Sh(F)^{\ast}(l(c))\cong l'(F(c))$. Now, the full subcategory of $\Sh({\cal D}, K)$ on such objects is equivalent to ${\cal D}_{F}^{s}$ by Corollary \ref{corhyperconnectedsieves}, so we have an equivalence $\xi:{\cal G}\to  \Sh({\cal D}_{F}^{s}, C_{K}^{s}|_{{\cal D}_{F}^{s}})$. The inverse image of the geometric morphism $p:\Sh({\cal D}, K)\to {\cal G}$ is precisely the inclusion functor of $\cal G$ into $\Sh({\cal D}, K)$, and the morphism of sites $i_{K}^{s}:({\cal D}, K)\to ({\cal D}_{K}^{s}, C_{K}^{s})$ yields an equivalence $\Sh({\cal D}_{K}^{s}, C_{K}^{s})\simeq \Sh({\cal D}, K)$.  So the canonical inclusion functor $i_{F}^{s}$ of ${\cal D}_{F}^{s}$ into ${\cal D}_{K}^{s}$ defines a morphism of sites $({\cal D}_{F}^{s}, C_{K}^{s}|_{{\cal D}_{F}^{s}})\to ({\cal D}_{K}^{s}, C_{K}^{s})$, since it induces the geometric morphism $\xi \circ p\circ \Sh(i_{K}^{s})$. On the other hand, the inverse image of the geometric morphism $p':{\cal G}\to \Sh({\cal C}, J)$ given by the localic part of the hyperconnected-localic factorization of $\Sh(F)$ is simply $\Sh(F)^{\ast}$, regarded as taking values in the subcategory $\cal G$ of $\Sh({\cal D}, K)$. Therefore the functor $F_{s}$ is a morphism of sites $({\cal C}, J)\to ({\cal D}_{F}^{s}, C_{K}^{s}|_{{\cal D}_{F}^{s}})$ as it induces the geometric morphism $p'\circ \xi^{-1}$.  
		This completes our proof. 
\end{proofs}
	
Let us now apply Theorem \ref{thmhyperconnectedlocalicfactorization} in the context of interpretations between geometric theories.

\begin{corollary}\label{corhyperclocalicintepretation}
	Let $I:{\cal C}_{\mathbb T}\to {\cal C}_{\mathbb S}$ be an interpretation of a geometric theory $\mathbb T$ into a geometric theory ${\mathbb S}$. Let ${\mathbb T}_{I}$ be the expansion of $\mathbb T$ whose signature is obtained by adding a predicate symbol $R_{\psi}$ of sorts $I(\vec{x})$ for any context $\vec{x}$ over the signature of $\mathbb T$ and any geometric formula $\psi$ in the context $I(\vec{x})$ over the signature of $\mathbb S$, and whose axioms are all the sequents whose image under the extension of $I$ sending each predicate symbol $R_{\psi}$ to the corresponding formula $\{I(\vec{x}). \psi\}$ is provable in $\mathbb S$, and let $I_{s}$ be the obvious interpretation of $\mathbb T$ into ${\mathbb T}_{I}$.
	Then the hyperconnected-localic factorization of the geometric morphism $\Set[{\mathbb S}]\to \Set[{\mathbb T}]$ induced by $I$ can be identified with the composite $\Sh(I')\circ \Sh(I_{s})$.
\end{corollary}
	
\begin{proofs}
	It is immediate to see that the factorization of $I$ as $I_{s}$ followed by $I'$ can be identified with the factorization of $I$, regarded as a morphism of sites $({\cal C}_{\mathbb T}, J_{\mathbb T})\to ({\cal C}_{\mathbb S}, J_{\mathbb S})$, provided by Theorem \ref{thmhyperconnectedlocalicfactorization} by using Proposition \ref{propcoverreflecting}, Corollary \ref{corinterpretationsfaithful} and Remark \ref{remsubcanonicalclosedsubobjects}.  
\end{proofs}		

\begin{remark}
	Corollary \ref{corhyperclocalicintepretation} generalizes the corresponding result for expansions of theories proved as Theorem 7.1.3 in \cite{OCBook}.
\end{remark}

\subsection{Equivalence of toposes}\label{sec:equivalencetoposes2}

From the above results we can obtain a criterion (alternative to Corollary \ref{flatequivalence}) for a $J$-continuous flat functor to induce an equivalence of toposes:

\begin{corollary}\label{corflatequivalence2}
		Let $({\cal C}, J)$ be a small-generated site, $\cal E$ a Grothendieck topos and $F:{\cal C}\to {\cal E}$ a $J$-continuous flat functor inducing a geometric morphism $f:{\cal E}\to \Sh({\cal C}, J)$. Then $f$ is an equivalence if and only if the following conditions are satisfied:
		\begin{enumerate}[(i)]
			\item $F$ is cover-reflecting;
			
			\item for every subobject $A\mono F(c)$ in $\cal E$ there exists a ($J$-closed) sieve $R$ on $c$ such that $A$ is the union of the images of the arrows $F(f)$ for $f\in R$;
			
			\item the objects of the form $F(c)$ for $c\in {\cal C}$ form a separating set for the topos $\cal E$.
		\end{enumerate} 
\end{corollary}

\begin{proofs}
	The functor $f^{\ast}$ is one half of an equivalence of categories if and only if it is full and faithful and essentially surjective. Since faithfulness and closure of its image under subobjects implies fullness for the functor $f^{\ast}$, we conclude that $f^{\ast}$ is one half of an equivalence if and only if it is faithful, essentially surjective and its image is closed under subobjects. Our thesis thus follows from Propositions \ref{propcoverreflecting}, \ref{prophypersubobjects} and \ref{propessentiallysurjective}. 
\end{proofs}

It is interesting to compare this characterization with the different, but equivalent, criterion provided by Corollary \ref{flatequivalence}: the difference lies in condition (ii) of Corollary \ref{corflatequivalence2} which is replaced by the ``local fullness'' condition of Corollary \ref{flatequivalence}.

This corollary can be notably applied in the context of interpretations between geometric theories to characterize those which induce Morita equivalences.  

\begin{corollary}\label{corsyntacticMorita}
	Let $I:{\cal C}_{\mathbb T}\to {\cal C}_{\mathbb S}$ be an interpretation of a geometric theory $\mathbb T$ into a geometric theory ${\mathbb S}$. Then $I$ induces a Morita-equivalence between $\mathbb T$ and $\mathbb S$ if and only if the following conditions are satisfied:
	
	\begin{enumerate}[(i)]
		\item for every sequent $\sigma$ over the signature of $\mathbb T$, if the sequent $I(\sigma)$ is provable in $\mathbb S$ then $\sigma$ is provable in $\mathbb T$;
		
		\item the image of $I$ is closed under subobjects (that is, for any sorts $A_{1}, \ldots, A_{n}$ of the signature of $\mathbb T$, every geometric formula $\psi$ over the signature of $\mathbb S$ in the context $I((x_{1}^{A_{1}}, \ldots, x_{n}^{A_{n}}))$ is $\mathbb T$-provably equivalent (in its context) to a formula of the form $I(\phi)$ for a geometric formula $\phi$ over the signature of $\mathbb T$ in the context $(x_{1}^{A_{1}}, \ldots, x_{n}^{A_{n}})$); 	 
		
		\item for every sort $B$ of the signature of $\mathbb S$, the object $\{x^{B}. \top\}$ is $J_{\mathbb S}$-covered by objects in the image of $I$.  
	\end{enumerate}  
\end{corollary}

\begin{proofs}
	Our thesis immediately follows from Corollaries \ref{corinterpretationsfaithful} and \ref{corsyntactichyperconnected}.
\end{proofs}

\begin{remark}
	Conditions (i) and (ii) of Corollary \ref{corsyntacticMorita} can be reformulated as the requirement that the functor $I$ should be an equivalence onto its image.
\end{remark}

\begin{example}
	Corollary \ref{corsyntacticMorita} can for instance be applied to the interpretation of the algebraic theory $\mathbb{MV}$ of MV-algebras into the theory ${\mathbb L}_{u}$ of lattice-ordered abelian groups with strong unit established in \cite{OC17}, which, as proved in $\cite{OC17}$, induces a Morita equivalence. Condition (iii) is satisfied since the formula $\{x. \top\}$ is covered in ${\cal C}_{{\mathbb L}_{u}}$ by the objects in the image of $I$. Indeed, by the axiom expressing the property of strong unit, the arrows $\{x'. 0\leq x'\leq u\}\to \{x. \top\}$ given by $x=nx'$ (for $n\in {\mathbb Z}$) generate a $J_{{\mathbb L}_{u}}$-covering sieve. One can also understand concretely why condition (ii) is satisfied; indeed, by using the construction of the $\ell$-group associated with an MV-algebra in terms of good sequences of elements of the latter, one can prove that for every atomic formula $\chi(\vec{x})$ over the signature of ${\mathbb L}_{u}$, there exists a Horn formula $\chi^{MV}(\vec{x})$ in the same context over the signature of $\mathbb{MV}$ such that the formula $I(\chi^{MV}(\vec{x}))$ is provably equivalent in ${\mathbb L}_{u}$ (in the context $\vec{x}$) to the formula $\chi \wedge (0\leq \vec{x} \leq u)$.
\end{example}

From Proposition \ref{propessentialsurjectivitysites} and Theorem \ref{Morphismsmaalgenerated}(a), we deduce the following criterion, alternative to Theorem \ref{thmweakdensityequivalence}, for a morphism of sites to induce an equivalence of toposes:
	
\begin{corollary}\label{corquivalencesites}
Let $F:({\cal C}, J)\to ({\cal D}, K)$ be a morphism of small-generated sites. Then the geometric morphism $\Sh(F):\Sh({\cal D}, K) \to \Sh({\cal C}, J)$ induced by $F$ is an equivalence if and only if $F$ is cover-reflecting, closed-sieve-lifting and satisfies condition (ii) of Proposition \ref{propweaklydensemorphism} (recall that, if $K$ is subcanonical, the latter condition is equivalent to $F$ being $K$-dense).
\end{corollary}	

\begin{proofs}
Since a geometric morphism $f$ is an equivalence if and only if it is hyperconnected and its inverse image is essentially surjective, the thesis follows from Proposition \ref{propessentialsurjectivitysites} and Theorem \ref{Morphismsmaalgenerated}(a).	
\end{proofs}

The following criterion, alternative to Corollary \ref{corinclusionequivalence}, for a geometric morphism to be an inclusion (resp. an equivalence) follows as an immediate consequence of Corollary \ref{corquivalencesites}.

\begin{corollary}\label{corinclusionequivalence2}
	Let $f:{\cal F}\to {\cal E}$ be a geometric morphism. Then
	
	\begin{enumerate}[(i)]
		\item $f$ is an inclusion if and only if $f^{\ast}$ is locally surjective and its image is closed under subobjects.

		\item $f$ is an equivalence if and only if $f^{\ast}$ is faithful, locally surjective and its image is closed under subobjects.		
	\end{enumerate}
\end{corollary}

\begin{proofs}
	As observed in the proof of Corollary \ref{corinclusionequivalence}, $f$ is an inclusion if and only if the morphism of sites $f^{\ast}:({\cal E}, ({J_{\cal E}^{\textup{can}}})_{f^{\ast}})\to ({\cal F}, J_{\cal F}^{\textup{can}})$ is dense (equivalently, weakly dense). Now, this morphism is $ J_{\cal F}^{\textup{can}}$-dense if and only if $f^{\ast}$ is locally surjective (in the sense of Corollary \ref{corinclusionequivalence}(i)), while it is closed-sieve-lifting if and only if the image of $f^{\ast}$ is closed under subobjects; indeed, for any Grothendieck topology $K$ on a Grothendieck topos $\cal G$ containing the canonical one, any $K$-closed sieve is equal to the principal sieve on a subobject. This proves the first part of the corollary; the second part follows from the first by using that a geometric morphism is an equivalence if and only if it is both a surjection and an inclusion.	
\end{proofs}

\subsection{An example}\label{sec:example}

Given a small-generated site $({\cal C}, J)$ and a presheaf $P:{\cal C}^{\textup{op}}\to \Set$, let us deduce, as a consequence of Corollary \ref{corflatequivalence2}, that we have an equivalence of toposes
\[
\Sh({\cal C}, J)\slash a_{J}(P)\simeq \Sh(\textstyle \int P, J_{P}),
\]
where $\int P$ is the category of elements of $P$ and $J_{P}$ is the Grothendieck topology on $\int P$ whose covering sieves are the sieves which are sent by the canonical projection $\pi_{P}:\int P \to {\cal C}$ to $J$-covering families. This equivalence, which was already established in Proposition III 5.4 \cite{grothendieck}, is induced by the (flat) $J_{P}$-continuous functor $F_{P}:\int P \to \Sh({\cal C}, J)\slash a_{J}(P)$ sending any object $(c, s)$ of $\int P$ to the image under $a_{J}$ of the arrow $y_{\cal C}(c)\to P$ induced by $(c, x)$ via the Yoneda lemma. First, we observe that, since we have an equivalence 
\[
[{\cal C}^{\textup{op}}, \Set]\slash P \simeq [({\textstyle \int P})^{\textup{op}}, \Set]
\]
(see, for instance, Proposition A1.1.7 \cite{El}) induced (by composition with the Yoneda embedding ${\int P}\to [({\int P})^{\textup{op}}, \Set]$) by the canonical functor $G_{P}:\int P\to [{\cal C}^{\textup{op}}, \Set]\slash P$ (sending any object $(c, x)\in \int P$ to the object $y_{\cal C}(c)\to P$ of $[{\cal C}^{\textup{op}}, \Set]\slash P$ given by the arrow corresponding to $x$ via the Yoneda Lemma), $G_{P}$ is flat. Now, $F_{P}$ is the composite of $G_{P}$ with the inverse image of the geometric morphism $\Sh({\cal C}, J)\slash a_{J}(P) \to [{\cal C}^{\textup{op}}, \Set]\slash P$ induced by the canonical geometric inclusion $\Sh({\cal C}, J)\hookrightarrow [{\cal C}^{\textup{op}}, \Set]$ as in Example \ref{exgeometricinclusion}, so it is flat as $G_{P}$ is. Notice that the Grothendieck topology $J_{P}$ is precisely the topology induced by $F_{P}$ (in the sense of Proposition \ref{reminducedtopology}).
Let us apply Corollary \ref{corflatequivalence2} to prove that $F_{P}$ induces an equivalence. First, we have to show that the objects in the image of the functor $F_{P}$ define a separating set of objects for the topos $\Sh({\cal C}, J)\slash a_{J}(P)$. Given two distinct arrows $\alpha, \beta:(\gamma:Q\to a_{J}(P))\to (\xi:Q'\to a_{J}(P))$ in $\Sh({\cal C}, J)\slash a_{J}(P)$, we have to show that there are an object $c$ of $\cal C$, an element $x\in P(c)$, corresponding to an arrow $P_{x}:l(c)\to a_{J}(P)$ in $\Sh({\cal C}, J)$, and an element $y\in Q(c)$, corresponding to an arrow $\overline{y}:l(c)\to Q$ in $\Sh({\cal C}, J)$, such that $P_{x}=\gamma \circ \overline{y}$ and $\alpha\circ \overline{y}\neq \beta\circ \overline{y}$. Consider the pullback $\tilde{\gamma}:\tilde{Q}\to P$ of $\gamma:Q\to a_{J}(P)$ along the unit $P\to a_{J}(P)$. Since the canonical arrow $z:\tilde{Q}\to Q$ is sent by $a_{J}$ to an isomorphism and $a_{J}$ is left  adjoint to the inclusion functor, $\alpha\circ z\neq \beta \circ z$. So there exists an object $c$ of $\cal C$ and an element $y'\in \tilde{Q}(c)$ such that $\alpha(z(c)(y'))\neq \beta(z(c)(y'))$. So, posing $x=\tilde{\gamma}(c)(y')$ and $y=z(c)(y')$, we obtain that $P_{x}=\gamma \circ \overline{y}$ and $\alpha\circ \overline{y}\neq \beta\circ \overline{y}$, as required. 

It now remains to show that condition (ii) is satisfied by $F_{P}$; but this is clear, since, for any object $(c, x)$ of $\int P$, any subobject $A$ of the object $P_{x}:l(c)\to a_{J}(P)$ in $\Sh({\cal C}, J)\slash a_{J}(P)$ corresponds to a subobject of $l(c)$ in $\Sh({\cal C}, J)$ and hence to a sieve $S$ on $c$, which clearly lifts to a sieve $S_{\pi}$ on the object $(c, x)$ in $\int P$ such that $A$ is isomorphic to the union of the images of the arrows $F_{P}(f)$ where $f\in S_{\pi}$.  

\subsection{A characterization of quotient interpretations}

Let us now apply Corollaries \ref{flatequivalence} and \ref{corflatequivalence2} in the context of interpretations of geometric theories, in order to characterize the interpretations of a geometric theory $\mathbb T$ which are isomorphic to the canonical interpretation of $\mathbb T$ into a quotient of it.

\begin{corollary}\label{corinterpretationinclusion}
	Let $I:{\cal C}_{\mathbb T}\to {\cal C}_{\mathbb S}$ be an interpretation of a geometric theory $\mathbb T$ into a geometric theory ${\mathbb S}$. Then $I$ is isomorphic (over ${\cal C}_{\mathbb T}$) to the canonical interpretation $I_{\mathbb T}^{{\mathbb T}'}:{\cal C}_{\mathbb T}\to {\cal C}_{{\mathbb T}'}$ of ${\mathbb T}$ into a quotient ${\mathbb T}'$ of ${\mathbb T}$ if and only if it satisfies the following conditions:
		\begin{enumerate}[(i)]		
		\item The image of $I$ is closed under subobjects (that is, for any sorts $A_{1}, \ldots, A_{n}$ of the signature of $\mathbb T$, every geometric formula $\psi$ over the signature of $\mathbb S$ in the context $I((x_{1}^{A_{1}}, \ldots, x_{n}^{A_{n}}))$ is $\mathbb T$-provably equivalent (in its context) to a formula of the form $I(\phi)$ for a geometric formula $\phi$ over the signature of $\mathbb T$ in the context $(x_{1}^{A_{1}}, \ldots, x_{n}^{A_{n}})$). 	 
		
		\item For every sort $B$ of the signature of $\mathbb S$, the object $\{x^{B}. \top\}$ is $J_{\mathbb S}$-covered by objects in the image of $I$.  
	\end{enumerate}
	Condition (i) can be equivalently replaced by the following condition : 
	\begin{enumerate}[(i)]		
		\item[(i)'] For every objects $\{\vec{x}. \phi\}$ and $\{\vec{y}. \psi\}$ of ${\cal C}_{\mathbb T}$ and any $\mathbb S$-provably functional formula $\theta(I(\vec{x}), I(\vec{y}))$ from $I(\{\vec{x}. \phi\})$ to $I(\{\vec{y}. \psi\})$ there exist a set of $\mathbb T$-provably functional formulae $\{\chi_{i}(\vec{x_i}, \vec{x}) \}:\{\vec{x_i}. \phi_{i}\}\to \{\vec{x}. \phi\}$ such that the sequent
		$(I(\phi) \vdash_{I(\vec{x})} \mathbin{\mathop{\textup{\huge $\vee$}}\limits_{i}}(\exists I(\vec{x_{i}}) I(\chi_{i}(\vec{x_i}, \vec{x})))$ is provable in $\mathbb S$ and $\mathbb T$-provably functional formulae $\xi_{i}(\vec{x_i}, \vec{y}):\{\vec{x_i}. \phi_{i}\} \to \{\vec{y}. \psi\}$ (for each $i$) such that $[\theta]\circ I([\chi_{i}])=I([\xi_{i}])$ in ${\cal C}_{\mathbb S}$ (i.e. the bi-sequent $((\exists I(\vec{x}))(\theta \wedge I([\chi_i])) \dashv \vdash_{I(\vec{x_{i}})} I([\xi_{i}]))$ is provable in $\mathbb S$). 
	\end{enumerate} 
\end{corollary}	  

\begin{proofs}
	By the duality theorem between quotients and subtoposes established in Chapter 3 of \cite{OCBook}, we have to characterize the interpretations $I$ which induce a geometric inclusion to the classifying topos of $\mathbb T$. The formulation of this condition as the conjunction of conditions (i) and (ii) follows from Corollary \ref{corflatequivalence2} by observing that the geometric morphism induced by $I$ is induced by the $(J_{\mathbb T})_{I}$-continuous flat functor ${\cal C}_{\mathbb T}\to \Sh({\cal C}_{\mathbb S}, J_{\mathbb S})$ given by the composite of the Yoneda embedding ${\cal C}_{\mathbb S}\hookrightarrow \Sh({\cal C}_{\mathbb S}, J_{\mathbb S})$ with $I$.
	
	On the other hand, conditions (i)' and (ii) are equivalent to the property of the interpretation $I$ to induce a dense morphism of sites $({\cal C}_{\mathbb T}, (J_{\mathbb T})_{I})\to ({\cal C}_{\mathbb S}, J_{\mathbb S})$, equivalently, by Theorem \ref{Morphismsmaalgenerated}(c), a geometric inclusion from the classifying topos of ${\mathbb S}$ to the classifying topos of $\mathbb T$. 	
\end{proofs}

\begin{remarks}
	\begin{enumerate}[(a)]
		\item The quotient ${\mathbb T}'$ of ${\mathbb T}$ corresponding to an interpretation $I$ satisfying the conditions of Corollary \ref{corinterpretationinclusion} is axiomatized by the sequents $\sigma$ over the signature of $\mathbb T$ such that $I(\sigma)$ is provable in $\mathbb S$ (cf. \ref{corsyntacticMorita}). 
		
		\item It is interesting to understand why property (i)' is true in the case of the canonical interpretation of $\mathbb T$ into a quotient ${\mathbb T}'$ of it.
		Given a ${\mathbb T}'$-provably functional formula $\theta(I(\vec{x}), I(\vec{y}))$ from $I(\{\vec{x}. \phi\})$ to $I(\{\vec{y}. \psi\})$, we can take $\chi$ to be the formula $(\theta \wedge \phi \wedge \psi)(\vec{x}, \vec{y}) \wedge \vec{x}=\vec{x'}$, which is ${\mathbb T}$-provably functional from $\{\vec{x}, \vec{y}. \theta \wedge \phi \wedge \psi\}$ to $\{\vec{x'}. \phi[\vec{x'}\slash \vec{x}]\}$, and $\xi$ to be the formula $(\theta \wedge \phi \wedge \psi)(\vec{x}, \vec{y}) \wedge \vec{y}=\vec{y'}$, which is $\mathbb T$-provably functional from $\{\vec{x}, \vec{y}. \theta \wedge \phi \wedge \psi\}$ to $\{\vec{y'}. \psi[\vec{y'}\slash \vec{y}]\}$.  		 
	\end{enumerate}
	\end{remarks}

\section{Properties of morphisms induced by comorphisms of sites}\label{sec:comorphismsofsites}

In this section we shall investigate the conditions on a comorphism of sites which correspond to the property of the corresponding geometric morphism to be a surjection (resp. an inclusion, hyperconnected, localic), and also describe the surjection-inclusion and hyperconnected-localic factorizations at the level of comorphisms of sites. 

\subsection{Surjections}\label{sec:surjectionscomorphisms}

\begin{proposition}\label{propsurjectioncomorphismofsites}
The geometric morphism $C_{F}:\Sh({\cal D}, K)\to \Sh({\cal C}, J)$ induced by a comorphism of sites $F:({\cal D}, K)\to ({\cal C}, J)$ is a surjection if and only if $J=K^{F}$ (in the sense of Proposition \ref{propimagetopology}), that is, if a sieve $S$ on an object $c\in {\cal C}$ satisfies the property that for every object $d$ of $\cal D$ and arrow $x:F(d)\to c$ in $\cal C$, there exists a $K$-covering sieve $T$ on $d$ such that $F(T)\subseteq x^{\ast}(S)$ then $S$ is $J$-covering. This condition implies that $F$ is $J$-dense and is equivalent to it if $F$ is cover-preserving.	
\end{proposition}

\begin{proofs}
	By Proposition \ref{propcoverreflecting}, $C_{F}$ is a surjection if and only the functor $A_{F}$ is cover-reflecting, that is, for any sieve $S$ on an object $c$, if the family of arrows $A_{F}(f)$ for $f\in S$ is epimorphic then $S$ is $J$-covering. Now, since a family of morphisms in $[{\cal D}^{\textup{op}},\Set]$ is sent by the associated sheaf functor $a_K:[{\cal D}^{\textup{op}},\Set] \to \Sh({\cal D}, K)$ to an epimorphic family if and only if it is $K$-locally surjective (cf. Corollary III.6 \cite{MM}), we obtain the following criterion: the family $\{A_{F}(f) \mid f\in S \}$ is epimorphic if and only if for every object $d$ of $\cal D$ and arrow $x:F(d)\to c$ there exists a $K$-covering sieve $T$ on $d$ such that for every $g\in T$, $x\circ F(g)\in S$, that is, if and only if $K^{F}\subseteq J$. From this, in light of Remark \ref{remcovlift}, the firt part of our thesis follows immediately. 
	
	It remains to show that the condition $J=K^{F}$ implies that $F$ is $J$-dense and is equivalent to it if $F$ is cover-preserving. 
	The fact that if $J=K^{F}$ then $F$ is $J$-dense follows from the fact that, for any $c\in {\cal C}$, the sieve $S_{c}$ on $c$ generated by the arrows from objects of the form $F(d)$ (for $d\in {\cal D}$) to $c$ is $K^{F}$-covering, by definition of $K^{F}$, and therefore $J$-covering. Conversely, let us show that if $F$ is cover-preserving and $J$-dense then $K^{F}\subseteq J$ (equivalently, $K^{F}=J$). Given a sieve $S$ on an object $c\in {\cal C}$ satisfying the property that for every object $d$ of $\cal D$ and arrow $x:F(d)\to c$ in $\cal C$, there exists a $K$-covering sieve $T$ on $d$ such that $F(T)\subseteq x^{\ast}(S)$, we want to prove that $S$ is $J$-covering. Since $F$ is cover-preserving, for any arrow $x:F(d)\to c$ in $\cal C$, $x^{\ast}(S)\in J(\textup{dom}(x))$; the fact that $S\in J(c)$ thus follows from the fact that $F$ is $J$-dense by the transitivity axiom for Grothendieck topologies.
\end{proofs}

\begin{remarks}\label{remsurjectioncomorphismofsites}
	\begin{enumerate}[(a)]
		\item If $F$ is $K$-full then the condition $J=K^{F}$ implies that $F$ is cover-preserving. Indeed, for any $R\in K(d')$, the sieve $\langle F(R)\rangle  $ is $J$-covering, since, by the $K$-fullness of $F$, for any arrow $x:F(d)\to F(d')$ in $\cal C$ there are a set $I$ and arrows $f_{i}:d_{i}\to d, g_{i}:d_{i}\to d'$ (for $i\in I$) such that the family $\{f_{i}:d_{i}\to d \mid i\in I\}$ is $K$-covering and $x\circ F(f_{i})=F(g_{i})$, and hence the sieve $T:=\{f_{i}\circ h \mid i\in I, h\in  g_{i}^{\ast}(R) \}$ is $K$-covering and satisfies the property $F(T)\subseteq x^{\ast}(S)$. 
		
		\item If $J$ and $K$ are the trivial topologies then $F$ is clearly cover-preserving and the condition that $F$ be $J$-dense is equivalent to the requirement that every object of $\cal D$ should be a retract of an object in the image of $F$. In fact, this was also proved with other means in Example A4.2.7(b) \cite{El}.	 
	\end{enumerate}
\end{remarks}

\subsection{Inclusions}

To characterize the comorphisms of sites which give rise to a geometric inclusion, we shall use the criterion provided by Corollary \ref{corinclusionflatfunctor}. 

The following lemma expresses certain properties of arrows in the image of the functor $A_{F}$ introduced above in terms of functional relations associated with them.
 
\begin{lemma}\label{lemprop}
	Let $c, c'\in {\cal C}$, $\xi:A_{F}(c)=a_{K}(\textup{Hom}_{\cal C}(F(-),c))\to A_{F}(c')=a_{K}(\textup{Hom}_{\cal C}(F(-),c'))$ an arrow in $\Sh({\cal D}, K)$ and $R_{\xi}$ the functional relation from $\textup{Hom}_{\cal C}(F(-),c)$ to $\textup{Hom}_{\cal C}(F(-),c')$ associated with it as in Remark \ref{remfunctionalrelations}(c). Then 
	\begin{enumerate}[(i)]
		\item For any arrows $f:c''\to c$ and $g:c''\to c'$, $\xi\circ A_{F}(f)=A_{F}(g)$ if and only if for every $x:F(e)\to c''$, $(f\circ x, g\circ x)\in R_{\xi}$. 
		
		\item There exist a family of arrows $f_{i}:c_{i}\to c$ in $\cal C$ which is sent by $A_{F}$ to an epimorphic family and a family of arrows $g_{i}:c_{i}\to c'$ in $\cal C$ such that $\xi\circ A_{F}(f_{i})=A_{F}(g_{i})$ for all $i$ if and only if there exists a family $(u_{j}, v_{j})$ of elements of $R_{\xi}$ such that the family $\{A_{F}(u_{j})\}$ is epimorphic and for any arrow $s:F(e)\to \textup{dom}(u_{j})=\textup{dom}(v_{j})$ in $\cal C$, $(u_{j}\circ s, v_{j}\circ s)\in R_{\xi}$.
	\end{enumerate}
\end{lemma}

\begin{proofs}
	(i) By Remark \ref{remfunctionalrelations}(c), $(f\circ x, g\circ x)\in R_{\xi}$ if and only if $\xi \circ \eta^{c} (f\circ x)=\eta^{c'}(g\circ x)$, where $\eta^{c}:\textup{Hom}_{\cal C}(F(-),c) \to a_{K}(\textup{Hom}_{\cal C}(F(-),c))$ and $\eta^{c'}:\textup{Hom}_{\cal C}(F(-),c') \to a_{K}(\textup{Hom}_{\cal C}(F(-),c'))$ are the canonical unit arrows. But $\xi \circ \eta^{c} (f\circ x)=(\xi\circ A_{F}(f)\circ \eta^{c''})(x)$ and $\eta^{c'}(g\circ x)=(A_{F}(g)\circ \eta^{c''})(x)$, where $\eta^{c''}:\textup{Hom}_{\cal C}(F(-),c'') \to a_{K}(\textup{Hom}_{\cal C}(F(-),c''))$ is the canonical unit arrow.
	
	(ii) Let us suppose that there exist a family of arrows $f_{i}:c_{i}\to c$ in $\cal C$ which is sent by $A_{F}$ to an epimorphic family and a family of arrows $g_{i}:c_{i}\to c'$ in $\cal C$ such that $\xi\circ A_{F}(f_{i})=A_{F}(g_{i})$ for all $i$. Then we may take as arrows $(u_{j}, v_{j})$ those of the form $(f_{i}\circ x, g_{i}\circ x)$ for $x$ an arrow $F(e)\to c_{i}$ for some $i\in I$. By (i), every pair of the form $(f_{i}\circ x, g_{i}\circ x)$ is in $R_{\xi}$, so it remains to show that the family of arrows $\{A_{F}(u_{j})\}$ is epimorphic (equivalently, for any arrow $u:F(e)\to c$ there exists a $K$-covering sieve $T$ on $e$ such that $u\circ F(t)$ factors through some $u_{j}$ for every $t\in T$). But this is clear since, the family $\{A_{F}(f_{i}) \mid i\in I \}$ being epimorphic, there is a $K$-covering sieve $T_{u}$ on $e$ such that, for every $t\in T_{u}$, $u\circ F(t)$ factors through some $f_{i}$, that is, is of the form $f_{i}\circ x$ for some arrow $x$. 
	
	Conversely, if there exists a family $(u_{j}, v_{j})$ of elements of $R_{\xi}$ such that the family $\{A_{F}(u_{j})\}$ is epimorphic and for any arrow $s:F(e)\to \textup{dom}(u_{j})=\textup{dom}(v_{j})$ in $\cal C$, $(u_{j}\circ s, v_{j}\circ s)\in R_{\xi}$ then by (i) $\xi\circ A_{F}(u_{j})=A_{F}(v_{j})$ for each $j$. So we may take as family of arrows $(f_{i}, g_{i})$ precisely the family of arrows $(u_{j}, v_{j})$.      
\end{proofs}

\begin{proposition}\label{propsatisfactionconditionii}
	Let  $F:({\cal D}, K)\to ({\cal C}, J)$ be a comorphism of sites. Then the flat functor $A_{F}:{\cal C}\to \Sh({\cal D}, K)$ satisfies condition (ii) of Corollary \ref{corinclusionflatfunctor} if and only if for every $K$-functional relation $R$ from $\textup{Hom}_{\cal C}(F(-), c)$ to $\textup{Hom}_{\cal C}(F(-), c')$ in $[{\cal D}^{\textup{op}}, \Set]$ there exists a family $(f_{i}, g_{i})$ of arrows $f_{i}:\textup{dom}(f_{i}) \to c$ and $g_{i}:\textup{dom}(g_{i}) \to c'$ with common domain $\textup{dom}(f_{i})=\textup{dom}(g_{i})$ such that
	\begin{itemize}
		\item for any arrow $x:F(e)\to c$ there exists a $K$-covering sieve $T$ on $e$ such that, for every $t\in T$, $x\circ F(t)$ factors through some $f_{i}$;
		
		\item for any arrow $x:F(e)\to \textup{dom}(f_{i})=\textup{dom}(g_{i})$ in $\cal C$, $(f_{i}\circ x, g_{i}\circ x)\in R$.
	\end{itemize}
In particular, this condition is satisfied if $F$ is $K$-full.   

If $F$ is cover-preserving then the conjunction of the two conditions above is equivalent to the requirement that for any $(f,g)\in R$ and any arrow $x:F(e)\to \textup{dom}(f)=\textup{dom}(g)$ in $\cal C$, $(f\circ x, g\circ x)\in R$.
\end{proposition}

\begin{proofs}
	Recall that condition (ii) of Corollary \ref{corinclusionflatfunctor} states that for every $c, c'\in {\cal C}$ and any arrow $\xi:A_{F}(c)\to A_{F}(c')$ in $\Sh({\cal D}, K)$, there exist a family of arrows $f_{i}:c_{i}\to c$ in $\cal C$ which is sent by $A_{F}$ to an epimorphic family and a family of arrows $g_{i}:c_{i}\to c'$ in $\cal C$ such that $\xi\circ A_{F}(f_{i})=A_{F}(g_{i})$ for all $i$.
	
	By Theorem \ref{thmfunctionalrelations}, the arrows $$\xi:A_{F}(c)=a_{K}(\textup{Hom}_{\cal C}(F(-),c))\to a_{K}(\textup{Hom}_{\cal C}(F(-),c'))=A_{F}(c')$$ in $\Sh({\cal D}, K)$ correspond to the $K$-functional relations $R_{\xi}$ from $\textup{Hom}_{\cal C}(F(-),c)$ to $\textup{Hom}_{\cal C}(F(-),c')$ in $[{\cal D}^{\textup{op}},\Set]$. 
	
	The first part of the proposition thus follows immediately from Lemma \ref{lemprop} by observing that a family of arrows $\{A_{F}(u_{j})\}$ is epimorphic if and only if for any arrow $x:F(e)\to c$ there exists a $K$-covering sieve $T$ on $e$ such that, for every $t\in T$, $x\circ F(t)$ factors through some $u_{j}$. 
	
	If $F$ is $K$-full then, for any arrow $\xi:A_{F}(c) \to A_{F}(c')$, the family of pairs of arrows in $R_{\xi}$ satisfies both conditions in the statement of the Proposition. The first condition is satisfied since, $R_{\xi}$ being a $K$-functional relation, the family of arrows $\{A_{F}(x) \mid x\in \pi_{1}(R_{\xi}) \}$, where $\pi_{1}$ is the canonical projection $\textup{Hom}_{\cal C}(F(-),c)\times \textup{Hom}_{\cal C}(F(-),c') \to \textup{Hom}_{\cal C}(F(-),c)$, is epimorphic in $\Sh({\cal D}, K)$, while the second holds by the following argument. Given $(f,g)\in R_{\xi}$ and $x:F(e)\to  \textup{dom}(f)=\textup{dom}(g)$, the $K$-fullness of $F$ implies that there exist a $K$-covering family $T$ on $e$ and for each $t\in T$ an arrow $s_{t}:\textup{dom}(t)\to e$ such that $x\circ F(t)=F(s_{t})$. Now, by the $K$-closure property of $R_{\xi}$, the condition $(f\circ x, g\circ x)\in R_{\xi}$ is equivalent to $(f\circ x\circ F(t),  g\circ x\circ F(t))\in R_{\xi}$ for every $t\in T$. But $(f\circ x\circ F(t),  g\circ x\circ F(t))=(f\circ F(s_{t}), g\circ F(s_{t}))$, and  $(f\circ F(s_{t}), g\circ F(s_{t}))\in R_{\xi}$ by the functoriality of $R_{\xi}$. 
	
	It remains to prove the last part of the proposition, namely that, provided that $F$ is cover-preserving, if $R$ satisfies the two conditions above then $R$ is closed under composition on the right with arbitrary arrows in $\cal C$. Note that the converse direction follows from the fact that if $R$ satisfies this property then, by taking the family $\{(f_i, g_i)\}$ to consist of the pairs in $R$, the first condition in the proposition is automatically satisfied by the $K$-functionality of $R$. So, let us suppose that $(f:F(e)\to c, g:F(e)\to c')\in R$ and that $x:F(a)\to F(e)$ is an arrow in $\cal C$; we want to show that $(f\circ x, g\circ x)\in R$. By the first condition, there is a $K$-covering sieve $T$ on $e$ such that, for each $t\in T$, $f\circ F(t)$ factors through some $f_{i_{t}}$, say, $f\circ F(t)=f_{i_{t}}\circ u_{i}^{t}$. Now, by the functoriality of $R$, $(f, g)\in R$ implies that $(f\circ F(t), g\circ F(t))\in R$. On the other hand, the second condition in the proposition implies that $(f_{i_{t}}\circ u_{i}^{t}, g_{i_{t}}\circ u_{i}^{t})\in R$. So, the $K$-functionality of $R$ implies that $g\circ F(t) \equiv_{K} g_{i_{t}}\circ u_{i}^{t}$. Since $F$ is cover-preserving, the sieve $\langle F(T)\rangle $ generated by the arrows $F(t)$ for $t\in T$ is $J$-covering; therefore its pullback along $x$ is also $J$-covering and hence by the covering-lifting property there is a $K$-covering sieve $H$ on $a$ such that for each $h\in H$, $x\circ F(h)=F(t_{h})\circ z_{h}$ for some arrow $t_{h}\in T$ and some arrow $z_{h}:F(\textup{dom}(h))\to F(\textup{dom}(t_{h}))$. Now, by the $K$-closure of $R$, $(f\circ x, g\circ x)\in R$ if and only if for every $h\in H$, $(f\circ x\circ F(h), g\circ x\circ F(h))\in R$. But $(f\circ x\circ F(h), g\circ x\circ F(h))=(f_{i_{t_{h}}}\circ u_{i}^{t_{h}}\circ z_{h}, g\circ F(t_{h})\circ  z_{h})$. Now, we have $g\circ F(t_{h}) \equiv_{K} g_{i_{t_{h}}}\circ u_{i}^{t_{h}}$, whence $g\circ F(t_{h}) \circ z_{h} \equiv_{K} g_{i_{t_{h}}}\circ u_{i}^{t_{h}}\circ z_{h}$. Since $R$ is a $K$-functional relation, the condition $(f_{i_{t_{h}}}\circ u_{i}^{t_{h}}\circ z_{h}, g\circ F(t_{h})\circ  z_{h})\in R$ is equivalent to $(f_{i_{t_{h}}}\circ u_{i}^{t_{h}}\circ z_{h}, g_{i_{t_{h}}}\circ u_{i}^{t_{h}}\circ z_{h})\in R$ (cf. Proposition \ref{propfuncionalassignment}), and this is satisfied by the second condition in the proposition.     
\end{proofs}

For investigating the satisfaction of condition (i) of Corollary \ref{corinclusionflatfunctor}, we need the following lemma. Notice that if $F:({\cal D}, K)\to ({\cal C}, J)$ is a comorphism of sites then for any $d\in {\cal D}$ we have an arrow $\chi_{d}:l'(d)\to A_{F}(F(d))$ corresponding, via the associated sheaf adjuntion and the Yoneda embedding, to the element of $A_{F}(F(d))(d)$ given by the identity arrow on $F(d)$.

\begin{lemma}\label{lemmacomorphisminclusion}
Let $F:({\cal D}, K)\to ({\cal C}, J)$ be a comorphism of sites. Then for any object $d$ of $\cal D$ and any arrow $g:d'\to d$ in $\cal D$, the arrows $\xi:A_{F}(F(d'))\to l'(d)$ such that $\xi\circ \chi_{d'}=l'(g)$ can be identified with the relations $R$ from $\textup{Hom}_{\cal C}(F(-),F(d'))$ to $\textup{Hom}_{\cal D}(-,d)$ which assign to each object $e$ of $\cal C$ a collection $R(e)$ of pairs of arrows $(x:F(e)\to F(d'), y:e\to d)$ in such a way that the following properties are satisfied:
		\begin{enumerate}[(i)]
		\item for any $\xi:e'\to e$ in $\cal D$, if $(x, y)\in R$ then $(x\circ F(\xi), y\circ \xi)\in R$;	
		
		\item  for any $e\in {\cal D}$ and any $(x, y)\in \textup{Hom}_{\cal C}(F(e), F(d'))\times \textup{Hom}_{\cal C}(e,d)$, if $\{\xi:e'\to e \mid\ (x\circ F(\xi), y\circ \xi)\in R(e')\}\in K(e)$ then $(x,y)\in R(e)$;
			
		\item for any $e\in {\cal D}$, if $(x, y), (x', y')\in R$ and $x=x'$ then $y\equiv_{K} y'$;  
				 			
		\item for any $e\in {\cal D}$ and any arrow $x:F(e)\to F(d')$ in $\cal C$, $\{t:e'\to e \mid \exists y:e'\to d, \,(x\circ F(t), y)\in R(e') \}\in K(e)$;
		
		\item for any arrow $f:e\to d'$ in $\cal D$, $(F(f), g\circ f)\in R(e)$.	
\end{enumerate}
In particular, if $F$ is $K$-full and $K$-faithful then for every $d\in {\cal D}$ the arrow $\chi_{d'}:l'(d)\to A_{F}(F(d))$ is an isomorphism.
\end{lemma}

\begin{proofs}
The first part of the lemma follows from the proof of Theorem \ref{thmfunctionalrelations}. It thus remains to show that if $F$ is $K$-full and $K$-faithful then there exists a relation which verifies the conditions in the statement of the lemma and is inverse to $\chi_{d'}$. We can define a relation $R_{F}$ on pairs of arrows $(x:F(e)\to F(d), y:e\to d)$ by stipulating that $(x,y)\in R_{F}$ if and only if $F(y)=x$. If $F$ is $K$-faithful and $K$-full then the relation $R_{F}$ clearly satisfies all the conditions of the lemma for the family of objects $\{c_{i} \mid i\in I\}$ given by the singleton $\{F(d)\}$ except for the closure condition (ii), which, by Remark \ref{remfunctionalrelations}(b), can be made to hold (without affecting the satisfaction of the other conditions) by taking its $K$-closure. The corresponding arrow $A_{F}(F(d))\to l'(d)$ (via the bijection of Theorem \ref{thmfunctionalrelations}) is actually inverse to the arrow $\chi_{d}$ (cf. Proposition \ref{proplocalfullnessandfaithfulness}), as required. 	
\end{proofs}

\begin{proposition}\label{propinclusioncomoprhismofsites}
	The geometric morphism $C_{F}:\Sh({\cal D}, K)\to \Sh({\cal C}, J)$ induced by a comorphism of sites $F:({\cal D}, K)\to ({\cal C}, J)$ is an inclusion if and only if $F$ satisfies the condition of Proposition \ref{propsatisfactionconditionii} and for every object $d$ of $\cal D$ there exist a $K$-covering family $\{g_{i}:d_{i}\to d \mid i\in I\}$ and arrows $\xi_{i}:A_{F}(F(d_{i}))\to l'(d)$ in $\Sh({\cal D}, K)$ such that $\xi_{i}\circ \chi_{d_{i}}=l'(g_{i})$ (equivalently, relations $R_{i}$ from $\textup{Hom}_{\cal C}(F(-),F(d_{i}))$ to $\textup{Hom}_{\cal D}(-,d)$  satisfying the conditions of Lemma \ref{lemmacomorphisminclusion}). 
	
	In particular, if $F$ is $K$-full and $K$-faithful then the geometric morphism $C_{F}$ is an inclusion. Conversely, if $F$ is cover-preserving and $C_{F}$ is an inclusion then $F$ is $K$-full.
	
	If $F$ is $(K, J)$-continuous (in the sense of Definition \ref{defcanonicalfunctor}) then $C_{F}$ is an inclusion if and only if $F$ is $K$-full and $K$-faithful.
\end{proposition} 
 
\begin{proofs}
By Corollary \ref{corinclusionflatfunctor}, $C_{F}$ is an inclusion if and only if conditions (i) and (ii) of the Corollary are satified by the flat functor $A_{F}$. Condition (ii) is equivalent to the condition of Proposition \ref{propsatisfactionconditionii}. On the other hand, condition (i) is satisfied by $A_{F}$ if and only if for each $d\in {\cal D}$, there exists a set $\{c_{j}\}$ of objects of $\cal C$ and an epimorphic family of arrows $\xi'_{j}:A_{F}(c_{j})\to l'(d)$. But, by using Proposition \ref{propfunctionalrelationmonoepi}(ii), one can easily see that this amounts to the existence of a $K$-covering family $\{g_{i}:d_{i}\to d \mid i\in I\}$ on $d$ and for each $i\in I$ an arrow $u_{i}:F(d_{i})\to c_{j_{i}}$ in $\cal C$ (for some $j_{i}$) such that $l'(g_{i})=\xi'_{j_{i}}\circ A_{F}(u_{i})\circ \chi_{d_{i}}$. So the composite of the arrow $\xi_{i}:=\xi'_{j_{i}}\circ A_{F}(u_{i}):A_{F}(F(d_{i}))\to l'(d)$ with $\chi_{d_{i}}$ is equal to $l'(g_{i})$.  

If $F$ is $K$-full and $K$-faithful then by Proposition \ref{propsatisfactionconditionii} $F$ satisfies the condition in it and by Lemma \ref{lemmacomorphisminclusion} for every $d\in {\cal D}$ the arrow $\chi_{d'}:l'(d)\to A_{F}(F(d))$ is an isomorphism whence we can take the trivial covering family on each $d$ to make the condition of the proposition satisfied. 

Let us now show that if $C_{F}$ is an inclusion and $F$ is cover-preserving then $F$ is $K$-full.
Let $f, f':e\to d$ be arrows in $\cal D$ such that $F(f)=F(f')$. Let $\{g_{i}:d_{i}\to d \mid i\in I\}$ be a $K$-covering family and $\xi_{i}:A_{F}(F(d_{i}))\to l'(d)$ be arrows in $\Sh({\cal D}, K)$, corresponding to $K$-functional relations $R_{i}$ from $\textup{Hom}_{\cal C}(F(-),F(d_{i}))$ to $\textup{Hom}_{\cal D}(-,d)$, such that $\xi_{i}\circ \chi_{d_{i}}=l'(g_{i})$. Let us first prove that, for any $i\in I$, we have $\chi_{d}\circ \xi_{i}=A_{F}(F(g_{i}))$. The $K$-functional relation corresponding to this arrow is the composite $R_{i}'$ of $R_{i}$ with the relation $R_{\chi_{d}}$ associated with $\chi_{d}$. Since $(1_{F(d_{i})}, g_{i})\in R_{i}$ (this follows from the fact that $R_{i}$ satisfies property (v) of Lemma \ref{lemmacomorphisminclusion} with respect to the arrow $g_{i}$) and $(g_{i}, F(g_{i})\in R_{\chi_{d}}$, we have $(1_{F(d_{i})}, F(g_{i}))\in R_{i}'$. So, by the last part of Proposition \ref{propsatisfactionconditionii}, for any arrow $x:F(e)\to F(d_{i})$, $(x, F(g_{i})\circ x)\in R_{i}'$; in other words, $\chi_{d}\circ \xi_{i}=A_{F}(F(g_{i}))$, as required. 

The sieve $\langle \{F(g_i) \mid i\in I\}\rangle $ generated by the arrows $F(g_i)$ for $i\in I$ is $J$-covering; so, for any arrow $y:F(e)\to F(d)$ in $\cal C$, its pullback along $y$ is also $J$-covering, whence by the covering-lifting property there is a $K$-covering sieve $H$ on $e$ such that for each $h\in H$, $y\circ F(h)=F(g_{i_{h}})\circ z_{h}$ for some $i_{h}\in I$ and some arrow $z_{h}:F(\textup{dom}(h))\to F(d_{i_{h}})$. Now, since $R_{i}$ is $K$-functional, there is a $K$-covering sieve $T_{h}$ on $\textup{dom}(h)$ and for each $t\in T_{h}$ an arrow $a_{t}:\textup{dom}(t)\to d$ such that $(z_{h}\circ F(t), a_{t})\in R_{i}$. Then, from the fact (proved above) that $\chi_{d}\circ \xi_{i}=A_{F}(F(g_{i}))$ it follows that $F(a_{t})$ is $K$-locally equal to $F(g_{i_{h}})\circ z_{h}\circ F(t)$, that is, there is a $K$-covering sieve $W_{t}$ on $\textup{dom}(t)$ such that $F(a_{t})\circ F(w)=F(g_{i_{h}})\circ z_{h}\circ F(t)\circ F(w)=y\circ F(h)\circ F(t)\circ F(w)$ for each $w\in W_{t}$. Therefore the $K$-covering sieve $U$ given by the multiomposite of $H$, the $T_{h}$'s (for $h\in H$) and the $W_{t}$'s (for $t\in T_{h}$) satisfies the property that for any $u\in U$ there is an arrow $a'_{u}:\textup{dom}(u)\to d$ such that $y\circ F(u)=F(a'_{u})$ for each $u\in U$ (indeed, if $u=h\circ t\circ w$ then we can take $a'_{u}=a_{t}\circ w$). So $F$ is $K$-full, as desired.   

It remains to prove the last part of the proposition. We have already shown that the condition on $F$ being $K$-faithful and $K$-full is a sufficient condition for $C_{F}$ to be an inclusion. We want to prove that this condition is also necessary if $F$ is $(K, J)$-continuous. But this follows from the fact that, by Corollary \ref{correstrictioncomorphism}, we have a commutative diagram

		$$
\xymatrix{
	{\mathcal C} \ar[r]^F \ar[d]_{l'} &{\cal D} \ar[d]^{l} \\
	\Sh({\cal D}, K) \ar[r]^{(C_{F})_{!}} & \Sh({\cal C}, J)
}
$$ 
in light of Proposition \ref{proplocalproperties}(i)-(ii) and Remark \ref{remcoveringliftingstrictequality}. 
\end{proofs}

\begin{remark}\label{remarkinclusion}
	If $J$ and $K$ are the trivial topologies then one may deduce from Proposition \ref{propinclusioncomoprhismofsites} the following criterion for the essential geometric morphism $C_{F}:[{\cal D}^{\textup{op}},\Set]\to [{\cal C}^{\textup{op}},\Set]$ induced by $F$ to be an inclusion: $C_{F}$ is an inclusion if and only if $F$ is full and faithful. 
	
	Indeed, the `if' direction follows immediately from the proposition, while the converse one can be proved as follows. The fullness of $F$ follows from the last part of the proposition, while its faithfulness follows from the fact that the condition in the proposition implies that, for every $d\in {\cal D}$, there is an arrow $\xi:A_{F}(F(d))\to l'(d)$ such that $\xi\circ \chi_{d}=1_{l'(d)}$. 
	
	An alternative way of proving that if $C_{F}$ is an inclusion then $F$ is full and faithful is to observe that, the inverse image functor of the geometric morphism $C_{F}$ having adjoints on both sides, the right adjoint is full and faithful if and only if the left adjoint is; but this latter condition implies that $F$ is full and faithful, since $F$ can be recovered from this functor as its restriction to the representables.
\end{remark}

The surjection-inclusion factorization of the geometric morphism induced by a comorphism of sites admits a particularly simple description when the latter is cover-preserving, as shown by the following proposition. 

\begin{proposition}\label{propsurinclfactorizationcomorphismofsites}
	Let $F:({\cal D}, K)\to ({\cal C}, J)$ be a comorphism of sites which is cover-preserving. Then the surjection-inclusion factorization of the geometric morphism $C_{F}:\Sh({\cal D}, K)\to \Sh({\cal C}, J)$ induced by $F$ can be identified with $C_{i}\circ C_{F'}$, where $F'$ is the functor $F$ regarded as a comorphism of sites from $({\cal D}, K)$ to the site $({\cal C}', J')$, where ${\cal C}'$ is the full subcategory of $\cal C$ on the objects in the image of $F$ and $J'$ is the smallest Grothendieck topology on ${\cal C}'$ making the inclusion $i$ of ${\cal C}'$ into ${\cal C}$ a comorphism of sites to $({\cal C}, J)$.	
\end{proposition}

\begin{proofs}
	Since $F$ is cover-preserving, by definition of the Grothendieck topology $J'$, $F'$ is also cover-preserving (cf. Corollary \ref{corliftcomorphismsofsites}). It is also surjective on objects, so, by Proposition \ref{propsurjectioncomorphismofsites}, $C_{F'}$ is a surjection. On the other hand, $i$ is full and faithful, so, by Proposition \ref{propinclusioncomoprhismofsites}, $C_{i}$ is an inclusion. By the uniqueness (up to equivalence) of the surjection-inclusion factorization, it follows that $C_{i}\circ C_{F'}$ is `the' surjection-inclusion factorization of $C_{F}$, as required.
\end{proofs}

\begin{remark}
	Proposition \ref{propsurinclfactorizationcomorphismofsites} generalizes the result (cf. Example A4.2.12(b) \cite{El}) describing the surjection-inclusion factorization of the essential geometric morphism induced by an arbitrary functor between small categories. 
\end{remark}

\subsection{Localic morphisms}

We shall now characterize the comorphisms of sites which induce localic geometric morphisms.

\begin{lemma}\label{lemmacomorphismlocalic}
	Let $F:({\cal D}, K)\to ({\cal C}, J)$ be a comorphism of sites. Then for any object $d$ of $\cal D$ and any arrow $g:d'\to d$ in $\cal D$, the arrows $\xi$ to $l'(d)$ whose domain is a subobject $s:S\mono A_{F}(F(d'))$ through which $\chi_{d'}$ factors as $\chi_{d'}=s\circ \overline{\chi_{d'}}$, and such that $\xi\circ \overline{\chi_{d'}}=l'(g)$, can be identified with the relations $R$ from $\textup{Hom}_{\cal C}(F(-),F(d'))$ to $\textup{Hom}_{\cal D}(-,d)$ satisfying all the conditions of Lemma \ref{lemmacomorphisminclusion} but condition (iv). 	

The functor $F$ is $K$-faithful if and only if for every object $d$ of $\cal D$, the arrow $\chi_{d}:l'(d)\to A_{F}(F(d))$ is a monomorphism. 
\end{lemma}

\begin{proofs}
	The first part of the lemma follows from the proof of Theorem \ref{thmfunctionalrelations}, so it remains to prove the second part. We notice that $\chi_{d}=a_{K}(\gamma_{d})$, where $\gamma_{d}$ is the arrow $y_{\cal D}(d)\to \textup{Hom}_{\cal C}(F(-), F(d))$ corresponding to the element $1_{F(d)}$ of $\textup{Hom}_{\cal C}(F(d), F(d))$ via the Yoneda embedding. So by Lemma \ref{lemmalift}(i) $\chi_{d}$ is a monomorphism if and only if for every $x,x':e\to d$ such that $F(x)=F(x')$, $x\equiv_{K} x'$, that is, if and only if $F$ is $K$-faithful.  
\end{proofs}

\begin{proposition}\label{proplocaliccomoprhismofsites}
	The geometric morphism $C_{F}:\Sh({\cal D}, K)\to \Sh({\cal C}, J)$ induced by a comorphism of sites $F:({\cal D}, K)\to ({\cal C}, J)$ is localic if and only if for every object $d$ of $\cal D$ there exist a $K$-covering sieve $\{g_{i}:d_{i}\to d \mid i\in I\}$ on $d$ and relations $R_{i}$ from $\textup{Hom}_{\cal C}(F(-),F(d_{i}))$ to $\textup{Hom}_{\cal D}(-,d)$  satisfying the conditions of Lemma \ref{lemmacomorphismlocalic}.
	
	In particular, if $F$ is $K$-faithful then the geometric morphism $C_{F}$ is localic.
\end{proposition} 

\begin{proofs}
	The proof uses Lemma \ref{lemmacomorphismlocalic} in a way analogous to that in which Lemma \ref{lemmacomorphisminclusion} is used in the proof of Proposition \ref{propinclusioncomoprhismofsites}.
	In fact, stating that, for each $d\in {\cal D}$, there is a set $\{c_{j}\}$ of objects of $\cal C$ and an epimorphic family of arrows $\{\xi_{j}:a_{K}(S_{j})\to l'(d)\}$, where each $S_{j}$ is a $K$-closed sieve on $c_{j}$ yielding a canonical subobject $s_{j}:a_{K}(S_{j})\mono A_{F}(c_{j})$, amounts to saying that there are a $K$-covering family $\{g_{i}:d_{i}\to d \mid i\in I\}$ and for each $i$ an element $j_{i}$
 	and an arrow $u_{j_{i}}:F(d_{i})\to c_{j_{i}}$ in $S_{j_{i}}$ such that $(u_{j_{i}}, g_{i})\in R_{j_{i}}$ (where $R_{j_{i}}$ is the functional relation corresponding to $\xi_{j_{i}}$). For each $i\in I$, by considering the pullback
 		$$
 	\xymatrix{
 		a_{K}(S'_{j_{i}}) \ar[rr]^{s'_{j_{i}}} \ar[d] & & A_{F}(F(d_{i})) \ar[d]^{A_{F}(u_{j_{i}})} \\
 		a_{K}(S_{j_{i}}) \ar[rr]^{s_{j_{i}}} & &  A_{F}(c_{j_{i}}) 
 	}
 	$$
 	 of each subobject $s_{j_{i}}:a_{K}(S_{j_{i}})\mono A_{F}(c_{j_{i}})$ along $A_{F}(u_{j_{i}})$, we obtain a subobject $s'_{j_{i}}:a_{K}(S'_{j_{i}})\mono A_{F}(F(d_{i}))$, where $S'_{j_{i}}$ is $K$-closed sieve on $F(d_{i})$, through which $\chi_{d_{i}}$ factors (this follows from the universal property of the above pullback, since $s_{j_{i}}\circ \widetilde{u_{j_{i}}}=A_{F}(u_{j_{i}})\circ \chi_{d_{i}}$, where $\widetilde{u_{j_{i}}}$ is the arrow $l'(d_{i})\to a_{K}(S_{j_{i}})$ corresponding to the element $u_{j_{i}}$ of $S_{j_{i}}$) and for which the condition $(u_{j_{i}}, g_{i})\in R_{j_{i}}$ can be reformulated as the condition $\xi_{j_{i}}\circ s'_{j_{i}} \circ \overline{\chi_{d_{i}}}=l'(g_{i})$, where $\overline{\chi_{d_{i}}}$ denotes the factorization of $\chi_{d_{i}}$ through $s'_{j_{i}}:a_{K}(S'_{j_{i}})\mono A_{F}(F(d_{i}))$.   
\end{proofs}

\subsection{Hyperconnected morphisms}

\begin{proposition}\label{prophyperconnectedcomoprhismofsites}
	The geometric morphism $C_{F}:\Sh({\cal D}, K)\to \Sh({\cal C}, J)$ induced by a comorphism of sites $F:({\cal D}, K)\to ({\cal C}, J)$ is hyperconnected if and only if $F$ satisfies the property of Proposition \ref{propsurjectioncomorphismofsites} and for every object $c$ of $\cal C$ and any set $A$ of arrows of the form $x:F(d)\to c$ (for an object $d$ of $\cal D$) which is functorial (in the sense that if $x\in A$ then $x\circ F(g)\in A$ for any arrow $g:d'\to d$ in $\cal D$) and $K$-closed (in the sense that for any $K$-covering sieve $T$ on $d$, if $x\circ F(t)\in A$ for every $t\in T$ then $x\in A$) there exists a ($J$-closed) sieve $S$ on $c$ such that 	
	\begin{equation*}\begin{split}
	A &= \{x:F(d)\to c \mid d\in {\cal D}\textup{, } \{t:\textup{dom}(t)\to d \mid x\circ F(t)\in S\}\in K(d)\}. 
	\end{split}\end{equation*}
	If particular, if $F$ is $K$-full this latter condition is satisfied.
\end{proposition} 

\begin{proofs}
The first part of the proposition follows immediately from Proposition \ref{prophyperconnected} in light of Remark \ref{remimageclosedsubobject}(a), so it remains to prove the second part. Given a set $A$ of arrows $F(d)\to c$ (for $d\in {\cal D}$) satisfying the property in Proposition \ref{prophyperconnectedcomoprhismofsites}, let us define $S$ to be the sieve on $c$ generated by it. We want to prove that $A= \{x:F(d)\to c \mid d\in {\cal D} \textup{ and } \{t:\textup{dom}(t)\to d \mid x\circ F(t)\in S\}\in K(d)\}.$ The inclusion $\subseteq$ is clear since $A\subseteq S$. To prove the converse one, since $A$ is functorial and $K$-closed, it clearly suffices to show that for any arrow $x:F(d)\to c$, if $x\in S$ then $x\in A$. If $x\in S$ then, by definition of $S$, there exists an arrow $y:F(d')\to c$ in $A$ and an arrow $\xi:F(d)\to F(d')$ such that $x=y\circ \xi$. Now, since $F$ is $K$-full, there exist a $K$-covering family $\{h:\textup{dom}(h)\to d \mid h\in R\}$ on $d$ and for each $h\in R$ an arrow $s_{h}:\textup{dom}(h)\to d'$ such that $\xi\circ F(h)=F(s_{h})$ for each $h\in R$. Since $A$ is $K$-closed, to show that $x\in A$ it suffices to show that $x\circ F(h)\in A$ for each $h\in R$. But $x\circ F(h)=y\circ \xi \circ F(h)=y\circ F(s_{h})$, which belongs to $A$ since $A$ is functorial.   
\end{proofs}

\begin{corollary}\label{corhyperconnectedpresheaves}
	Let $f:{\cal D}\to {\cal C}$ be a functor between essentially small categories. Then the geometric morphism $C_{F}:[{\cal D}^{\textup{op}}, \Set]\to [{\cal C}^{\textup{op}}, \Set]$ is hyperconnected if and only if $F$ is full and every object of $\cal D$ is a retract of an object in the image of $F$.
\end{corollary}

\begin{proofs}
	By Proposition \ref{prophyperconnectedcomoprhismofsites} and Remark \ref{remsurjectioncomorphismofsites}(b), $C_{F}$ is hyperconnected if and only if every object of $\cal D$ is a retract of an object in the image of $F$ and for every functorial set $A$ of arrows $x:F(d)\to c$ there exists a sieve $S$ on $c$ such that for any arrow $x:F(d)\to c$ in $\cal D$, $x\in A$ if and only if $x\in S$. Let us prove that this latter condition is satisfied if and only if $F$ is full.
	If $F$ is full, then we can get the condition satisfied by defining $S$ to be the sieve on $c$ generated by the arrows in $A$. Indeed, if $x\in A$ then $x\in S$ since $A\subseteq S$. Conversely, suppose that $x:F(d)\to c$ is in $S$; then we can write $x=y\circ g$, where $y:F(d')\to c$ is in $A$. But by the fullness of $F$, $g=F(\xi)$ for some arrow $\xi:d\to d'$ in $\cal D$ and hence $x=y\circ F(\xi)\in A$ by the functoriality of $A$. Conversely, supposing that $F$ satisfies the condition, we want to prove that it is full. Given $d\in {\cal D}$, let $A_{d}$ be the collection of all the arrows of the form $F(\xi)$ for an arrow $\xi$ of $\cal D$ with codomain $d$. If $A_{d}$ satisfies the condition then there is a sieve $S_{d}$ on $F(d)$ such that an arrow $F(d')\to F(d)$ lies in $S_{d}$ if and only if it lies in $A_{d}$. It follows that the identity $1_{F(d)}$ lies in $S_{d}$, so every arrow $F(d')\to F(d)$ in $\cal C$ (where $d'$ is an object of $\cal D$) lies in $S_{d}$ and hence in $A_{d}$; that is, it is of the form $F(\xi)$ for some $\xi$. 
\end{proofs}

\begin{remark}
Corollary \ref{corhyperconnectedpresheaves} generalizes and strengthens Example A4.6.9 \cite{El}, where it is observed that if $F$ is bijective on objects and full then the geometric morphism $C_{F}$ is hyperconnected.
\end{remark}

The following result describes the hyperconnected-localic factorization of the geometric morphism induced by a cover-preserving comorphism of sites; it generalizes the well-known result for presheaf toposes (Example A4.6.9 \cite{El}).

\begin{proposition}\label{prophyperlocfactorizationcomorphismofsites}
	Let $F:({\cal D}, K)\to ({\cal C}, J)$ be a comorphism of sites which is cover-preserving. Then the hyperconnected-localic factorization of the geometric morphism $C_{F}:\Sh({\cal D}, K)\to \Sh({\cal C}, J)$ induced by $F$ can be identified with $C_{\tilde{F}}\circ C_{\pi}$, where $\tilde{F}$ is the functor $F$ regarded as a comorphism of sites from the site $({\cal E}, L)$ whose underlying category $\cal E$ is the quotient of the category $\cal D$ by the congruence induced by $F$ and whose Grothendieck topology $L$ has as covering sieves the sieves whose inverse image under the canonical projection functor $\pi:{\cal D}\to {\cal E}$ is $K$-covering.	
\end{proposition}

\begin{proofs}
	Since the functor $\pi$ is cover-preserving and full on objects and arrows, by Propositions \ref{prophyperconnectedcomoprhismofsites} and \ref{propsurjectioncomorphismofsites}, $C_{\pi}$ is hyperconnected. On the other hand, since $\tilde{F}$ is faithful then $C_{\tilde{F}}$ is localic by Proposition \ref{proplocaliccomoprhismofsites}. Since we clearly have $F=\tilde{F}\circ \pi$, our thesis follows. 
\end{proofs}

\subsection{Equivalences of toposes}

Finally, we can deduce from Propositions \ref{proplocaliccomoprhismofsites} and \ref{prophyperconnectedcomoprhismofsites} a criterion for a comorphism of sites to induce an equivalence of toposes (by observing that a geometric morphism is an equivalence if and only if it is both localic and hyperconnected):

\begin{proposition}
The geometric morphism $C_{F}:\Sh({\cal D}, K)\to \Sh({\cal C}, J)$ induced by a comorphism of sites $F:({\cal D}, K)\to ({\cal C}, J)$ is an equivalence if and only if the following conditions are satisfied:
\begin{enumerate}[(i)]
	\item for every object $d$ of $\cal D$ there exist a $K$-covering sieve $\{g_{i}:d_{i}\to d \mid i\in I\}$ on $d$ and  relations $R_{i}$ from $\textup{Hom}_{\cal C}(F(-),F(d_{i}))$ to $\textup{Hom}_{\cal D}(-,d)$ satisfying all the conditions of Lemma \ref{lemmacomorphisminclusion} but condition (iv); 
	
	\item $J=K^{F}$, that is, if a sieve $S$ on an object $c\in {\cal C}$ satisfies the property that for every object $d$ of $\cal D$ and arrow $x:F(d)\to c$ there exists a $K$-covering sieve $T$ on $d$ such that for every $g\in T$, $x\circ F(g)\in S$ then $S$ is $J$-covering;  
	
	\item for every object $c$ of $\cal C$ and any set $A$ of arrows of the form $x:F(d)\to c$ for an object $d$ of $\cal D$ which is functorial (in the sense that if $x\in A$ then $x\circ F(g)\in A$ for any arrow $g:d'\to d$ in $\cal D$) and $K$-closed (in the sense that for any $K$-covering sieve $T$ on $d$, if $x\circ F(t)\in A$ for every $t\in T$ then $x\in A$) there exists a ($J$-closed) sieve $S$ on $c$ such that 	
	\begin{equation*}\begin{split}
	A &= \{x:F(d)\to c \mid d\in {\cal D}\textup{, } \{t:\textup{dom}(t)\to d \mid x\circ F(t)\in S\}\in K(d)\}. 
	\end{split}\end{equation*}
\end{enumerate}	
	In particular, if $F$ is $K$-full and $K$-faithful and $J$-dense then $C_{F}$ is an equivalence. 
\end{proposition}\qed

As shown by the following proposition, the condition on $F$ being $K$-full, $K$-faithful and $J$-dense is not only sufficient but also necessary for $C_{F}$ to be an equivalence if $F$ is $(K, J)$-continuous:

\begin{proposition}
	Let $F$ be a $(K, J)$-continuous (in the sense of Definition \ref{defcanonicalfunctor}(b)) comorphism of sites $({\cal D}, K)\to ({\cal C}, J)$. Then the geometric morphism $C_{F}:\Sh({\cal D}, K)\to \Sh({\cal C}, J)$ is an equivalence if and only if $F$ is $K$-full, $K$-faithful and $J$-dense. 
\end{proposition}	

\begin{proofs}
	By Proposition \ref{propcharJcanonical}(iii), every $(K, J)$-continuous functor is cover-preserving.
	
	Now, a geometric morphism is an equivalence if and only if it is both an inclusion and a surjection. Our thesis thus follows from Propositions \ref{propsurjectioncomorphismofsites} and \ref{propinclusioncomoprhismofsites}. 
\end{proofs}

The following result provides a criterion for a functor which is both a morphism and a comorphism of sites to induce an equivalence.

\begin{proposition}\label{propmorphismcomorphismequivalence}
	Let $F:{\cal D}\to {\cal C}$ be a functor which is both a comorphism of sites and a morphism of sites $({\cal D}, K)\to ({\cal C}, J)$, where $J$ and $K$ are Grothendieck topologies respectively on $\cal C$ and $\cal D$. Then the following conditions are equivalent:
	
	\begin{enumerate}[(i)]
		\item The geometric morphism $C_{F}$ is an equivalence with quasi-inverse $\Sh(F)$.
		
		\item $F$ is $K$-full and $J$-dense (equivalently, a dense morphism of sites $({\cal D}, K)\to ({\cal C}, J)$).
	\end{enumerate} 
\end{proposition}

\begin{proofs}
	(i) $\imp$ (ii) If $F$, as a morphism of sites $({\cal D}, K)\to ({\cal C}, J)$, induces an equivalence then, since $F$ has the covering-lifting property, by Corollary \ref{corequivalencecovliftproperty} $F$ is $K$-full and $J$-dense. 
	
	(ii) $\imp$ (i) Since the construction of the geometric morphism induced by a comorphism of sites is functorial, and any canonical geometric inclusion $\Sh({\cal A}, Z) \hookrightarrow [{\cal A}^{\textup{op}}, \Set]$ is induced by the identity functor on $A$, regarded as comorphism of sites $({\cal A}, Z)\to ({\cal A}, T)$ where $T$ is the trivial topology on $\cal A$, we have a commutative square
	$$
	\xymatrix{
		\Sh({\cal D}, K) \ar[rr]^{C_{F}} \ar[d] & & \Sh({\cal C}, J) \ar[d] \\
		[{\cal D}^{\textup{op}}, \Set] \ar[rr]^{E(F)} & &  [{\cal C}^{\textup{op}}, \Set], 
	}
	$$
	where $E(F)$ is the geometric morphism induced by $F$ (regarded as a comorphism of trivial sites) and the vertical arrows are the canonical geometric inclusions. Let $\Sh(F):\Sh({\cal C}, J) \to \Sh({\cal D}, K)$ be the geometric morphism induced by $F$ as a morphism of sites $({\cal D}, K)\to ({\cal C}, J)$. Since $F$ is $J$-dense and has the covering-lifting property, by Proposition \ref{proprecoverytopology}, $J=K^{F}$ and hence, by Lemma \ref{lemmafullfaithfulcoverreflecting}, $F$ is cover-reflecting. So $F$ is a dense morphism of sites (cf. Remark \ref{remmorphismsofsitesfaithfulness}) and hence $\Sh(F)$ is an equivalence. To prove that $g$ is a quasi-inverse to $C_{F}$, it is clearly enough to show that $C_{F}^{\ast}\circ {\Sh(F)}^{\ast}\circ l_{\cal D}\cong l_{\cal D}$. Now, $\Sh(F)^{\ast}\circ l_{\cal D}\cong l_{\cal C}\circ F$, so $C_{F}^{\ast}\circ {\Sh(F)}^{\ast}\circ l_{\cal D}\cong C_{F}^{\ast}\circ l_{\cal C}\circ F$. But $C_{F}^{\ast}\circ  l_{\cal C}\cong C_{F}^{\ast} \circ a_{J}\circ y_{\cal C}\cong a_{K}\circ E(F)^{\ast}\circ y_{\cal C}$, where the second isomorphism follows from the commutativity of the above square. Now, by Proposition \ref{proplocalfullnessandfaithfulness}, since $F$ is $K$-full and $K$-faithful, the canonical arrow $l_{\cal D}\cong a_{K}\circ y_{\cal D} \to a_{K}\circ E(F)^{\ast}\circ y_{\cal C}\circ F$ is an isomorphism, so $C_{F}^{\ast}\circ g^{\ast}\circ l_{\cal D}\cong l_{\cal D}$, as required.  
\end{proofs}

\begin{remark}
	For a comorhism of sites $F:({\cal D}, K)\to ({\cal C}, J)$ to induce an equivalence $C_{F}$ it is necessary that $F$ be $J$-dense (by Proposition \ref{propsurjectioncomorphismofsites}). So inducing an equivalence is not a sufficient condition on a morphism of sites  $F:({\cal D}, K)\to ({\cal C}, J)$ for it to be also a comorphism between the same sites (since there are weakly dense morphisms which are not dense); in fact, this is the case if and only if $F$ is $J$-dense and $K$-full (cf. Corollary \ref{corequivalencecovliftproperty}). 
\end{remark}

\subsection{Local morphisms}\label{sec:localmorphisms}

Recall that a geometric morphism $f:{\cal F}\to {\cal E}$ is said to be \emph{local} if $f_{\ast}$ has a fully faithful right adjoint. For different characterizations of the property of a geometric morphism to be local see, for instance, Theorem C3.6.1 \cite{El}.

Note that this definition makes sense more generally for a weak morphism of toposes (in the sense of section \ref{sec:weakmorphisms}).
 
\begin{theorem}
	Let $F:{\cal D}\to {\cal C}$ be a continuous comorphism of sites (also regarded as a weak morphism of sites) $({\cal D}, K)\to ({\cal C}, J)$. Then:
	\begin{enumerate}[(i)]
		\item The geometric morphism $C_{F}:\Sh({\cal D}, K)\to \Sh({\cal C}, J)$ is essential, and 
		\[
		(C_{F})_{!}\cong \Sh(F)^{\ast} \dashv \Sh(F)_{\ast}\cong (C_{F})^{\ast}=D_{F}:=(-\circ F^{\textup{op}}) \dashv (C_{F})_{\ast}
		\]
		
		\item The weak morphisms $\Sh(F)$ and $C_{F}$ form an adjoint pair in the $2$-category of Grothendieck toposes, weak morphisms and geometric transformations between them ($\Sh(F)$ being the right adjoint and $C_{F}$ the left adjoint); indeed, $\Sh(F)^{\ast}$ is left adjoint to $(C_{F})^{\ast}$. 
		
		In particular, if $F$ is a morphism of sites then the geometric morphisms $\Sh(F)$ and $C_{F}$ form an adjoint pair in the $2$-category $\mathfrak{Btop}$ of Grothendieck toposes, geometric morphisms and geometric transformations. 
		
		\item The weak morphism $\Sh(F):\Sh({\cal C}, J) \to \Sh({\cal D}, K)$ is local if and only if $C_{F}$ is an inclusion, that is, if and only if $F$ is $K$-faithful and $K$-full. 
		
		\item The canonical geometric transformation 
		\[
		1_{\Sh({\cal D}, K)} \to \Sh(F)\circ C_{F}
		\] 
		(given by the unit of the adjunction between $\Sh(F)$ and $C_{F}$  as in (ii)) is an isomorphism if (and only if) $F$ is $K$-faithful and $K$-full. In this case, if $F$ is moreover a morphism of sites $({\cal D}, K)\to ({\cal C}, J)$, the morphisms $C_{F}$ and $\Sh(F)$ realize the topos $\Sh({\cal D}, K)$ as a (coadjoint\footnote{Recall that a retract $(i, r)$ in a $2$-category (in the sense of section B1.1 of \cite{El}), where $r\circ i\cong 1$, is said to be \emph{coadjoint} if $i$ is left adjoint to $r$.}) retract of $\Sh({\cal C}, J)$ in the $2$-category $\mathfrak{Btop}$.		
	\end{enumerate} 
\end{theorem}

\begin{proofs}
(i) The fact that $C_{F}$ is essential follows from Corollary \ref{correstrictioncomorphism}. By definition of $(C_{F})^{\ast}$, for any $J$-sheaf $P$ on $\cal C$, $(C_{F})^{\ast}(P)=a_{K}(P\circ F^{\textup{op}})$. But, since $F$ is $(K, J)$-continuous, $P\circ F^{\textup{op}}$ is already a $K$-sheaf, so $(C_{F})^{\ast}=(-\circ F^{\textup{op}})$, and by definition of the weak morphism of toposes $\Sh(F)$ associated with a weak morphism of sites $F$, we have $\Sh(F)_{\ast}=D_{F}:=(-\circ F^{\textup{op}})$. The other isomorphism in (i) follow from this one in light of the uniqueness (up to isomorphism) of (right or left) adjoints to a given functor. 

(ii) This follows immediately from (i). 

(iii) The right adjoint $(C_{F})_{\ast}$ of $\Sh(F)_{\ast}$ is full and faithful if and only if $C_{F}$ is a geometric inclusion. But this holds, by Proposition \ref{propinclusioncomoprhismofsites}, precisely when $F$ is $K$-faithful and $K$-full.

(iv) We have a canonical geometric transformation $\alpha:1_{\Sh({\cal D}, K)} \to \Sh(F)\circ C_{F}$, given by the unit of the adjunction between $\Sh(F)$ and $C_{F}$. We want to show that this transformation is an isomorphism if and only if $F$ is $K$-faithful and $K$-full. Let $l:{\cal C}\to \Sh({\cal C}, J)$ and $l':{\cal D}\to \Sh({\cal D}, K)$ be the canonical functors. Note that $\alpha^{\ast}:1_{\Sh({\cal D}, K)}\to (C_{F})^{\ast}\circ \Sh(F)^{\ast}$ is determined by its values at the objects of the form $l'(d)$ (for $d\in {\cal D}$), since all the functors involved preserve arbitrary colimits. Now, for any $d\in {\cal D}$, $\alpha^{\ast}(d)$ is the natural transformation 
\[
l'(d)\to (C_{F})^{\ast}(\Sh(F)^{\ast}(l'(d)))=(C_{F})^{\ast}(l(F(d)))\cong l(F(d))\circ F^{\textup{op}}
\]
which corresponds to the element of $(l(F(d))\circ F^{\textup{op}})(d)=l(F(d))(F(d))$ associated with the identity on $F(d)$. Now, by Proposition \ref{propcompositionsheafificationcomorphism}, we have an isomorphism
\[
a_{K}(y_{\cal D}(F(d))\circ F^{\textup{op}})\cong l(F(d))\circ F^{\textup{op}},
\]
and it is immediate to see that the natural transformation $\alpha^{\ast}(d)$ corresponds under this isomorphism to the natural transformation
\[
l'(d)\to a_{K}(y_{\cal D}(F(d))\circ F^{\textup{op}})
\]  
considered in Proposition \ref{proplocalfullnessandfaithfulness}. It thus follows from that proposition that $\alpha^{\ast}(d)$ is an isomorphism for every $d\in {\cal D}$ if and only if $F$ is $K$-faithful and $K$-full, as required. 
\end{proofs}

\bigskip
\textbf{Acknowledgements:} We gratefully acknowledge MIUR for the support in the form of a ``Rita Levi Montalcini'' position, and IH\'ES, where a significant part of this work has been written. We also thank Laurent Lafforgue and Riccardo Zanfa for carefully reading and commenting on a preliminary draft of this text.

\vspace{1cm}

\textsc{Olivia Caramello} 

\vspace{0.3cm}
{\small \textsc{Dipartimento di Scienza e Alta Tecnologia, Universit\`a degli Studi dell'Insubria, via Valleggio 11, 22100 Como, Italy.}\\
	\emph{E-mail address:} \texttt{olivia.caramello@uninsubria.it}}

\vspace{0.5cm}

{\small \textsc{Institut des Hautes \'Etudes Scientifiques, 35 Route de Chartres,
		91440 Bures-sur-Yvette, France.}\\
	\emph{E-mail address:} \texttt{olivia@ihes.fr}}

\end{document}